\documentclass[11pt, amssymb]{article}
\usepackage{amsthm, amsmath, amssymb}
\newtheorem{theorem}{Theorem}[section]
\newtheorem{mtheorem}{Main Theorem}
\newtheorem{proposition}{Proposition}[section]

\newtheorem{hyp}{Hypothesis}[section]
\newtheorem{lemma}{Lemma}[section]
\newcommand{\bb}{\bar{b}}

\newcommand{\bd}{\bar{d}}
\newcommand{\bc}{\bar{c}}

\begin{document}
\title{On Normalized Table Algebras Generated by a Faithful  Non-real Element of Degree 3-II}
\author{\small Zvi Arad\\
\small Department of Mathematics and Computer Science,\\
\small Bar-Ilan University, Ramat Gan 52900,Israel\\
\small and\\
\small Netanya Academic College,
\small Kiryat Rabin, Netanya, Israel\\
\small Guiyun Chen\\
\small  Department of Mathematics,
\small Southwest University, \\
\small 400715, Chongqing, P. R. China. Email:gychen@swu.edu.cn\\
\small Arisha Haj Ihia Hussam,\\
\small Department of Mathematics,\\
\small Alqasemi Academic College of Education,\\
\small P. O. Box 124, Baqa El-Gharbieh 30100, Israel\\
\small and\\
\small Department of Mathematics,
\\
\small Beit Berl Academic College,\\
\small Beit Berl, Israel, 44905}\maketitle \baselineskip 16pt

\section{Introduction}

\par The concept of ``table algebra'' was introduced by Z. Arad and H. Blau
in \cite{ab} in order to study in a uniform way properties of
products of
 conjugacy classes and of irreducible characters of a finite group.
\\
\par {\bf DEFINITION.} A table algebra (A, B) is a finite dimensional,
commutative algebra A with identity element 1 over the complex
numbers
 $C$, and a distinguished base $B=\{b_1=1,\ b_2,\ ...,\ b_k\}$ such that
  the following properties
hold ( where $(b_i,\ a)$ denotes the coefficient of $b_i$ in $a\in
A$, $a$ written as a linear combination of B; and where $R^*$
denotes $R^+\cup\{0\}$,
 the set of non-negative real numbers):
\par (I) For all i, j, m, $b_ib_j=\sum_m\delta_{ijm}b_m$ with $\delta_{ijm}$
a nonnegative real number.
\par (II) There is an algebra automorphism (denoted by $^-$) of A,
whose order divides 2, such that $b_i\in B$ implies that
$\bar{b_i}\in B$.
 ( Then $\bar{i}$ is defined by $\bar{b_i}=b_{\bar{i}}$).
\par (III) Hypothesis (II) holds and there is a function $g$:$B\times B\rightarrow R^+$(the positive reals) such that
$$(b_m,b_ib_j)=g(b_i,b_j)\cdot (b_i,\bar{b_j}b_m),$$
where $g(b_i, b_m)$ is independent of j for all i, j, m.
\par B is called the {\bf table basis} of (A, B). Clearly $1\in B$, we
always use $b_1$ to denote base element 1, and $B^{\sharp}$ to
denote $B\setminus\{b_1\}$. The elements of B are called the {\bf
irreducible components} of (A, B), and nonzero nonnegative  linear
combinations of elements of B with coefficients in $R^+$ are
called components. If $a=\sum^k_{m=1}\lambda_mb_m$ is a
{$component$ ($\lambda_m\in R^*)$ then
$Supp\{a\}=\{b_m|\lambda_m\not=0\}$ is called the set of {\bf
irreducible constituents} of a. An element $a\in A$ is
 called a {\bf  real element} if $a=\bar{a}$.
\par Two table algebras $(A,\ B)$ and $(A',\ B')$ are called isomorphic
(denoted $B\cong B'$) when there exists an algebra isomorphism
$\psi:\ A\rightarrow A'$ such that $\psi (B)$ is a rescaling of
$B'$; and the algebras are called $exactly\ isomorphic$(denoted
$B\cong_xB'$) when $\psi(B)=B'$. So $B\cong_xB'$ means that $B$
and $B'$ yield the same structure constants.
\par Proposition 2.2 of \cite{ab} shows that if (A, B) is a table algebra,
then there exists a basis $B'$,  which consists of suitable
positive real
 scalar multiples $b_i'$ of the elements $b_i$ of $B$, $g(b'_i,b'_j)=1$ for
  any $b'_i,$ $b'_j\in B'$. Such a basis $B'$ is called a {\bf normalized
  basis}. Now $Supp(b'_{i_1}b'_{i_2}\cdots b'_{i_t})$ consists of the
  corresponding scalar multiples of $Supp(b_{i_1}b_{i_2}\cdots b_{i_t})$,
  for any sequence $i_1,\ i_2,\ ...,\ i_t$ of indices. So in the proof of
  any proposition which identifies the irreducible constituents of a
  product of basis elements, we may assume that B is normalized.
\par Suppose that B is normalized. It follows from (III) and \cite{ab} {\bf Sect. 2} that A has a positive definite Hermitian form, with B as an orthonormal basis, and such that
$$(a,bc)=(a\bb,c)$$
for all $a,\ b,\ c\in \mathcal{R}B$.

\par It follows from Section 2 in \cite{ab}  that there exists an algebra
homomorphism f from A to $\mathit{C}$ such that $f(B)\subseteq
R^+$. For an element $b\in B$, $f(b)$ is called the degree of b.
\par Each finite group G yields two natural
table algebras: the table algebra of conjugacy classes (denoted by
$(\mathit{ZC}(G)$, $\mathit{Cl}(G))$) and the table algebra of
generalized characters (denoted by ($\mathit{Ch}(G),\
\mathit{Irr}(G))$.
\par Both table algebras arising from group theory have an additional
property: their structure constants and degrees are non-negative
integers. Such algebras were defined in \cite{hb3} as {\bf
integral table algebras} (abbrev. as ITA).
\par Each of the elements of a table algebra is contained in a unique table
subalgebra which may be considered as a table subalgebra generated
by this element. So it is natural to start the study of integral
table algebras from
 those which are generated by a single element. Normalized integral table
 algebras (abbrev. as NITA) generated by
 an element of degree 2 were completely classified by H. Blau in \cite{hb}.

\par The finite linear groups in dimension $n\leq 7$ have been completely
classified. See Feit \cite{feit}. Deep theory and properties of
finite groups are heavily used for the proofs. For n=2, 3, 4 see
Blichfeldt \cite{blic}. For n=5 see Brauer \cite{brauer}. For n=6
see Lindsey \cite{lindsey}. For n=7 see Wales \cite {wales}. In
order to generalize these results to normalized integral table
algebra, the authors
 began to classify
normalized integral table algebras  $(A, B)$ generated by a
faithful non-real element of degree 3 and without nonidentity
irreducible
elements of degree 1 or 2 in \cite{cgyarad}, and came to the following theorem:\\
\begin{theorem}\label{mthm1} Let $(A,\ B)$ be NITA generated by a non-real element $b_3\in B$ of degree 3 and without non-identity basis element of degree 1 or 2. Then  $b_3\bb_3=1+b_8$ and one of the following holds:
\par (1) $(A,\ B)\cong_x(\mathit{Ch}$$(PSL(2, 7)),\ Irr$$(PSL(2, 7))$;
\par (2) $b^2_3=b_4+b_5$, where $b_4\in B$ and $b_5\in B$.
\par (3) $b^2_3=\bb_3+b_6$, where $b_6\in B$ nonreal;
\par (4) $b^2_3=c_3+b_6$, where $c_3,\ b_6\in B$ and $c_3\not=b_3,\ \bb_3$;

\end{theorem}

\par The purpose of this paper is to investigate  NITA generated by a non-real element $b_3$  and satisfies (2),  (3) and (4) in {\bf  Theorem \ref{mthm1}}. The following theorem is proved:
\newpage
\begin{mtheorem}\nonumber\label{mthm3} Let $(A,\ B)$ be NITA generated by a non-real element $b_3\in B$ of degree 3 and without non-identity basis element of degree 1 or 2. Then $b_3\bb_3=1+b_8$, $b_8\in B$ and one of the following hold:
\par (1) There exists a real element $b_6\in B$ such that $b^2_3=\bb_3+b_6$ and $(A,\ B)\cong_x(\mathit{Ch}$$(PSL(2, 7))$,\ $Irr$$(PSL(2, 7))$;
\par (2) There exist $b_4,\ b_5\in B$ such that $b^2_3=b_4+b_5$, $(b_3b_8,\ b_3b_8)=3$ and $(b^2_4,\ b^2_4)=3$.
\par (3) There exist $b_6,\  b_{10}, b_{15}\in B$, where $b_6$ is non-real such that
\[
\begin{array}{lll}
b^2_3=\bb_3+b_6, &\bb_3b_6=b_3+b_{15}, &b_3b_6=b_8+b_{10}.
\end{array}
\]
and $(b_3b_8,\ b_3b_8)=3$. Moreover if $b_{10}$ is real, then
 $(A,\ B)\cong_x (\mathit{CH}$$(3\cdot A_6),\ \mathit{Irr}$$(3\cdot A_6))$;
\par (4)  There exist  $c_3,\ b_6\in B$, $c_3\not=b_3,\ \bb_3$,  such that $b^2_3=c_3+b_6$ and either
$(b_3b_8,\ b_3b_8)=3$ or $4$.
\par When $(b_3b_8,\ b_3b_8)=3$. If $c_3$
is non-real, then $(A,\ B)$ is exactly  isomorphic to $(A(3\cdot
A_6\cdot 2),\ B_{32})$ (see Theorem 7.2). If $c_3$ is real, then
$(A,\ B)$ is exactly isomorphic to the NITA $(A(7\cdot 5\cdot
10),\ B_{22})$ of dimension 22.(see  Theorem 8.1).
\end{mtheorem}
\par {\bf Proof.} The theorem follows from {\bf Theorems 1.1, 3.1, 5.1, 7.2}.

\section{Two NITA Generated by a Non-real element of Degree 3 not derived From a Group and Lemmas}
\par Here are two  examples of NITA not induced from group theory, which has either 32 or 22 basis elements. The NITA of dimension 32 contains a subalgebra that is strictly isomorphic to the algebra of dimension 17 of characters of the group $3\cdot A_6$, we denote it as $(A(3\cdot A_6\cdot 2),\ B_{32})$, where its basis is $B_{32}$.

For $B_{32}$ there are three table subsets: $\{1\}\subseteq
C\subseteq D\subseteq B$, where
\[
\begin{array}{rll}
C&=&\{1,\ b_8,\ x_{10},\ b_5,\ c_5,\ c_8,\ x_9\},\\
D&=&C\cup\{c_3,\ \bc_3,\ d_3,\ \bar{d}_3,\ c_9,\ \bc_9,\ b_6,\ \bb_6,\ y_{5},\ \bar{y}_{15}\},\\
B_{32}&=&D\cup\{b_3,\ \bb_3,\ r_3,\ x_6,\ \bar{x}_6,\ s_6,\ b_9,\
\bb_9,\ x_{15},\ \bar{x}_{15},\ t_{15},\ d_9,\ y_3,\ z_3,\
\bar{z}_3\}.
\end{array}
\]
\par From the below equations one can check that $B_{32}$ has the following table subsets: $\{1\}\subseteq\ C\subseteq E=C\cup\{r_3,\ s_6,\ t_{15},\ d_9,\ y_3\}\subseteq B_{32}$. The table subsets $E$ and $D$ are two maximal table subsets of $B_{32}$. \\

\par The structure algebra constants of $C$.

\[
\begin{array}{rll}
b^2_5&=&1+x_9+x_{10}+b_5,\\
b_5c_5&=&c_8+b_8+x_9,\\
b_5c_8&=&c_5+c_8+b_8+x_9+x_{10},\\
b_5b_8&=&c_5+c_8+b_8+x_9+x_{10},\\
b_5x_9&=&c_8+c_5+b_8+b_5+x_9+x_{10},\\
b_5x_{10}&=&c_8+b_8+x_9+b_5+2x_{10},\\
\hline
\end{array}
\]
\[
\begin{array}{rll}
c^2_5&=&1+x_9+x_{10}+c_5,\\
c_5c_8&=&b_5+c_8+b_8+x_9+x_{10},\\
c_5b_8&=&b_5+c_8+b_8+x_9+x_{10},\\
c_5x_9&=&b_8+c_8+b_5+c_5+x_9+x_{10},\\
c_5x_{10}&=&b_8+c_8+x_9+c_5+2x_{10},\\
\hline
\end{array}
\]
\[
\begin{array}{rll}
b^2_8&=&1+b_5+c_5+c_8+2b_8+x_9+2x_{10},\\
b_8x_9&=&b_5+c_5+2c_8+b_8+2x_9+2x_{10},\\
b_8x_{10}&=&b_5+c_5+2c_8+2b_8+2x_9+2x_{10}.\\
\hline
\end{array}
\]
\[
\begin{array}{rll}

c^2_8&=&1+b_5+c_5+2c_8+b_8+x_9+2x_{10},\\
c_8b_8&=&b_5+c_5+c_8+b_8+2x_9+2x_{10},\\
c_8x_9&=&b_5+c_5+c_8+2b_8+2x_9+2x_{10},\\
c_8x_{10}&=&b_5+c_5+2c_8+2b_8+2x_9+2x_{10},\\
\hline
\end{array}
\]
\[
\begin{array}{rll}

x^2_9&=&1+b_5+c_5+2c_8+2b_8+2x_9+2x_{10},\\
x_9x_{10}&=&b_5+c_5+2c_8+2b_8+2x_9+3x_{10},\\
\hline
x^2_{10}&=&1+2b_5+2c_5+2c_8+2b_8+3x_9+2x_{10},\\
\end{array}
\]
\par The structure algebra constants of $D$(Since $C\subseteq D$ we will not repeat the equations of $C$ as shown above):
\[
\begin{array}{rll}
c_3\bc_3&=&1+b_8,\\
c^2_3&=&\bc_3+\bb_6,\\
c_3d_3&=&x_9,\\
c_3\bd_3&=&\bar{c}_9,\\
c_3b_6&=&\bc_3+\bar{y}_{15},\\
c_3\bb_6&=&b_8+x_{10},\\
c_3c_9&=&d_3+\bar{c}_9+\bar{y}_{15},\\
c_3\bc_9&=&c_8+x_9+x_{10},\\
c_3y_{15}&=&\bb_6+2\bar{y}_{15}+\bar{c}_9,\\
c_3\bar{y}_{15}&=&b_5+c_5+c_8+b_8+x_9+x_{10},\\
c_3b_8&=&c_3+b_6+y_{15},\\
c_3x_{10}&=&b_6+c_9+y_{15},\\

c_3b_5&=&y_{15},\\
c_3c_5&=&y_{15},\\
c_3c_8&=&y_{15}+c_9,\\
c_3x_9&=&y_{15}+c_9+\bar{d}_3,\\
\hline
\end{array}
\]
\[
\begin{array}{rll}
d_3\bd_3&=&1+c_8,\\
d^2_3&=&\bd_3+b_6,\\
d_3b_6&=&c_8+x_{10},\\
d_3\bb_6&=&\bar{d}_3+y_{15},\\

d_3c_9&=&b_8+x_9+x_{10},\\
d_3\bc_9&=&c_3+y_{15}+c_9,\\
d_3y_{15}&=&b_5+c_5+c_8+b_8+x_9+x_{10},\\
d_3\bar{y}_{15}&=&2y_{15}+b_6+c_9,\\
d_3b_8&=&\bar{c}_9+\bar{y}_{15},\\
d_3x_{10}&=&\bar{y}_{15}+\bb_6+\bc_9,\\

d_3b_5&=&\bar{y}_{15},\\
d_3c_5&=&\bar{y}_{15},\\
d_3c_8&=&\bar{y}_{15}+b_6+d_3,\\
d_3x_9&=&\bar{y}_{15}+\bc_3+\bc_9,\\
\hline
\end{array}
\]
\[
\begin{array}{rll}
b_6\bb_6&=&1+b_5+c_5+c_8+b_8+x_9,\\
b^2_6&=&2\bb_6+\bar{y}_{15}+\bc_9,\\
b_6c_9&=&2\bar{y}_{15}+2\bc_9+\bb_6,\\
b_6\bc_9&=&b_5+c_5+c_8+b_8+2x_9+x_{10},\\
b_6\bar{y}_{15}&=&b_5+c_5+2c_8+2b_8+2x_9+3x_{10},\\
b_6y_{15}&=&4\bar{y}_{15}+\bb_6+\bc_3+d_3+2\bc_9,\\
b_6b_5&=&b_6+y_{15}+c_9,\\
b_6c_5&=&y_{15}+b_6+c_9,\\
b_6c_8&=&\bar{d}_3+b_6+2y_{15}+c_9,\\
b_6b_8&=&c_3+b_6+c_9+2y_{15},\\
b_6x_9&=&2y_{15}+b_6+2c_9,\\
b_6x_{10}&=&3y_{15}+c_3+\bar{d}_3+c_9,\\
\hline
\end{array}
\]
\[
\begin{array}{rll}
y_{15}\bar{y}_{15}&=&1+3b_5+3c_5+5c_8+5b_8+6x_9+6x_{10},\\
y^2_{15}&=&9\bar{y}_{15}+4\bb_6+2\bc_3+2d_3+6\bc_9,\\
y_{15}b_5&=&3y_{15}+b_6+c_3+\bd_3+2c_9,\\
y_{15}c_5&=&3y_{15}+b_6+c_3+\bd_3+2c_9,\\
y_{15}c_8&=&5y_{15}+c_3+\bd_3+2b_6+3c_9,\\
y_{15}b_8&=&5y_{15}+c_3+\bd_3+2b_6+3c_9,\\
y_{15}x_9&=&6y_{15}+2b_6+\bd_3+c_3+3c_9,\\
y_{15}x_{10}&=&6y_{15}+3b_6+c_3+\bd_3+4c_9,\\
\hline
\end{array}
\]
\[
\begin{array}{rll}
c_9\bc_9&=&1+2c_8+2b_8+2x_9+2x_{10}+c_5+b_5,\\
c_9^2&=&3\bar{y}_{15}+d_3+\bar{c}_3+2\bb_6+2\bc_9,\\
c_9y_{15}&=&6\bar{y}_{15}+2\bb_6+\bc_3+d_3+3\bc_9,\\
c_9\bar{y}_{15}&=&2b_5+2c_5+3c_8+3b_8+3x_9+4x_{10},\\
c_9b_5&=&2y_{15}+b_6+c_9,\\
c_9c_5&=&2y_{15}+b_6+c_9,\\
c_9b_8&=&3y_{15}+\bar{d}_3+b_6+2c_9,\\
c_9c_8&=&3y_{15}+b_6+c_3+2c_9,\\
c_9x_9&=&c_3+2c_9+3y_{15}+\bar{d}_3+2b_6,\\
c_9x_{10}&=&c_3+2c_9+4y_{15}+\bar{d}_3+b_6,\\
\hline
\end{array}
\]
\par From the equations one can see that $C$ and $D$ are table subsets of $B_{32}$. Furthermore the NITA subalgebra $(<D>,\ D)$ of dimension 17 is strictly isomorphic to $(CH(3\cdot A_6),\ Irr(3\cdot A_6))$ while our NITA of dimension 32 is not induced from a finite group $G$.
\par The following equations give out  other products of basis elements in $B_{32}$:
\[
\begin{array}{l}
b_3\bb_3=1+b_8,\\
b^2_3=c_3+b_6,\\
b_3r_3=\bc_3+\bb_6,\\
b_3x_6=c_3+y_{15},\\
b_3\bar{x}_6=b_8+x_{10},\\
b_3s_6=\bc_3+\bar{y}_{15},\\
b_3b_9=y_{15}+c_9+\bar{d}_3,\\
b_3\bb_9=x_9+c_8+x_{10},\\
b_3x_{15}=b_6+2y_{15}+c_9,\\
b_3\bar{x}_{15}=b_5+c_5+b_8+c_8+x_9+x_{10},\\
b_3t_{15}=\bb_6+\bc_9+2\bar{y}_{15},\\
b_3d_9=\bc_9+\bar{y}_{15}+d_3,\\
b_3y_3=\bc_9,\\
b_3z_3=x_9,\\
b_3\bar{z}_3=c_9,\\
b_3b_8=b_3+x_6+x_{15},\\
b_3x_{10}=x_6+x_{15}+b_9,\\
b_3b_5=x_{15},\\
b_3c_5=x_{15},\\
b_3c_8=x_{15}+b_9,\\
b_3x_9=\bar{z}_3+b_9+x_{15},\\
b_3c_3=r_3+s_6,\\
b_3\bc_3=\bb_3+\bar{x}_6,\\
b_3d_3=\bb_9, \\
b_3\bar{d}_3=d_9,\\
b_3c_9=t_{15}+d_9+y_3,\\
b_3\bc_9=\bb_9+\bar{x}_{15}+z_3,\\
b_3b_6=r_3+t_{15},\\
b_3\bb_6=\bb_3+\bar{x}_{15},\\
b_3y_{15}=s_6+2t_{15}+d_9, \\
b_3\bar{y}_{15}=\bar{x}_6+2\bar{x}_{15}+\bb_9.\\
\hline
\end{array}
\]
\[
\begin{array}{l}
r^2_3=1+b_8,\\
r_3x_6=\bc_3+\bar{y}_{15},\\
r_3s_6= b_8+x_{10},\\
r_3b_9=d_3+\bc_9+\bar{y}_{15},\\
r_3x_{15}=\bar{b}_6+\bar{c}_9+2\bar{y}_{15},\\
r_3t_{15}=b_5+c_5+b_8+c_8+x_9+x_{10},\\
r_3d_9=c_8+x_9+x_{10},\\
r_3y_3=x_9,\\
r_3z_3=c_9,\\
r_3b_8=s_6+r_3+t_{15},\\
r_3x_{10}=s_6+t_{15}+d_9,\\
r_3b_5=t_{15}, \\
r_3c_5=t_{15},\\
r_3c_8=t_{15}+d_9,\\
r_3x_9=t_{15}+d_9+y_3,\\
r_3c_3=\bb_3+\bar{x}_6,\\
r_3d_3=b_9,\\
r_3c_9=\bb_9+\bar{x}_{15}+z_3,\\
r_3b_6=\bar{x}_{15}+\bb_3,\\
r_3y_{15}=\bar{x}_6+2\bar{x}_{15}+\bb_9,\\
\hline
\\
x^2_6=2b_6+y_{15}+c_9,\\
x_6\bar{x}_6=1+b_8+b_5+c_5+c_8+x_9, \\
x_6s_6=2\bb_6+\bc_9+\bar{y}_{15},\\
x_6b_9=\bb_6+2\bc_9+2\bar{y}_{15},\\
x_6\bb_9=b_5+c_5+x_{10}+b_8+c_8+2x_9,\\
x_6x_{15}=4y_{15}+b_6+2c_9+d_3+c_3,\\
x_6\bar{x}_{15}=b_5+c_5+2b_8+3x_{10}+2c_8+2x_9,\\
x_6t_{15}=\bc_3+d_3+2\bc_9+4\bar{y}_{15}+\bb_6\\
x_6d_9=\bb_6+2\bc_9+2\bar{y}_{15},\\
x_6y_3=d_3+\bar{y}_{15},\\
x_6z_3=x_{10}+c_8.\\
x_6\bar{z}_3=\bar{d}_3+y_{15},\\
x_6b_8=b_3+x_6+2x_{15}+b_9,\\
x_6x_{10}=b_9+b_3+\bar{z}_3+3x_{15},\\
\end{array}
\]
\[
\begin{array}{l}
x_6b_5=x_6+b_9+x_{15},\\
x_6c_5=x_6+b_9+x_{15},\\
x_6c_8=x_6+b_9+2x_{15}+\bar{z}_3,\\
x_6x_9=2b_9+x_6+2x_{15},\\
x_6c_3=r_3+t_{15},\\
x_6\bar{c}_3=\bar{b}_3+\bar{x}_{15},\\
x_6d_3=z_3+\bar{x}_{15},\\
x_6\bar{d}_3=y_3+t_{15},\\
x_6c_9=2t_{15}+2d_9+s_6,\\
x_6\bar{c}_9=2\bar{x}_{15}+\bar{x}_6+2\bb_9,\\
x_6b_6=2s_6+t_{15}+d_9\\
x_6\bar{b}_6=2\bar{x}_6+\bar{b}_9+\bar{x}_{15},\\
x_6y_{15}=r_3+4t_{15}+s_6+2d_9+y_3,\\
x_6\bar{y}_{15}=\bb_3+z_3+\bar{x}_6+4\bar{x}_{15}+2\bb_9, \\

\hline
\end{array}
\]
\[
\begin{array}{l}

s^2_6=1+b_5+c_5+b_8+c_8+x_9,\\
s_6b_9=\bb_6+2\bc_9+2\bar{y}_{15},\\
s_6x_{15}=\bc_3+d_3+\bb_6+2\bc_9+4\bar{y}_{15},\\
s_6t_{15}=b_5+c_5+2b_8+2c_8+2x_9+3x_{10},\\
s_6d_9=b_5+c_5+2x_9+b_8+c_8+x_{10},\\
s_6y_3=c_8+x_{10},\\
s_6z_3=\bar{d}_3+y_{15},\\
s_6b_8=r_3+s_6+2t_{15}+d_9,\\
s_6x_{10}=r_3+y_3+d_9+3t_{15},\\
s_6b_5=s_6+d_9+t_{15},\\
s_6c_5=s_6+t_{15}+d_9,\\
s_6c_8=y_3+s_6+d_9+2t_{15},\\
s_6x_9=s_6+2d_9+2t_{15},\\
s_6c_3=\bar{b}_3+\bar{x}_{15},\\
s_6d_3=\bar{z}_3+x_{15},\\
s_6c_9=\bar{x}_6+2\bb_9+2\bar{x}_{15},\\
s_6b_6=2\bar{x}_6+\bb_9+\bar{x}_{15},\\
s_6y_{15}=\bb_3+z_3+4\bar{x}_{15}+\bar{x}_6+2\bb_9.\\

\hline
\end{array}
\]
\[
\begin{array}{l}
b_9\bb_9=1+b_5+c_5+2b_8+2c_8+2x_9+2x_{10},\\
b^2_9=\bar{d}_3+c_3+2b_6+2c_9+3y_{15},\\
b_9x_{15}=6y_{15}+3c_9+\bar{d}_3+2b_6+c_3,\\
b_9\bar{x}_{15}=2b_5+2c_5+3c_8+3b_8+4x_{10}+3x_9,\\
b_9t_{15}=d_3+\bc_3+6\bar{y}_{15}+3\bc_9+2\bb_6,\\
b_9d_9=\bc_3+d_3+2\bb_6+2\bc_9+3\bar{y}_{15},\\
b_9y_3=\bc_3+\bar{y}_{15}+\bc_9,\\
b_9z_3=b_8+x_9+x_{10},\\
b_9\bar{z}_3=c_9+c_3+y_{15},\\
b_9b_8=\bar{z}_3+x_6+2b_9+3x_{15},\\
b_9x_{10}=\bar{z}_3+b_3+2b_9+x_6+4x_{15},\\
b_9b_5=x_6+b_9+2x_{15},\\
b_9c_5=b_9+2x_{15}+x_6,\\
b_9c_8=r_3+t_{15}+2b_9+x_6+2x_{15},\\
\hline
\end{array}
\]
\[
\begin{array}{l}
b_9x_9=2b_9+3x_{15}+\bar{z}_3+b_3+2x_6\\
b_9c_3=t_{15}+d_9+\bar{y}_3,\\
b_9\bc_3=z_3+\bb_9+\bar{x}_{15},\\
b_9d_3=\bb_3+\bb_9+\bar{x}_{15},\\
b_9\bar{d}_3=r_3+t_{15}+d_9\\
b_9c_9=r_3+y_3+2s_6+2d_9+3t_{15},\\
b_9\bar{c}_9=\bar{b}_3+z_3+2\bar{x}_6+3\bar{x}_{15}+2\bar{b}_9,\\
b_9b_6=s_6+2d_9+2t_{15},\\
b_9\bb_6=2\bar{b}_9+2\bar{x}_{15}+\bar{x}_6\\
b_9\bar{y}_{15}=6\bar{x}_{15}+\bar{b}_3+z_3+2\bar{x}_6+3\bar{b}_9,\\
b_9y_{15}=r_3+y_3+2s_6+6t_{15}+3d_9\\
\hline
\\
x^2_{15}=2c_3+2\bar{d}_3+4b_6+6c_9+9y_{15}\\
x_{15}\bar{x}_{15}=1+6x_9+3b_5+3c_5+5b_8+5c_8+6x_{10}\\
x_{15}t_{15}=4\bb_6+6\bc_9+9\bar{y}_{15}+2\bc_3+2d_3\\
x_{15}d_9=\bc_3+d_3+2\bb_6+3\bc_9+6\bar{y}_{15}\\
x_{15}y_3=\bb_6+\bc_9+2\bar{y}_{15},\\
x_{15}\bar{z}_3=b_6+c_9+2y_{15},\\
x_{15}z_3=b_5+c_5+b_8+c_8+x_9+x_{10}\\
\end{array}
\]
\[
\begin{array}{l}

x_{15}b_8=b_3+\bar{z}_3+2x_6+5x_{15}+3b_9,\\
x_{15}x_{10}=b_3+\bar{z}_3+3x_6+4b_9+6x_{15}\\
x_{15}b_5=\bar{z}_3+b_3+x_6+2b_9+3x_{15}\\
x_{15}c_5=b_3+\bar{z}_3+x_6+2b_9+3x_{15}\\
x_{15}c_8=b_3+\bar{z}_3+2x_6+3b_9+5x_{15}\\
x_{15}x_9=b_3+\bar{z}_3+2x_6+3b_9+6x_{15}\\
x_{15}c_3=2t_{15}+s_6+d_9,\\
x_{15}\bar{c}_3=2\bar{x}_{15}+\bar{x}_6+\bar{b}_9\\
x_{15}d_3=\bar{x}_6+2\bar{x}_{15}+\bb_9\\
x_{15}\bar{d}_3=s_6+2t_{15}+d_9,\\
x_{15}c_9=r_3+y_3+2s_6+6t_{15}+3d_9\\
x_{15}\bc_9=\bb_3+z_3+3\bb_9+2\bar{x}_6+6\bar{x}_{15},\\
x_{15}b_6=r_3+s_6+y_3+2d_9+4t_{15},\\
x_{15}\bb_6=\bb_3+z_3+2\bb_9+4\bar{x}_{15}+\bar{x}_6,\\
x_{15}y_{15}=2y_3+2r_3+4s_6+6d_9+9t_{15}\\
x_{15}\bar{y}_{15}=2z_3+9\bar{x}_{15}+2\bb_3+4\bar{x}_6+6\bb_9.\\
\hline
\end{array}
\]
\[
\begin{array}{ll}

t^2_{15}=1+5b_8+5c_8+3b_5+3c_5+6x_{10}+6x_9\\
t_{15}d_9=2b_5+2c_5+3b_8+3c_8+3x_9+4x_{10}\\
t_{15}y_3=b_5+c_5+b_8+c_8+x_9+x_{10}\\
t_{15}z_3=b_6+c_9+2y_{15}\\
t_{15}b_8=r_3+y_3+5t_{15}+2s_6+3d_9\\
t_{15}x_{10}=r_3+y_3+4d_9+3s_6+6t_{15}\\
t_{15}b_5=r_3+y_3+s_6+2d_9+3t_{15}\\
t_{15}c_5=r_3+y_3+s_6+2d_9+3t_{15}\\
t_{15}c_8=r_3+y_3+2s_6+3d_9+5t_{15}\\
t_{15}x_9=r_3+y_3+3d_9+2s_6+6t_{15}\\
t_{15}c_3=\bar{x}_6+\bb_9+2\bar{x}_{15},\\
t_{15}d_3=x_6+b_9+2x_{15}\\
t_{15}c_9=\bb_3+z_3+3\bb_9+2\bar{x}_6+6\bar{x}_{15}\\
t_{15}b_6=\bb_3+z_3+\bar{x}_6+2\bb_9+4\bar{x}_{15}\\
t_{15}y_{15}=2\bb_3+2z_3+4\bar{x}_6+6\bb_9+9\bar{x}_{15}\\
\hline
\end{array}
\]
\[
\begin{array}{l}

d^2_9=1+b_5+c_5+2b_8+2c_8+2x_9+2x_{10}\\
d_9y_3=x_9+b_8+x_{10},\\
d_9z_3=c_3+c_9+y_{15},\\
d_9b_8=y_3+s_6+2d_9+3t_{15},\\
d_9x_{10}=y_3+r_3+s_6+2d_9+4t_{15}\\
d_9b_5=s_6+d_9+2t_{15}\\
d_9c_5=s_6+d_9+2t_{15}\\
d_9c_8=r_3+s_6+2d_9+3t_{15}\\
d_9x_9=r_3+y_3+2s_6+2d_9+3t_{15}\\
d_9c_3=\bb_9+z_3+\bar{x}_{15},\\
d_9d_3=b_3+b_9+x_{15}\\
d_9c_9=\bb_3+z_3+2\bar{x}_6+2\bb_9+3\bar{x}_{15}\\
d_9b_6=\bar{x}_6+3\bb_9+2\bar{x}_{15}, \\
d_9y_{15}=\bb_3+z_3+2\bar{x}_6+3\bb_9+6\bar{x}_{15}\\

\hline
\end{array}
\]
\[
\begin{array}{ll}

y^2_3=1+c_8,\\
y_3z_3=\bar{d}_3+b_6,\\
y_3b_8=d_9+t_{15},\\
y_3x_{10}=s_6+t_{15}+d_9,\\
y_3b_5=t_{15},\\
y_3c_5=t_{15},\\
y_3c_8=s_6+t_{15}+y_3,\\
y_3x_9=t_{15}+d_9+r_3,\\
y_3c_3=\bb_9,\\
y_3d_3=x_6+\bar{z}_3,\\
y_3c_9=\bb_3+\bb_9+\bar{x}_{15},\\
y_3b_6=\bar{x}_{15}+z_3,\\
y_3y_{15}=\bar{x}_6+\bb_9+2\bar{x}_{15}.\\
\hline
\\
z_3\bar{z}_3=1+c_8,\\
z^2_3=d_3+\bb_6,\\
z_3b_8=\bb_9+\bar{x}_{15},\\
z_3x_{10}=\bar{x}_6+\bb_9+\bar{x}_{15},\\
z_3b_5=\bar{x}_{15},\\
z_3c_5=\bar{x}_{15},\\
z_3c_8=\bar{x}_6+\bar{x}_{15}+z_3,\\

\end{array}
\]
\[
\begin{array}{ll}

z_3x_9=\bb_3+\bb_9+\bar{x}_{15}\\
z_3c_3=b_9\\
z_3\bc_3=d_9,\\
z_3d_3=s_6+y_3,\\
z_3\bar{d}_3=\bar{z}_3+x_6,\\
z_3c_9=b_3+b_9+x_{15},\\
z_3\bc_9=r_3+d_9+t_{15},\\
z_3b_6=\bar{z}_3+x_{15}\\
z_3\bb_6=y_3+t_{15},\\
z_3y_{15}=x_6+b_9+2x_{15},\\
z_3\bar{y}_{15}=s_6+2t_{15}+d_9.\\

\hline
\end{array}
\]

\par Now let us see the following  example of NITA of dimension 22, which satifies $b^2_3=r_3+s_6$, $r_3$ is  real, denoted as $(A(7\cdot 5\cdot 10),\ B_{22})$, where $B_{22}=\{1, b_8,\ x_{10},\ b_5,\ c_5,\ c_8,\ x_9,\ r_3,\ y_3,\ s_6,\ t_{15},\ d_9,\ b_3,\ \bb_3,\ t_6,\ \bar{t}_6,\ b_{15},\ \bb_{15},\ y_9,\ \bar{y}_9,\ x_3,\bar{x}_3\}$. The interesting thing is that this example also contains table subsets:
$$\{1\}\subseteq C\subseteq E\subseteq B_{22},$$
where $C$ and $E$ are the same as defined in the previous example
of dimension 32. In particular both examples of dimensions 22 and
32 contain each the same sub-table algebra of dimension 12
generated by the basis $E$. Also $E$ is a maximal table subset of
$B_{22}$. The equations of products of the basis elements of $C$
and $E$ are as in the example of $A(3\cdot A_6\cdot 2,\ 32)$
described above..
\[
\begin{array}{l}
b_3\bb_3=1+b_8,\\
b^2_3=r_3+s_6,\\
b_3t_6=b_8+x_{10},\\
b_3\bar{t}_6=r_3+t_{15},\\
b_3b_{15}=2t_{15}+s_6+d_9,\\
b_3\bar{b}_{15}=b_8+x_{10}+b_5+c_5+c_8+x_9,\\
b_3y_9=t_{15}+d_9+y_3,\\
b_3\bar{y}_9=c_8+x_9+x_{10},\\
b_3x_3=d_9,\\
b_3\bar{x}_3=x_9,\\
\end{array}
\]
\[
\begin{array}{l}

b_3b_8=b_3+\bar{t}_6+b_{15},\\
b_3x_{10}=b_{15}+\bar{t}_6+y_9,\\
b_3b_5=b_{15},\\
b_3c_5=b_{15},\\
b_3c_8=b_{15}+y_9,\\
b_3x_9=b_{15}+y_9+x_3,\\
b_3r_3=\bb_3+t_6,\\
b_3y_3=\bar{y}_9,\\
b_3s_6=\bb_3+\bb_{15},\\
b_3t_{15}=2\bb_{15}+t_6+\bar{y}_9,\\
b_3d_9=\bb_{15}+\bar{y}_9+\bar{x}_3.\\
\hline
 \\
t_6\bar{t}_6=1+b_5+c_5+b_8+c_8+x_9,\\
t^2_6=2s_6+t_{15}+d_9,\\
t_6b_{15}=2b_8+2c_8+b_5+c_5+3x_{10}+2x_9,\\
t_6\bb_{15}=c_3+b_6+4t_{15}+y_3+2d_9,\\
t_6y_9=b_8+b_5+c_5+c_8+2x_9+x_{10},\\
t_6\bar{y}_9=s_6+2d_9+2t_{15},\\
t_6x_3=x_{10}+c_8,\\
t_6\bar{x}_3=t_{15}+y_3,\\
t_6b_8=\bb_3+2\bb_{15}+t_6+\bar{y}_9,\\
t_6x_{10}=3\bb_{15}+\bar{y}_9+\bb_3+\bar{x}_3,\\

\end{array}
\]
\[
\begin{array}{l}
t_6b_5=\bb_{15}+t_6+\bar{y}_9,\\
t_6c_5=\bb_{15}+t_6+\bar{y}_9,\\
t_6c_8=2\bb_{15}+t_6+\bar{y}_9+\bar{x}_3,\\
t_6x_9=2\bb_{15}+t_6+2\bar{y}_9,\\
t_6r_3=b_{15}+b_3,\\
t_6y_3=b_{15}+x_3,\\
t_6s_6=2\bar{t}_6+b_{15}+y_9,\\
t_6t_{15}=2y_9+x_3+b_3+\bar{t}_6+4b_{15},\\
t_6d_9=2y_9+2b_{15}+\bar{t}_6.\\
\hline
\\
b_{15}\bb_{15}=1+3b_5+3c_5+5b_8+5c_8+6x_{10}+6x_9,\\
b^2_{15}=2r_3+2y_3+4s_6+6d_9+9t_{15},\\
b_{15}y_9=r_3+y_3+2s_6+3d_9+6t_{15},\\
b_{15}\bar{y}_9=2b_5+2c_5+3b_8+3c_8+3x_9+4x_{10},\\
b_{15}x_3=s_6+2t_{15}+d_9,\\
b_{15}\bar{x}_3=c_5+c_8+x_9+b_8+x_{10}+b_5,\\
b_{15}b_8=5b_{15}+3y_9+x_3+b_3+2\bar{t}_6,\\
\end{array}
\]
\[
\begin{array}{l}
b_{15}x_{10}=6b_{15}+4y_9+b_3+x_3+3\bar{t}_6,\\
b_{15}b_5=b_3+x_3+\bar{t}_6+2y_9+3b_{15},\\
b_{15}c_5=3b_{15}+2y_9+b_3+x_3+\bar{t}_6,\\
b_{15}c_8=5b_{15}+3y_9+x_3+b_3+2\bar{t}_6,\\
b_{15}x_9=6b_{15}+3y_9+x_3+b_3+2\bar{t}_6,\\
b_{15}r_3=2\bb_{15}+t_6+\bar{y}_9,\\
b_{15}y_3=2\bb_{15}+\bar{y}_9+t_6,\\
b_{15}s_6=4\bb_{15}+2\bar{y}_9+\bb_3+\bar{x}_3+t_6,\\
b_{15}t_{15}=9\bb_{15}+6\bar{y}_9+2\bb_3+2\bar{x}_3+4t_6,\\
b_{15}d_9=6\bb_{15}+3\bar{y}_9+\bb_3+\bar{x}_3+2t_6.\\
\hline
\\
y_9\bar{y}_9=1+b_5+c_5+2b_8+2c_8+2x_9+2x_{10},\\
y^2_9=r_3+y_3+2s_6+2d_9+3t_{15},\\
y_9x_3=r_3+d_9+t_{15},\\
y_9\bar{x}_3=b_8+x_{10}+x_9,\\
y_9b_8=\bar{t}_6+2y_9+3b_{15}+x_3,\\
y_9x_{10}=\bar{t}_6+2y_9+b_3+x_3+4b_{15},\\
y_9b_5=\bar{t}_6+y_9+2b_{15},\\
y_9c_5=\bar{t}_6+y_9+2b_{15},\\

\end{array}
\]
\[
\begin{array}{l}
y_9c_8=3b_{15}+b_3+2y_9+\bar{t}_6,\\
y_9x_9=2\bar{t}_6+2y_9+3b_{15}+x_3+b_3,\\
y_9r_3=\bb_{15}+\bar{y}_9+\bar{x}_3,\\
y_9y_3=\bb_{15}+\bb_3+\bar{y}_9,\\
y_9s_6=2\bb_{15}+t_6+2\bar{y}_9,\\
y_9t_{15}=6\bb_{15}+2t_6+\bb_3+3\bar{y}_9+\bar{x}_3,\\
y_9d_9=3\bb_{15}+2t_6+\bb_3+2\bar{y}_9+\bar{x}_3.\\
\hline
\end{array}
\]
\[
\begin{array}{l}
x_3\bar{x}_3=1+c_8,\\
x^2_3=s_6+y_3,\\
x_3b_8=b_{15}+y_9,\\
x_3x_{10}=b_{15}+\bar{t}_6+y_9,\\
x_3b_5=b_{15},\\
x_3c_5=b_{15},\\
x_3c_8=x_3+b_{15}+\bar{t}_6,\\
x_3x_9=b_3+y_9+b_{15},\\
x_3r_3=\bar{y}_9,\\
x_3y_3=\bar{x}_3+t_6,\\
x_3s_6=\bb_{15}+\bar{x}_3,\\
x_3t_{15}=2\bb_{15}+t_6+\bar{y}_9,\\
x_3d_9=\bb_{15}+\bb_3+\bar{y}_9.\\
\hline
\end{array}
\]

\begin{table}
    \caption{NITA $(A(3\cdot A_6\cdot 2), B_{32})$}
    \begin{center}
\begin{tabular}{|c|c|c|}
\hline\mbox{The powers of $b_3$}&\mbox{contituents of the powers}&\mbox{Basis  of quotient of $B_{32}/C$}\\
\hline &\\
$b_3^2$&$c_3,\ b_6$&$C=\{1,\ b_8,\ x_{10},\ b_5,\ c_5,\ c_8,\ x_9\}$\\
$b_3^3$&$r_3,\ s_6,\ t_{15}$&$Cb_3=\{b_3,\ x_6,\ x_{15},\ b_9,\ \bar{z}_3\}$\\
$b_3^4$&$\bc_3,\ \bb_6,\ \bar{y}_{15},\ \bc_9$&$C\bb_3=\{\bb_3,\ \bar{x}_6,\ \bar{x}_{15},\ \bb_9,\ z_3\}$\\
$b_3^5$&$\bb_3,\ \bar{x}_6,\ \bar{x}_{15},\ \bb_9,\ z_3$&$Cr_3=\{r_3,\ t_{15},\ d_9,\ y_3,\ s_6\}$\\
$b_3^6$&$1,\ b_8,\ x_{10},\ b_5,\ c_5,\ c_8,\ x_9$&$Cc_3=\{c_3,\ b_6,\ y_{15},\ c_9,\ \bar{d}_3\}$\\
$b_3^7$&$b_3,\ x_6,\ x_{15},\ b_9,\ \bar{z}_3$&$C\bar{c}_3=\{\bar{c}_3,\ \bb_6,\ \bar{y}_{15},\ \bc_9,\ d_3\}$\\
$b_3^8$&$c_3,\ b_6,\ y_{15},\ c_9,\ \bar{d}_3$&\\
$b_3^9$&$r_3,\ s_6,\ t_{15},\ d_9,\ y_3$&\\
$b_3^{10}$&$\bc_3,\ \bb_6,\ \bar{y}_{15},\ \bc_9,\ d_3$&\\
\hline
\end{tabular}
\end{center}
\end{table}
\par If we define in a table algebra $(A,\ B)$ for $R,\ S\subseteq B$,
$$R*S=\{\cup Supp(R\cdot S)|r\in R,\ s\in S\},$$
 then:
$Cb_3*C\bb_3=Cr^2_3=Cc_3*C\bc_3=C,\ C\mbox{ is identity.}$
$(Cb_3)^2=C\bc_3$,\ $Cb_3*Cc_3=Cr_3$,\ $Cb_3*C\bc_3=C\bb_3$ and
$Cb_3Cr_3=C\bc_3$. For the definitions of abelian table algebra,
quotient of table algebra and its basis and the sub-table algebra
generated by an element $b\in B$ denoted by $<b>=B_b$ see
\cite{ab} and \cite{hb3}. Also the definitions of linear and
faithful element $b\in B$ are found there.
\par One can derived from our table that $B_{32}=\cup^{10}_{i=1}Bb^i_3=\cup^{\infty}_{i=0}Bb_3^i$ and $b_3$ is a faithful element of $B_{32}$ generates $B_{32}$. Also the quotient table algebra with basis $B_{32}/C$ is strictly isomorphic the table algebra induced by the cyclic group $Z_6$ of order 6.
\begin{table}
    \caption{NITA $(A(7\cdot 5\cdot 10),\ B_{22})$}
    \begin{center}
\begin{tabular}{|c|c|c|}
\hline\mbox{The powers of $b_3$}&\mbox{new contituents of the powers}&\mbox{Basis  of quotient of $B_{32}/C$}\\
\hline
$b_3^2$&$r_3,\ s_6$&$C=\{1,\ b_8,\ x_{10},\ b_5,\ c_5,\ c_8,\ c_9\}$\\
$b_3^3$&$\bb_3,\ t_6,\ \bb_{15}$&$C\bb_3=\{\bb_3,\ t_6,\ \bb_{15},\ \bar{y}_9,\ \bar{x}_3\}$\\
$b_3^4$&$1,\ b_8,\ x_{10},\ b_5,\ c_5,\ c_8,\ x_9$&$Cr_3=\{r_3,\ s_6,\ t_{15},\ d_9,\ y_3\}$\\
$b_3^5$&$b_3,\ \bar{t}_6,\ b_{15},\ y_9,\ x_3$&$Cb_3=\{b_3,\ \bar{t}_6,\ b_{15},\ y_9,\ x_3\}$\\
$b_3^6$&$r_3,\ s_6,\ t_{15},\ d_9,\ y_{15},\ y_3$&\\
$b_3^7$&$\bb_3,\ t_6,\ \bb_{15},\ \bar{y}_9,\ \bar{x}_3$&\\
\hline
\end{tabular}
\end{center}
\end{table}
\newpage
\par One can derived from this table that $B_{22}=\cup^7_{i=1}Bb_3^i=\cup^{\infty}_{i=0}Bb^i_3$ and $b_3$ is a faithful element of $B_{22}$ generates $B_{22}$. Also the quotient table algebra with basis $B_{22}/C$ is strictly isomorphic to the table algebra induced by the cyclic group $Z_n$  of order 4.\\
\par In order to prove Main Theorems, the following fact is frequently used:
\begin{lemma}\label{la1}
 Let $e_m,\ f_m,\ u_n,\ v_n\ \in B$ such that $e_m\bar{e}_m=f_m\bar{f}_m$. Then $(e_mu_n,\ e_mu_n)=(\bar{e}_mu_n,\ \bar{e}_mu_n)=(f_mv_n,\ f_mv_n)=(\bar{f}_mv_n,\ \bar{f}_mv_n).$
\end{lemma}
\par {\bf Proof.} (1) follows from $(e_mu_n,\ e_mu_n)=(e_m\bar{e}_m,\ u_n\bar{u}_n)=(f_m\bar{f}_m,\ u_n\bar{u}_n)=(f_mu_n,\ f_mu_n).$
\par Notice that $t_3\bar{t}_3,\ b_8),\ (s_4\bar{s}_4,\ b_8)\leq 1$, we come to (2).

\begin{lemma}\label{t} Let $b_3,\ t_3,\ s_4\in B$, $b_3\bb_3=1+b_8$. Then
$(b_3t_3,\ b_3t_3)$, ($b_3s_4,\ b_3s_4)\leq 2$. If $(b_3t_3,\
b_3t_3)=2$, then $t_3\bar{t}_3=1+b_8$.
\end{lemma}
\par {\bf Proof.}  If $(b_3x_4,\ b_3x_4)\geq 3$, then $(x_4\bar{x}_4,\ b_8)\geq 2$ for $b_3\bb_3=1+b_8$, a contradiction. (1) follows.
\par If $(b_3c_3,\ b_3c_3)=2$, then $(b_3\bb_3,\ c_3\bc_3)=2$. (2) follows from $b_3\bb_3=1+b_8$.
\\

\section{NITA Generated by $b_3$ and Satisfying $b^2_3=b_4+b_5$}
\par In this section, we shall investigate NITA generated by $b_3$ and satisfying $b^2_3=b_4+b_5$ and come to the following theorem:
\begin{theorem}\label{mthm2} Let $(A,\ B)$ be NITA generated by a non-real element
$b_3\in B$ of degree 3 and without non-identity basis element of
degree 1 or 2. Then $b_3\bb_3=1+b_8$, $b_8\in B$. Assume that
$b^2_3=b_4+b_5$, $b_4\in B$ and $b_5\in B$. Then  $(\bb_5b_5,
b_8)=1$,  $(b^2_4,\ b^2_4)=3$ and $(b_3b_8,\ b_3b_8)=3$.
\end{theorem}
Now we start to investigate NITA satifying
\begin{eqnarray}
\label{1.1}b_3\bb_3&=&1+b_8,\\
\label{1.2}b^2_3&=&b_4+b_5.
\end{eqnarray}
The Proof of this theorem is in the rest of this section.

Based on the above equations, one may set that
\begin{eqnarray}\label{1.3}
\bb_3b_4=b_3+b_9,\ \mbox{some}\ b_9\in\mathit{N}^*B.
\end{eqnarray}
Since  $(b_4\bb_4,\ b_8)\leq 1$, we have that $(\bb_3b_4,\
\bb_3b_4)=2$, $b_9\in B$ and
\begin{eqnarray}\label{1.4}
b_4\bb_4=1+b_8+c_7, \mbox{ where } c_7\in N^*B.
\end{eqnarray}
\par Collecting the above equations, checking the associative laws for basis elements and making necessary assumptions, we have the following equations:
\begin{lemma}\label{be} The following equations always hold:
\[
\begin{array}{rlll}
\bb_3b_4&=&b_3+b_9,&\\
b_4\bb_4&=&1+b_8+c_7, c_7\in N^*B&\\
\bb_3b_5&=&b_3+x_{12},\ \ x_{12}\in\mathit{N}^*B& \mbox{assumption},\\
b_3b_8&=&b_3+b_9+x_{12}\ &\mbox{calculating}\ b^2_3\bb_3=b_3(b_3\bb_3),\\
b_3c_7&=&b_9+f_{12},\  f_{12}\in\mathit{N}^*B,& \mbox{assumption},\\
b_4\bb_9&=&b_3+b_9+x_{12}+f_{12},&\mbox{calculating}\ b_3(\bb_4b_4)=(b_3\bb_4)b_4,\\
b_4b_8&=&b_5+b_3b_9&\mbox{calculating}\ b_3(\bb_3b_4)=(b_3\bb_3)b_4,\\
\bb_4b_5&=&b_8+z_{12},\ z_{12}\in\mathit{N}^*B&\mbox{assumption},\\
\bb_3b_9&=&c_7+b_8+\bar{z}_{12}&\mbox{calculating}\ \bb^2_3b_4=\bb_3(\bb_3b_4),\\
b_8^2&=&1+c_7+b_8+\bar{z}_{12}+\bb_3x_{12}&\mbox{calculating}\ (b_3\bb_3)b_8=(b_3b_8)\bb_3,\\
b_8c_7&=&b_8+\bar{z}_{12}+\bb_3f_{12}&\mbox{calculating}\ (b_3\bb_3)c_7=(b_3b_7)\bb_3,\\
b_5\bb_9&=&b_9+b_3z_{12},&calculating\ (\bb_3b_4)\bb_5=\bb_3(b_4\bb_5),\\
b_9\bb_9&=&1+b_8-z_{12}+\bar{z}_{12}+\bb_3f_{12}+\bb_3x
_{12}&\mbox{caculating}\ (b_4\bb_4)b_8=\bb_4(b_4b_8).
\end{array}
\]
\end{lemma}
\par {\bf Proof.} Based on the known equations, checking associative laws one by one and making necessary assumptions, one can easily  come to the equations above. Here is the proof of the last equation. It follows from
\[
\begin{array}{rll}
(b_4\bb_4)b_8&=&b_8+b^2_8+b_8c_7\\
&=&b_8+1+c_7+b_8+\bar{z}_{12}+\bb_3x_{12}+b_8+\bar{z}_{12}+\bb_3f_{12},\\
\bb_4(b_4b_8)&=&\bb_4b_5+(\bb_4b_3)b_9\\
&=&\bb_4b_5+\bb_3b_9+b_9\bb_9\\
&=& b_8+z_{12}+c_7+b_8+\bar{z}_{12}+b_9\bb_9.
\end{array}
\]
\\
\par Obviously $(b_5\bb_5,\ b_8)\leq 3$,  $c_7$ is either
$c_3+c_4$ or irreducible.
\par If $c_7\not\in B$, then there exist $c_3,\ c_4\in B$ such that $c_7=c_3+c_4$. At first we have the following lemma:
\begin{lemma}\label{lla1} There exits no NITA
 generated by $b_3$ and satisfying $b_3\bb_3=1+b_8$,
 $b^2_3=b_4+b_5$ and $c_7=c_3+c_4$.
\end{lemma}
\par {\bf Proof.}  By the assumption, we have that
 \begin{eqnarray}\label{1.5}
 b_4\bb_4=1+b_8+c_3+c_4.
 \end{eqnarray}
 Then
\[ \begin{array}{rll}
b^2_3\bb_4&=&b_4\bb_4+b_5\bb_4\\
&=&1+c_3+c_4+b_8+\bb_4b_5,\\
b_3(b_3\bb_4)&=&b_3\bb_3+b_3\bb_9\\
&=&1+b_8+b_3\bb_9.
\end{array}
\]
Hence $(c_3,\ b_3\bb_9)=(c_4,\ b_3\bb_9)=1$. Therefore

\begin{eqnarray}
\label{1.6}b_3c_3&=&b_9,\\
\label{1.7}b_3c_4&=&d_3+b_9,\ \mbox{some}\ d_3\in B,\ d_3\neq b_3.
\end{eqnarray}
Now one has $(b_3\bar{d}_3,\ c_4)=1$. Then
\begin{eqnarray}\label{1.8}b_3\bar{d}_3=c_4+d_5,\ \mbox{some}\ d_5\in B,
\end{eqnarray}
from which $(b_3\bb_3,\ d_3\bar{d}_3)=2$ follows. Hence
\begin{eqnarray}\label{1.9}
d_3\bar{d}_3=1+b_8.
\end{eqnarray}
\par By (\ref{1.5}), we have that $(b_4c_3,\ b_4)=1$. Let
\begin{eqnarray}\label{1.10}
b_4c_3=b_4+y_8,\ y_8\in\mathit{N}^*B.
\end{eqnarray}
\par Since $(b_3\bb_3)c_3=c_3+c_3b_8,\ \bb_3(b_3c_3)=\bb_3b_9,$
one has that $\bb_3b_9=c_3+c_3b_8$. But $(c_4,\ \bb_3b_9)=1$. Then
$(c_4,\ c_3b_8)=1$ and so $(c_4c_3,\ b_8)=1$. Thus
\begin{eqnarray}
\label{1.11}c_3c_4=x_4+b_8, \mbox{some}\ x_4\in B.
\end{eqnarray}
 By (\ref{1.10}) we have that
 $$2\leq(b_4c_3,\ b_4c_3)=(b_4\bb_4,\ c^2_3)=(1+c_3+c_4+b_8,\ c^2_3),$$
which implies that $c^2_3$ equals one of
$$1+b_8,\ 1+c_3+x_5,\ 1+c_4+y_4,\ 1+2c_4.$$
\par By (\ref{1.11}), we know that $c^2_3$ cannot be $1+c_4+y_4$ or $1+2c_4$.
And by (\ref{1.6}), $c^2_3\neq 1+b_8$. Hence
\begin{eqnarray}\label{1.12}
c^2_3=1+c_3+x_5.
\end{eqnarray}
Therefore $(b_4c_3,\ b_4c_3)=(c^2_3,\ b_4\bb_4)=2$, which implies
that $y_8\in B$.
\par By {\bf Lemma} \ref{be}, (\ref{1.8}) and (\ref{1.9}), we have
\[\begin{array}{rl}
b_3(b_3\bb_3)&=b_3+b_3b_8=2b_3+b_9+x_{12},\\
\bb_3(b_3^2)&=\bb_3b_4+\bb_3b_5=b_3+b_9+\bb_3b_5,\\
(b_3\bar{d}_3)d_3&=c_4d_3+d_5d_3.
\end{array}\]
\par Hence $b_3b_8=b_9+\bb_3b_5,\
\mbox{and}\ b_9\in Supp\{c_4d_3+d_5d_3\}.$
\par Assume that $b_9\in Supp\{d_3d_5\}$. By (\ref{1.8}), it follow that
\[
\begin{array}{rll}
d_3d_5&=&b_9+b_3+l_3,\ l_3\in B,\\
\bb_3b_5&=&c_4d_3+l_3,\\
\end{array}
\]
so that $b_3b_8=b_9+l_3+c_4d_3,$ which implies that $b_8\in
Supp\{b_3\bar{l}_3\}$.
 Hence $l_3=b_3$. So $(d_3d_5,\ b_3)=3$, which implies that $(b_4\bar{d}_3,\ d_5)=2$,  a contradiction to (\ref{1.8}). Thus $b_9\in Supp\{c_4d_3\}$ and
 \begin{eqnarray}
 \label{1.13}c_4d_3=b_3+b_9.
 \end{eqnarray}
Thus $(d_3\bar{d}_3,\ c^2_4)=2$, we may assume that
\begin{eqnarray}
\label{1.16} c^2_4=1+b_8+f_7,\ \mbox{some}\ f_7\in \mathit{N}^*B.
\end{eqnarray}
Then $c_4\not\in Supp\{c_4x_5\}$. Otherwise, $c_4\in
Supp\{c_4x_5\}$, which implies that $x_5\in Supp\{c^2_4\}$, and so
$f_7=x_5+x_2$, a contradiction.
\par Since
\[
\begin{array}{rll}
c^2_3c_4&=&c_4+c_3c_4+x_5c_4\\
&=&c_4+x_4+b_8+c_4x_5,\\
c_3(c_3c_4)&=&c_3b_8+c_3x_4,
\end{array}
\]
we have that
$$c_3b_8+c_3x_4=c_4+x_4+b_8+c_4x_5.$$
But left side of the above equation contains two $c_4$ by
(\ref{1.11}). Hence $c_4=x_4$ by $c_4\not\in Supp\{c_4x_5\}$.
Therefore
\begin{eqnarray}\label{x1}c_3c_4=b_8+c_4,\end{eqnarray}
which implies that $(c^2_4,\ c_3)=1$ and
$$
c^2_4=1+b_8+c_3+e_4, \mbox{some}\ e_4\in B.
$$
Now we have following equations
\[
\begin{array}{rll}
b_3c^2_4&=&b_3+b_3b_8+b_3c_3+b_3e_4\\
&=&2b_3+2b_9+x_{12}+b_3e_4,\\
\mbox{By (7) and (13)\ \ \ }c_4(b_3c_4)&=&c_4b_9+c_4d_3\\
&=&c_4b_9+b_9+b_3.
\end{array}
\]
Then $c_4b_9=b_3+b_9+x_{12}+b_3e_4.$ Since $d_3\in
Supp\{c_4b_9\}$, one has that $d_3\in Supp\{x_{12}\}$ or $b_3e_4$.
\par If $d_3\in Supp\{x_{12}\}$, then $d_3\in Supp\{b_3b_8\}$, which implies that $b_3\bar{d}_3=1+b_8$.
Thus $b_3=d_3$, which implies that $b_3\bb_3=c_4+d_5$, a
contradiction. Hence $d_3\in Supp\{b_3e_4\}$, $e_4\in
Supp\{b_3\bar{d}_3\}$. Therefore by (8) $e_4=c_4$ and
\[
\begin{array}{rll}
\mbox{by (7)\ \ }c^2_4&=&1+c_3+c_4+b_8=b_4\bb_4,\\
c_4b_9&=&b_3+d_3+2b_9+x_{12}.
\end{array}
\]
Furthermore
\[
\begin{array}{rll}
b_3c^2_4&=&b_3+b_3b_8+b_3c_4+b_3c_3\\
&=&b_3+b_3+b_9+x_{12}+d_3+b_9+b_9,\\
(b_3\bb_4)b_4&=&\bb_3b_4+\bb_9b_4\\
&=&b_3+b_9+b_4\bb_9,\\
(b_3c_4)c_4&=&c_4b_9+d_3c_4\\
&=&c_4b_9+b_3+b_9.
\end{array}
\]
Thus
$$b_4\bb_9=b_3+d_3+2b_9+x_{12}=c_4b_9.$$
\par Hence by (6) $(b_3c_3)c_4=b_9c_4=b_3+d_3+2b_9+x_{12}.$ But by (7)
$(b_3c_4)c_3=c_3b_9+c_3d_3.$

Since $d_3\not\in Supp\{c_3d_3\}$, we have that $d_3\in
Supp\{c_3b_9\}$

\begin{eqnarray}
c_3d_3&=&b_9, \nonumber \\
\label{1.19}c_3b_9&=&b_3+d_3+b_9+x_{12}.
\end{eqnarray}
Therefore $(b_3c_3,\ c_3d_3)=1$, so $(b_3\bar{d}_3,\ c^2_3)=1$. By
(\ref{1.8}) and (\ref{1.12}), we have that $(b_3\bar{d}_3,\
x_5)=1$, which implies that $x_5=d_5$. Consequently $d_5$ is real.
Then
\[
\begin{array}{rll}
b_3\bar{d}_3&=&c_4+d_5=\bb_3d_3,\\
c^2_3&=&1+c_3+d_5.
\end{array}
\]
Moreover $(b^2_3,\ d^2_3)=(b_3\bar{d}_3,\ \bb_3d_3)=2.$ So
$$d^2_3=1+2b_4\ \mbox{or}\ b_4+b_5.$$
Since (\ref{1.9}) implies that $d_3$ is non-real. So
$$d^2_3=b_4+b_5.$$
 \par By (\ref{1.8}), we may set that $$\bb_3d_5=\bar{d}_3+u_{12},\ u_{12}\in\mathit{N}^*B.$$
 Then
 \[
\begin{array}{rll}
\bb_3c^2_3&=&\bb_3+\bb_3c_3+\bb_3d_5\\
&=&\bb_3+\bb_9+\bar{d}_3+u_{12},\\
(\bb_3c_3)c_3&=&c_3\bb_9.
\end{array}
\]
Thus
$$c_3b_9=b_3+d_3+b_9+\bar{u}_{12}.$$
Hence $x_{12}=\bar{u}_{12}$ by (\ref{1.19}) and so
$b_3d_5=d_3+x_{12}$. Since $c_4c^2_3=c_4+c_3c_4+d_5c_4$ and
$c_3(c_3c_4)=c_3c_4+c_3b_8$ by (\ref{x1}), we have that
$$c_3b_8=c_4+c_4d_5.$$
Now, by $c^2_4=b_4\bb_4$, we have that
\[
\begin{array}{rll}
c_3c^2_4&=&c_3+c^2_3+c_3c_4+c_3b_8\\
&=&c_3+1+c_3+d_5+c_4+b_8+c_4+c_4d_5,\\
(c_3b_4)\bb_4&=&b_4\bb_4+y_8\bb_4\\
&=&1+c_3+c_4+b_8+\bb_4y_8.
\end{array}
\]
So
$$\bb_4y_8=c_3+c_4+d_5+c_4d_5.$$

\par We have known that $y_8$ is irreducible, then $(c_4b_4,\ y_8)=1$.
But $(b_4c_4,\ b_4c_4)=4$. Hence $(c_4b_4-y_8,\ c_4b_4-y_8)=3$, which is impossible for $L_1(B)=\{1\}$ and $L_2(B)=\emptyset$. The lemma follows.\\

\par Now we begin to investigate NITA such that $c_7\in B$ based on the value of $(b_5\bb_5,\ b_8)$.

\begin{lemma}\label{lla2}There exits no NITA such that $c_7\in B$ and $(b_5\bb_5,\ b_8)=3$.\end{lemma}

\par {\bf Proof.} It is easy to see that $(b_5\bb_5,\ b_8)=3$ implies that
$(\bb_3b_5,\ \bb_3b_5)=4$. By {\bf Lemma \ref{be}}, we may assume
that
\[
\begin{array}{rll}
x_{12}&=&x+y+z,\ x,\ y,\ z\in B,\\
\bb_3b_5&=&b_3+x+y+z,\\
b_3b_8&=&b_3+b_9+x+y+z.
\end{array}
\]
By the last equation above, one has that $b_8$ is contained in
$\bb_3x$, $\bb_3y$ and $\bb_3z$. Hence no one of $x$, $y$ and $z$
has degree 3. Therefore
\[
\begin{array}{rll}
x_{12}&=&x_4+y_4+z_4,\\
\bb_3x_4&=&r_4+b_8,\ \mbox{some}\ r_4\in B,\\
\bb_3y_4&=&s_4+b_8,\ \mbox{some}\ s_4\in B,\\
\bb_3z_4&=&t_4+b_8,\ \mbox{some}\ t_4\in B.
\end{array}
\]
Furthermore $b_3x_4$, $b_3y_4$ and $b_3z_4$ have exactly two
constituents. Let
\[
\begin{array}{rll}
b_3x_4&=&b_5+r_7,\ \mbox{some}\ r_7\in B,\\
b_3y_4&=&b_5+s_7,\ \mbox{some}\ s_7\in B,\\
b_3z_4&=&b_5+t_7,\ \mbox{some}\ t_7\in B.
\end{array}
\]
\par Since
\[
\begin{array}{rll}
\bb_3(b_3x_4)&=&\bb_3b_5+\bb_3r_7\\
&=&b_3x_4+y_4+z_4+\bb_3r_7,\\
b_3(\bb_3x_4)&=&b_3b_8+b_3r_4\\
&=&b_3+b_9+x_4+y_4+z_4+b_3r_4,
\end{array}
\]
we have that $\bb_3r_7=b_9+b_3r_4$. By the same reasoning we have
that $\bb_3s_7=b_9+b_3s_4$ and $\bb_3t_7=b_9+b_3t_4$. Hence $r_7,\
s_7,\ t_7\in Supp\{b_3b_9\}$.
\par If $r_7$, $s_7$ and $t_7$ are distinct, then
$$b_3b_9=b_4+r_7+s_7+t_7+e_2,\ \mbox{some}\ e_2\in B,$$
a contradiction. Hence at least two of $r_7$, $s_7$ and $t_7$ are
equal, without loss of generality, let $r_7=s_7$. Then $(b_3x_4,\
b_3y_4)=(\bb_3x_4,\ \bb_3y_4)=2$. Thus $r_4=s_4$. BY {\bf Lemma
\ref{be}}, we have
$$
b^2_8=1+4b_8+c_7+z_{12}+2r_4+t_4.$$ So $\bar{z}_{12}=2r_4+t_4.$
{\bf Lemma \ref{be}} implies that $$\bb_3b_9=c_7+b_8+2r_4+t_4.$$
Then $(\bb_3b_9,\ r_4)=(b_9,\ b_3r_4)\geq 2$, a contradiction,
which lemma follows
from.\\
\begin{lemma} There exits no NITA such that $c_7\in B$ and $(b_5\bb_5,\ b_8)=2$.\end{lemma}
\par {\bf Proof.} Suppose that $(b_5\bb_5,\ b_8)=2$. Let
\[
b_5\bb_5=1+2b_8+d_8,\ \mbox{some}\ d_8\in \mathit{N}^*B.
\]
Since $(\bb_3b_5,\ \bb_3b_5)=(b_3\bb_3,\ b_5\bb_5)=3$, by {\bf
Lemma \ref{be}}, we have that $x_{12}=x+y$ and
\[
\begin{array}{rll}
\bb_3b_5&=&b_3+x+y,\\
b_3b_8&=&b_3+b_9+x+y.
\end{array}
\]

Thus $b_8\in Supp\{\bb_3x\}$, $Supp\{\bb_3y\}$, which implies that
degrees of $x$
 and $y$ are $\geq 4$. We have three possibilities:
$(f(x),\ f(y))=(4,\ 8),\ (5,\ 7),\ (6,\ 6)$.
\par We shall divide the proof into three propositions. And the lemma follows from
{\bf Proposition \ref{prop3.1}, \ref{prop3.2},
 \ref{prop3.3}.}\\
\begin{proposition}\label{prop3.1}
No NITA such that $b_5\bb_5=1+2b_8+d_8$ and $(f(x),\ f(y))=(4,\
8)$.
\end{proposition}
\par {\bf Proof.} Assume $(f(x),\ f(y))=(4,\ 8)$. Then $x_{12}=x_4+y_8$ and we may set

\begin{eqnarray}
\nonumber\bb_3x_4&=&b_8+x'_4,\ x'_4\in B,\\
\label{x2}\bb_3y_8&=&b_8+y_{16},\ y_{16}\in \mathit{N}^*B,
\end{eqnarray}
Then
$$b^2_8=1+c_7+3b_8+x'_4+\bar{z}_{12}+y_{16}.$$
Hence
\[
\begin{array}{rll}
(b_3\bb_3)^2&=&2+5b_8+c_7+x'_4+\bar{z}_{12}+y_{16},\\
b^2_3\bb^2_3&=&b_4\bb_4+b_4\bb_5+\bb_4b_5+b_5\bb_5\\
&=&1+b_8+c_7+b_4\bb_5+\bb_4b_5+1+2b_8+d_8,
\end{array}
\]
so $z_{12}+d_8=x'_4+y_{16}.$
\par Now two cases appear, $x'_4\in Supp\{d_8\}$ or $x'_4\not\in Supp\{d_8\}$.\\
\par {\bf Case I}. Suppose that $x'_4\in Supp\{z_{12}\}$. Then there exists $s_8$ such that
\[
\begin{array}{rll}
z_{12}&=&x'_4+s_8,\\
y_{16}&=&d_8+s_8.
\end{array}
\]
\par Since $(\bb_4b_5,\ \bb_4b_5)=(b_4\bb_4,\ b_5\bb_5)=3$, we have that $s_8\in B$ and
$$\bb_3b_9=c_7+b_8+\bar{x'}_4+\bar{s}_8.$$

Hence $b_9\in Supp\{b_3x'_4\}$, which means that $(b_3\bar{x}'_4,\
b_3\bar{x}'_4)=2$. So $(b_3x'_4,\ b_3x'_4)=2$. Since
$\bb_3x_4=b_8+x'_4$, we may set $b_3\bar{x}'_4=a_3+b_9$ and
$b_3x'_4=x_4+e_8,\ \mbox{some}\ a_3,\ e_8\in B$.
\par Suppose that $b_3\bar{a}_3=x'_4+l_5,\ \mbox{some}\ l_5\in B$. Then
\[
\begin{array}{rll}
b_3(\bb_3x'_4)&=&b_3\bb_9+b_3\bar{a}_3\\
&=&c_7+b_8+x'_4+s_8+x'_4+l_5,\\
\bb_3(b_3x'_4)&=&\bb_3x_4+\bb_3e_8\\
&=&b_8+x'_4+\bb_3e_8.
\end{array}
\]
Hence
$$\bb_3e_8=c_7+s_8+x'_4+l_5.$$
Therefore $e_8\in Supp\{b_3c_7\}$. Let $b_3c_7=b_9+e_8+e_4,\
e_4\in B$. By $(b_3\bb_3)\bar{a}_3=\bb_3(b_3\bar{a}_3)$, we come
to
$$\bar{a}_3b_8=\bb_9+\bb_3l_5.$$
\par Since $(\bar{a}_3b_8,\ \bar{a}_3b_8)=(a_3\bar{a}_3,\ b^2_8)=(b_3\bar{b}_3,\ b^2_8)=4$, we have that $(\bb_3l_5,\ \bb_3l_5)=3$. So $(b_3l_5,\ b_3l_5)=3$. Let
$$b_3l_5=e_8+\alpha_3+\beta_4,$$
where $\alpha_3,\ \beta_4\in B.$ Obviously $\alpha_3\not=b_3$.
\[
\begin{array}{rll}
b_3(\bb_3e_8)&=&b_3c_7+b_3b_8+b_3x'_4+b_3l_5\\
&=&2b_9+3e_8+e_4+b_3+2x_4+y_8+\alpha_3+\beta_4,\\
(b_3\bb_3)e_8&=&e_8+e_8b_8.
\end{array}
\]
 Then
$$e_8b_8=2b_9+2e_8+e_4+b_3+2x_4+y_8+\alpha_3+\beta_4.$$
 Hence $e_8\in Supp\{b_3b_8\}$, which implies that $e_8=y_8$. Hence $b_3c_7=b_9+y_8+e_4$, so that $c_7\in Supp\{\bb_3y_8\}$, which implies that $c_7\in Supp\{y_{16}\}=\{d_8,\ s_8\}$, a contradiction.\\
\par {\bf Case II.} NITA satisfying $x'_4\in Supp\{ d_8\}$.\\
\par If $x'_4\in d_8$, then $d_8=x'_4+\bar{x}'_4$ ($x'_4$ maybe real). Then
\[
\begin{array}{rll}
b_5\bb_5&=&1+2b_8+x'_4+\bar{x}'_4.
\end{array}
\]
Hence $y_{16}=\bar{x}'_4+z_{12}$ and
\begin{eqnarray}
\label{x3}b^2_8&=&1+c_7+3b_8+x'_4+\bar{x}'_4+z_{12}+\bar{z}_{12},\\
\label{x4}\bb_3y_8&=&b_8+\bar{x}'_4+z_{12},\\
\label{x5}b_9\bb_9&=&1+3b_8+x'_4+z_{12}+\bar{x}'_4+\bb_3f_{12},
\end{eqnarray}
by (\ref{x2}). Therefore $(y_8,\ b_3\bar{x}'_4)=1$. Let
\begin{eqnarray}\label{x7}b_3\bar{x}'_4=y_8+t_4,\ b_3x'_4=x_4+j_8,\ \mbox{where}\ t_4,\ j_8\in B.\end{eqnarray}
Since
\[
\begin{array}{rll}
b_3(b_5\bb_5)&=&b_3+2b_3b_8+b_3x'_4+b_3\bar{x}'_4\\
&=&3b_3+3b_9+2x_4+3y_8+t_4+x_4+j_8,\\
b_5(b_3\bb_5)&=&b_5\bb_3+b_5\bar{x}_4+b_5\bar{y}_8\\
&=&b_3+x_4+y_8+b_5\bar{x}_4+b_5\bar{y}_8.
\end{array}
\]
we have that
$$\bar{x}_4b_5+\bar{y}_8b_5=2b_3+2b_9+2x_4+2y_8+t_4+j_8.$$
Since
\[
\begin{array}{rll}
x'_4(\bb_3b_4)&=&b_3x'_4+b_9x'_4\\
&=&x_4+j_8+b_9x'_4,\\
(\bb_3x'_4)b_4&=&\bar{y}_8b_4+\bar{t}_4b_4,
\end{array}
\]
we have that $x_4\in Supp\{b_4\bar{y}_8\}$ or
$Supp\{\bar{t}_4b_4\}$.
\par We assert that $x_4\in Supp\{b_4\bar{t}_4\}$. Otherwise  $x_4\in Supp\{b_4\bar{y}_8\}$, and $y_8\in Supp\{b_4\bar{x}_4\}.$
But $b_9\in Supp\{b_4\bar{x}_4\}$ by {\bf Lemma \ref{be}}, a
contradiction. Hence $x_4\in Supp\{b_4\bar{t}_4\}$,
 which implies that $t_4\in Supp\{b_4\bar{x}_4\}$. Therefore
\begin{eqnarray}\label{x6}b_4\bar{x}_4=t_4+b_9+r_3,\ r_3\in B.\end{eqnarray}
Since
\[
\begin{array}{rll}
\bb_3(b_4b_8)&=&\bb_3b_5+\bb_3(b_3b_9)\\
&=&b_3+x_4+y_8+b_9+b_3(\bb_3b_9)\\
&=&b_3+x_4+y_8+b_9+b_3c_7+b_3b_8+b_3\bar{z}_{12}\\
&=&b_3+x_4+y_8+b_9+b_9+f_{12}+b_3+x_4+y_8+b_9+b_3\bar{z}_{12},\\
b_4(\bb_3b_8)&=&b_4\bb_3+b_4\bar{x}_4+b_4\bar{y}_8+b_4\bb_9\\
&=&b_3+b_9+t_4+b_9+r_3+b_4\bar{y}_8+b_3+b_9+x_4+y_8+f_{12},\\
(\bb_3b_4)b_8&=&b_3b_8+b_8b_9\\
&=&b_3+x_4+y_8+b_9+b_8b_9,
\end{array}
\]
 we have that
\[
\begin{array}{rll}
x'_4b_8&=&b_8+z_{12}+b_3\bar{t}_4,\\
b_8b_9&=&b_3+x_4+y_8+b_9+y_{12}+b_3\bar{z}_{12},\\
b_4\bar{y}_8+b_9+t_4+r_3&=&x_4+y_8+b_3\bar{z}_{12}
\end{array}
\]
We continue the proof based on $t_4=x_4$ or $t_4\neq x_4$.
\\
\par {\bf Subcase I.} No NITA such that $t_4=x_4$.\\
\par If $t_4=x_4$, then $(x_4,\ b_3\bar{x}'_4)=1$ by $b_3\bar{x}'_4=y_8+t_4$, which implies that
$x'_4\in Supp\{b_3\bar{x}_4\}=\{b_8,\ \bar{x}'_4\}$. And
$\bb_3t_4=\bb_3x_4=b_8+x'_4$. Hence $x'_4$ is real and
\begin{eqnarray}
b_5\bb_5&=&1+2b_8+2x'_4,\\
b_3x'_4&=&x_4+y_8.
\end{eqnarray}
Thus
\[
\begin{array}{rll}
\bb_3(b_3x'_4)&=&\bb_3y_8+\bb_3x_4\\
&=&b_8+x'_4+\bar{z}_{12}+b_8+x'_4,\\
(b_3\bb_3)x'_4&=&x'_4+x'_4b_8.
\end{array}
\]
So
\[
x'_4b_8=x'_4+2b_8+z_{12}.
\]
Since \[
\begin{array}{rll}
b_3(b_5\bb_5)&=&b_3+2b_3b_8+2b_3x'_4\\
&=&b_3+2b_3+2b_9+2x_4+2y_8+2x_4+2y_8,\\
b_5(b_3\bb_5)&=&\bb_3b_5+\bar{x}_4b_5+\bar{y}_8b_5\\
&=&b_3+x_4+y_8+\bar{x}_4b_5+\bar{y}_8b_5,
\end{array}
\]
we have that $\bar{x}_4b_5+\bar{y}_8b_5=2b_3+2b_9+3x_4+3y_8.$
Because $b_3\in Supp\{\bar{x}_4b_5\}$ and $b_3\in
Supp\{\bar{y}_8b_5\}$, there exactly two possibilities:
\[
\begin{array}{ll}

\bar{x}_4b_5=b_3+b_9+y_8,& \bar{y}_8b_5=b_3+b_9+3x_4+2y_8,\\
\bar{x}_4b_5=b_3+b_9+2x_4,& \bar{y}_8b_5=b_3+b_9+x_4+3y_8.
\end{array}
\]
\par If $\bar{x}_4b_5=b_3+b_9+y_8$ and $\bar{y}_8b_5=b_3+b_9+3x_4+2y_8$, then $(\bar{y}_8b_5,\ x_4)=3$, which implies that $(b_5\bar{x}_4,\ y_8)=3$, a contradiction.
\par If $\bar{x}_4b_5=b_3+b_9+2x_4$ and $\bar{y}_8b_5=b_3+b_9+x_4+3y_8.$ Then
$(\bar{x}_4b_5,\ x_4)=3$. So $(x^2_4,\ b_5)=2$. It follow that $(x^2_4,\ x^2_4)\geq 5$. It is impossible for either $x_4\bar{x}_4=1+b_8+r_7,\ r_7\in B$ or $x_4\bar{x}_4=1+b_8+r_4+r_3,\ r_3,\ r_4\in B$. Hence {\bf Subcase I} follows.\\
\par {\bf Subcase II.} No NITA such that $t_4\neq x_4$.\\
\par If $t_4\neq x_4$, by $b_4\bar{y}_8+b_9+t_4+r_3=x_4+y_8+b_3\bar{z}_{12}$, we have that
$ r_3+t_4+b_9$ is a part of $b_3\bar{z}_{12}$, $x_4+y_8$ is a part
of $b_4\bar{y}_8$ and
$b_3+x_4+y_8+b_9+y_{12}+r_3+t_4+b_9$ is a part of $b_8b_9$.\\
\par  If $b_9\in Supp\{y_{12}\}$, then there exists $h_3\in B$ such that
$y_{12}=h_3+b_9$. So $b_3c_7=h_3+2b_9$,\ $b_3\bar{h}_3$ contains
$c_7$, which is impossible for $L_2(B)=\emptyset$. Hence
$b_9\not\in Supp\{ y_{12}\}$.
\par Since $(b_9\bb_9,\ b_8)=3$, so $(b_8b_9,\ b_9)=3$, which implies that $(b_3\bar{z}_{12},\ b_9)=2$ and there exists  $x_{11}\in N^*B$ of degrees 11 such that
\[
\begin{array}{rll}
b_3\bar{z}_{12}&=&2b_9+t_4+r_3+x_{11},\\
b_4\bar{y}_8&=&x_4+y_8+b_9+x_{11}.
\end{array}
\]
\par At the following we'll discuss based on what $z_{12}$ is like. Since $(z_{12},\ z_{12})=2$, we have that $z_{12}=\alpha+\beta$,
where $(f(\alpha),\ f(\beta))=(3,\ 9)$, $(4,\ 8)$, $(5,\ 7)$ or $(6,\ 6)$.\\
\par {\bf Step 1.} There exists no NITA such that $z_{12}=\alpha_3+\beta_9$.\\
\par If $z_{12}=\alpha_3+\beta_9$. Then $b_3(\bar{\alpha}_3+\bar{\beta}_9)=r_3+t_4+b_9+x_{11}$, $b_3\bar{\alpha}_3=b_9$ and
\[
\begin{array}{rll}
\bb_3(b_3\bar{\alpha}_3)&=&c_7+b_8+\bar{\alpha}_3+\bar{\beta}_9,\\
(\bb_3b_3)\bar{\alpha}_3&=&\bar{\alpha}_3+\bar{\alpha}_3b_8.
\end{array}
\]
Thus
$$\alpha_3b_8=c_7+b_8+\beta_9,$$
which implies that $(\alpha_3b_8,\ \alpha_3b_8)=3$. Hence $(\alpha_3\bar{\alpha}_3,\ b^2_8)=3$.
It follows that $\alpha_3\bar{\alpha}_3=1+x'_4+\bar{x}'_4$ by the expression of $b^2_8$ (see (\ref{x3}), which contradicts Theorem 1.1.\\
\par {\bf Step 2.} There exists no NITA such that  $z_{12}=\alpha_4+\beta_8$.\\
\par If $z_{12}=\alpha_4+\beta_8$. Then $b_3(\bar{\alpha}_4+\bar{\beta}_8)=r_3+t_4+2b_9+x_{11}$. Since $b_9\in Supp\{b_3\bar{\alpha}_4\}$, we have that
$$b_3\bar{\alpha}_4=r_3+b_9,\ b_3\bar{\beta}_8=t_4+b_9+x_{11}.$$
\par Let $\bb_3r_3=\bar{\alpha}_4+e_5$. Then
\[
\begin{array}{rll}
\bb_3(b_3\bar{\alpha}_4)&=&\bb_3b_9+\bb_3r_3\\
&=&c_7+b_8+\bar{\alpha}_4+\bar{\beta}_8+\bar{\alpha}_4+e_5,\\
\bar{\alpha}_4(b_3\bb_3)&=&\bar{\alpha}_4+\bar{\alpha}_4b_8,
\end{array}
\]
which means that $\alpha_4b_8=c_7+b_8+\beta_8+e_5.$ Since
$\bb_3y_8=\alpha_4+\beta_8+\bar{x}'_4+b_8$, we have that
$b_3\alpha_4=y_8+\gamma_4,\ \gamma_4\in B.$
\par Since
\[
\begin{array}{rll}
\bb_3(b_3\alpha_4)&=&\bb_3y_8+b_3\gamma_4\\
&=&b_8+\alpha_4+\beta_8+\bar{x}'_4+b_3\gamma_4,\\
(\bb_3b_3)\alpha_4&=&\alpha_4+\alpha_4b_8\\
&=&2\alpha_4+c_7+b_8+\beta_8+e_5,
\end{array}
\]
we have that $\bar{x}'_4+b_3\gamma_4=\alpha_4+c_7+e_5$, which
implies that $\bar{x}'_4=\alpha_4$. But by (21) $b_3x'_4=j_8+x_4$,
$j_8\in B$ and
$b_3\bar{\alpha}_4=b_9+r_3$, a contradiction.\\
\par {\bf Step 3.} There exits no NITA such that $z_{12}=\alpha_5+\beta_7$.\\
\par If $z_{12}=\alpha_5+\beta_7$. Then
$$b_3(\bar{\alpha}_5+\bar{\beta}_7)=2b_9+r_3+t_4+2b_9+x_{11}.$$
We need a sum of degree 6 to make up $b_3\bar{\alpha}_5$ for
$b_9\in Supp\{b_3\bar{\alpha}_5\}$.
\par If $r_3\in Supp\{b_3\bar{\beta}_7\}$, then $\bb_3r_3$ contains an element of degree 2, a
contradiction. Hence $r_3\in Supp\{b_3\bar{\alpha}_5\}$. We may
set
$$b_3\bar{\alpha}_5=b_9+r_3+\delta_3,\ \delta_3\in Supp\{x_{11}\}.$$
Furthermore $\bb_3r_3=\bar{\alpha}_5+d_4,\ \mbox{some}\ d_4\in B.$
It follows that $r_3\bar{r}_3=1+b_8$. Now we have that
$(\alpha_5\bar{\alpha}_5,\ b_8)=2$ and
$$2b_3+x_4+y_8+b_9=b_3(b_3\bb_3)=r_3(b_3\bar{r}_3)=r_3\alpha_5+r_3\bar{d}_4,$$
 from which the following holds:
$$r_3\alpha_5+r_3\bar{d}_4=2b_3+x_4+y_8+b_9.$$
Comparing the two sides, the following equations hold:
\[
\begin{array}{rll}
r_3\bar{d}_4&=&b_3+b_9,\\
r_3\alpha_5&=&b_3+x_4+y_8.\\
\end{array}
\]
Therefore
\[
\begin{array}{rll}
\bar{b}_3(r_3\alpha_5)&=&b_3\bb_3+x_4\bb_3+y_8\bb_3\\
&=&1+b_8+b_8+x'_4+b_8+\alpha_5+\beta_7+\bar{x}'_4,\\
\alpha_5(\bb_3r_3)&=&\alpha_5\bar{\alpha}_5+d_4\alpha_5.
\end{array}
\]
\par Since $(\alpha_5\bar{\alpha}_5,\ b_8)=2$, we can only have that
\[
\begin{array}{rll}
\alpha_5\bar{\alpha}_5&=&1+2b_8+x'_4+\bar{x}'_4,\\
d_4\alpha_5&=&b_8+\alpha_5+\beta_7.
\end{array}
\]
Hence $d_4=x'_4$ or $\bar{x}'_4$. So $\bb_3r_3=\alpha_5+x'_4$ or
$\alpha_5+\bar{x}'_4$.
\par If $\bb_3r_3=\alpha_5+\bar{x}'_4$, then $r_3\in Supp\{b_3\bar{x}'_4\}=\{y_8,\ t_4\}$ by (\ref{x7}), a contradiction. Hence
 $\bb_3r_3=\alpha_5+x'_4$, which implies that $r_3\in Supp\{b_3x'_4\}=\{x_4,\ j_8\}$,  a contradiction.\\
\par {\bf Step 4.} There exists no NITA such that $z_{12}=\alpha_6+\beta_6$.\\
\par If  $z_{12}=\alpha_6+\beta_6$, then
$$b_3(\bar{\alpha}_6+\bar{\beta}_6)=2b_9+t_4+r_3+x_{11}.$$
Hence $\alpha_6\neq \beta_6$. Without loss of generality, let
$r_3\in Supp\{b_3\alpha_6\}$. We may set that
$\bb_3r_3=\bar{\alpha}_6+d_3, $ which implies that
$r_3\bar{r}_3=1+b_8$. Therefore we have the following equation:
\[
\begin{array}{rll}
\alpha_6r_3+\bar{d}_3r_3&=&(b_3\bar{r}_3)r_3\\
&=&(b_3\bb_3)b_3\\
&=&2b_3+b_9+x_4+y_8.
\end{array}
\]
It is easy to see that we cannot sum up $\bar{d}_3r_3$, a contradiction. The {\bf Proposition 1} follows.\\
\begin{proposition}\label{prop3.2}
No NITA such that $(f(x),\ f(y))=(5,\ 7)$.
\end{proposition}
\par {\bf Proof.} If $(f(x),\ f(y))=(5,\ 7)$, that is to say, $x_{12}=x_5+y_7$, and
\[
\begin{array}{rll}
\bb_3b_5&=&b_3+x_5+y_7,\\
b_3b_8&=&b_3+b_9+x_5+y_7.
\end{array}
\]
So we may set that
\[
\begin{array}{rll}
\bb_3x_5&=&b_8+e_7,\ e_7\in \mathit{N}^*B,\\
b_3x_5&=&b_5+e_{10},\ e_{10}\in \mathit{N}^*B,\\
\bb_3y_7&=&b_8+y_{13},\ y_{13}\in\mathit{N}^*B,\\
b_3y_7&=&b_5+y_{16},\ y_{16}\in \mathit{N}^*B.
\end{array}
\]
Hence
\[
\begin{array}{rll}
\bb_3(b_3b_8)&=&\bb_3b_3+\bb_3b_9+\bb_3x_5+\bb_3y_7\\
&=&1+b_8+c_7+b_8+\bar{z}_{12}+b_8+e_7+b_8+y_{13},\\
(b_3\bb_3)b_8&=&b_8+b^2_8,\\
(\bb_3b_5)b_3&=&b^2_3+b_3x_5+b_3y_7\\
&=&b_4+3b_5+e_{10}+y_{16},\\
(b_3\bb_3)b_5&=&b_5+b_5b_8,
\end{array}
\]
one has that
\[
\begin{array}{rll}
b^2_8&=&1+c_7+e_7+3b_8+\bar{z}_{12}+y_{13},\\
b_5b_8&=&b_4+2b_5+e_{10}+y_{16}.
\end{array}
\]
\par On the other hand,
\[
\begin{array}{rll}
(b_3\bb_3)^2&=&1+2b_8+1+c_7+e_7+3b_8+\bar{z}_{12}+y_{13},\\
b^2_3\bb^2_3&=&1+b_8+c_7+b_8+z_{12}+b_8+\bar{z}_{12}+1+2b_8+d_8,
\end{array}
\]
which implies that $z_{12}+d_8=e_7+y_{13}.$\\
\par {\bf Step 1.} There exists no NITA such that $e_7\not\in Supp\{z_{12}\}$.\\
\par Suppose that $e_7\not\in Supp\{z_{12}\}$. Since $e_7\not\in Supp\{d_8\}$, there exist
$\alpha_3\in Supp\{z_{12}\}$, $g_4\in Supp\{d_8\}$ or $\alpha_4\in
Supp\{z_{12}\}$, $g_3\in Supp\{d_8\}$ such that
$$e_7=\alpha_3+g_4\ \mbox{or}\ \alpha_4+g_3.$$
\par {\bf Substep 1.} There exists no NITA such that $e_7=\alpha_3+g_4$.\\
\par If $e_7=\alpha_3+g_4$, then there exist $z_9\in \mathit{N}^*B$ and $h_4\in B$ such that  $z_{12}=\alpha_3+z_9$ and $d_8=g_4+h_4$. Further
\[
\begin{array}{rll}
\bb_4b_5&=&\alpha_3+z_9+b_8,\\
\bb_3b_9&=&c_7+b_8+\bar{\alpha}_3+\bar{z}_9,\\
\bb_3x_5&=&b_8+\alpha_3+g_4,\\
\end{array}
\]
Hence
\[
\begin{array}{rll}
b_3\alpha_3&=&x_5+i_4,\ i_4\in B,\\
b_3\bar{\alpha}_3&=&b_9,
\end{array}
\]
which is impossible for $(b_3\bar{\alpha}_3,\ b_3\bar{\alpha}_3)=1$ and $(b_3\alpha_3,\ b_3\alpha_3)=2$.\\
\par {\bf Substep 2.} No NITA such that $e_7=\alpha_4+g_3$.\\
\par If $e_7=\alpha_4+g_3$, then there exist $z_8\in Supp\{\mathit{N}^*B\}$ and
 $h_5\in B$ such that $\bar{z}_{12}=\alpha_4+z_8$, $d_8=g_3+h_5$. Hence
\[
\begin{array}{rll}
\bb_3x_5&=&g_3+\alpha_4+b_8,\\
b_5\bb_5&=&1+2b_8+g_3+h_5,\\
\bb_3y_7&=&b_8+z_8+h_5.
\end{array}
\]
So there exists $l_4\in B$ and $n_8\in\mathit{N}^*B$ such that
$b_3g_3=x_5+l_4$ and $b_3h_5=y_7+n_8.$ Hence $g_3\bar{g}_3=1+b_8$
and $\bb_3l_4=g_3+m_9$. Here $m_9\in B$, otherwise $\bar{l}_4l_4$
will contains two $b_8$, which is a contradiction.
\par Since \[
\begin{array}{rll}
b_3(b_3\bb_3)&=&2b_3+b_9+x_5+y_7,\\
(b_3g_3)\bar{g}_3&=&\bar{g}_3x_5+\bar{g}_3l_4,
\end{array}
\]
it follows that $\bar{g}_3l_4=b_3+b_9$ and
\begin{eqnarray}\label{yy}\bar{g}_3x_5=b_3+x_5+y_7=\bb_3b_5.
\end{eqnarray}
\par We assert that $Supp\{n_8\}\cap\{b_3,\ x_5,\ y_7,\ l_4\}=\emptyset.$ Obviously $y_7\not\in Supp\{n_8\}$. By the expressions of $b_3h_5$ and $b_3\bb_3$, we have that $b_3\not\in Supp\{n_8\}$. If $x_5\in Supp\{n_8\}$, then $h_5\in Supp\{\bb_3x_5\}$, a contradiction. If $l_4\in Supp\{n_8\}$, then $h_5\in Supp\{\bb_3l_4\}$, a contradiction.
Since
\[
\begin{array}{rll}
b_3(b_5\bb_5)&=&b_3+2b_3b_8+b_3g_3+b_3h_5\\
&=&3b_3+2b_9+2x_5+2y_7+x_5+l_4+y_7+n_8,\\
(b_3\bb_5)b_5&=&b_5(\bb_3+\bar{x}_5+\bar{y}_7)\\
&=&b_3+x_5+y_7+\bar{x}_5b_5+\bar{y}_7b_5,
\end{array}
\]
It holds that
\begin{eqnarray}\label{yyy}
b_5\bar{x}_5+b_5\bar{y}_7=2b_3+2b_9+2x_5+2y_7+l_4+n_8.
\end{eqnarray}

\par Now we assert that $(b_5\bar{y}_7,\ b_9)=1$. In fact, it follows from that $b_5\bb_9=b_9+b_3\alpha_4+b_3z_8$ by {\bf Lemma \ref{be}}  and
$(b_3\alpha_4,\ y_7)=0$, $(b_3z_8,\ y_7)=1$. Therefore $(b_5\bar{x}_5,\ b_3)=(b_5\bar{x}_5,\ b_9)=1$ and $(b_5\bar{y}_7,\ b_3)=(b_5\bar{y}_7,\ b_9)=1$.
\par By (\ref{yy}), $(\bar{g}_3x_5,\ \bb_3b_5)=3$. So $(b_3\bar{g}_3,\ b_5\bar{x}_5)=3$. By (\ref{yyy}), $b_5\bar{x}_5$ cannot contains
three multiple of a constituent of degree 3 or two multiple of
constituent of degree 6. Then the common constituents of
$b_3\bar{g}_3$ and $b_5\bar{x}_5$ are degree 5 or 4 for
$b_3\bar{g}_3$ can only have constituents of degree 5 and 4 at
this moment. From $(b_3\bar{g}_3,\ b_5\bar{x}_5)=3$ and
(\ref{yyy}), we have that $b_5\bar{x}_5$ must contains $2x_5$ and
some element $v_4$ of degree 4 from $l_4+n_8$. It leads to that
$$b_5\bar{x}_5=b_3+b_9+2x_5+v_4,$$
a contradiction for left side having degree 26.\\

\par {\bf Step 2.} There exists no NITA such that $e_7\in Supp\{z_{12}\}$.\\
\par If  $e_7\in Supp\{z_{12}\}$, then there exists $w_5\in B$ such that  $z_{12}=e_7+w_5$. So $y_{13}=d_8+w_5$ and
\begin{eqnarray}
\label{x8}\bb_3x_5=b_8+e_7,\\
\label{x9}\bb_3y_7&=&b_8+d_8+w_5,\\
\label{x10}\bb_4b_5&=&w_5+e_7+b_8,\\
\label{x11}\bb_3b_9&=&\bar{w}_5+c_7+\bar{e}_7+b_8.
\end{eqnarray}
\par Since $(\bb_4b_5,\ \bb_4b_5)=3$, we have that $e_7$, $w_5\in B$. By {Lemma \ref{be}}
\begin{eqnarray}\label{x17}b_5\bb_9=b_9+b_3e_7+b_3w_5.\end{eqnarray}
\par Since $y_7\in Supp\{b_3w_5\}$, $b_9\in Supp\{b_3\bar{w}_5\}$,\ and $x_5\in Supp\{b_3e_7\}$, we may set that
\begin{eqnarray}
\label{x12}b_3w_5&=&y_7+v_8,\  v_8\in\mathit{N}^*B,\\
\label{x13}b_3\bar{w}_5&=&b_9+d_6,\  d_8\in\mathit{N}^*B,\\
\label{x14}b_3c_7&=&b_9+f_{12},f_{12}\in N^*B\\
\label{x15}b_3e_7&=&x_5+h_{16},h_{16}\in N^*B\\
\label{x16}b_3\bar{e}_7&=&b_9+g_{12}, g_{12} \in N^*B.
\end{eqnarray}

\par {\bf Case 1.}There exists no NITA such that either $v_8$ or $d_6$ reducible.\\

\par It is easy to see that one of $v_8$ and $d_6$ is reducible will imply another one is reducible. Let $d_6=d_3+g_3$. Then
$$b_3\bar{w}_5=d_3+g_3+b_9.$$
It is easy to show that $d_3\neq g_3$.  Hence we may set that
\[
\begin{array}{rll}
b_3\bar{d}_3&=&w_5+\delta_4,\\
b_3\bar{g}_3&=&w_5+\psi_4,\\
\end{array}
\]
which implies that $d_3\bar{d}_3=g_3\bar{g}_3=1+b_8$. If
$\delta_4=\psi_4$, then $b_3\bar{d}_3=b_3\bar{g}_3$,
$(b_3\bar{d}_3,\ b_3\bar{g}_3)=2$, which implies that
$d_3\bar{g}_3=1+b_8$. Thus $g_3=d_3$ and $b_3\bar{w}_5=b_9+2g_3$,
so that $(b_3\bar{g}_3,\ w_5)=2$, a contradiction. Hence
$\delta_4\not=\psi_4$. Let's calculate
$b^2_3\bb_3=b_3d_3\bar{d}_3=b_3g_3\bar{g}_3$, we have that
\[
\begin{array}{rll}
(b_3\bb_3)b_3&=&2b_3+b_9+x_5+y_7,\\
(b_3\bar{d}_3)d_3&=&w_5d_3+\delta_4d_3,\\
(b_3\bar{g}_3)g_3&=&w_5g_3+\psi_4g_3.
\end{array}
\]
It is easy to see that
\[
\begin{array}{rll}
\delta_4d_3&=&\psi_4g_3=b_3+b_9,\\
w_5d_3&=&w_5g_3=b_3+x_5+y_7.
\end{array}
\]
\par On the other hand, we have that
\[
\begin{array}{rll}
w_5\bar{w}_5&=&1+2b_8+u_8,\\
\bb_3(g_3\psi_4)&=& \bb_3(d_3\delta_4)=\bb_3b_3+\bb_3b_9\\
&=&1+b_8+c_7+\bar{e}_7+\bar{w}_5+b_8,\\
(\bb_3g_3)\psi_4&=&\bar{w}_5\psi_4+\bar{\psi}_4\psi_4,\\
 (\bb_3d_3)\delta_4&=&\bar{w}_5\delta_4+\bar{\delta}_4\delta_4,
\end{array}
\]
so that
\[
\begin{array}{rll}
\bar{w}_5\psi_4&=&\bar{w}_5\delta_4=\bar{e}_7+\bar{w}_5+b_8,\\
\psi_4\bar{\psi}_4&=&\delta_4\bar{\delta}_4=1+b_8+c_7.
\end{array}
\]
Hence $\psi_4$, $\delta_4\in Supp\{w_5\bar{w}_5\}$. Therefore
$$w_5\bar{w}_5=1+2b_8+\delta_4+\psi_4.$$
But $(g_3\bar{d}_3,\ w_5\bar{w}_5)=(g_3w_5,\ d_3w_5)=3$. We have
that
  $1\in Supp\{g_3\bar{d}_3\}$, which implies that $g_3=d_3$, a contradiction. \\
\par {\bf Case 2.} There exists no NITA such that $v_8$ and $d_6$ irreducible. \\
\par Suppose $v_8$ and $d_6$ irreducible. We assert that $f_{12}$ is irreducible.
Otherwise, let $f_{12}=r+s$. Suppose that $c_7\in Supp\{\bb_3r\}$,
which implies that $f(r)\geq 4$.
\par By {\bf Lemma \ref{be}}, one has that
\[
\begin{array}{rll}
b_8c_7&=&e_7+w_5+b_8+b_3\bar{f}_{12},\\
\bb_9b_9&=&1+d_8+3b_8+e_7+w_5+b_3\bar{f}_{12}.
\end{array}
\]

\par {\bf Subcase 1.} $w_5$ is not real.\\
\par If $w_5$ is real, then $y_7+v_8=b_3w_5=b_3\bar{w}_5=b_9+d_6$ by (\ref{x12}) and (\ref{x13}), a contradiction.
\par {\bf Subcase 2.} $f_{12}\not\in B$.\\
\par By the expression of $b_8c_7$, we have that $\bar{w}_5\in Supp\{b_3\bar{f}_{12}\}$. If $f_{12}$ is irreducible, then $f_{12}\in Supp\{b_3w_5\}=\{y_7,\ v_8\}$ by (\ref{x12}), a contradiction. So $f_{12}$ is reducible.\\
\par {\bf Subcase 3.} $e_7$ is not real.
\par If $e_7$ is real, then there exists an element $u_7\in \mathit{N}^*B$ by (\ref{x8}) and (\ref{x11}), such that
$$b_3e_7=x_5+b_9+u_7.$$
Hence $(e_7^2,\ b_8)\geq 2$. So $(e_7b_8,\ e_7)\geq 2$, By the
expressions of $b_8c_7$ and $b^2_8$, we have that $(e_7b_8,\
c_7)=1$, $(e_7b_8,\ b_8)\geq 2$. Thus there exists an element
$x_{19}\in\mathit{N}^*B$ of degree 19 such that
$$e_7b_8=2e_7+2b_8+c_7+x_{19}.$$
\par On ther other hand, from (\ref{x8}) and (\ref{x11}), we have that
\[
\begin{array}{rll}
\bb_3(b_3e_7)&=&\bb_3x_5+\bb_3b_9+\bb_3u_7\\
&=&b_8+e_7+c_7+e_7+\bar{w}_5+b_8+\bb_3u_7,\\
(\bb_3b_3)e_7&=& e_7+b_8e_7\\
&=&3e_7+2b_8+c_7+x_{19}.
\end{array}
\]
Then $\bar{w}_5+\bb_3u_7=e_7+x_{19}$. Since right side is real, we
have that $w_5\in Supp\{\bb_3u_7\}$ by {\bf Substep 1}, which
implies that there exists a constituent $t$ of $u_7$ such that
$w_5\in Supp\{\bb_3t\}$. Then $t\in Supp\{b_3w_5\}$, so $t=y_7$ or
$v_8$. It follows that $u_7=y_7$. By (\ref{x9}) implies that
$$w_5+\bar{w}_5+b_8+d_8=\bar{w}_5+\bb_3y_7=e_7+x_{19}, $$
 a contradiction.\\
\par {\bf Subcase 4.} {\bf Case 2} follows.\\
\par By {\bf Subcases 1, 2, 3,} we have that $e_7$ and $w_5$ are not real, and $f_{12}$ is reducible. Notice $b_8c_7$ is real, by the expression of $b_8c_7$, we have that $b_3\bar{f}_{12}$ contains $\bar{w}_5$ and $\bar{e}_7$. Hence there exists a constituent $r$ of $f_{12}$ such that
$b_3\bar{r}$ contains $\bar{w}_5$, which implies that $r\in
Supp\{b_3w_5\}$. Therefore $r=y_7$ or $v_8$.
\par If $r=y_7$, then  $b_3\bar{y}_7=b_8+d_8+\bar{w}_5$ by (\ref{x9}). But (\ref{x14}) means that $b_3\bar{y}_7$ contains $c_7$, a contradiction.
\par If $r=v_8$, then there exists $s_4$ such that $f_{12}=v_8+s_4$. Since $c_7\in Supp\{b_3\bar{s}_4\}$ by (\ref{x14}), we have that $\bar{e}_7\not\in Supp\{b_3\bar{s}_4\}$. So $\bar{e}_7\in Supp\{b_3\bar{v}_8\}$. Hence there exists a basis element $a_5$ of degree 5 such that
\begin{eqnarray}
\label{x23} b_3\bar{v}_8&=&\bar{w}_5+\bar{e}_7+c_7+a_5,\\
\label{x21}b_3\bar{s}_4&=&c_7+c_5,\\
\label{x24} b_3c_7&=&v_8+s_4+b_9.
\end{eqnarray}

\par Checking $b_3(b_5\bb_5)=b_5(b_3\bb_5)$, one implies that
\begin{eqnarray}\label{x20}
b_5\bar{x}_5+b_5\bar{y}_7=2b_3+2b_9+x_5+y_7+b_3d_8.\end{eqnarray}
\par We assert that $d_8\in B$. In fact, if $d_8\not\in B$, then $d_8=2x_4$ or $x_4+y_4$.
\par If $d_8=2x_4$, then $\bb_3y_7=b_8+2x_4+w_5$ by (\ref{x9}). Hence $(b_3x_4,\ y_7)=2$. It is impossible.
\par If $d_8=x_4+y_4$, then $\bb_3y_7=b_8+x_4+y_4+w_5$ by (\ref{x9}). Let $b_3x_4=y_7+p_5$ and $b_3y_4=y_7+q_5$. By (\ref{x20}), it holds that
$$b_5\bar{x}_5+b_5\bar{y}_7=2b_3+2b_9+x_5+3y_7+p_5+q_5.$$
But $(b_5\bar{x}_5,\ b_3)=(b_5\bar{x}_5,\ b_9)=1$. So it is
necessary to find  elements from the left side of above equation
such that their sum having degree 13 to make up $b_5\bar{x}_5$. It
is impossible. Therefore $d_8\in B$.
\par Since
\[
\begin{array}{rll}
\bb_3(b_3w_5)&=&\bb_3y_7+\bb_3v_8\\
&=&b_8+d_8+w_5+c_7+w_5+e_7+\bar{a}_5,\\
b_3(\bb_3w_5)&=&b_3\bb_9+b_3\bar{d}_6\\
&=&w_5+c_7+e_7+b_8+b_3\bar{d}_6,
\end{array}
\]
then
\begin{eqnarray}\label{x22}b_3\bar{d}_6=d_8+w_5+\bar{a}_5.\end{eqnarray}
\par Now we assert that $c_5$ and $e_5$ are nonreal.  Otherwise, $d_6,\ v_8,\ \in Supp\{b_3a_5\}$ by (\ref{x23}) and (\ref{x22}).
So $1\in Spp\{b_3a_5\}$, a contradiction.Hence $c_5$ and $a_5$ are
nonreal, and $c_5=\bar{a}_5$ by the expression of $b^2_8$. Hence
$b_3\bar{d}_6=d_8+w_5+c_5$. Associated with(\ref{x21}), we have
that $(b_3\bc_5,\ d_6)=(b_3\bc_5,\ s_4)=1$, which implies that
$(b_3\bc_5,\ b_3\bc_5)=3$. Furthermore $(b_3c_5,\ b_3c_5)=3$. Let
$$b_3c_5=v_8+\delta_3+\lambda_4.$$
Thus we may set that $\bb_3\delta_3=c_5+\Delta_4$. Then
$\bar{\delta}_3\delta_3=1+b_8$. And
\[
\begin{array}{rll}
b_3(\delta_3\bar{\delta}_3)&=&b_3+b_3b_8\\
&=&2b_3+b_9+x_5+y_7,\\
(b_3\bar{\delta}_3)\delta_3&=&\bc_5\delta_3+\bar{\Delta}_4\delta_3.
\end{array}
\]
We  have that
\begin{eqnarray}\label{x23}\bar{\Delta}_4\delta_3=b_3+b_9\ \mbox{and}\  \bc_5\delta_3=b_3+x_5+y_7.\end{eqnarray}
Further
\[
\begin{array}{rll}
\bb_3(\bar{\Delta}_4\delta_3)&=&\bb_3b_3+\bb_3b_9\\
&=&1+2b_8+c_7+\bar{w}_5+\bar{e}_7,\\
(\bb_3\delta_3)\bar{\Delta}_4&=&c_5\bar{\Delta}_4+\Delta_4\bar{\Delta}_4.
\end{array}
\]
It follow that
$$\Delta_4\bar{\Delta}_4=1+b_8+c_7\ \mbox{and}\ c_5\bar{\Delta}_4=b_8+\bar{w}_5+\bar{e}_7=b_4\bb_5.$$
Hence
$(b_3c_5)\bar{\Delta}_4=b_3(c_5\bar{\Delta}_4)=b_3(b_4\bb_5)=b_4(b_3\bb_5)$.
But
\[
\begin{array}{rll}
b_4(b_3\bb_5)&=&b_4\bb_3+b_4\bar{x}_5+b_4\bar{y}_7\\
&=&b_3+b_9+b_4\bar{x}_5+b_4\bar{y}_7,\\
(b_3c_5)\bar{\Delta}_4=v_8\bar{\Delta}_4+\delta_3\bar{\Delta}_4+\lambda_4\bar{\Delta}_4.
\end{array}
\]
So $b_3$ is a constituent of one of $v_8\bar{\Delta}_4$,
$\delta_3\bar{\Delta}_4$ and $\lambda_4\bar{\Delta}_4$. Then one
of $v_8$, $\delta_3$ and $\lambda_4$ is a constituent of
$b_3\Delta_4$. On the other hand $\delta_3\in Supp\{b_3\Delta_4\}$
by (\ref{x23}). But $b_3\Delta_4$ contains at most two
constituents, a contradiction. {\bf Subcase 4} follows, which
concludes {\bf Step 2.} Consequently the proposition follows.

\par \begin{proposition}\label{prop3.3}No NITA such that $(f(x),\ f(y))=(6,\ 6).$
\end{proposition}
\par {\bf Proof.} If $x_{12}=x_6+y_6$, then
\[
\begin{array}{rll}
\bb_3b_5&=&b_3+x_6+y_6,\\
b_3b_8&=&b_3+b_9+x_6+y_6.
\end{array}
\]
Assume that
\[
\begin{array}{rlll}
\bb_3x_6&=&b_8+h_{10},\ \mbox{for some}\ h_{10}\in\mathit{N}^*B,\\
\bb_3y_6&=&b_8+i_{10},\ \mbox{for some}\ i_{10}\in\mathit{N}^*B,\\
b_3x_6&=&b_5+g_{13},\ \mbox{for some}\ g_{13}\in\mathit{N}^*B,\\
b_3y_6&=&b_5+f_{13},\ \mbox{for some}\ f_{13}\in\mathit{N}^*B.
\end{array}
\]
Therefore
$$b^2_8=1+c_7+3b_8+h_{10}+i_{10}+\bar{z}_{12}.$$
\par By the above equations and {\bf Lemma \ref{be}},  one has that
\[
\begin{array}{rll}
(\bb_3b_3)^2&=&2+c_7+5b_8+h_{10}+i_{10}+\bar{z}_{12},\\
\bb^2_3b^2_3&=&\bb_4b_4+\bb_4b_5+\bb_5b_4+\bb_5b_5\\
&=&2+5b_8+c_7+\bar{z}_{12}+z_{12}+d_8.
\end{array}
\]
Hence $$d_8+z_{12}=h_{10}+i_{10}.$$
\par We assert that $d_8$ cannot be totally contained in one of $h_{10}$ and $i_{10}$.
Otherwise, without loss of generality, assume $d_8$ is contained
in $h_{10}$. Then $h_{10}$ contains an element of degree 2, a
contradiction. Thus $d_8$ must be reducible and cannot be
contained in one of $h_{10}$ and $i_{10}$. For $d_8$, there are
four possibilities:
\begin{center}
$d_8=c_3+d_5$, $c_4+d_4$, $c_4+\bar{c}_4$ and $2c_4.$
\end{center}

\par {\bf Step 1.} No NITA such that $d_8=c_3+d_5$.\\
\par If $d_8=c_3+d_5$, then $c_3+d_5+z_{12}=h_{10}+i_{10}$. Without loss of generality, let
$$h_{10}=c_3+h_7,\ i_{10}=d_5+i_5.$$
Then $z_{12}=h_7+i_5$ and

\begin{eqnarray}
\label{y1}\bb_3x_6&=&c_3+h_7+b_8,\\
\label{y2}\bb_3y_6&=&d_5+i_5+b_8,\\
\label{y3}\label{3.1}\bb_4b_5&=&b_8+i_5+h_7,\\
\label{y4}\bb_3b_9&=&c_7+b_8+\bar{i}_5+\bar{h}_7,\\
\label{y5}b^2_8&=&1+c_7+3b_8+c_3+h_7+d_5+i_5+\bar{h}_7+\bar{i}_5.
\end{eqnarray}

Hence $c_3$ is real. By (\ref{3.1}) and $(\bb_4b_5,\ \bb_4b_5)=3$,
we have that  $h_7\in B$. Set $b_3c_3=x_6+u_3,\ u_3\in B$ by {\bf
Lemma \ref{t} }.
 Thus $(b_3c_3,\ b_3c_3)=2$ so that
$$c^2_3=1+b_8.$$
Moreover
\[
\begin{array}{rll}
(b_3c_3)c_3&=&x_6c_3+u_3c_3,\\
b_3c^2_3&=&b^2_3\bb_3=b_3+b_3b_8\\
&=&2b_3+b_9+x_6+y_6.
\end{array}
\]
Hence
\[
\begin{array}{rll}
c_3u_3&=&b_3+x_6\ \mbox{or}\ b_3+y_6,\\
c_3x_6&=&b_3+b_9+y_6\ \mbox{or}\ b_3+b_9+x_6.
\end{array}
\]
\par If $c_3u_3=b_3+x_6,\ c_3x_6=b_3+b_9+y_6$, then
$c_3\in Supp\{\bar{u}_3x_6\}$. But
\[
\begin{array}{rll}
(b_3\bb_3)(c^2_3)&=&2+5b_8+c_7+c_3+d_5+h_7+\bar{h}_7+i_5+\bar{i}_5,\\
(b_3c_3)(\bb_3c_3)&=&x_6\bar{x}_6+u_3\bar{x}_6+\bar{u}_3x_6+1+b_8.
\end{array}
\]
Notice that $ u_3\bar{x}_6+\bar{u}_3x_6$ contains two $c_3$, we
come to a contradiction. Therefore $c_3u_3=b_3+y_6,\
c_3x_6=b_3+b_9+x_6$. Let $\bb_3u_3=c_3+w_6$. Then
$u_3\bar{u}_3=1+b_8$ by {\bf Lemma \ref{t}} and
\[
\begin{array}{rll}
(c_3b_3)\bb_3&=&\bb_3x_6+\bb_3u_3\\
&=&b_8+c_3+h_7+c_3+w_6
\end{array}
\]
\par Since  $c_3(b_3\bb_3)$ is real, it follows that $h_7$ and $w_6$ are real.
Since
\[
\begin{array}{rll}
(\bb_3u_3)c_3&=&c^2_3+c_3w_6\\
&=&1+b_8+c_3w_6,\\
(\bb_3c_3)u_3&=&\bar{x}_6u_3+\bar{u}_3u_3\\
&=&1+b_8+\bar{x}_6u_3,\\
\bb_3(c_3u_3)&=&\bb_3b_3+\bb_3y_6\\
&=&1+b_8+b_8+d_5+i_5,
\end{array}
\]
we have that
$$c_3w_6=\bar{x}_6u_3=d_5+i_5+b_8,$$
Since $c_3$, $w_6$ are real,  $i_5$ is real too by above equation.
By (\ref{y1}), (\ref{y2}) and (\ref{y4}), one has that
\[
\begin{array}{rll}
b_3i_5&=&y_6+b_9,\\
b_3h_7&=&b_9+x_6+l_6,\ \mbox{for some}\ l_6\in \mathit{N}^*B.
\end{array}
\]
Considering $c_3^2u_3\bar{u}_3$, we have that
\[
\begin{array}{rll}
(c_3u_3)(c_3\bar{u}_3)&=&b_3\bb_3+\bb_3y_6+b_3\bar{y}_6+y_6\bar{y}_6\\
&=&1+b_8+2b_8+2d_5+2i_5+y_6\bar{y}_6,\\
c^2_3(u_3\bar{u}_3)&=&2+5b_8+c_7+c_3+2i_5+2h_7+d_5.
\end{array}
\]
Comparing the number of elements of degree 5, we come to a
contradiction.
{\bf Step 1} follows.\\
\par {\bf Step 2.} There exists no NITA such that $d_8=2c_4$.\\
\par If $d_8=2c_4$, then $2c_4+z_{12}=h_{10}+i_{10}$. There exist $h_6$ and $i_6$ such that $z_{12}=h_6+i_6$.
Since $(\bb_4b_5,\bb_4b_5)\le 4$, we obtain that  $h_6,\ i_6\in
B$. Furthermore
$$\bb_3x_6=b_8+c_4+h_6,\ \bb_3y_6=b_8+c_4+i_6.$$
Thus $b_3c_4=x_6+y_6$. Then we can set
$$c_4^2=1+b_8+g_7,\ g_7\in\mathit{N}^*B.$$
\par Since  $b_5\bb_5=1+2b_8+2c_4$,  $(c_4b_5,\ b_5)=2$. Further $(c_4b_5,\ c_4b_5)\geq 5$. Hence
$$c_4^2=1+b_8+c_4+e_3$$
and $(c_4b_5,\ c_4b_5)=5$. There exists $t_{10}\in B$ such that
$$c_4b_5=2b_5+t_{10}.$$
\par Since
\[
\begin{array}{rll}
b_3(b_5\bb_5)&=&b_3+2b_3b_8+2b_3c_4\\
&=&3b_3+2b_9+4x_6+4y_6,\\
b_5(b_3\bb_5)&=&b_5(\bb_3+\bar{x}_6+\bar{y}_6)\\
&=&b_3+x_6+y_6+\bar{x}_6b_5+\bar{y}_6b_5,
\end{array}
\]
we have that $b_5\bar{x}_6+b_5\bar{y}_6=2b_3+2b_9+3x_6+3y_6.$ On
the other hand
\[
\begin{array}{rll}
\bb_3(c_4b_5)=2\bb_3b_5+\bb_3t_{10}\\
&=&2b_3+2x_6+2y_6+\bb_3t_{10},\\
(\bb_3c_4)b_5&=&\bar{x}_6b_5+\bar{y}_6b_5\\
&=&2b_3+2b_9+3x_6+3y_6.
\end{array}
\]
Hence $\bb_3t_{10}=2b_9+x_6+y_6$, from which it follows that $(b_3b_9,\ t_{10})=2$. Therefore $(b_3b_9,\ b_3b_9)\geq 5$. But $(\bb_3b_9,\ \bb_3b_9)=4$ by {\bf Lemma \ref{be}} and $h_6$, $i_6\in B$,  a contradiction.\\
\par {\bf Step 3.} There exists no NITA such that $d_8=c_4+d_4$.\\
\par If $d_8=c_4+d_4$, then there exist $h_6$, $i_6\in\mathit{N}^*B$ such that
$h_{10}=c_4+h_6$, $i_{10}=d_4+i_6$,
 $z_{12}=h_6+i_6$ and

\begin{eqnarray}
\label{y6}\bb_3x_6&=&c_4+h_6+b_8,\\
\label{y7}\bb_3y_6&=&d_4+i_6+b_8,\\
\label{3.2}\bb_4b_5&=&b_8+i_6+h_6,\\
\label{y9}\bb_3b_9&=&c_7+b_8+\bar{i}_6+\bar{h}_6,\\
\label{y10}b^2_8&=&1+c_7+3b_8+c_4+h_6+d_4+i_6+\bar{h}_6+\bar{i}_6.
\end{eqnarray}
\par By (\ref{3.2}) and $(\bb_4b_5,\ \bb_4b_5)=3$, we have that
$i_6,\ h_6\in B$. Now we may set that
\begin{eqnarray}
\label{y11}b_3c_4&=&x_6+z_6,\ z_6\in B,\\
\label{y12}b_3d_4&=&y_6+p_6,\ p_6\in B.
\end{eqnarray}

\par {\bf Substep 1.} $h_6$ and $i_6$ are non-real.\\

\par If $h_6$ is real, then there exist $e_3$ and $f_3$ such that
\[
\begin{array}{rll}
\bb_3h_6&=&\bb_9+\bar{x}_6+e_3,\\
b_3e_3&=&f_3+h_6.
\end{array}
\]
 Thus $e_3\bar{e}_3=1+b_8$ bu {\bf Lemma \ref{t}}.
Since
\[
\begin{array}{rll}
b_3(e_3\bar{e}_3)&=&2b_3+b_9+x_6+y_6,\\
(b_3e_3)\bar{e}_3&=&h_6\bar{e}_3+f_3\bar{e}_3.
\end{array}
\]
Hence
\[
\begin{array}{rll}
f_3\bar{e}_3&=&b_3+x_6\ \mbox{or}\ b_3+y_6,\\
h_6\bar{e}_3&=&b_3+b_9+y_6\ \mbox{or}\ b_3+b_9+x_6.
\end{array}
\]
Therefore
\[
\begin{array}{rll}
\bb_3(f_3\bar{e}_3)&=&\bb_3b_3+\bb_3x_6\ \mbox{or}\ \bb_3b_3+\bb_3y_6\\
&=&1+2b_8+c_4+h_6\  \mbox{or}\ 1+2b_8+d_4+i_6,\\
f_3(\bb_3\bar{e}_3)&=&f_3(h_6+\bar{f}_3)\\
&=&1+b_8+f_3h_6.
\end{array}
\]
Thus
 $$f_3h_6=b_8+c_4+h_6\  \mbox{or}\ b_8+d_4+i_6.$$
\par On the other hand, one has the following equations:
 \[
\begin{array}{rll}
(b_3\bb_3)h_6&=&h_6+h_6b_8,\\
(b_3h_6)\bb_3&=&(\bar{e}_3+x_6+b_9)\bb_3\\
&=&\bb_3b_9+\bb_3x_6+\bb_3\bar{e}_3\\
&=&c_7+2b_8+3h_6+\bar{i}_6+c_4+\bar{f}_3.
\end{array}
\]
\ Then
\begin{eqnarray}
\label{3.4}h_6b_8=c_7+c_4+2b_8+2h_6+\bar{i}_6+\bar{f}_3,
\end{eqnarray}

which implies that $i_6$ and $f_3$ are real.

\par Now let us calculate $(b_3\bb_3)(e_3\bar{e}_3)$:
\[
\begin{array}{rll}
(b_3\bb_3)(e_3\bar{e}_3)&=&2+5b_8+c_7+c_4+d_4+2h_6+2i_6,\\
(b_3e_3)(\bb_3\bar{e}_3)&=&(h_6+f_3)^2\\
&=&h^2_6+2f_3h_6+f^2_3\\
&=&1+h^2_6+3b_8+2c_4+2h_6\ \mbox{or}\ 1+h^2_6+3b_8+2d_4+2i_6.
\end{array}
\]
Comparing the number of elements of degree 4, we have that
$c_4=d_4$, a contradiction.
\par  Symmetricly, we can prove that $i_6$ is non-real.\\
\par {\bf Substep 2.} There exists no NITA such that $d_8=c_4+d_4$ and $h_6$, $i_6$ are non-real.\\
\par Since $(b_3\bb_3)c_4=c_4+c_4b_8$ and $\bb_3(b_3c_4)=\bb_3x_6+\bb_3z_6=c_4+h_6+b_8+\bb_3z_6$, we have that
$$c_4b_8=b_8+h_6+\bb_3z_6, $$
Because $h_6$ is non-real and $c_4b_8$ is real, $\bar{h}_6\in
Supp\{\bb_3z_6\}$. By (\ref{y11}), there exists $r_8\in
\mathit{N}^*B$ such that
$$\bb_3z_6=c_4+\bar{h}_6+r_8.$$
By (\ref{y9}), there exists $\alpha_3$ such that
$$b_3\bar{h}_6=b_9+z_6+\alpha_3.$$
Thus there exists $\gamma_3$ such that
$b_3\bar{\alpha}_3=h_6+\gamma_3$. Hence
$\alpha_3\bar{\alpha}_3=1+b_8$ by {\bf Lemma \ref{t}}.
\par Since
\[
\begin{array}{rll}
b_3(\alpha_3\bar{\alpha}_3)&=&2b_3+b_9+x_6+y_6,\\
(b_3\alpha_3)\bar{\alpha}_3&=&h_6\bar{\alpha}_3+\gamma_3\bar{\alpha}_3,
\end{array}
\]
we have that
\[
\begin{array}{rll}
\gamma_3\bar{\alpha}_3&=&b_3+x_6\ \mbox{or}\ b_3+y_6,\\
h_6\bar{\alpha}_3&=&b_3+b_9+y_6\ \mbox{or}\ b_3+b_9+x_6.
\end{array}
\]
Hence
\[
\begin{array}{rll}
\bb_3(\gamma_3\bar{\alpha}_3)&=&\bb_3b_3+\bb_3x_6\ \mbox{or}\ \bb_3b_3+\bb_3y_6\\
&=&1+2b_8+c_4+h_6\ \mbox{or}\ 1+2b_8+d_4+i_6,\\
\gamma_3(\bb_3\bar{\alpha}_3)&=&\gamma_3\bar{h}_6+\gamma_3\bar{\gamma}_3\\
&=&1+b_8+\gamma_3\bar{h}_6.
\end{array}
\]
Therefore $$\gamma_3\bar{h}_6=b_8+c_4+h_6\ \mbox{or}\
b_8+d_4+i_6.$$
\par Since
\[
\begin{array}{rll}
(b_3\bb_3)(\alpha_3\bar{\alpha}_3)&=&2+5b_8+c_7+c_4+d_4+2h_6+2i_6,\\
(b_3\alpha_3)(\bb_3\bar{\alpha}_3)&=&(h_6+\gamma_3)(\bar{h}_6+\bar{\gamma}_3)\\
&=&h_6\bar{h}_6+\gamma_3\bar{h}_6+\bar{\gamma}_3h_6+\gamma_3\bar{\gamma}_3\\
&=&1+3b_8+2c_4+h_6+\bar{h}_6+h_6\bar{h}_6\ \mbox{or}\
1+3b_8+2d_4+i_6+\bar{i}_6+h_6\bar{h}_6.
\end{array}
\]
\par Comparing the number of degree 4, we come to $c_4=d_4$, a contradiction. So {\bf Substep 2} follows.\\
\par {\bf Step 4.} There exists no NITA such that $d_8=c_4+\bar{c}_4$.\\
\par If $d_8=c_4+\bar{c}_4$, then $b_5\bb_5=1+2b_8+c_4+\bar{c}_4.$ And (\ref{y6}),   (\ref{y7}),\ (\ref{3.2}),\ (\ref{y9}), (\ref{y10}), (\ref{y11}) and (\ref{y12}) still hold with $d_4=\bc_4$. Thus
\[
\begin{array}{rll}
\bb_3(b_3c_4)&=&\bb_3x_6+\bb_3z_6,\\
&=&c_4+h_6+b_8+\bb_3z_6,\\
b_3(\bb_3c_4)&=&b_3\bar{y}_6+b_3\bar{p}_6\\
&=&c_4+\bar{i}_6+b_8+b_3\bar{p}_6.
\end{array}
\]
Then
\begin{eqnarray}\label{z1}h_6+\bb_3z_6=\bar{i}_6+b_3\bar{p}_6.\end{eqnarray}\\
\par {\bf Substep 1.} There exists no NITA such that $h_6=\bar{i}_6$.\\
\par If $h_6=\bar{i}_6$, then it follows from (\ref{y6}) and (\ref{y9}) that
there is $\alpha_3\in B$ such that
\[
\begin{array}{rll}
b_3h_6&=&b_9+x_6+\alpha_3,\\
b_3\bar{h}_6&=&b_9+y_6+\beta_3.
\end{array}
\]
It is easy to see that $\alpha_3\neq \beta_3$. Since
$z_{12}=h_6+\bar{h}_6$ right now, it means by {\bf Lemma \ref{be}}
that
$$b_5\bb_9=x_6+y_6+3b_9+\alpha_3+\beta_3.$$
Thus $(b_5\bb_9,\ b_5\bb_9)=13.$
\par By {\bf Lemma \ref{be}}, we have that
$$b_9\bb_9=1+3b_8+c_4+\bc_4+h_6+\bar{h}_6+\bb_3f_{12}.$$
But $b_5\bb_5=1+2b_8+c_4+\bc_4$. Since $(b_9\bb_9,\ b_8)=3$, we
have that
$$13=(b_5\bb_9,\ b_5\bb_9)=9+(c_4,\ \bb_3f_{12})+(\bc_4,\ \bb_3f_{12}).$$
It concludes that $(c_4,\ \bb_3f_{12})=(\bc_4,\ \bb_3f_{12})=2$.
Thus
$$(b_3c_4,\ f_{12})=(b_3\bc_4,\ f_{12})=2.$$
From  (\ref{y11}) and (\ref{y12}), it follows that
$$(x_6+z_6,\ f_{12})=(y_6+p_6,\ f_{12})=2.$$
Since $x_6\neq y_6$ by $(b_3\bb_3,\ b_5\bb_5)=3$ and $z_6\neq
x_6$, $y_6\neq p_6$ by {\bf Lemma \ref{t}}, it follows that
$f_{12}=2z_6$ and $z_6=p_6$. Hence $b_3c_7=b_9+2z_6$ by {\bf Lemma
\ref{be}}. Therefore $(b_3c_7,\ b_3c_7)=5$ and $(\bb_3z_6,\
c_7)=2$. So $\bb_3z_6=2c_7+c_4$.  But
$$b_8c_7=b_8+h_6+\bar{h}_6+2\bb_3z_6=b_8+h_6+\bar{h}_6+4c_7+2c_4.$$
It is impossible for $b_8c_7$ is real and $c_4$ not real.\\
\par {\bf Substep 2.} There exists no NITA such that $h_6\neq \bar{i}_6$.\\
\par  Suppose that $h_6\neq\bar{i}_6$. By (\ref{z1}), it follows that
$(b_3\bar{p}_6,\ h_6)=(p_6,\ b_3\bar{h}_6)=(\bar{i}_6,\
\bb_3z_6)=(b_3\bar{i}_6,\ z_6)=1$. But $b_9\in
Supp\{b_3\bar{h}_6\}$ by (\ref{y9}) and $(c_4,\
b_6\bar{p}_6)=(c_4,\ \bb_3z_6)=1$ by (\ref{y11}) and (\ref{y12}).
By (\ref{z1}) , (\ref{y9}) and (\ref{y12}),  there exist
$\delta_3$,\ $z_3$ and $v_8\in \mathit{N}^*B$ such that
\begin{eqnarray}
\label{z5}b_3\bar{h}_6&=&b_9+p_6+\delta_3,\\
\label{z6}\bb_3z_6&=&\bar{i}_6+c_4+v_8,\\
\label{z7}b_3\bar{p}_6&=&h_6+c_4+v_8.
\end{eqnarray}
Hence $\bb_3\delta_3=\bar{h}_6+\gamma_3$, some $\gamma_3\in B$.
Therefore $(b_3\gamma_3,\ b_3\gamma_3)\geq 2$. It follows that
$$\delta_3\bar{\delta}_3=\gamma_3\bar{\gamma}_3=1+b_8$$ by {\bf
Lemma \ref{t}}.
\par Since
\[
\begin{array}{rll}
(b_3\bb_3)c_4&=&c_4+c_4b_8,\\
b_3(\bb_3c_4)&=&b_3\bar{y}_6+b_3\bar{p}_6,\\
&=&c_4+\bar{i}_6+b_8+h_6+c_4+v_8,
\end{array}
\]
we have that $c_4b_8=c_4+h_6+\bar{i}_6+v_8+b_8$. Then $(c_4b_8,\
c_4b_8)\geq 5$. Let $c_4\bc_4=1+b_8+r_7$, $r_\in \mathit{N}^*B$
and $r_7$ is real. By (\ref{y10}) and $c_4$ is not real, we have
that $r_7=c_7$ and then $(c_4b_8,\ c_4b_8)= 5$. Consequently
$v_8\in B$.
 By $b_3(\delta_3\bar{\delta}_3)=(b_3\bar{\delta}_3)\delta_3$, we have that
\[
\begin{array}{rll}
\delta_3\bar{\gamma}_3&=&b_3+x_6\ \mbox{or}\ b_3+y_6,\\
\delta_3\bar{h}_6&=&b_3+b_9+y_6\ \mbox{or}\ b_3+b_9+x_6.
\end{array}
\]
Without loss of generality, let $\delta_3\bar{\gamma}_3=b_3+x_6$,
then
\[
\begin{array}{rll}
1+2b_8+b_8^2=(\delta_3 \bar{\delta}_3)(\gamma_3 \bar{\gamma}_3)=
(\delta_3\bar{\gamma}_3)(\bar{\delta}_3\gamma_3)&=&b_3\bb_3+b_3\bar{x}_6+\bb_3x_6+x_6\bar{x}_6\\
&=&1+3b_8+c_4+\bc_4+h_6+\bar{h}_6+x_6\bar{x}_6.
\end{array}
\]
 It holds
by (\ref{y10}) that
$$x_6\bar{x}_6=1+2b_8+c_7+i_6+\bar{i}_6.$$
Furthermore $(b_5\bb_5,\ x_6\bar{x}_6)=5$. On the other hand, the
following equations hold:
\[
\begin{array}{rll}
b_3(b_5\bb_5)&=&b_3+2b_3b_8+b_3c_4+b_3\bc_4,\\
&=&3b_3+2b_9+2x_6+2y_6+x_6+z_6+y_6+p_6,\\
(b_3\bb_5)b_5&=&\bb_3b_5+\bar{x}_6b_5+\bar{y}_6b_5\\
&=&b_3+x_6+y_6+\bar{x}_6b_5+\bar{y}_6b_5.
\end{array}
\]
Hence
\begin{eqnarray}\label{y14}\bar{x}_6b_5+\bar{y}_6b_5=2b_3+2b_9+2x_6+2y_6+z_6+p_6.\end{eqnarray}
It is proved that $(b_5\bb_5,\ x_6\bar{x}_6)=5$. By above equation
and $(b_5\bar{x}_5,\ b_3)=(b_5\bar{y}_5,\ b_3)=1$, we can only
have that
\[
\begin{array}{rll}
b_5\bar{x}_6&=&b_3+b_9+x_6+y_6+z_6\ \mbox{or}\   b_3+b_9+x_6+y_6+p_6,\\
b_5\bar{y}_6&=&b_3+b_9+x_6+y_6+p_6\ \mbox{or}\
b_3+b_9+x_6+y_6+z_6.
\end{array}
\]
Hence
\[
\begin{array}{rll}
\bb_3(b_5\bar{x}_6)&=&\bb_3b_3+\bb_3b_9+\bb_3x_6+\bb_3y_6+\bb_3z_6\ \mbox{or}\ \bb_3b_3+\bb_3b_9+\bb_3x_6+\bb_3y_6+\bb_3p_6\\
&=&1+b_8+c_7+b_8+\bar{h}_6+\bar{i}_6+c_4+h_6+b_8+\bc_4+i_6+b_8+\bar{i}_6+c_4+v_8\\
&\mbox{or}&1+b_8+c_7+b_8+\bar{h}_6+\bar{i}_6+c_4+h_6+b_8+\bc_4+i_6+b_8+\bar{h}_6+\bc_4+\bar{v}_8,\\
(\bb_3b_5)\bar{x}_6&=&b_3\bar{x}_6+x_6\bar{x}_6+y_6\bar{x}_6\\
&=&\bc_4+\bar{h}_6+b_8+1+2b_8+c_7+i_6+\bar{i}_6+y_6\bar{x}_6.
\end{array}
\]
Therefore

$$y_6\bar{x}_6=2c_4+h_6+b_8+\bar{i}_6+v_8\ \mbox{or}\ c_4+\bc_4+h_6+b_8+\bar{h}_6+\bar{v}_8.$$
\par  If $y_6\bar{x}_6=2c_4+h_6+b_8+\bar{i}_6+v_8$, then $(c_4x_6,\ y_6)=2$. Hence
$(c_4x_6,\ c_4x_6)\geq 5$, a contradiction. Therefore
\begin{eqnarray}\label{z2}y_6\bar{x}_6=c_4+\bc_4+h_6+b_8+\bar{h}_6+\bar{v}_8.\end{eqnarray}
Thus $(\bc_4x_6,\ y_6)=1$.
\par Considering the following equations:
\[
\begin{array}{rll}
b_3(c_4\bc_4)&=&b_3+b_3b_8+b_3c_7\\
&=&2b_3+b_9+x_6+y_6+b_9+f_{12},\\
(b_3c_4)\bc_4&=&\bc_4x_6+\bc_4z_6,\\
(b_3\bc_4)c_4&=&c_4y_6+c_4p_6.
\end{array}
\]
By the expression of $x_6\bar{x}_6$, we have that $(x_6,\
\bc_4x_6)=0$. Consequently  $(x_6,\ \bc_4z_6)=1$ by above
equations. Because $b_3c_7=b_9+f_{12}$ and $L_1(B)=\{1\}$ and
$L_2(B)=\emptyset$, so the constituents of $f_{12}$ are of degrees
larger than 3. Since $(\bc_4x_6,\ b_3)=(\bc_4z_6,\ b_3)=1$ by
(\ref{y11}), $(\bc_4x_6,\ \bc_4x_6)=4$ and it is proved that
$(\bc_4x_6,\ y_6)=1$, we have that
$$\bc_4x_6=b_3+y_6+b_9+g_6,\ g_6\in B.$$
Consequently $\bc_4z_6=b_3+x_6+b_9+f_6,$ where $f_{12}=f_6+g_6$
and $f_6\in\mathit{N}^*B$.
\par Let $\bb_3g_6=c_7+x_{11}$ and
$\bb_3f_6=c_7+y_{11}$, $x_{11},y_{11} \in N^*(B)$. Then
\[
\begin{array}{rll}
\bb_3(\bc_4x_6)&=&\bb_3b_3+\bb_3y_6+\bb_3b_9+\bb_3g_6\\
&=&1+b_8+\bc_4+i_6+b_8+c_7+b_8+\bar{h}_6+\bar{i}_6+c_7+x_{11},\\
\bc_4(\bb_3x_6)&=&\bc_4c_4+\bc_4h_6+\bc_4b_8\\
&=&1+b_8+c_7+\bc_4h_6+\bc_4+\bar{h}_6+i_6+\bar{v}_8+b_8;\\
\bb_3(\bc_4z_6)&=&\bb_3b_3+\bb_3x_6+\bb_3b_9+\bb_3f_6\\
&=&1+b_8+c_4+h_6+b_8+c_7+b_8+\bar{h}_6+\bar{i}_6+c_7+y_{11},\\
\bc_4(\bb_3z_6)&=&\bc_4c_4+\bc_4\bar{i}_6+\bc_4v_8\\
&=&1+b_8+c_7+\bc_4\bar{i}_6+\bc_4v_8.
\end{array}
\]
Hence
\begin{eqnarray}
\label{z3}\bc_4h_6+\bar{v}_8&=&b_8+\bar{i}_6+c_7+x_{11}\\
\label{z4}\bc_4\bar{i}_6+\bc_4v_8&=&c_4+h_6+2b_8+\bar{i}+\bar{h}_6+c_7+y_{11}.
\end{eqnarray}

\par We assert that $v_8\neq b_8$. Otherwise, if $v_8=b_8$, then $(c_4b_8,\  b_8)=2$. So $(b^2_8,\ c_4)=2$, a contradiction by (\ref{y10}).
Hence $\bar{v}_8\in Supp\{x_{11}\}$ by (\ref{z3}). Let
$x_{11}=\bar{v}_8+x_3$, then it follows by (\ref{z3}) that
$$\bc_4h_6=b_8+\bar{i}_6+c_7+x_3.$$
\par By (\ref{z4}), we have that $\bar{h}_6\in Supp\{\bc_4\bar{i}_6+\bc_4v_8\}$.
So  $(\bc_4h_6,\ i_6)\geq 1$ or $(\bc_4h_6,\ \bar{v}_8)\geq 1$.
The latter case will leads to $v_8=b_8$, a contradiction. Hence $(\bc_4h_6,\ i_6)\geq 1$. Thus $i_6=\bar{i}_6$ and $(\bc_4h_6,\ i_6)=1$. Furthermore
$$x_6\bar{x}_6=1+2b_8+c_7+2i_6.$$
And it follows by (\ref{y7}) and (\ref{y9}) that
$$b_3i_6=b_9+y_6+z_3,\ \mbox{some}\ z_3\in B.$$
Thus
\[
\begin{array}{rll}
b_3(x_6\bar{x}_6)&=&b_3+2b_3b_8+b_3c_7+2b_3i_6\\
&=&3b_3+2b_9+2x_6+2y_6+b_9+g_6+f_6+2b_9+2y_6+2z_3,\\
x_6(b_3\bar{x}_6)&=&x_6\bc_4+x_6\bar{h}_6+x_6b_8\\
&=&b_3+y_6+b_9+g_6+x_6\bar{h}_6+x_6b_8.
\end{array}
\]
Then $x_6\bar{h}_6+x_6b_8=2b_3+2b_9+2x_6+3y_6+f_6+2b_9+2z_3$
\par By the expression of $y_6\bar{x}_6$ and $v_8\neq b_8$, we have that
$(b_8x_6,\ y_6)=1$. Hence $(x_6\bar{h}_6,\ y_6)=2$ and so
$(\bar{x}_6y_6,\ \bar{h}_6)=2$. Again by the expression of
$y_6\bar{x}_6$, we have that $h_6=\bar{h}_6$. Now (\ref{y6})
implies that $x_6\in Supp\{b_3h_6\}$. Then $p_6=x_6$ by
(\ref{z5}). Consequently $b_3\bc_4=x_6+y_6$, from which we have
$(bc_4,\ \bb_3x_6)=1$. Thus $c_4=\bc_4$ by (\ref{y6}).
{\bf Step 4} follows. This complete the proof of the lemma 3.4. \\
\par
Therefore $(b_5\bb_5, b_8)=1$ and $(b_3\bb_3,b_5\bb_5)=2$.
Consequently by lemma 3.1 we have $x_{12} \in B$ and
$(b_3b_8,b_3b_8)=3$. Thus Theorem 3.1 holds. Furthermore, in Lemma
3.1 we proved that $c_7 \in B$.

 \section{General Information of NITA Generated by $b_3$ and  Satisfying $b_3^2=\bb_3+b_6$ and $b^2_3=c_3+b_6$}

\par Since  the structure of NITA generated by $b_3$ and  satisfying $b_3^2=\bb_3+b_6$
and $b^2_3=c_3+b_6$ are related, we investigate them
simultaneously. Some results are written in the same proposition.
The purpose of this section is to obtain (3) and (4) in the {\bf
Main Theorem}. The following lemmas are useful during our
investigation.

\begin{lemma}\label{la2} (1) If $b_3^2=\bb_3+b_6$, then $b_3b_8=b_6+\bb_3b_6$;
\par (2) If $b^2_3=c_3+b_6$, then there exists an element $x_6\in B$ such that $\bb_3c_3=b_3+x_6$ and $c_3\bc_3=1+b_8$. Further $b_3b_8=x_6+\bb_3b_6$ and $c_3b_8=b_6+b_3x_6$.
\end{lemma}
\par {\bf Proof.} (1) follows from $(b_3\bb_3)b_3=b^2_3\bb_3$.
\par Since $b_3\in Supp\{\bb_3c_3\}$, $(b_3\bb_3,\ c_3\bar{c}_3)\geq 2$.   But $(b_3\bb_3,\ c_3\bc_3)\leq 2$, so $(b_3\bb_3,\ c_3\bc_3)=2$ and
$$c_3\bar{c}_3=1+b_8,\  \bb_3c_3=b_3+x_6,\ x_6\in B.$$
Checking $(b_3\bb_3)b_3=b^2_3\bb_3$ and
$c_3(c_3\bc_3)=b_3(\bb_3c_3)$, one has (2).

\begin{lemma}\label{la3} If $(b_3b_8,\ b_3b_8)\geq 4$, then any constituent of $b_3b_8$ and $\bb_3b_6$ is either  $b_3$ or  of degree $>3$ and $(b_3b_8,\ b_3b_8)=4,\ 5$.
\end{lemma}
\par {\bf Proof.} Since $(b_3b_8,\ b_3)=1$, the degree of the sum of the remaining constituents of $b_3b_8$ must be 21.
\par On the other hand, any constituent $y_n$ of $b_3b_8$ must have degree$\geq 4$. Otherwise, $n=3$ and $b_3\bar{y}_3=1+b_8$, which implies that $b_3=y_3$, a contradiction. \\

\begin{theorem}\label{tthm1}
There is no NITA generated by $b_3$ and  satisfying
$b_3^2=\bb_3+b_6$ or $b^2_3=c_3+b_6$ and $(b_3b_8,\ b_3b_8)=2$.
\end{theorem}
\par {\bf Proof.} Since $b_3\bb_3=1+b_8$ and $(b_3b_8,\ b_3b_8)=2$, we may assume that $b_3b_8=b_3+b_{21}$, where $b_{21}\in B$. By Lemma \ref{la2}, one has that
$$b_3+b_{21}=b_6+\bb_3b_6\ \mbox{or}\ x_6+\bb_3b_6.$$
It is impossible.\\

\begin{theorem}\label{tthm2} There is no  NITA generated by $b_3$ and  satisfying
$b_3^2=\bb_3+b_6$ or $b^2_3=c_3+b_6$ and $(b_3b_8,\ b_3b_8)=5$.
\end{theorem}
\par {\bf Proof.} (1) there is no NITA satisfying $b^2_3=\bb_3+b_6$ and $(b_3b_8,\ b_3b_8)=5$.
\par By $b^2_3\bb_3=(b_3\bb_3)b_3$, we have that $b_3b_8=b_6+\bb_3b_6$,
which implies that $b_8\in Supp\{\bb_3b_6\}$. But $b_3\in Supp\{\bb_3b_6\}$. We need a sum with degree 7 of two constituents to make up $\bb_3b_6$, which is impossible by Lemma \ref{la2}.\\

\par (2) there is no NITA satisfying $b^2_3=c_3+b_6$ and $(b_3b_8,\ b_3b_8)=5$.\\

\par By Lemma \ref{la2} and $(b_3b_8,\ b_3b_8)=5$, we may assume that
\[
\begin{array}{rlll}
(A)&b_3b_8&=&b_3+x_6+x_4+y_4+z_7,\\
(B)&b_3b_8&=&b_3+x_6+x_4+y_5+z_6,\\
(C)&b_3b_8&=&b_3+x_6+x_5+y_5+z_5.
\end{array}
\]

\par Since
\[
\begin{array}{rll}
\bb_3(\bb_3c_3)&=&b_3\bb_3+x_6\bb_3\\
\bb_3^2c_3&=&c_3\bar{c}_3+c_3\bb_6\\
&=&1+b_8+c_3\bb_6,
\end{array}
\]
we have that $\bb_3x_6=c_3\bb_6.$
\par Since  $(c_3b_8,\ c_3b_8)=(b_3b_8,\ b_3b_8)=5$ by {\bf Lemma } \ref{la2} and  $c_3b_8=b_6+b_3x_6$ by {\bf Lemma} \ref{la2}, we have  $(b_3x_6,\ b_3x_6)=4.$ Further $(\bb_3x_6,\ \bb_3x_6)=4$. But $(\bb_3x_6,\ b_8)=1$. We may assume
$$c_3\bb_6=\bb_3x_6=b_8+m_3+n_3+v_4,$$
where  $m_3,\ n_3,\ v_4\in B$ and $m_3,\ n_3,\ v_4$ are distinct.
We may assume
\[
\begin{array}{rll}
c_3\bar{m}_3&=&b_6+m^*_3,\ m^*_3\in B,\\
c_3\bar{n}_3&=&b_6+n^*_3.
\end{array}
\]
Hence $(c_3\bar{m}_3,\ c_3\bar{m}_3)=2$. So $m_3\bar{m}_3=1+b_8$.
Moreover
\[
\begin{array}{rll}
(c_3\bc_3)\bar{m}_3&=&\bar{m}_3+\bar{m}_3b_8,\\
(c_3\bar{m}_3)\bc_3&=&\bc_3b_6+\bc_3m^*_3\\
&=&b_8+\bar{m}_3+\bar{n}_3+\bar{v}_4+\bc_3m^*_3
\end{array}
\]
Then $$\bar{m}_3b_8=b_8+\bar{n}_3+\bar{v}_4+\bc_3m^*_3,$$
which implies that $m_3\bar{n}_3=1+b_8$. So $m_3=n_3$, a contradiction.\\

\begin{theorem}\label{tthm3} There exists no NITA generated by $b_3$ and satisfying
 $b^2_3=\bb_3+b_6$, $(b_3b_8,\ b_3b_8)=4.$
\end{theorem}
\par {\bf Proof.} Because $b_3b_8=b_6+\bb_3b_6$ by Lemma \ref{la2}
and by assumption that
  $(b_3b_8,\ b_3b_8)=4$, then $b_8\in Supp\{\bb_3b_6\}$ and $(\bb_3b_6,\ \bb_3b_6)=3$. But $b_3\in Supp\{\bb_3b_6\}$. Thus there exists $p_7\in B$ such that
$$\bb_3b_6=b_3+b_8+p_7,\ b_3b_8=b_3+b_6+b_8+p_7.$$
Hence
\[
\begin{array}{rll}
(b_3\bb_3)b_8&=&b_8+b^2_8\\
\bb_3(b_3b_8)&=&\bb_3b_3+\bb_3b_6+\bb_3b_8+\bb_3p_7\\
&=&1+b_8+b_3+b_8+p_7+\bb_3+\bb_6+\bb_8+\bar{p}_7+\bb_3p_7.
\end{array}
\]
Consequently
$b^2_8=1+b_3+\bb_3+p_7+\bar{p}_7+2b_8+\bb_6+\bb_3p_7.$ On the
other hand,
\[
\begin{array}{rll}
(b_3\bb_3)^2&=&1+2b_8+b^2_8,\\
b^2_3\bb^2_3&=&b_3\bb_3+\bb_3b_6+b_3\bb_6+b_6\bb_6\\
&=&1+b_8+b_3+b_8+p_7+\bb_3+b_8+\bar{p}_3+b_6\bb_6.
\end{array}
\]
We have that $b^2_8=b_3+\bb_3+p_7+\bar{p}_7+b_8+b_6\bb_6.$ Hence
$b_6\bb_6=1+\bb_6+b_8+\bb_3p_7.$
\par Since
\[
\begin{array}{rll}
\bb^2_3b_6&=&b_3b_6+\bb_6b_6\\
&=&b_3b_6+1+b_8+\bb_6+\bb_3p_7,\\
\bb_3(\bb_3b_6)&=&\bb_3b_3+\bb_3b_8+\bb_3p_7\\
&=&1+b_8+\bb_3+\bb_6+b_8+\bar{p}_7+\bb_3p_7,
\end{array}
\]
it holds that $b_3b_6=\bb_3+\bar{p}_7+b_8.$ Now the following
equations follow:
\[
\begin{array}{rll}
b_3(\bb_3b_6)&=&b^2_3+b_3b_8+b_3p_7\\
&=&\bb_3+b_6+b_3+b_6+b_8+p_7+b_3p_7,\\
(b_3b_6)\bb_3&=&\bb_3^2+\bb_3\bar{p}_7+b_8\bb_3\\
&=&b_3+\bb_6+\bb_3+\bb_6+b_8+\bar{p}_7+\bb_3\bar{p}_7.
\end{array}
\]
Then $2b_6+p_7+b_3p_7=2\bb_6+\bar{p}_7+\bb_3\bar{p}_7.$  Since
$b_6$ is nonreal, we have that $(b_3p_7,\bb_6)=2$, a contradiction
to the expression of $\bb_3b_6$. The theorem follows.
\begin{theorem}\label{tthm4} There exists no NITA satisfying $b^2_3=c_3+b_6$,
$(b_3b_8,\ b_3b_8)=4$ and $c^2_3=r_4+s_5.$
\end{theorem}
\par {\bf Proof.} Let investigate Sub-NITA generated by $c_2$. By assumptions, we have that $(b_3\bb_3,\ b_3\bb_3)=(c_3\bc_3,\ c_3\bc_3)$, hence $b_3\bb_3=c_3\bc_3$.
Therefore $(c_3\bc_3,\ b_8^2)=(b_3\bb_3,\ b_8^2)=4$, thus $(c_3b_8,\
c_3b_8)=4$, a contradiction to {\bf Theorem \ref{mthm2}}.

\section{NITA Generated by $b_3$ Satisfying $b_3^2=\bar{b}_3+b_6$ and $b_6$ non-real and $b_{10}\in B$ is real}}
\par
 Now we discuss NITA satisfying the following hypothesis.
\begin{hyp}\label{hyp1}Let $(A,\ B)$  be a $NITA$ generated by an non-real
element $b_{3}\in B$ satisfies $b_{3}\overline{b}_{3}=1+c_{8}$ and
$b_{3}^{2}=\overline{b}_{3}+b_{6}$,\  where $b_{6}\in B$ is
non-real. Here we still assume that $L(B)=1$ and $L_{2}(B)=\O$. In
this section,\  we use the symbols $b_{i},\ c_{i},\ d_{i}$ denote
elements of $B$ with degree $i$,\  where $i\geq 2$.
\end{hyp}
\par By the associative law and {\bf Hypothesis \ref{hyp1}}, we have
$b_{3}(b_{3}\overline{b}_{3})=(b_{3}^{2})\overline{b}_{3}$,
$b_{3}+c_{8}b_{3}=\overline{b}_{3}^{2}+b_{6}\overline{b}_{3}$, so
$$ c_{8}b_{3}=\overline{b}_{6}+b_{6}\overline{b}_{3}.$$
 Hence $(b_{3}b_{6},\ c_{8})=(b_{3}c_{8},\ \overline{b}_{6})=1$,\  so
$$b_{3}b_{6}=c_{8}+b_{10}, $$
therefore $b_{10}\in B$ since $(b_3c_8,\ b_3c_8)=3$. \par  {\bf
Now we shall deal with the case $b_{10}=\overline{b}_{10}\in B$.}

\begin{lemma}\label{lemma4.1}Let $(A,\ B)$ satisfies {\bf Hypothesis \ref{hyp1}} and $\overline{b}_{10}=b_{10}\in B$,\ then

1)  $b_{3}b_{6}=c_{8}+b_{10}$;

2)  $\overline{b}_{3}b_{6}=b_{3}+b_{15}$,\ where $b_{15}\in B$;

3)  $b_{3}c_{8}=b_{3}+\overline{b}_{6}+b_{15}$;

4)  $b_{3}b_{10}=b_{15}+\overline{b}_{6}+\overline{c}_{9},\
c_{9}\in B$;

5)  $b_{3}b_{15}=2\overline{b}_{15}+b_{6}+c_{9}$;

6)  $b_{6}^{2}=2\overline{b}_{6}+\overline{c}_{9}+b_{15}$;

7)  $b_{6}c_{8}=\overline{b}_{3}+2\overline{b}_{15}+b_{6}+c_{9}$;

8)  $\overline{b}_{6}b_{6}=1+c_{8}+x+y+z+w$,\  where $x,\ y,\ z,\
w\in B,\ \mid x+y+z+w\mid =27$,\  and the degrees of $x,\ \ y,\ \
z$ and $w \geq 5$. Moreover $\overline{x+y+z+w}=x+y+z+w$;

9)  $ c_{8}^{2}=1+2c_{8}+2b_{10}+x+y+z+w$;

10)  $b_{3}\overline{b}_{15}=b_{10}+c_{8}+x+y+z+w$;

11)  $(\overline{c}_{9}b_{3},\ c_{9})= 1$ or $(b_{6}c_{9},\
\overline{c}_{9})\geq 1$.
\end{lemma}
\begin{proof}We have seen that
$$ c_{8}b_{3}=\overline{b}_{6}+b_{6}\overline{b}_{3}.\eqno(1) $$
 and
$$b_{3}b_{6}=c_{8}+b_{10},\eqno(2)
$$
where $b_{10}$ is real.
 So we get $(
\overline{b}_{3}b_{6},\ \overline{b}_{3}b_{6}) =2$ and
consequently we get

$$\overline{b}_{3}b_{6}=b_{3}+b_{15},\ \eqno(3)$$ where $b_{15}\in B$.
Hence by (1)
$$c_{8}b_{3}=b_{3}+\overline{b}_{6}+b_{15}  \eqno(4)$$
By the associative law and (2),\ (3), we have
$b_{3}(\overline{b}_{3}b_{6})=(b_{3}b_{6})\overline{b}_{3}$
 and $b_{3}^{2}+b_{15}b_{3}=c_{8}\overline{b}_{3}+b_{10}\overline{b}_{3}$, so by (4) and {\bf Hypothesis
 \ref{hyp1}} one gets
$$b_{15}b_{3}=\overline{b}_{15}+\overline{b}_{3}b_{10}   \eqno(5)$$

If $x\in B$ and $(b_{15}b_{3},\ x)=(b_{15},\ \overline{b}_{3}x)$,\
so $\mid x\mid \geq 5,\  \forall x\in B.$ So every constituent of
$ \overline{b}_{3}b_{10}$ have degrees $\geq 5$. Since $b_{10}$ is
real also constituents of $b_{3}b_{10}$ have degrees $\geq 5$,\ by
the associative law and Hypothesis \ref{hyp1} and (4),\ (2):
$$(\overline{b}_{3}b_{3})c_{8}=(b_{3}c_{8})\overline{b}_{3},\ $$
$$c_{8}+c_{8}^{2}=b_{3}\overline{b}_{3}+\overline{b}_{6}\overline{b}_{3}+b_{15}\overline{b}_{3},\ $$
$$c_{8}^{2}=1+c_{8}+b_{10}+\overline{b}_{3}b_{15}  \eqno(6)$$
 By the associative law and Hypothesis \ref{hyp1} and (2):
 $$(b_{3}\overline{b}_{3})(b_{3}\overline{b}_{3})=b_{3}^{2}\overline{b}_{3}^{2},\ $$
 $$1+2c_{8}+c_{8}^{2}=\overline{b}_{3}b_{3}+\overline{b}_{3}\overline{b}_{6}+b_{6}b_{3}+b_{6}\overline{b}_{6},\ $$
 $$c_{8}^{2}=2b_{10}+c_{8}+b_{6}\overline{b}_{6} \eqno(7)$$
 so we get $$1=(\bb_{3}b_{15},\ b_{10})=(b_{15},\ b_{3}b_{10}).$$
 By (2),\
 $$(b_{3}b_{10},\ \overline{b}_{6})=(b_{3}b_{6},\ b_{10})=1.$$
 $b_{3}b_{10}$ have constituents of degree $\geq 5$ as we proves.
Therefore
$$b_{3}b_{10}=b_{15}+\overline{b}_{6}+\overline{c}_{9}  \eqno(8)$$
where $c_{9}\in B$. So by (5),\
$$b_{3}b_{15}=2\overline{b}_{15}+b_{6}+c_{9} . \eqno(9)$$
By the associative law and Hypothesis \ref{hyp1} and (2),\ (3),\
(4):
$$(b_{3}^{2})b_{6}=(b_{3}b_{6})b_{3},\ $$
$$\overline{b}_{3}b_{6}+b_{6}^{2}=c_{8}b_{3}+b_{10}b_{3},\ $$
consequently $$b_6^2=2\bb_6+b_{15}+\bar{c}_9.$$ Now
$$b_3(\bb_3b_6)=(b_3\bb_3)b_6,$$
$$b_3^2 +b_{15}b_3=b_6+b_6b_8.$$ Thus
$$c_{8}b_{6}=\overline{b}_{3}+2\overline{b}_{15}+b_{6}+c_{9}.$$
So$$(\overline{b}_{6}b_{6},\ c_{8})=(b_{6},\ b_{6}c_{8})=1.$$
Since
$$( b_{6}^{2},\ b_{6}^{2})=6=(
\overline{b}_{6}b_{6},\ \overline{b}_{6}b_{6}),\ $$
 hence
 $$ b_{6}\overline{b}_{6}=1+c_{8}+x+y+z+w$$
 where $x,\ \ y,\ \ z,\ \ w\in B$. Therefore by (7),\ (6)
 $$ c_{8}^{2}=1+2c_{8}+2b_{10}+x+y+z+w,\ $$
 and $$b_{3}\overline{b_{5}}=b_{10}+c_{8}+x+y+z+w.$$
  So $(b_{3}\overline{b}_{15},\ \ m)=(b_{3}\bar{m},\ \ b_{15})$,\  so $\mid m\mid
 \geq5$,\  hence $$ \mid x\mid,\ \ \mid y\mid,\ \ \mid z\mid,\ \ \mid w\mid \geq
 5.$$
 By (8),\
 $$ (b_{3}c_{9},\ b_{10})=(b_{3}b_{10},\ \overline{c}_{9})=1,\ $$
 hence
 $$b_{3}c_{9}=b_{10}+\mathbf{\alpha},\ $$
 where $\alpha \in N^*B$. By the associative law and Hypothesis \ref{hyp1}:
 $$(b_{3}^{2})c_{9}=(b_{3}c_{9})b_{3},\ $$
 $$c_{9}\overline{b}_{3}+c_{9}b_{6}=b_{10}b_{3}+\alpha
 b_{3}=b_{15}+\overline{b}_{6}+\overline{c}_{9}+\alpha b_{3}.$$
 hence,\ $(c_{9}\overline{b}_{3},\ \overline{c}_{9})\geq 1$,\ or
 $(b_{6}c_{9},\ \overline{c}_{9})\geq 1$.
By (9) we conclude that either $(c_9\bb_3,\bar{c}_9)=1$ or
$(b_6c_9, \bar{c}_9)\ge 1$.
\par
We start to investigate the case $(\bar{c}_9b_3,c_9)=1$ of Lemma
\ref{lemma4.1}.

\end{proof}

\begin{lemma}\label{lemma4.2}Let $(A,\ B)$ satisfies Hypothesis \ref{hyp1}. Assume that
$\overline{b}_{10}=b_{10}\in B$ and $( \overline{c}_{9}b_{3},\
c_{9})= 1.$  Then we get the following:

1)  $b_{3}\overline{c}_{9}=\overline{b}_{15}+c_{9}+c_{3},\ c_{3},\
c_{9}\in B$;

2) $b_{6}b_{10}=3\overline{b}_{15}+\bar{b}_{3}+c_{3}+c_{9},\
b_{15}\in B$;

3)
$c_{8}b_{15}=5b_{15}+b_{3}+\overline{c}_{3}+2\overline{b}_{6}+3\overline{c}_{9}$;

4) $b_{3}\overline{c}_{3}=c_{9}$;

5) $b_{3}c_{3}=b_{9},\ \overline{b}_{9}=b_{9}\in B$;

6) $b_{3}c_{9}=b_{9}+b_{10}+b_{8},\ b_{9}=\overline{b}_{9},\
\overline{b}_{8}=b_{8}\in B$;

7) $ \overline{c}_{3}b_{6}=b_{10}+b_{8}$;

8)
$b_{6}\overline{b}_{15}=4b_{15}+\overline{b}_{6}+b_{3}+\overline{c}_{3}+2\overline{c}_{9}$;

9)  $c_{3}\overline{c}_{3}=1+b_{8}$;

10) $ \overline{b}_{6}b_{6}=1+c_{8}+b_{9}+b_{8}+c_{5}+b_{5},\
\overline{c_{5}+b_{5}}=c_{5}+b_{5}$;

11) $c_{8}^{2}=1+2c_{8}+2b_{10}+b_{9}+b_{8}+b_{5}+c_{5}$;

12) $b_{3}\overline{b}_{15}=b_{10}+c_{8}+b_{9}+b_{8}+b_{5}+c_{5}$;

13) $ b_{6}b_{15}=2c_{8}+3b_{10}+2b_{9}+2b_{8}+b_{5}+c_{5}$;

14) $ c_{8}b_{10}=2c_{8}+2b_{10}+2b_{9}+2b_{8}+b_{5}+c_{5}$;

15) $b_{6}\overline{c}_{9}=b_{10}+c_{8}+2b_{9}+b_{8}+b_{5}+c_{5}$;

16) $b_{10}^{2}=1+2c_{8}+2b_{10}+3b_{9}+2b_{8}+2b_{5}+2c_{5}$.
\end{lemma}
\begin{proof}By 5) in {\bf Lemma \ref{lemma4.1}}
,\ $(\overline{b}_{3}c_{9},\ b_{15})=(c_{9},\ b_{3}b_{15})=1$. So
by the assumption in the lemma,\ $(b_{3}\overline{c}_{9},\
c_{9})=1$,\ hence
$$b_{3}\overline{c}_{9}=\overline{b}_{15}+c_{9}+c_{3},\  \eqno(1)$$
where $ c_{3}\in B$. By the associative law and Hypothesis
\ref{hyp1} and 4) in {\bf Lemma \ref{lemma4.1}}
$$(b_{3}^{2})b_{10}=(b_{3}b_{10})b_{3},\ $$
$$
\overline{b}_{3}b_{10}+b_{6}b_{10}=b_{15}b_{3}+\overline{b}_{6}b_{3}+\overline{c}_{9}b_{3},\
$$ by 2),\ 4),\ 5) in {\bf Lemma \ref{lemma4.1}} and (1),\
$$b_{6}b_{10}=3\overline{b}_{15}+\overline{b}_{3}+c_{3}+c_{9}.\eqno(2)$$
By the associative law and Hypothesis \ref{hyp1} and 5) in {\bf
Lemma \ref{lemma4.1}}:
$$(\overline{b}_{3}b_{3})b_{15}=(b_{3}b_{15})\overline{b}_{3},\ $$
$$c_{8}b_{15}+b_{15}=2\overline{b}_{15}\overline{b}_{3}+\overline{b}_{3}b_{6}+c_{9}\overline{b}_{3},\ $$
by (1) and 2),\ 5) in {\bf Lemma \ref{lemma4.1}}:
$$c_{8}b_{15}=5b_{15}+b_{3}+\overline{c}_{3}+2\overline{b}_{6}+3\overline{c}_{9}.
\eqno(3)$$ By (1),\ $$( \overline{c}_{3}b_{3},\
c_{9})=(b_{3}\overline{c}_{9},\ c_{3})=1,\ $$ hence $$
\overline{c}_{3}b_{3}=c_{9}    \eqno(3').$$ So $$( c_{3}b_{3},\
c_{3}b_{3})=( b_{3}\overline{c}_{9},\ c_{3})=1,\ $$ hence
$$c_{3}b_{3}=b_{9}. \eqno(4)$$ By the associative law and Hypothesis
\ref{hyp1} and (3)
$$\overline{c}_{3}(b_{3}^{2})=(\overline{c}_{3}b_{3})b_{3},\ $$
$$
\overline{c}_{3}\overline{b}_{3}+\overline{c}_{3}b_{6}=c_{9}b_{3},\
$$ hence $$ c_{9}b_{3}=\overline{b}_{9}+\overline{c}_{3}b_{6}
\eqno(5)$$ (1) and (4) in {\bf Lemma \ref{lemma4.1}},\
$(c_{9}b_{3},\ b_{10})=(b_{3}b_{10},\ \overline{c}_{9})=1,\ (
b_{3}c_{9},\ b_{3}c_{9})=( b_{3}\overline{c}_{9},\
b_{3}\overline{c}_{9})=3$,\ hence
$$c_{9}b_{3}=\overline{b}_{9}+b_{10}+b_{8} ,\ \eqno(6)$$
where $b_{8}\in B$. So
$$\overline{c}_{3}b_{6}=b_{10}+b_{8}   \eqno(7)$$
By the associative law and Hypothesis \ref{hyp1} and 4) in {\bf
Lemma \ref{lemma4.1}}:
$$(\overline{b}_{3}b_{3})b_{10}=(b_{3}b_{10})\overline{b}_{3}$$
$$b_{10}+c_{8}b_{10}=b_{15}\overline{b}_{3}+\overline{b}_{6}\overline{b}_{3}+\overline{c}_{9}\overline{b}_{3},\ $$
by 1),\ 10) in {\bf Lemma \ref{lemma4.1}} and (6)
$$c_{8}b_{10}=c_{8}+\overline{x}+\overline{y}+\overline{z}+\overline{w}+\overline{b}_{8}.
\eqno(8)$$ $ b_{10},\ c_{8}$ are reals ,\  so if
$\overline{b}_{8}\neq b_{8}$,\  then $(c_{8}b_{10},\ b_{8})\neq
0$,\ and without lose of generality  $\overline{x}=b_{8}$,\   so
$x=b_{8}$ and by 8) in {\bf Lemma \ref{lemma4.1}} and without lose
of generality $ y=b_{8}$,\  so by (8) $(c_{8}b_{10},\
\overline{b}_{8})=2$,\  hence $(c_{8}b_{10},\ b_{8})=2$,\  so
without lose of generality  $ \overline{z}=b_{8}$ ,\   hence
$z=\overline{b}_{8}$ and by 8) in {\bf Lemma \ref{lemma4.1}}
without lose of generality $w=b_{8}$ then $(c_{8}b_{10},\
b_{8})=3$ a contradiction. Hence $b_{8}=\overline{b}_{8}$.  In the
same way $b_{9}=\overline{b}_{9}$.

By the associative law and 2),\ 6) in {\bf Lemma \ref{lemma4.1}}:
$$(b_{6}^{2})\overline{b}_{3}=(b_{6}\overline{b}_{3})b_{6}$$
$$2\overline{b}_{6}\overline{b}_{3}+\overline{c}_{9}\overline{b}_{3}+b_{15}\overline{b}_{3}=b_{3}b_{6}+b_{15}b_{6}.$$
and by {\bf Lemma \ref{lemma4.1}} and (6)
$$b_{6}b_{15}=2c_{8}+3b_{10}+b_{8}+b_{9}+x+y+w+z. \eqno(9)$$
By the associative law and Hypothesis and 4) in {\bf Lemma
\ref{lemma4.1}}:
$$b_{10}(b_{3}\overline{b}_{3})=(\overline{b}_{3}b_{10})b_{3},\ $$
$$b_{10}+c_{8}b_{10}=b_{3}\overline{b}_{15}+b_{3}b_{6}+b_{3}c_{9},\ $$
by (6) and 1) ,\ 10) in {\bf Lemma \ref{lemma4.1}}
$$c_{8}b_{10}=2c_{8}+2b_{10}+x+y+z+w+b_{8}+b_{9}.   \eqno(10)$$
By the associative law and 2),\ 3),\ 6) in {\bf Lemma
\ref{lemma4.1}} and (3):
$$c_{8}(\overline{b}_{3}b_{6})=(\overline{b}_{3}c_{8})b_{6},\ $$
$$b_{3}c_{8}+b_{15}c_{8}=\overline{b}_{3}b_{6}+b_{6}^2+\overline{b}_{15}b_{6},\ $$
$$b_{6}\overline{b}_{15}=4b_{15}+\overline{b}_{6}+b_{3}+\overline{c}_{3}+2\overline{c}_{9}.$$
By the associative law and (1),\ (4),\ (6) and Hypothesis1:
$$(b_{3}^2)\overline{c}_{9}=(b_{3}\overline{c}_{9})b_{3},\ $$
$$\overline{b}_{3}\overline{c}_{9}+b_{6}\overline{c}_{9}=\overline{b}_{15}b_{3}+c_{9}b_{3}+c_{3}b_{3},\ $$
$$b_{6} \overline{c}_{9}=b_{10}+c_{8}+x+y+z+w+b_{9} . \eqno(11)$$
By the associative law and 1)in the {\bf Lemma \ref{lemma4.1}} and
(2):
$$(b_{3}b_{6})b_{10}=(b_{6}b_{10})b_{3}$$
$$c_{8}b_{10}+b_{10}^2=3\overline{b}_{15}b_{3}+\overline{b}_{3}b_{3}+c_{3}b_{3}+c_{9}b_{3},\ $$
by (4),\ (6),\ (10) and Hypothesis \ref{hyp1}
$$b_{10}^2=1+2c_{8}+b_{9}+2b_{10}+2x+2y+2z+2w.  \eqno(12)$$
By (2),\ $( b_{6}b_{10},\ b_{6}b_{10})=12$.  But by 8) in {\bf
Lemma \ref{lemma4.1}} and (12) $( b_{6} \overline{b}_{6},\
b_{10}^2)=10$,\ so we get without lose of generality
$$x=b_{9}. \eqno(13)$$
By (4) and that $b_{9}=\overline{b}_{9},\ ( c_{3}b_{3},\
\overline{c}_{3}\overline{b}_{3})=( b_{3}^2,\
\overline{c}_{3}^2)=1$,\  but by (3') $(c_{3}^2,\ b_{3})=(c_{3}
\overline{b}_{3},\ \overline{c}_{3})=0$,\  so by hypothesis
\ref{hyp1}
$$(\overline{c}_{3}^2,\ b_{6})=1.\eqno(14)$$
Hence  $$( c_{3} \overline{c}_{3},\ c_{3} \overline{c}_{3})=2,\ $$
so
$$c_{3} \overline{c}_{3}=1+t_{8}. \eqno(15)$$ By (4)$$(
c_{3}b_{3},\ c_{3}b_{3})=( c_{3}\overline{c}_{3},\
b_{3}\overline{b}_{3})=1,\ $$  therefore by hypothesis $t_{8}\neq
c_{8}$.  By (12),\ (13) and 9) in {\bf Lemma \ref{lemma4.1}}
$$( b_{10}^2,\ c_{8}^2)=18=(
c_{8}b_{10},\ c_{8}b_{10}).$$ So by (8) we conclude that without
lose of generality
$$y=b_{8}. \eqno(16)$$
So by (13),\ (16) and 8) in lemma 1 we obtain $\mid z\mid=\mid
w\mid=5$,\   set $z=b_{5}$ and $w=c_{5}$.  By (7),\
$$(
\overline{c}_{3}b_{6},\ \overline{c}_{3}b_{6})=(
\overline{c}_{3}c_{3},\ \overline{b}_{6}b_{6})=2,\ $$ so by (15)
and 8) in {\bf Lemma \ref{lemma4.1}} and by (16) $t_{8}=b_{8}$.
So $c_{3}\overline{c}_{3}=1+b_{8}$. By 8),\ 9),\ 10) in {\bf Lemma
\ref{lemma4.1}} and (10),\ (11),\ (12),\ (13),\ (16),\ (9)
$$ \overline{b}_{6}b_{6}=1+c_{8}+b_{9}+b_{8}+c_{5}+b_{5},\ $$
$$c_{8}^2=1+2c_{8}+2b_{10}+b_{9}+b_{8}+b_{5}+c_{5},\ $$
$$b_{3}\overline{b}_{15}=b_{10}+c_{8}+b_{9}+b_{8}+b_{5}+c_{5},\ $$
$$b_{6}b_{15}=2c_{8}+3b_{10}+2b_{9}+2b_{8}+b_{5}+c_{5},\ $$
$$c_{8}b_{10}=2c_{8}+2b_{10}+2b_{9}+2b_{8}+b_{5}+c_{5},\ $$
$$b_{6}\overline{c}_{9}=b_{10}+c_{8}+2b_{9}+b_{8}+b_{5}+c_{5},\ $$
$$b_{10}^2=1+2c_{8}+2b_{10}+3b_{9}+2b_{8}+2b_{5}+2c_{5}.$$

\end{proof}

\begin{lemma}\label{lemma4.3}Let $(A,\ B)$ satisfy Hypothesis \ref{hyp1}. Assume that $
\overline{b}_{10}=b_{10}\in B$ and $(\overline{c}_{9}b_{3},\
c_{9})= 1$.  Then we get the following equations:

1)  $b_{3}b_{9}=b_{15}+\overline{c}_{9}+\overline{c}_{3}$;

2)  $c_{3}c_{8}=c_{9}+\overline{b}_{15}$;

3)  $b_{3}b_{8}=b_{15}+\overline{c}_{9}$;

4)  $b_{3}c_{5}=b_{15}$;

5)  $b_{3}b_{5}=b_{15}$;

6)  $ b_{6}b_{5}=\overline{b}_{15}+b_{6}+c_{9}$;

7)  $b_{6}c_{9}=2b_{15}+2\overline{c}_{9}+\overline{b}_{6}$;

8)  $c_{3}b_{6}=\overline{c}_{3}+b_{15}$;

9)  $b_{6}b_{9}=2\overline{b}_{15}+b_{6}+2c_{9}$;

10)  $b_{6}c_{5}=\overline{b}_{15}+b_{6}+c_9$;

11)  $b_{6}b_{8}=2\overline{b}_{15}+b_{6}+c_{3}+c_{9}$;

12)  $c_3b_{10}=\overline{b}_{15}+b_6+c_9$;

13)  $c_3c_9=b_{15}+\overline{c}_9+b_3$;

14)  $c_3b_8=c_3+b_6+\overline{b}_{15}$;

15)  $c_3b_9=\overline{b}_{15}+\overline{b}_3+c_9$;

16)  $c_3\overline{b}_{15}=2b_{15}+\overline{b}_6+\overline{c}_9$;

17)
$b_8\overline{b}_{15}=5\overline{b}_{15}+\overline{b}_3+c_3+2b_6+3c_9$;

18)  $b_8c_9=3\overline{b}_{15}+2c_9+b_6+\overline{b}_3$.
\end{lemma}
\begin{proof}By 5),\ 6),\ 12),\ in Lemma \ref{lemma4.2}
$$(b_3b_9,\ \overline{c}_3)=(b_3c_3,\ b_9)=1,\ $$
$$(b_3b_9,\ \overline{c}_9)=(b_3c_9,\ b_9)=1,\ $$
$$(b_3b_9,\ b_{15})=(b_3\overline{b}_{15},\ b_9)=1,\ $$
then
$$b_3b_9=b_{15}+\overline{c}_9+\overline{c}_3.\eqno(1)$$
By the associative law and 1),\ 4) in {\bf Lemma \ref{lemma4.2}}
and Hypothesis \ref{hyp1}
$$( \overline{b}_3c_3)b_3=(\overline{b}_3b_3)c_3,\ $$
$$ \overline{c}_9b_3=c_3+c_3c_8,\ $$
so
$$c_3c_8=c_9+\overline{b}_{15}.   \eqno(2)$$
By 6),\ 12) in {\bf Lemma \ref{lemma4.2}},\
$$(b_3b_8,\ b_{15})=(b_3\overline{b}_{15},b_8)=1,\ $$
$$(b_3b_8,\ \overline{c}_9)=(b_3c_9,\ b_8)=1,\ $$
then $$b_3b_8=b_{15}+\overline{c}_9.\eqno(3)$$ By the associative
law and 2) in {\bf Lemma \ref{lemma4.1}} and 10) in {\bf Lemma
\ref{lemma4.2}}:
$$(b_6\overline{b}_6)b_3=(b_3\overline{b}_6)b_6,\ $$
$$b_3+c_8b_3+b_9b_3+b_8b_3+c_5b_3+b_5b_3=\overline{b}_3b_6+\overline{b}_{15}b_6.$$
by 2),\ 3)in {\bf Lemma \ref{lemma4.1}} and (3),\ (1) and 8) in
{\bf Lemma \ref{lemma4.2}}:
$$c_5b_3+b_5b_3=2b_{15},\ $$
so $$c_5b_3=b_{15}   ,\  \eqno(4)$$
$$b_5b_3=b_{15}.\eqno(5)$$
By the associative law and Hypothesis \ref{hyp1} and (5) and 5) in
Lemma \ref{lemma4.1}
$$(b_{3}^2)b_5=(b_3b_5)b_3=b_3b_{15},\ $$
$$ \overline{b}_3b_5+b_6b_5=2\overline{b}_{15}+b_6+c_9,\ $$
by (5)
$$(\overline{b}_3b_5,\ \overline{b}_3b_5)=(b_5b_3,\ b_5b_3)=1.$$
So $\overline{b}_3b_5=\overline{b}_{15}$,\   therefore
$$b_6b_5=\overline{b}_{15}+b_6+c_9.  \eqno(6)$$
By the associative law and 4),\ 7) in {\bf Lemma \ref{lemma4.2}}:
$$(b_3\overline{c}_3)b_6=( \overline{c}_3b_6)b_3,\ $$
$$c_9b_6=b_{10}b_3+b_8b_3,\ $$
by 4) in {\bf Lemma \ref{lemma4.1}} and by (3)
$$b_6c_9=2b_{15}+2\overline{c}_9+\overline{b}_6.\eqno(7)$$
By the associative law and Hypothesis \ref{hyp1} and 5) in {\bf
Lemma \ref{lemma4.2}}:
$$(b_{3}^2)c_3=(b_3c_3)b_3,\ \overline{b}_3c_3+b_6c_3=b_9b_3,\ $$
so  by (1) and 4) in {\bf Lemma \ref{lemma4.2}}
$$b_6c_3=\overline{c}_3+b_{15}. \eqno(8)$$
10),\ 13),\ 15) in {\bf Lemma \ref{lemma4.2}},\
$$(b_6b_9,\ \overline{b}_{15})=(b_6b_{15},\ b_9)=2,\ $$
$$(b_6b_9,\ b_6)=(b_9,\ \overline{b}_6b_6)=1,\ $$
$$(b_6b_9,\ c_9)=(b_6\overline{c}_9,\ b_9)=2,\ $$
so $$b_6b_9=2\overline{b}_{15}+b_6+2c_9. \eqno(9)$$ By the
associative law and Hypothesis \ref{hyp1} and (4):
$$(b_{3}^2)c_5=(b_3c_5)b_3,\ \overline{b}_3c_5+b_6c_5=b_{15}b_3,\ $$
by 5) in {\bf Lemma
\ref{lemma4.1}}
$$\overline{b}_3c_5+b_6c_5=2\overline{b}_{15}+b_6+c_9,\
$$
$$(c_5\overline{b}_3,\ c_5\overline{b}_3)=(c_5b_3,\ c_5b_3)=1,\ $$
so $\overline{b}_3c_5=\overline{b}_{15}$,\   hence by 7), 10), 13),
and 15) in Lemma \ref{lemma4.2}
$$b_6c_5= \bb_{15}+b_6+c_9 \eqno(10)$$
$$(b_6b_8,\ c_3)=(b_6\overline{c}_3,\ b_8)=1,\ (b_6b_8,\ b_6)=(b_8,\ \overline{b}_6b_6)=1,\ $$
$$(b_6b_8,\ c_9)=(b_8,\ \overline{b}_6c_9)=1,\ (b_6b_8,\ \overline{b}_{15})=(b_6b_{15},\ b_8)=2,\ $$
so $$b_6b_8=2\overline{b}_{15}+b_6+c_3+c_9.  \eqno(11)$$ By the
associative law and 1) in {\bf Lemma \ref{lemma4.1}} and (8)
$$c_3(b_3b_6)=(c_3b_6)b_3,\ c_3c_8+c_3b_{10}=\overline{c}_3b_3+b_{15}b_3,\ $$
by (2),\ 4) in {\bf Lemma \ref{lemma4.2}} and 5) in {\bf Lemma
\ref{lemma4.1}}
$$c_3b_{10}=\overline{b}_{15}+b_6+c_9. \eqno(12)$$
By the associative law and 4),\ 9) in {\bf Lemma \ref{lemma4.2}}:
$$(c_3\overline{c}_3)b_3=(
\overline{c}_3b_3)c_3,\ b_3+b_8b_3=c_9c_3,\ $$ by (3)
$$c_3c_9=b_{15}+\overline{c}_9+b_3. \eqno(13)$$
By the associative law and 6) in {\bf Lemma \ref{lemma4.2}} and
(13)
$$(b_3c_9)c_3=(c_3c_9)b_3,\ $$
$$b_9c_3+b_{10}c_3+b_8c_3=b_{15}b_3+\overline{c}_9b_3+b_{3}^2,\ $$
by (12) and 5) in {\bf Lemma \ref{lemma4.1}} and 1) in {\bf Lemma
\ref{lemma4.2}} and Hypothesis \ref{hyp1}
$$b_9c_3+b_8c_3+b_{15}+b_6+c_9=2\overline{b}_{15}b_6+c_9+\overline{b}_{15}+c_9+c_3+\overline{b}_3+b_6,\ $$
$$b_9c_3+b_8c_3=2\overline{b}_{15}+c_9+c_3+\overline{b}_3+b_6,\ $$
by 5),\ 7),\ 9) in {\bf Lemma \ref{lemma4.2}}
$$(b_8c_3,\ \overline{b}_3)=(c_3b_3,\ b_8)=0,\ (b_8c_3,\ c_3)=(b_8,\ \overline{c}_3c_3)=1,\ (b_8c_3,\ b_6)=(b_8,\ \overline{c}_3b_6)=1,\ $$
so $$b_8c_3=c_3+b_6+\overline{b}_{15},\  \eqno(14)$$ hence
$$b_9c_3=\overline{b}_{15}+\overline{b}_3+c_9.\eqno(15)$$
By the associative law and 2) in {\bf Lemma \ref{lemma4.1}}and 7)
in {\bf Lemma \ref{lemma4.2}}:
$$c_3(b_3\overline{b}_6)=(c_3\overline{b}_6)b_3,\ c_3\overline{b}_3+c_3\overline{b}_{15}=b_{10}b_3+b_8b_3,\ $$
by 4) in {\bf Lemma \ref{lemma4.2}} and 4) in {\bf Lemma
\ref{lemma4.1}} and (3):
$$c_3\overline{b}_{15}=2b_{15}+\overline{b}_6+\overline{c}_9.
\eqno(16)$$ By the associative law and 9) in {\bf Lemma
\ref{lemma4.2}}:
$$(c_3\overline{b}_{15})\overline{c}_3=(c_3\overline{c}_3)\overline{b}_{15},\
2b_{15}\overline{c}_3+\overline{b}_6\overline{c}_3+\bar{c}_9\bar{c}_3=\overline{b}_{15}+b_8\overline{b}_{15},\
$$ by (16) and (8) ,\ (13)
$$b_8\overline{b}_{15}=5\overline{b}_{15}+\overline{b}_3+c_3+2b_6+3c_9.\eqno(17)$$
By the associative law and (13) and 9) in {\bf Lemma
\ref{lemma4.2}}:
$$ (c_3\overline{c}_3)c_9=(c_3c_9)\overline{c}_3,\
c_9+b_8c_9=b_{15}\overline{c}_3+\overline{c}_9\overline{c}_3+b_3\overline{c}_3,\
$$ by (13),\ (16) and 4) in {\bf Lemma \ref{lemma4.2}}:
$$c_9+b_8c_9=2\overline{b}_{15}+b_6+c_9+\overline{b}_{15}+c_9+\overline{b}_3+c_9,\ $$
then
$$b_8c_9=3\overline{b}_{15}+2c_9+b_6+\overline{b}_3.\eqno(18)$$

\end{proof}

\begin{lemma}\label{lemma4.4}Let $(A,\ B)$ satisfies Hypothesis \ref{hyp1}. Assume that $
\overline{b}_{10}=b_{10}\in B$ and $(\overline{c}_{9}b_{3},\
c_{9})= 1$.  Then we get the following equations:

1)  $c_3\overline{c}_9=b_{10}+b_9+c_8$;

2)  $c_3b_{15}=b_5+c_5+c_8+b_8+b_9+b_{10}$;

3)  $c_{3}^2=\overline{c}_3+\overline{b}_6$;

4)  $\overline{b}_{15}b_{15}=1+3b_5+3c_5+5c_8+5b_8+6b_6+6b_{10}$;

5)  $b_{15}^2=9\overline{b}_{15}+4b_6+2\overline{b}_3+2c_3+6c_9$;

6)
$b_{10}b_{15}=6b_{15}+3\overline{b}_6+b_3+\overline{c}_3+4\overline{c}_9$;

7)  $c_9b_{15}=2b_5+2c_5+3c_8+3b_8+3b_9+4b_{10}$;

8)
$c_9\overline{b}_{15}=6b_{15}+2\overline{b}_6+b_3+\overline{c}_3+3\overline{c}_9$;

9)
$c_{9}^2=3b_{15}+b_3+\overline{c}_3+2\overline{b}_6+2\overline{c}_9$;

10)  $c_9\overline{c}_9=1+2c_8+2b_9+2b_{10}+c_5+b_5+2b_8$;

11)  $c_8c_9=3\overline{b}_{15}+b_6+c_3+2c_9$;

12)  $b_9c_9=\overline{b}_3+2c_9+3\overline{b}_{15}+c_3+2b_6$;

13)  $c_9b_{10}=\overline{b}_3+2c_9+3\overline{b}_{15}+c_3+2b_6$;

14)  $b_5c_9=2\overline{b}_{15}+b_6+c_9$;

15)  $c_5c_9=2\overline{b}_{15}+b_6+c_9$;

16)  $b_5=\overline{b}_5,\ c_5=\overline{c}_5$;

17)  $c_3b_5=\overline{b}_{15}$;

18)  $c_3c_5=\overline{b}_{15}$.

\end{lemma}
\begin{proof}By 2),\ 12),\ 15) in {\bf Lemma \ref{lemma4.3}},\
$(c_3\overline{c}_9,\ b_{10})=(c_3b_{10},\ c_9)=1,\
(c_3\overline{c}_9,\ b_9)=(c_3b_9,\ c_9)=1,\ (c_3\overline{c}_9,\
c_8)=(c_3c_8,\ c_9)=1$,\   then
$$c_3\overline{c}_9=b_{10}+b_9+c_8. \eqno(1)$$
By the associative law and 3) and 14) in {\bf Lemma
\ref{lemma4.3}}:
$$(b_3b_8)
c_3=(c_3b_8)b_3,\
b_{15}c_3+\overline{c}_9c_3=c_3b_3+b_6b_3+\overline{b}_{15}b_3,\
$$ by 5),\ 12) in  {\bf Lemma \ref{lemma4.2}} and (1) and 1) in
{\bf Lemma \ref{lemma4.1}}
$$c_3b_{15}=b_5+c_5+c_8+b_8+b_9+b_{10}.\eqno(2)$$
By the associative law and (1) and 1) {\bf Lemma \ref{lemma4.2}}:
$$(b_3\overline{c}_9)c_3=(c_3\overline{c}_9)b_3,\ \ \overline{b}_{15}c_3+c_9c_3+c_{3}^2=b_{10}
b_3+b_9b_3+c_8b_3,\ $$ by 1),\ 13),\ 16) in {\bf Lemma
\ref{lemma4.3}} and 3),\ 4) in {\bf Lemma \ref{lemma4.1}}
$$c_{3}^2=\overline{c}_3+\overline{b}_6.\eqno(3)$$
By the associative law and 2) in {\bf Lemma \ref{lemma4.1}} and 8)
in {\bf Lemma \ref{lemma4.2}}:
$$(b_3\overline{b}_6)b_{15}=(\overline{b}_6b_{15})b_3,\ $$
$$
\overline{b}_3b_{15}+\overline{b}_{15}b_{15}=4\overline{b}_{15}b_3+b_6b_3+\overline{b}_3b_3+c_3b_3+2c_9b_3.\
$$
By 6), 5), 12) in Lemma \ref{lemma4.2} and 1) in Lemma
\ref{lemma4.1} and Hypothesis \ref{hyp1} we obtain that
$$
\overline{b}_{15}b_{15}=1+3b_5+3c_5+5c_8+5b_8+6b_6+6b_{10}.\eqno(4)$$
By the associative law  and 2) in {\bf Lemma \ref{lemma4.1}} and
13) in {\bf Lemma \ref{lemma4.2}}:
$$
(b_3\overline{b}_6)\overline{b}_{15}=b_3(\overline{b}_6\overline{b}_{15}),\
$$
$$
\overline{b}_3\overline{b}_{15}+\overline{b}_{15}^2=2c_8b_3+3b_{10}b_3+2b_9b_3+2b_8b_3+b_5b_3+c_5b_3,\
$$ by 3),\ 4) in {\bf Lemma \ref{lemma4.1}} and 1) , 3),\ 4),\ 5) in {\bf
Lemma \ref{lemma4.3}}
$$
b_{15}^2=9\overline{b}_{15}+4b_6+2\overline{b}_3+3c_3+6c_9.\eqno(5)$$
By  4),\ 12),\ 13) in {\bf Lemma \ref{lemma4.2}} and 12) in {\bf
Lemma \ref{lemma4.3}} and 5) in {\bf Lemma \ref{lemma4.1}}:
$$(b_{10}b_{15},\ b_3)=(b_{10},\ \overline{b}_{15}b_{15})=6,\
(b_{10}b_{15},\ \overline{b}_{6})=(b_{15}b_6,\ b_{10})=3,\ $$
$$(b_{10}b_{15},\ b_3)=(\overline{b}_{15}b_3,\ b_{10})=1,\
(b_{10}b_{15},\ \overline{c}_3)=(b_{10}c_3,\ \bb_{15})=1,\
$$
$$ (b_{10}b_{15},\ \overline{b}_3c_3)=(b_{15}b_3,\ b_{10}c_3)=4,\ $$
so $ (b_{10}b_{15},\ \overline{c}_9)=4$.  hence
$$b_{10}b_{15}=6b_{15}+3\overline{b}_6+b_3+\overline{c}_3+4\overline{c}_9.\eqno(6)$$
By the associative law  and 4) in {\bf Lemma \ref{lemma4.2}} and
16) in {\bf Lemma \ref{lemma4.3}}:
$$(b_3\overline{c}_3)b_{15}=b_3(\overline{c}_3b_{15}),\
 c_9b_{15}=2\overline{b}_{15}b_3+b_3b_6+b_3c_9,\ $$ by 6),\ 12) in
 {\bf Lemma \ref{lemma4.2}} and 1) in {\bf Lemma \ref{lemma4.1}}:
 $$ c_9b_{15}=2b_5+2c_5+3c_8+3b_8+3b_9+4b_{10}.\eqno(7)$$
 By 4) and Lemma \ref{lemma4.2} and (2) we get:
 $$(b_3\overline{c}_3)\overline{b}_{15}=(\overline{c}_3\overline{b}_{15})b_3,\ $$
 $$c_9\overline{b}_{15}=b_5b_3+c_5b_3+c_8b_3+b_8b_3+b_9b_3+b_{10}b_3,\ $$
 by 1),\ 3),\ 4),\ 5) in {\bf Lemma \ref{lemma4.3}} and 3),\ 4) in {\bf Lemma \ref{lemma4.1}},\
 $$c_9\overline{b}_{15}=6b_{15}+2\overline{b}_6+b_3+\overline{c}_3+3\overline{c}_9.\eqno(8)$$
By the associative law  and 4) in {\bf Lemma \ref{lemma4.2}} and
(1):
$$(b_3\overline{c}_3)c_9=b_3(\overline{c}_3c_9),\ c_{9}^2=b_{10}b_3+b_9b_3+c_8b_3,\ $$
by 3),\ 4) in {\bf Lemma \ref{lemma4.1}} and 1) in {\bf Lemma
\ref{lemma4.3}}
$$c_{9}^2=3b_{15}+b_3+\overline{c}_3+2\overline{b}_6+2\overline{c}_9.\eqno(9)$$
By the associative law  and 4) in {\bf Lemma \ref{lemma4.2}} and
13) in {\bf Lemma \ref{lemma4.3}}:
$$(b_3\overline{c}_3)\overline{c}_9=b_3(\overline{c}_3\overline{c}_9),\
c_9\overline{c}_9=b_3\overline{b}_{15}+b_3c_9+b_3\overline{b}_3,\
$$ by 6),\ 12) in {\bf Lemma \ref{lemma4.2}} and Hypothesis
\ref{hyp1}:
$$c_9\overline{c}_9=1+2c_8+2b_9+2b_{10}+c_5+b_5+2b_8.\eqno(10)$$
By the associative law  and 4) in {\bf Lemma \ref{lemma4.2}} and
2) in {\bf Lemma \ref{lemma4.3}}:
$$ (b_3\overline{c}_3)c_8=( \overline{c}_3c_8)b_3,\
c_9c_8=\overline{c}_9b_3+b_{15}b_3,\ $$ by 1) in {\bf Lemma
\ref{lemma4.2}} and 5) in {\bf Lemma \ref{lemma4.1}}:
$$ c_9c_8=3\overline{b}_{15}+b_6+c_3+2c_9.\eqno(11)$$
By the associative law  and 5) in {\bf Lemma \ref{lemma4.2}} and
13) in {\bf Lemma \ref{lemma4.3}}:
$$ (
b_3c_3)c_9=(c_3c_9)b_3,\
b_9c_9=b_{15}b_3+\overline{c}_9b_3+b_{3}^2,\ $$ by 5) in {\bf
Lemma \ref{lemma4.1}} and Hypothesis \ref{hyp1} and 1) in lemma
\ref{lemma4.2}
$$ b_9c_9=\overline{b}_3+2c_9+3\overline{b}_{15}+c_3+2b_6.\eqno(12)$$
By the associative law  and 4) in {\bf Lemma \ref{lemma4.2}} and
12) in {\bf Lemma \ref{lemma4.3}}:
$$(b_3\overline{c}_3)b_{10}=(\overline{c}_3b_{10})b_3,\
c_9b_{10}=b_{15}b_3+\overline{b}_6b_3+\overline{c}_9b_3,\ $$ by 2)
and 5) in {\bf Lemma \ref{lemma4.1}} and 1) in {\bf Lemma
\ref{lemma4.2}}
$$ c_9
b_{10}=\overline{b}_3+2c_9+4\overline{b}_{15}+c_3+b_6.\eqno(13)$$
By (2),\
$$( \overline{c}_3b_5,\ b_{15})=(b_5,\ c_3b_{15})=1,\ (
\overline{c}_3c_5,\ b_{15})=(c_{5},\ c_3b_{15})=1.$$ So
$$ \overline{c}_3b_5=b_{15},\ \eqno(14)$$
$$ \overline{c}_3c_5=b_{15}.\eqno(15)$$
Hence by 4),\ 5) in {\bf Lemma \ref{lemma4.2}} and the associative
law:
$$ (
\overline{c}_3b_5)b_3=(b_3\overline{c}_5)b_5,\
(\overline{c}_3c_5)b_3=(b_3\overline{c}_3)c_5,\ $$
$$ b_{15}b_3=c_9b_5,\ b_{15}b_3=c_9c_5.\ $$
By 5) in Lemma \ref{lemma4.1} we get:
$$ b_5c_9=2\overline{b}_{15}+b_6+c_9,\ \eqno(16)$$
$$ c_5c_9=2\overline{b}_{15}+b_6+c_9. \eqno(17)$$
By the associative law  and 12) in {\bf Lemma \ref{lemma4.2}} and
5) in {\bf Lemma \ref{lemma4.3}}:
$$ (b_3\overline{b}_{15})b_5=(b_3b_5)\overline{b}_{15}=b_{15}
\overline{b}_{15},\ $$
$$
b_{10}b_5+c_8b_5+b_9b_5+b_8b_5+b_{5}^2+c_5b_5=b_{15}\overline{b}_{15},\
$$ so $b_5$ is real and so $c_5$ is real.  So by (14),\ (15)
$c_3b_5=\overline{b}_{15},\ \ c_3b_5=b_{15}$.
\end{proof}

\begin{lemma}\label{lemma4.5}Let $(A,\ B)$ satisfies Hypothesis \ref{hyp1}. Assume that $
\overline{b}_{10}=b_{10}\in B$ and $(\overline{c}_{9}b_{3},\
c_{9})= 1$.  Then we get the following equations:

1)  $b_5c_8=c_5+c_8+b_8+b_9+b_{10}$;

2)  $b_5b_8=c_5+c_8+b_8+b_9+b_{10}$;

3)  $b_5b_9=b_{10}+c_8+b_9+b_8+b_5+c_5$;

4)  $b_5b_{10}=c_8+b_8+b_9+b_5+2b_{10}$;

5)  $b_{5}^2=1+b_{10}+b_9+b_5$;

6)  $c_5b_5=b_8+c_8+b_9$;

7)  $c_8b_8=b_5+c_5+c_8+b_8+2b_9+2b_{10}$;

8)  $c_8b_9=2b_8+2b_9+2b_{10}+b_5+c_5+c_8$;

9)  $b_{9}^2=1+b_5+c_5+2c_8+2b_8+2b_9+2b_{10}$;

10)  $b_{8}^2=1+b_5+c_5+c_8+2b_8+b_9+2b_{10}$;

11)  $b_8b_9=b_5+c_5+2c_8+b_8+2b_9+2b_{10}$;

12)  $b_8b_{10}=b_5+c_5+2c_8+2b_8+2b_9+2b_{10}$;

13)  $c_5c_8=b_5+c_8+b_8+b_9+b_{10}$;

14)  $c_5b_8=b_5+c_8+b_8+b_9+b_{10}$;

15)  $c_5b_9=b_{10}+c_8+b_9+b_8+b_5+c_5$;

16)  $c_5b_{10}=b_8+c_8+b_9+c_5+2b_{10}$;

17)  $c_{5}^2=1+b_9+b_{10}+c_5$;

18)
$b_5b_{15}=3b_{15}+\overline{b}_6+b_3+\overline{c}_3+2\overline{c}_9$;

19)
$c_5b_{15}=3b_{15}+\overline{b}_6+b_3+\overline{c}_3+2\overline{c}_9$;

20)
$b_9b_{15}=6b_{15}+3\overline{c}_9+b_3+2\overline{b}_6+\overline{c}_3$;

21)  $b_9b_{10}=b_5+c_5+2c_8+2b_8+2b_9+3b_{10}$.

\end{lemma}

\begin{proof}By the associative law and Hypothesis \ref{hyp1} and 12) in
{\bf Lemma \ref{lemma4.2}}:
$$b_5(b_3\overline{b}_3)=( \overline{b}_3b_5)b_3=\bb_{15}b_3,\ $$
$$b_5+b_5c_8=\bb_{15}b_3=b_{10}+c_8+b_9+b_8+b_5+c_5.$$
Thus
$$c_8b_5=c_5+c_8+b_8+b_9+b_{10}.\eqno(1)$$
By the associative law and 9) in {\bf Lemma \ref{lemma4.2}} and
18) in {\bf Lemma \ref{lemma4.4}}:
$$(c_3\overline{c}_3)b_5=(\overline{c}_3b_5)c_3,\ \
b_5+b_5b_8=b_{15}c_3,\ $$ by 2) in lemma \ref{lemma4.4}
$$b_5b_8=c_5+c_8+b_8+b_9+b_{10}.\eqno(2)$$
By the associative law and 5),\ 12) in {\bf Lemma \ref{lemma4.2}}
and 18) in {\bf Lemma \ref{lemma4.4}}:
$$ (b_3c_3)b_5=(c_3b_5)b_3,\ $$
$$ b_9b_5=\overline{b}_{15}b_3=b_{10}+c_8+b_9+b_8+b_5+c_5.\eqno(3)$$
By the associative law and 1) in {\bf Lemma \ref{lemma4.1}} and 5)
in {\bf Lemma \ref{lemma4.3}}:
$$(b_3b_6)b_5=(b_3b_5)b_6,\ \ c_8b_5+b_{10}b_5=b_{15}b_6,\ $$
by (1) and 13) in {\bf Lemma \ref{lemma4.2}}:
$$b_5b_{10}=c_8+b_8+b_9+b_5+2b_{10}.\eqno(4)$$
By the associative law and 12) in {\bf Lemma \ref{lemma4.2}} and
5) in {\bf Lemma \ref{lemma4.3}}:
$$(b_3\overline{b}_{15})b_5=(b_3b_5)\overline{b}_{15},\ $$
$$b_{10}b_5+c_8b_5+b_9b_5+b_8b_5+b_{5}^2+c_5b_5=b_{15}\overline{b}_{15},\ $$
so by 4) in {\bf Lemma \ref{lemma4.4}} and (1),\ (2),\ (3),\ (4),\
$$b_{5}^2+c_5b_5=1+b_5+2b_9+b_{10}+b_8+c_8,\ $$
by (2),\ (4),\ (3),\ (1),\
$$(b_{5}^2,\ b_{10})=(b_5,\ b_5b_{10})=1,\ \
(b_{5}^2,\ b_9)=(b_5,\ b_5b_9)=1,\ $$ so
$$b_{5}^2=1+b_{10}+b_9+b_5,\ \eqno(5)$$
$$ c_5b_5=b_8+c_8+b_9.\eqno(6)$$
By the associative law and Hypothesis \ref{hyp1} and 3) in {\bf
Lemma \ref{lemma4.3}}:
$$(b_3\overline{b}_3)b_8=( \overline{b}_3b_8)b_3,\ \
b_8+c_8b_8=\overline{b}_{15}b_3+c_9b_3,\ $$
 6),\ 12) in {\bf Lemma \ref{lemma4.2}} and
 $$c_8b_8=b_5+c_5+c_8+b_8+2b_9+2b_{10}.\eqno(7)$$
By the associative law and Hypothesis \ref{hyp1} and 1) in {\bf
Lemma \ref{lemma4.3}}:
$$(b_3\overline{b}_3)b_9=( \overline{b}_3b_9)b_3,\ \
b_9+c_8b_9=\overline{b}_{15}b_3+c_9b_3+c_3b_3.$$ 5),\ 6),\ 12) in
lemma \ref{lemma4.2}
$$b_9+c_8b_9=b_{10}+c_8+b_9+b_8+b_5+c_5+b_9+b_{10}+b_8+b_9,\ $$
By 5) in Lemma \ref{lemma4.2} and 15) in Lemma \ref{lemma4.3} we
obtain
$$c_8b_9=2b_8+2b_9+2b_{10}+b_5+c_5+c_8.\eqno(8)$$
$$(b_3c_3)b_9=(c_3b_9)b_3,\ \
b_{9}^2=\overline{b}_{15}b_3+\overline{b}_3b_3+c_9b_3,\ $$ by 12)
in {\bf Lemma \ref{lemma4.2}} and Hypothesis \ref{hyp1} and 6) in
{\bf Lemma \ref{lemma4.2}},\
$$b_{9}^2=1+b_5+c_5+2c_8+2b_8+2b_9+2b_{10}.\eqno(9)$$
By the associative law  and 9) in {\bf Lemma \ref{lemma4.2}} and
14) in {\bf Lemma \ref{lemma4.3}}:
$$(c_3\overline{c}_3)b_8=c_3( \overline{c}_3b_8),\ \
b_8+b_{8}^2=c_3\overline{c}_3+c_3\overline{b}_6+c_3b_{15},\ $$ by
7) ,\ 9) in {\bf Lemma \ref{lemma4.2}} and 2) in {\bf Lemma
\ref{lemma4.4}}:
$$b_{8}^2=1+b_5+c_5+c_8+2b_8+b_9+2b_{10}.\eqno(10)$$
By the associative law and 9) in {\bf Lemma \ref{lemma4.2}} and
15) in {\bf Lemma \ref{lemma4.3}}:
$$(c_3\overline{c}_3)b_9=(\overline{c}_3b_9)c_3,\ \
b_9+b_8b_9=b_{15}c_3+b_3c_3+\overline{c}_9c_3,\ $$ by 1),\ 2) in
{\bf Lemma \ref{lemma4.4}} and 5) in {\bf Lemma \ref{lemma4.2}}:
$$ b_9+b_8b_9=b_5+c_5+c_8+b_8+b_9+b_{10}+b_9+b_{10}+b_9+c_8,\ $$
$$b_8b_9=b_5+c_5+2c_8+b_8+2b_9+2b_{10}.\eqno(11)$$
By the associative law and 9) in {\bf Lemma \ref{lemma4.2}},\ 12)
in {\bf Lemma \ref{lemma4.3}}:
$$(c_3\overline{c}_3)b_{10}=( \overline{c}_3b_{10})c_3,\ \
b_{10}+b_8b_{10}=b_{15}c_3+\overline{b}_6c_3+\overline{c}_9c_3,\
$$ so by 1),\ 2) in {\bf Lemma \ref{lemma4.4}} and 7) in {\bf
Lemma \ref{lemma4.2}},\
$$b_8b_{10}=b_5+c_5+2c_8+2b_8+2b_9+2b_{10}.\eqno(12)$$
By the associative law and Hypothesis \ref{hyp1} and 4) in {\bf
Lemma \ref{lemma4.3}} 12) in {\bf Lemma \ref{lemma4.2}}:
$$ (b_3\overline{b}_3)c_5=(\overline{b}_3c_5)b_3,\ $$
$$ c_5+c_5c_8=\overline{b}_{15}b_3=b_{10}+c_8+b_9+b_8+b_5+c_5,\ $$
$$c_5c_8=b_{10}+c_8+b_9+b_8+b_5.\eqno(13)$$
By the associative law and 9) in {\bf Lemma \ref{lemma4.2}} and
2),\ 18) in {\bf Lemma \ref{lemma4.4}}:
$$(c_3\overline{c}_3)c_5=( \overline{c}_3c_5)c_3,\ $$
$$c_5+c_5b_8=b_{15}c_3=b_5+c_5+c_8+b_8+b_9+b_{10},\ $$
$$c_5b_8=b_5+c_8+b_8+b_9+b_{10}.\eqno(14)$$
By the associative law and 5),\ 12) in {\bf Lemma
\ref{lemma4.2}},\ 18) in {\bf Lemma \ref{lemma4.4}} :
$$(b_3c_3)c_5=(c_3c_5)b_3,\ $$
$$b_9c_5=\overline{b}_{15}b_3=b_{10}+c_8+b_9+b_8+b_5+c_5.\eqno(15)$$
By the associative law and 1) in {\bf Lemma \ref{lemma4.1}} and
10) in {\bf Lemma \ref{lemma4.3}}
$$(b_3b_6)c_5=(b_6c_5)b_3,\ $$
$$c_8c_5+b_{10}c_5=\overline{b}_{15}b_3+b_3b_6+b_3c_9.\ $$
By 6), 12) in Lemma \ref{lemma4.2} and (13) and 1) in Lemma
\ref{lemma4.1} we get:
$$b_{10}+c_8+b_9+b_8+b_5+b_{10}c_5=b_{10}+c_8+b_9+b_8+b_5+c_5+c_8+b_{10}+b_9+b_{10}+b,\ $$
$$b_{10}c_5=c_5+c_8+2b_{10}+b_9+b_8.\eqno(16)$$
By the associative law and 10) in {\bf Lemma \ref{lemma4.2}} ,\
10) in {\bf Lemma \ref{lemma4.3}}:
$$(b_6\overline{b}_6)c_5=\overline{b}_6(b_6c_5),\ $$
$$c_5+c_5c_8+b_9c_5+b_8c_5+c_{5}^2+b_5c_5=\overline{b}_6\overline{b}_{15}+\overline{b}_6b_6,\ $$
by (13),\ (15),\ (14),\ (6) and 10) ,\ 13) in {\bf Lemma
\ref{lemma4.2}}
$$c_{5}^2=1+b_9+b_{10}+c_5.\eqno(17)$$
By the associative law and 18) in {\bf Lemma \ref{lemma4.4}} and
(5):
$$c_3b_{5}^2=(c_3b_5)b_5,\ $$
$$c_3+b_{10}c_3+b_9c_3+b_5c_3=\overline{b}_{15}b_5,\ $$
so by 12),\ 15) in {\bf Lemma \ref{lemma4.3}} and 17) in {\bf
Lemma \ref{lemma4.4}},\
$$b_5\overline{b}_{15}=3\overline{b}_{15}+b_6+\overline{b}_3+c_3+2c_9.\eqno(18)$$
By the associative law and 18) in {\bf Lemma \ref{lemma4.4}} and
(17):
$$(c_{5}^2)c_3=(c_5c_3)c_5,\ \
c_3+b_9c_3+b_{10}c_3+c_5c_3=\overline{b}_{15}c_5,\ $$ so by 12),\
15) in {\bf Lemma \ref{lemma4.3}} and 18) in {\bf Lemma
\ref{lemma4.4}}:
$$\overline{b}_{15}c_5=3\overline{b}_{15}+b_6+\overline{b}_3+c_3+2c_9.\eqno(19)$$
By the associative law and 4) in {\bf Lemma \ref{lemma4.3}}and
(15)
$$(b_3c_5)b_9=b_3(c_5b_9),\ $$
$$b_{15}b_9=b_3b_8+c_8b_3+b_5b_3+c_5b_3+b_9b_3+b_{10}b_3,\ $$
by 3),\ 4) in {\bf Lemma \ref{lemma4.1}} and 1),\ 3),\ 4),\ 5) in
{\bf Lemma \ref{lemma4.3}}:
$$b_{15}b_9=b_{15}+\overline{c}_9+b_3+\overline{b}_6+b_{15}+2b_{15}+2b_{15}+2\overline{c}_9+\overline{c}_3+\overline{b}_6,\ $$
$$b_{15}b_9=6b_{15}+3\overline{c}_9+b_3+2\overline{b}_6+\overline{c}_3.$$
By the associative law and 5) in {\bf Lemma \ref{lemma4.2}} and 4)
in {\bf Lemma \ref{lemma4.1}}:
$$b_{10}(b_3c_3)=(b_{10}b_3)c_3,\ $$
$$b_{10}b_9=b_{15}c_3+\overline{b}_6c_3+\overline{c}_9c_3,\ $$
$$b_{10}b_9=b_5+c_5+c_8+b_8+b_9+b_{10}+b_8+c_8+b_9+b_{10},\ $$
$$b_{10}b_9=b_5+c_5+2c_8+2b_8+2b_9+3b_{10},\ $$

\end{proof}
So we finished the investigation of 11) in Lemma \ref{lemma4.1}
subcase $(\bar{c}_9b_3,c_9)=1$. Thus we have
 in the remaining case that $(b_3\overline{c}_{9},\ c_{9})=0$.  In the
rest of this section let us deal with the investigation of 11)
Lemma \ref{lemma4.1} subcase $(b_6c_9,\bar{c}_9)\ge 1$.\par

\begin{lemma}\label{lemma4.6}Let $(A,\ B)$ satisfies Hypothesis \ref{hyp1}. Assume that $
\overline{b}_{10}=b_{10}\in B$  and $(b_6c_9,\ \overline{c}_9)\geq
1$. Then we get:

1)  $b_3c_9=b_{10}+\overline{x}+s,\ s\in B,\ \mid
\overline{x}+s\mid=17,\ \mid \overline{s}\mid=\mid s\mid$;

2)  $(c_8b_{15},\ \overline{c}_9)=(c_8c_9,\ \overline{b}_{15})=2$;

3)  $
\overline{b}_6\overline{c}_9=3c_9+\overline{b}_{15}+b_6+t_6,\
t_6\in N^*B$;

4)  $c_8b_{10}=2c_8+2b_{10}+2x+y+z+w+s$;

5)  $b_6b_{15}=3b_{10}+2c_8+x+y+z+w+\overline{x}+s$;

6)  $\overline{b}_3 c_9=b_{15}+\varepsilon+\theta,\ \varepsilon,\
\theta \in B,\ \mid \varepsilon+\theta\mid=12$;

7)
$b_6b_{10}=\overline{b}_3+3\overline{b}_{15}+\overline{\varepsilon}+\overline{\theta}$;

8)
$c_8b_{15}=5b_{15}+2\overline{b}_6+2\overline{c}_9+b_3+\varepsilon+\theta$.
 \end{lemma}

 \begin{proof}By the associative law and Hypothesis \ref{hyp1} and 10) in
 {\bf Lemma \ref{lemma4.1}}:
 $$(b_{3}^2)\overline{b}_{15}=(b_3\overline{b}_{15})b_3,\ $$
 $$
 \overline{b}_3\overline{b}_{15}+b_6\overline{b}_{15}=b_{10}b_3+c_8b_3+b_3(x+y+z+w).$$
 So by 5),\ 4),\ 3) in {\bf Lemma \ref{lemma4.1}},\
 $$ b_6\overline{b}_{15}=b_3+\overline{b}_6+b_3(x+y+z+w).\eqno(1)$$
By the associative law and Hypothesis \ref{hyp1} and 4) in
 {\bf Lemma \ref{lemma4.1}}:
 $$(b_{3}^2)b_{10}=(b_3b_{10})b_3,\ $$
 $$ \overline{b}_3b_{10}+b_6b_{10}=b_{15}b_3+\overline{b}_{6}b_3
+\overline{c}_9b_3,\ $$ by 4),\ 5),\ 2) in {\bf Lemma
\ref{lemma4.1}}:
$$
b_6b_{10}=\overline{b}_3+2\overline{b}_{15}+b_3\overline{c}_9.\eqno(2)$$
By the associative law and(2) and 4) in
 {\bf Lemma \ref{lemma4.1}}:
 $$ (b_3b_{10})\overline{b}_6=(\overline{b}_6b_{10}
)b_3,\ $$
$$
\overline{b}_6^2+b_{15}\overline{b}_6+\overline{b}_6\overline{c}_9=b_{3}^2+2b_{15}
b_3+( \overline{b}_3b_3) c_9.$$ so by 5),\ 6) in {\bf Lemma
\ref{lemma4.1}} and Hypothesis \ref{hyp1}:
$$2b_6+c_9+\overline{b}_{15}+\overline{b}_3+b_6+\overline{b}_3(x+y+z+w)+\overline{b}_6\overline{c}_9=\overline{b}_3+b_6
+4\overline{b}_{15}+2b_6+2c_9+c_9+c_8c_9,\ $$
$$
\overline{b}_6\overline{c}_9+\overline{b}_3(x+y+z+w)=3\overline{b}_{15}+2c_9+c_8c_9.\eqno(3)$$
By the associative law and 6),\ 7) in {\bf Lemma \ref{lemma4.1}}:
$$(b_{6}^2)c_8=(b_6c_8)b_6,\ $$
$$2\overline{b}_6c_8+\overline{c}_9c_8+b_{15}
c_8=\overline{b_3}b_6+2\overline{b}_{15}b_6+b_{6}^2+c_9b_6,\ $$ by
7),\ 2),\ 6) in {\bf Lemma \ref{lemma4.1}} and (3),\ (1)
$$c_8b_{15}=b_3+2\overline{b}_6+\overline{c}_9+b_{15}+b_3(x+y+z+w).\eqno(4)$$
By the associative law and Hypothesis \ref{hyp1} :
$$ (b_3\overline{b}_3)c_9=(\overline{b}_3c_9)b_3,\ \ c_9+c_9c_8=(
\overline{b}_3c_9)b_3,\ $$ but by 5) in {\bf Lemma
\ref{lemma4.1}},\
$$ (\overline{b}_3c_9,\ b_{15}
)=(c_9,\ b_3b_{15})=1,\ $$ so
$$ c_9+c_9c_8=b_{15}b_3+\alpha_{1}b_3,$$
hence by 5) in Lemma \ref{lemma4.1} we obtain that:
$$c_9c_8=2\overline{b}_{15}+b_6+\alpha_{1}b_3.\eqno(5)$$ Hence
$$ (c_8c_9,\ \overline{b}_{15})\geq 2.\eqno(6)$$
By 10) in {\bf Lemma \ref{lemma4.1}},\ $$ ( \overline{b}_3x,\
\overline{b}_{15})=1,\ \  ( \overline{b}_3y,\ \overline{b}_{15})=1,\
( \overline{b}_3z,\ \overline{b}_{15})=1,\ ( \overline{b}_3w,\
\overline{b}_{15})=1.$$ So by (3),\ $$ (
\overline{b}_6\overline{c}_9,\ \overline{b}_{15})\geq 1.\eqno(7)$$
By (6),\ $ (c_8b_{15} ,\ \overline{c}_9)\geq 2$,\  so by (4),\ $$
(b_3(x+y+z+w),\ \overline{c}_9)\neq 0,\ $$ hence one of $b_3x$,
$b_3y$, $b_3z$, and $b_3w$ contains $\bar{c}_9$, thus W. L. O. G, we
may assume that $$ (b_3x,\ \bar{c}_9)\neq 0.\eqno(8)$$ By the
associative law and Hypothesis \ref{hyp1} and 4) in {\bf Lemma
\ref{lemma4.1}}:
$$(b_3\overline{b}_3)b_{10}=(b_3b_{10})\overline{b}_3,\ $$
$$
b_{10}+c_8b_{10}=b_{15}\overline{b}_3+\overline{b}_6\overline{b}_3+\overline{c}_9\overline{b}_3,\
$$ by 10),\ 1) in {\bf Lemma \ref{lemma4.1}}
$$
c_8b_{10}=2c_8+b_{10}+x+y+z+w+\overline{b}_3\overline{c}_9.\eqno(9)$$
By (8) and 4) in {\bf Lemma \ref{lemma4.1}},\ $ (b_3c_9,\
b_{10})=(b_3b_{10},\ \overline{c}_9)=1,\ \ (b_3c_9,\
\overline{x})\neq 0$. So
$$b_3c_9=b_{10}+\overline{x}+\alpha,\ \alpha \in N^*B.\eqno(10)$$
Then $ (b_3x,\ \overline{c}_9)=( b_3c_9,\ \overline{x}) \neq 0$. By
(9) since $ c_8,\ b_{10},\ x+y+z+w$ are reals then $ b_3c_9$ is
real. By 10) in \textbf{Lemma 5.1}, $ (b_3x,\ b_{15})=1$,\ therefore
$$ b_3x=\overline{c}_9+b_{15}+\beta ,\ \eqno(11)$$
where $\beta \in N^*B$.  So $ \mid x\mid \geq 8$. Hence by (10)
$$ 3\leq (b_3c_9,\ b_3c_9)\leq 5.$$
Thus $$3\leq (b_3\overline{b}_3,\ c_9\overline{c}_9)\leq 5.$$ So
by Hypothesis $ 2\leq (c_9\overline{c}_9,\ c_8)\leq 4.$   Thus
$$ 2\leq (c_9,\ c_9c_8)\leq 4.\eqno(12)$$
\par If $ (c_8b_{15},\ \overline{c}_9)\geq 4$  then by (4), we have
that $(b_3(x+y+z+w),\ \overline{c}_9)\geq 3$ .  Since $1+(c_8,\
b_{15}\bb_{15})=(1+c_8,\ b_{15}\bb_{15})=(b_3\bb_3,\
b_{15}\bb_{15})=(b_3b_{15},\ b_3b_{15})=6$ by 10) in\textbf{ Lemma
5.1}, we have that $(c_8b_{15},\ b_{15})=5$. Hence we get
 without lose of generality by (11) that
$$ b_3y=\overline{c}_9+b_{15}+\beta_1,\  b_3z=\overline{c}_9+b_{15}+\beta_2.$$ So the degrees of $x$, $y$ and $z$ are all
larger or equal to 8, but $|x+y+z+w|=27$, hence $\mid x\mid=\mid
y\mid=|z|=8$ by (11) and $ (b_3c_9,\ x)=(b_3c_9,\ y)=(b_3c_9,\
z)=1$, a contradiction to (10). If $ (c_8b_{15},\ \bc_9)=3$. Then by
(5) and (12)
$$c_8c_9=3\overline{b}_{15}+b_6+2c_9+c_3.\eqno(13)$$
where $ c_3\in B$.  By (4) ,\ $(b_3(x+y+z+w),\ \overline{c}_9)=2$,\
so by (3) $( \overline{b}_6\overline{c}_9,\ c_9)=2$ and by 10),\  6)
in {\bf Lemma \ref{lemma4.1}},\ $ ( \overline{b}_6\overline{c}_9,\
b_6)=1,\ ( \overline{b}_3(x+y+z+w),\ \overline{b}_{15})=4$,\   so by
(3),\ (13) $(\overline{b}_6 \overline{c}_9,\ \overline{b}_{15})=2$.
Hence
$$ b_6c_9=\overline{b}_6+2\overline{c}_9+2b_{15}.\eqno(14)$$ So by
(3),\
$$ \overline{b}_3(x+y+z+w)=4\overline{b}_{15}+2c_9+c_3.$$
Then we get without lose of generality the only possibility we set
$$ x=b_8,\ y=b_9,\ z=b_5,\ w=c_5,\ $$
when $ b_8,\ b_9,\ b_5,\ c_5$ are reals. So
$$ \overline{b}_3b_8=\overline{b}_{15}+c_9,\ \eqno(15)$$
$$ \overline{b}_3b_9=\overline{b}_{15}+c_9+c_3,\ \eqno(16)$$
$$ \overline{b}_3b_{5}=\overline{b}_{15},\ \eqno(17)$$
$$ \overline{b}_3c_5=\overline{b}_{15},\ \eqno(18)$$
So $ (b_3c_9,\ b_8)=1$,\   and by 4) in {\bf Lemma
\ref{lemma4.1}},\ $(b_3c_9,\ b_{10})=1$. Thus $$
b_3c_9=b_{10}+b_8+b_9.\eqno(19)$$ So by the associative law and
Hypothesis \ref{hyp1}
$$(b_{3}^2)c_9=(b_3c_9)b_3,\ $$
$$c_9\overline{b}_3+c_9b_6=b_{10}b_3+b_8b_3+b_9b_3,\ $$
so by (15),\ (16),\ (14) $
c_9\overline{b}_3=b_{15}+\overline{c}_9+\overline{c}_3$  which
contradicts the assumption in the lemma ,\   we conclude that $
(c_8b_{15},\ \overline{c}_9)\leq 2$ and by (8) and (4),\
$$(c_8b_{15},\ \overline{c}_9)=2.\eqno(20)$$
By (4) $$ (b_3(x+y+z+w),\ c_9)=1,\ \eqno(21)$$ and by (3),\ (12)
we get $ 3\leq (\overline{b}_6 \overline{c}_9,\ c_9)\leq 5$,\ by
(20) $ (c_8c_9,\ \overline{b}_{15})=2$.  So by (3) and that $ (
\overline{b}_3(x+y+z+w),\ \overline{b}_{15})=4$ we get $(
\overline{b}_6 \overline{c}_9,\ \overline{b}_{15})=1$.  And by 6)
in {\bf Lemma \ref{lemma4.1}},\ $(\overline{b}_6 \overline{c}_9,\
{b}_6)=1$.  So $( \overline{b}_6 \overline{c}_9,\ c_9)=3$. Then by
(3),\ (21),\ $(c_8c_9,\ c_9)=2$.  So $(c_8,\
\overline{c}_9c_9)=2$,\  and by Hypothesis \ref{hyp1} $
(\overline{c}_9c_9,\ \overline{b}_3b_3)=3$.  Thus $(b_3c_9,\
b_3c_9)=3$.  So by 4) in {\bf Lemma \ref{lemma4.1}} and (10) we
get that
$$b_3c_9=b_{10}+\overline{x}+s,\ \eqno(22)$$
$s\in B$.  Also $$\overline{b}_6
\overline{c}_9=3c_9+\overline{b}_{15}+b_6+t_6,\ \eqno(23)$$ where
$ t_6\in N^*B,\ t_6\neq b_6,\ \overline{b}_6$.  By (9) and (22),\
$$ c_8b_{10}=2c_8+2b_{10}+2x+y+z+w+\overline{S}.\eqno(24)$$
By the associative law and 5) in {\bf Lemma \ref{lemma4.1}}:
$$(b_{3}^2)b_{15}=(b_3b_{15})b_3,\ $$
$$
\overline{b}_3b_{15}+b_6b_{15}=2\overline{b}_{15}b_3+b_6b_3+c_9b_3,\
$$ by (22) and 1),\ 10) in {\bf Lemma \ref{lemma4.1}}
$$b_6b_{15}=3b_{10}+2c_8+x+y+z+w+\overline{x}+S.\eqno(25)$$
By the associative law and Hypothesis \ref{hyp1} and 5) in {\bf
Lemma \ref{lemma4.1}}: \par Since
$b_{15}(b_3\bb_3)=\bb_3(b_{15}b_3),$ then
$$c_8b_{15}=4b_{15}+2\overline{b}_6+2\overline{c}_9+b_3+\overline{b}_3c_9.\eqno(26)$$
By (22) and 5) in {\bf Lemma \ref{lemma4.1}}
$$ \overline{b}_3c_9=b_{15}+\varepsilon+\theta,\ \eqno(27)$$,\ where
$\varepsilon,\ \theta \in B,\ \mid \varepsilon +\theta \mid=12$,\
so by (26)
$$c_8b_{15}=5b_{15}+2\overline{b}_6+2\overline{c}_9+b_3+\varepsilon+\theta.$$
By the associative law and Hypothesis \ref{hyp1} and 4) in {\bf
Lemma \ref{lemma4.1}}:
$$(b_{3}^2)b_{10}=b_3(b_3b_{10}),\ $$
$$b_{10}\overline{b}_3+b_{10}b_6=b_{15}b_3+\overline{b}_6b_3+\overline{c}_9b_3,\ $$
so by (27) and 5),\ 2),\ 4) in {\bf Lemma \ref{lemma4.1}}:
$$b_6b_{10}=\overline{b}_3+3\overline{b}_{15}+\overline{\varepsilon}+\overline{\theta}.$$\end{proof}

\begin{lemma}\label{lemma4.7}Let $(A,\ B)$ satisfies Hypothesis \ref{hyp1}. Assume that $
\overline{b}_{10}=b_{10}\in B$ and $(b_3\overline{c}_{9},\
c_{9})=0$. Then $(b_6c_9,\ \overline{c}_9)\geq 1$ is impossible.
\end{lemma}

\begin{proof}We assume that our assumption is possible.  By
Hypothesis \ref{hyp1} and 10) in {\bf Lemma \ref{lemma4.1}} and
the associative law:
$$ (b_3\overline{b}_3)b_{15}=b_3(\overline{b}_3b_{15}),\ $$
$$
b_{15}+c_8b_{15}=b_3b_{10}+b_3c_8+b_3\overline{x}+b_3\overline{y}+b_3\overline{z}+b_3\overline{w},\
$$ by 3),\ 4) in {\bf Lemma \ref{lemma4.1}},\ 8) in {\bf Lemma
\ref{lemma4.6}}
$$b_{15}+5b_{15}+2\overline{b}_6+2\overline{c}_9+b_3+\varepsilon+\theta=b_{15}+\overline{b}_6+\overline{c}_9+b_3+
\overline{b}_6+b_{15}+b_3\overline{x}+b_3\overline{y}+b_3\overline{z}+b_3\overline{w},\
$$
$$
4b_{15}+\overline{c}_9+\varepsilon+\theta=b_3\overline{x}+b_3\overline{y}+b_3\overline{z}+b_3\overline{w}.\eqno(1)$$
by 10) in {\bf Lemma \ref{lemma4.1}},\ and 1) in {\bf Lemma
\ref{lemma4.6}}:
$$ (b_3\overline{x},\ b_{15})=(\overline{b}_{15}b_3,\ x)=1,\ \
(b_3\overline{y},\ b_{15})=(b_3\overline{b}_{15},\ y)=1,\ $$
$$ (b_3\overline{z},\ b_{15})=(\overline{b}_{15}b_3,\ z)=1,\ \
(b_3\overline{w},\ b_{15})=(b_3\overline{b}_{15},\ w)=1,\ $$
$(b_3\overline{x},\ \overline{c}_9)=(b_3c_9,\ x)=1$,\   so
$x=\overline{x}$.  If otherwise two of $ x,\ y,\ z,\ w$ have
degrees $\geq 8$,\   impossible.Then $$ b_3x=\overline{c}_9+b_{15}
,\ \eqno(2)$$ then $ x=x_8$.  So $ \mid
\overline{y}+\overline{z}+\overline{w}\mid=19$,\  and also we get
$ 5\leq \mid \overline{y}\mid ,\ \mid \overline{z}\mid,\  \mid
\overline{w}\mid \leq 9$  and also without lose of generality
$\mid \overline{y}\mid=5$,\   hence $
\overline{y}=\overline{y}_5$,\  so
$$ b_3\overline{y_5}=b_{15}.\eqno(3)$$ So $$ \mid
\overline{z}+\overline{w}\mid=14.\eqno(4)$$ By 1) in {\bf Lemma
\ref{lemma4.6}} and that $ x=x_8$ we get that
$$ b_3c_9=b_{10}+x_8+s_9.\eqno(5)$$
So by the associative law and Hypothesis \ref{hyp1} and (5):
$$ (b_{3}^2)c_9=b_3(b_3c_9),\ $$
$$ c_9\overline{b}_3+b_6c_9=b_3b_{10}+b_3x_8+s_9b_3,\ $$
then by (2) and (6) ,\ 3) in {\bf Lemma \ref{lemma4.6}},\ 4) in
{\bf Lemma \ref{lemma4.1}}:
$$b_{15}+\varepsilon+\theta+3\overline{c}_9+b_{15}+\overline{b}_6+\overline{t}_6=
b_{15}+\overline{b}_6+\overline{c}_9+b_{15}+\overline{c}_9+s_9b_3,\
$$
$$ s_9b_3=\varepsilon+\theta+\overline{c}_9+\overline{t}_6.\eqno(6)$$
If $ \mid z\mid=5$,\   then by (4) $ \mid w\mid=9$,\   hence $
z=\overline{z}_5,\ w=\overline{{w}_9}$.  So
$$b_3z_5=b_{15},\ $$
$$b_3{w}_9=b_{15}+\varepsilon+\theta.\eqno(7)$$
So by (6)
$$(\overline{b}_3\varepsilon,\ s_9)=(\varepsilon,\ b_3s_9)=1,\
(\overline{b}_3\theta,\ s_9)=(\theta,\ b_3s_9)=1,\ $$
$$(\overline{b}_3\varepsilon,\ {{w}_9})=(\varepsilon,\ b_3{w}_9)=1,\
(\overline{b}_3\theta,\ {w}_9)=(\theta,\ b_3{w}_9)=1.$$ hence
$$ \overline{b}_3\varepsilon =s_9+{w}_9,\ \eqno(8)$$
$$\overline{b}_3\theta =s_9+{w}_9,\   \eqno(9)$$
hence $ \varepsilon=\varepsilon_6,\ \theta=\theta_6$,\   by 4) in
the {\bf Lemma \ref{lemma4.6}} $ s_9=\overline{s}_9$.  So by the
associative law and Hypothesis \ref{hyp1}:
$$
\overline{\varepsilon}_6(b_{3}^2)=b_3(\overline{\varepsilon}_6b_3),\
$$
$$
\overline{\varepsilon}_6\overline{b}_3+\overline{\varepsilon}_6b_6=s_9b_3+{w}_9b_3,\
$$ so by (6),\ (7)
$$\overline{\varepsilon}_6\overline{b}_3+\overline{\varepsilon}_6b_6=\varepsilon_6+\theta_6+\overline{c}_9+\overline{t}_6
+b_{15}+\varepsilon_6+\theta_6,\ $$
$$
\overline{\varepsilon}_6\overline{b}_3+\overline{\varepsilon}_6b_6=2\varepsilon_6+2\theta_6+\overline{c}_9+\overline{t}_6
+b_{15},\ $$ by 6) in {\bf Lemma \ref{lemma4.6}}
$(\overline{\varepsilon}_6\overline{b}_3,\
\overline{c}_9)=(\varepsilon_6,\ \overline{b}_3c_9)=1$ and
$(\overline{\varepsilon}_6\overline{b}_3,\ \delta)=1$,\   where $
\delta=\varepsilon_6$ or $ \theta_6$,\   so $
\overline{\varepsilon}_6\overline{b}_3=\delta+\overline{c}_9+m_3$,\
so $ \overline{t}_6=m_3+k_3$,\   and $(b_3m_3,\
\overline{\varepsilon}_6)=(\overline{b}_3\overline{\varepsilon}_6,\
m_3)=1$. Therefor
$$ (b_3m_3,\ b_3m_3)=(b_3\overline{b}_3,\ m_3\overline{m}_3)>1.$$
In the (6) we get $$ (\overline{b}_3m_3,\ s_9)=(m_3,\ b_3s_9)=1,\
$$ then
$$(\overline{b}_3m_3,\ \overline{b}_3m_3)=(\overline{b}_3b_3,\ \overline{m}_3m_3)=1,\ $$
a contradiction.

If $ \mid \overline{z}\mid=6$ then by (4) $ \mid
\overline{w}\mid=8$,\   hence $z=\overline{z},\ w=\overline{w}$.
So
$$b_3z_6=b_{15}+\varepsilon_3,\ b_3w_8=b_{15}+\theta_9.$$
So $ ( \varepsilon_3\overline{b}_3,\ z_6)=( \varepsilon_3,\
b_3z_6)=1$,\  hence $ (\varepsilon_3\overline{b}_3,\
\varepsilon_3\overline{b}_3)>1$.   By (6) ,\ $(s_9b_3,\
\varepsilon_3)=(s_9,\ \overline{b}_3z_6)=1$,\   hence $
\overline{b}_3\varepsilon_3=s_9$,\   so
$(\varepsilon_3\overline{b}_3,\ \varepsilon_3\overline{b}_3)=1$  a
contradiction.  If $ \mid \overline{z}\mid=7$ then by (4) $ \mid
\overline{w}\mid=7$,\   so $$ z=\overline{z},\ w=\overline{w} \
or\  \overline{z}=w,\ $$
$$ b_3\overline{z}_7=b_{15}+\varepsilon_6,\ \eqno(10)$$
$$ b_3\overline{w}_7=b_{15}+\theta_6.\eqno(11)$$
So $(\overline{b}_3\varepsilon_6,\ \overline{z}_7)=(\varepsilon_6,\
b_3\overline{z}_7)=1$. By (6) $(\overline{b}_3\varepsilon_6,\
s_9)=(\varepsilon_6,\ b_3s_9)=1$,\ hence
$$ \overline{b}_3\varepsilon_6=\overline{z}_7+s_9+m_2.$$
a contradiction.
\end{proof}
Consequently 11) in Lemma \ref{lemma4.1} subcase
$(b_6c_9,\bar{c}_9)\ge 1$ is impossible and we have only the case
$(b_3\bar{c}_9, c_9)=1$. We can state now Theorem
\ref{theorem4.1}.
\begin{theorem} \label{theorem4.1}Let $(A,\ B)$  be a $NITA$ generated by an non-real
element $b_{3}\in B$ satisfies $b_{3}\overline{b}_{3}=1+c_{8}$ and
$b_{3}^{2}=\overline{b}_{3}+b_{6}$,\  where $b_{6}\in B$ is non-real
and satisfying Hypothesis \ref{hyp1}. Then $b_3b_6=c_8+b_{10} $ and
if $b_{10}=\bar{b}_{10}$ then $(A,\ B)\cong(CH(3A_6),\ Irr(3A_6))$.
$(A,\ B)$ is table algebra of dimension 17: $B=\{1,\ b_3,\
\bar{b_3},\ c_3,\ \bar{c_3},\ b_6,\ \bar{b}_6,\ c_9,\ \bar{c}_9,\
b_{15},\ \bar{b}_{15},\ b_5,\ c_5,\ b_8,\ c_8,\ b_9,$ $\ b_{10}\}$
and $B$ has an increasing series of table subsets
$\{b_1\}\subseteq\{b_1,\ b_5,\ c_5,\ b_8,\ c_8,\ b_9,$ $\
b_{10}\}\subseteq B$ defined by

1)  $b_3\bar{b_3}=1+c_8$;

2)  $b_{3}^2=\overline{b_3}+b_6$;

3)  $b_3c_3=b_9$;

4)  $b_3\bar{c_3}=c_9$;

5)  $b_3b_6=c_8+b_{10}$;

6)  $b_3\overline{b}_6=\bar{b_3}+\bar{b}_{15}$;

7)  $b_3c_9=b_9+b_{10}+b_8$;

8)  $b_3\bar{c}_9=\bar{b}_{15}+c_9+c_3$;

9)  $b_3b_{15}=2\bar{b}_{15}+b_6+c_9$;

10)  $b_3\bar{b}_{15}=b_{10}+c_8+b_9+b_8+b_5+c_5$;

11)  $b_3b_5=b_{15}$;

12)  $b_3c_5=b_{15}$;

13)  $b_3b_8=b_{15}+\bar{c}_9$;

14)  $b_3c_8=b_3+\bar{b}_6+b_{15}$;

15)  $b_3b_9=b_{15}+\bar{c}_9+\bar{c_3}$;

16)  $b_3b_{10}=b_{15}+\bar{b}_6+c_9$;

18)  $c_{3}^2=\bar{c_3}+\bar{b}_6$;

19)  $c_3b_6=\bar{c_3}+b_{15}$;

20)  $c_3\bar{b}_6=b_{10}+b_8$;

21)  $c_3c_9=b_3+b_{15}+\bar{c}_9$;

22)  $c_3\bar{c}_9=c_8+b_9+b_{10}$;

23)  $c_3b_{15}=b_5+c_5+c_8+b_8+b_9+b_{10}$;

24)  $c_3{\bb_{15}}=2b_{15}+\bar{b}_6+\bar{c}_9$;

25)  $c_3b_5={\bb_{15}}$;

26)  $c_3c_5={\bb_{15}}$;

27)  $c_3b_8=c_3+b_6+{\bb_{15}}$;

28)  $c_3c_8=c_9+{\bb_{15}}$;

29)  $c_3b_9=b_{15}+\bar{b_3}+c_9$;

30)  $c_3b_{10}=\bar{b}_{15}+b_6+c_9$;

31)  $b_6\bar{b}_6=1+c_8+b_9+b_8+c_5+b_5$;

32)  $b_6c_9=2b_{15}+2\bar{c}_9+\bar{b}_6$;

33)  $b_6\bar{c}_9=b_{10}+c_8+2b_9+b_8+b_5+c_5$;

34)  $b_6b_{15}=2c_8+3b_{10}+2b_9+2b_8+b_5+c_5$;

35)  $b_6{\bb_{15}}=4b_{15}+\bar{b}_6+b_3+\bar{c_3}+2\bar{c}_9$;

36)  $b_6b_5={\bb_{15}}+b_6+c_9$;

37)  $b_6c_5={\bb_{15}}+b_6+c_9$;

38)  $b_6b_8=2{\bb_{15}}+b_6+c_3+c_9$;

39)  $b_6c_8=\bar{b_3}+2{\bb_{15}}+b_6+c_9$;

40)  $b_6b_9=b_6+2c_9+2\bb_{15}$;

41)  $b_6b_{10}=3\overline{b}_{15}+\bar{b_3}+c_3+c_9$;

42)  $b_{6}^2=2\bar{b}_6+\bar{c}_9+b_{15}$;

43)  $c_9\bar{c}_9=1+2c_8+2b_9+2b_{10}+c_5+b_5+2b_8$;

44)  $c_9b_{15}=2b_5+2c_5+3c_8+3b_8+3b_9+4b_{10}$;

45)  $c_9\bar{b}_{15}=6b_{15}+2\bar{b}_6+b_3+\bar{c_3}+3\bar{c}_9$;

46)  $c_9b_5=2\bar{b}_{15}+b_6+c_9$;

47)  $c_9c_5=2\bar{b}_{15}+b_6+c_9$;

48)  $c_9b_8=3\bar{b}_{15}+\overline{b_3}+b_6+2c_9$;

49)  $c_9c_8=3\bar{b}_{15}+c_3+b_6+2c_9$;

50)  $c_9b_9=\bar{b}_3+2c_9+3\bar{b}_{15}+c_3+2b_6$;

51)  $c_9b_{10}=\bar{b}_3+2c_9+4\bar{b}_{15}+c_3+b_6$;

52)  $ c_{9}^2=3b_{15}+b_3+\bar{c_3}+2\bar{b}_6+2\bar{c}_9$;

53)  $ b_{15}\bar{b}_{15}=1+3b_5+3c_5+5c_8+5b_8+6b_9+6b_{10}$;

54)  $b_{15}^2=9\bar{b}_{15}+4b_6+2\bar{b}_3+2c_3+6c_9$;

55)  $b_{15}b_5=3b_{15}+\bar{b}_6+b_3+\bar{c_3}+2\bar{c}_9$;

56)  $b_{15}c_5=3b_{15}+\bar{b}_6+b_3+\bar{c_3}+2\bar{c}_9$;

57)  $b_{15}b_8=5b_{15}+b_3+\bar{c_3}+2\bar{b}_6+3\bar{c}_9$;

58)  $b_{15}c_8=5b_{15}+b_3+\bar{c_3}+2\bar{b}_6+3\bar{c}_9$;

59)  $b_{15}b_9=6b_{15}+3\bar{c}_9+b_3+2\bar{b}_6+\bar{c_3}$;

60)  $b_{15}b_{10}=6b_{15}+3\bar{b}_6++b_3+\bar{c_3}+4\bar{c}_9$;

61)  $ b_{5}^2=1+b_{10}+b_9+b_5$;

62)  $b_5c_5=b_8+c_8+b_9$;

63)  $b_5b_8=c_5+c_8+b_8+b_9+b_{10}$;

64)  $b_5c_8=c_5+b_8+b_9+b_{10}+c_8$;

65)  $b_5b_9=c_8+c_5+b_8+b_5+b_9+b_{10}$;

66)  $b_5b_{10}=c_8+b_8+b_9+b_5+2b_{10}$;

67)  $c_{5}^2=1+b_9+b_{10}+c_5$;

68)  $c_5b_8=b_5+c_8+b_8+b_9+b_{10}$;

69)  $c_5c_8=b_5+c_8+b_8+b_9+b_{10}$;

70)  $c_5b_9=b_8+c_8+b_5+c_5+b_9+b_{10}$;

71)  $c_5b_{10}=b_8+c_8+b_9+c_5+2b_{10}$;

72)  $b_{8}^2=1+b_5+c_5+c_8+2b_8+b_9+2b_{10}$;

73)  $b_{8}c_8=b_5+c_5+c_8+b_8+2b_9+2b_{10}$;

74)  $b_{8}b_9=b_5+c_5+2c_8+b_8+2b_9+2b_{10}$;

75)  $b_{8}b_{10}=b_5+c_5+2c_8+2b_8+2b_9+2b_{10}$;

76)  $c_{8}^2=1+2c_8+2b_{10}+b_9+b_8+b_5+c_5$;

77)  $c_8b_9=b_5+c_5+c_8+2b_8+2b_9+2b_{10}$;

78)  $c_8b_{10}=2c_8+2b_{10}+2b_9+2b_8+b_5+c_5$;

79)  $b_9b_{10}=b_5+c_5+2c_8+2b_8+2b_9+3b_{10}$;

80)  $b_{9}^2=1+b_5+c_5+2c_8+2b_8+2b_9+2b_{10}$;

81)  $b_{10}^2=1+2c_8+2b_{10}+3b_9+2b_8+2b_5+2c_5$.

\begin{proof}
 The Theorem follows by Lemmas \ref{lemma4.1}---\ref{lemma4.7}.
\end{proof}

\end{theorem}

\par In this section, we don't have information about NITA
satisfying \textbf{} and $b_{10}$ is non-real. But we
conjecture:\\
\par \textbf{Conjecture 5.1} \emph{There exists no NITA satisfying  \textbf{Hypothesis 5.1} and
$b_{3}b_{6}=c_{8}+b_{10}$, where $b_{10}$ is non-real.}
\section{NITA Generated by $b_3$ Satisfying $b^2_3=c_3+b_6$, $c_3\neq b_3,\ \bb_3$, $b_6$ non-real, $(b_3b_8,\ b_3b_8)=4$
and $c_3^2=r_3+s_6$ }
 \par Now we discuss NITA satisfying the following hypothesis.
\begin{hyp}\label{hypx1}Let $(A,\ B)$  be a $NITA$ generated by an non-real
element $b_{3}\in B$ satisfies $b^2_3=c_3+b_6$, $c_3\neq b_3,\
\bb_3$, $b_6$ non-real, $(b_3b_8,\ b_3b_8)=4$ and $c_3^2=r_3+s_6$.
Here we still assume that $L(B)=1$ and $L_{2}(B)=\O$.
\end{hyp}
\par In this section, we cannot classify NITAs satisfying \textbf{Hypothesis
 \ref{hypx1}}, but we have following theorem.
\begin{theorem}
 There exists no $NITA$ satisfying \textbf{Hypothesis
 \ref{hypx1}} and $c_3$ generates a sub-NITA isomorphic to one of
 four known NITAs in \textbf{Main Theorem}.
 \end{theorem}

\par {Proof}. By assumptions we have that (1) $b_3\bb_3=1+b_8=c_3\bc_3$, (2)
$b_3^2=c_3+b_6$, (3) $(b_3b_8,\ b_3b_8)=4$ and (4) $c_3^2=r_3+s_6$.

\par \textbf{Step 1} There is no NITA satisfying hypothesis (1) to (4) and
containing NITA strictly isomorphic to $(Ch(PSL(2,\ 7)),
Irr(PSL(2,\
7))$.\\
\par By \cite{cgyarad}, $(Ch(PSL(2,\
7)), Irr(PSL(2,\ 7))$ the sub-NITA generated by $c_3$ is of
dimension 6 and has base elements: 1,\ $c_3,\ \bc_3,\ s_6,\ b_7,\
b_8$. Moreover, we have following equations:
\par (a1) $c_3b_8=c_3+s_6+b_7+b_8$,
\par (a2) $s_6^2=1+2s_6+b_7+b_8$,
\par (a3) $c_3s_6=\bc_3+b_7+b_8$.
\par It follows by (1) that $(b_3\bc_3,\ b_3\bc_3)=(b_3\bb_3,\
c_3\bc_3)=2$, hence we get equation (a4): $b_3\bc_3=\bb_3+t_6$,
where $t_6\in B$. By associative law and equations (a1) and (a4),
we have
\par $(b_3\bb_3)\bc_3=\bc_3+\bc_3b_8=\bc_3+\bc_3+s_6+b_7+b_8,$
\par $(b_3\bc_3)\bb_3=\bb_3^2+\bb_3t_6=\bc_3+\bb_6+\bb_3t_6.$

Then we get $\bb_6=s_6$, so $b_6$ is real. Moreover we have
equation (a5): $\bb_3t_6=\bc_3+s_6+b_7+b_8$. By
$(b_3\bb_3)b_3=b_3^2\bb_3$, we arrive at
$b_3b_8=\bar{t}_6+\bb_3s_6.$ Now by (3) we obtain equation
 (a6): $\bb_3s_6=b_3+x+y$, where distinct $b_3,\ x$ and $y$ such that $x+y$ is of degree
  15.
 Hence we have equation (a7): $b_3b_8=b_3+\bar{t}_6+x+y$.
Therefore by (a2),
 (a3), (a5) and (a6)
\par
$b_3^2s_6=(c_3+s_6)s_6=c_3s_6+s_6s_6=\bc_3+b_7+b_8+1+2s_6+b_7+b_8,$
\par $b_3(b_3s_6)=b_3(\bb_3+\bar{x}+\bar{y})=1+b_8+b_3\bar{x}+b_3\bar{y}.$\\
So $b_3\bar{x}+b_3\bar{y}=\bc_3+2s_6+2b_7+b_8$, but  by (a7), the
left side of this equation contains two $b_8$, right side only
one,
, a contradiction.\\
 \par \textbf{Step 2. }
 There is no NITA satisfying hypothesis (1) to (4) and containing a sub-NITA generated
 by $c_3$ is a NITA of dimension 17, 32 or 22.\\

\par By \cite{cgyarad}, we see that for the NITAs of dimensions 17, 32 and 22 one can
find new base elements: $b_5,\ c_5,\ s_6,\ x_8,\ b_9$ and $b_{10}$
such that following equations hold:
\par  (a1) $c_3b_8=c_3+\bar{s}_6+b_{15}$,
\par  (a2) $s_6\bar{s}_6=1+b_5+c_5+b_8+x_8+b_9,$
\par (a3)  $c_3s_6=b_8+b_{10}.$\\
 \par Suppose $b_3\bc_3=\bb_3+t_6$, where $t_6\in B$. Then by hypothesis and equation (a1)
 it follows that $(b_3\bb_3)\bc_3=\bc_3+\bc_3b_8=\bc_3+\bc_3+s_6+\bb_{15}$,
 $(b_3\bc_3)\bb_3=\bb_3^2+\bb_3t_6=\bc_3+\bb_6+\bb_3t_6$.
 we get $s_6=\bb_6$ and
 equation (a4): $\bb_3t_6=\bc_3+\bb_{15}$. Hence $b_3+b_3b_8=(b_3\bb_3)b_3=b_3^2\bb_3=\bb_3c_3+\bb_3\bar{s}_6$,
  we get equation:
 $b_3b_8=\bar{t}_6+\bb_3\bar{s}_6.$ Now by (3) we obtain equation
 (a6): $\bb_3\bar{s}_6=b_3+x+y$, where distinct $b_3,\ x$ and $y$ such that $x+y$ is of degree
  15. Hence we
 have equation (a7): $b_3b_8=b_3+\bar{t}_6+x+y$. Therefore
\par
$b_3^2s_6=(c_3+\bar{s}_6)s_6=c_3s_6+s_6\bar{s}_6=b_8+b_{10}+1+b_5+c_5+b_8+x_8+b_9,$
\par $ b_3(b_3s_6)=b_3(\bb_3+\bar{x}+\bar{y})=1+b_8+b_3\bar{x}+b_3\bar{y}.$\\
So $b_3\bar{x}+b_3\bar{y}=1+b_5+c_5+b_8+x_8+b_9$, which implies
$x_8=b_8$ by (a7), a contradiction. This concludes the theorem.\\

\section{Structure of NITA generated by $b_3$ and satisfying
$b^2_3=c_3+b_6$, $c_3\neq b_3,\ \bb_3$, $(b_3b_8,\ b_3b_8)=3$ and
$c_3$ non-real}

\par
In this section we classify NITA generated by $b_3$ and satisfying
$b^2_3=c_3+b_6$, $c_3\neq b_3,\ \bb_3$ and $c_3$ non-real and
$(b_3b_8,\ b_3b_8)=3$. At first, we state Hypothesis \ref{hyp7.1}.
\par
\begin{hyp}\label{hyp7.1}
Let $(A,B)$ be a NITA generated by a non-real element $b_3 \in B$
such that $b_3\bb_3=1+b_8$ and $b_3^2=c_3+b_6, \ c_3,\ b_6 \in B$,
$c_3 \ne \bb_3$, $c_3 \ne \bc_3$ and $(b_3b_8,b_3b_8)=3$.
\end{hyp}
\begin{lemma}\label{s1} Let $(A,B)$ be a NITA satisfies Hypothesis
\ref{hyp7.1}.Then $c_3\bc_3=b_3\bb_3=1+b_8$ and
$b_3b_8=b_3+x_6+y_{15}$ and $c_3b_8=c_3+y_6+z_{15}$, where $b_8,\
x_6,\ y_{15}\in B,$ and there exist $r_3,\ s_6,\ u_3,\ v_6, w_3,\
z_6\in B$ such that
\[
\begin{array}{rll}
b^2_3=r_3+s_6,& \bb_3r_3=b_3+x_6,& \bb_3s_6=b_3+y_{15},\\
b_3c_3=u_3+v_6,&\bb_3u_3=c_3+y_6,&\bb_3v_6=c_3+z_{15},\\
\bb_3c_3=w_3+z_6,&b_3w_3=c_3+y_6,&b_3z_6=c_3+z_{15}.
\end{array}
\]
\end{lemma}
\par {\bf Proof.} It is obvious that $(b_3^2,\ b_3^2)=2$. So $b^2_3=r_3+s_6$ or $r_4+s_5$. Hence
\[
\begin{array}{rll}
b^2_3\bb_3&=&\bb_3r_3+\bb_3s_6\ \mbox{or}\ \bb_3r_4+\bb_3s_6,\\
(b_3\bb_3)b_3&=&b_3+b_3b_8\\
&=&2b_3+x_6+y_{15}.
\end{array}
\]
Since $b_3\in Supp\{\bb_3r_3\}$ and $Supp\{\bb_3s_6\}$ or $b_3\in
Supp\{\bb_3r_4\}$ and $Supp\{\bb_3s_5\}$, the first part of the
lemma follows.
\par Since $b_3\bb_3=c_3\bc_3=1+b_8$, it follow that $b_3c_3=u_3+v_6\ \mbox{or}\ u_4+v_5$ and $\bb_3c_3=w_3+z_6$ or $w_4+z_5$. The second part and the third part of the lemma follows from $(b_3\bb_3)c_3=(b_3c_3)\bb_3=(\bb_3c_3)b_3$.

\begin{lemma}\label{s2} Let $(A,B)$ satisfies Hypothesis
\ref{hyp7.1}. Then there are non-real base elements $x_6,\
y_{15},\ x_{15}$ and real basis elements $r_3$,  $s_6,\ d_9,\
x_{10},\ t_{15}$ such that the following equations hold:
\[
\begin{array}{ll}
b_3\bb_3=1+b_8,&b^2_3=c_3+b_6,\\
b_3\bc_3=\bb_3+\bar{x}_6,&b_3x_6=c_3+y_{15},\\
b_3b_8=b_3+x_6+x_{15},&b_3\bar{x}_6=b_8+x_{10},\\
b_3\bb_6=\bb_3+\bar{x}_{15},&b_3c_6=\bc_3+\bar{y}_{15},\\
b_3c_3=r_3+s_6,&
b_3r_3=\bc_3+\bb_6,\\
b_3s_6=\bc_3+\bar{y}_{15},&b_3x_{10}=x_6+x_{15}+b_9,\\
b_3\bar{y}_{15}=\bar{x}_6+2\bar{x}_{15}+\bb_9,&b_3y_{15}=s_6+2t_{15}+d_9, \\
b_3b_6=r_3+t_{15}, &b_3x_{15}=b_6+2y_{15}+c_9,\\
b_3t_{15}=\bb_6+\bc_9+2\bar{y}_{15},&\\
\\
c_3\bar{c}_3=1+b_8,&c_3r_3=\bb_3+\bar{x}_6,\\
c_3s_6=\bb_3+\bar{x}_{15},&c^2_3=\bc_3+\bb_6,\\
c_3b_8=c_3+b_6+y_{15},&c_3\bar{x}_6=b_3+x_{15},\\
c_3x_6=r_3+t_{15},&c_3x_{10}=b_6+y_{15}+c_9,\\
c_3y_{15}=\bb_6+\bc_9+2\bar{y}_{15},&c_3b_6=\bc_3+\bar{y}_{15}, \\
c_3t_{15}=\bar{x}_6+\bb_9+2\bar{x}_{15},&c_3x_{15}=2t_{15}+s_6+d_9,\\
c_3\bar{x}_{15}=2x_{15}+x_6+b_9,&c_3\bb_6=b_8+x_{10},\\
\\
r_3\bc_3=c_3+x_6,&r_3\bb_3=c_3+b_6,\\
r_3\bar{x}_6=c_3+y_{15},&r_3y_{15}=\bar{x}_6+2\bar{x}_{15}+\bb_9,\\
r_3b_6=\bar{x}_{15}+\bb_3,&r_3b_8=s_6+r_3+t_{15},\\
r_3\bar{x}_{15}=b_6+c_9+2y_{15},&\\
r_3x_{10}=s_6+t_{15}+d_9,&\\
\\
x^2_6=2b_6+y_{15}+c_9,&x_6b_8=b_3+x_6+2x_{15}+b_9,\\
\\
b_6^2=2\bb_6+\bar{y}_{15}+\bc_9,&b_6\bar{x}_6=2x_6+b_9+x_{15},\\
b_6x_6=2s_6+t_{15}+d_9,&b_6b_8=c_3+b_6+c_9+2y_{15},\\
\\
s_6\bb_6=2x_6+x_{15}+b_9,&s_6b_8=r_3+s_6+2t_{15}+d_9,\\
s_6\bc_3=b_3+x_{15}&\\
\\
b_8y_{15}=c_3+2b_6+2c_9+4y_{15}+\bb_3d_9.\\
\\
\end{array}
\]
\end{lemma}
\par {\bf Proof.} By Hypothesis \ref{hyp7.1} $(b_3b_8,\ b_3b_8)=3$. Then the following equations hold by {\bf Lemma \ref{la2}}:\\
\[
\begin{array}{ll}
b_3\bb_3=1+b_8,& c_3\bar{c}_3=1+b_8,\\
b^2_3=c_3+b_6,&\bb_3c_3=b_3+x_6.
\end{array}
\]
\par Let $b_3x_6=c_3+y_{15}$, then
\[
\begin{array}{rll}
(b_3\bb_3)c_3&=&c_3+c_3b_8,\\
(\bb_3c_3)b_3&=&b^2_3+b_3x_6\\
&=&c_3+b_6+c_3+y_{15}.
\end{array}
\]
Hence
$$c_3b_8=c_3+b_6+y_{15}.$$
Therefore $y_{15}\in B$ by {\bf Lemma \ref{la2}}.
\par Let $c_3\bar{x}_6=b_3+x_{15}$, $x_{15}\in \mathit{N}^*$$B$. Since  $(c_3\bc_3)b_3=(b_3\bc_3)c_3$, it follows that $x_{15}\in B$ and
$$b_3b_8=b_3+x_6+x_{15}$$
by $b^2_3\bb_3=(b_3\bb_3)b_3$. Because $y_{15}\in B$, so
$(\bb_3x_6,\ \bb_3x_6)=(b_3x_6,\ b_3x_6)=2$. Thus
$$\bb_3x_6=b_8+x_{10},\  x_{10}\in B.$$
\par Since $b^2_3\bc_3=(b_3\bc_3)b_3$ and $b^2_3\bb_3=(b_3\bb_3)b_3$, one has that
\[
\begin{array}{ll}
c_3\bb_6=b_8+x_{10},&\bb_3b_6=b_3+x_{15}.
\end{array}
\]
\par By {\bf Lemma \ref{s1}}, there exist $w_3$, $z_6$, $r_3$, $s_6\in B$ such that
\[
\begin{array}{rll}
c^2_3=w_3+z_6,&w_3\bc_3=c_3+b_6,&z_6\bc_3=c_3+y_{15},\\
b_3c_3=r_3+s_6,&b_3r_3=w_3+z_6,&b_3s_6=w_3+z_{15}.
\end{array}
\]
\par By {\bf Lemma \ref{s1}} again, it follow that
\[
\begin{array}{ll}
r_3\bc_3=b_3+x_6, &\bc_3s_6=b_3+x_{15},\\
r_3\bb_3=c_3+b_6,&\bb_3s_6=c_3+y_{15}.
\end{array}
\]
Thus $(\bc_3r_3,\ \bc_3r_3)=2$, which implies that
$r_3\bar{r}_3=1+b_8$.
\par Since
\[
\begin{array}{rll}
(b_3r_3)\bc_3&=&w_3\bc_3+z_6\bc_3\\
&=&2c_3+b_6+y_{15},\\
(b_3\bc_3)r_3&=&r_3\bb_3+r_3\bar{x}_6\\
&=&c_3+b_6+r_3\bar{x}_6,
\end{array}
\]
which implies that $r_3\bar{x}_6=c_3+y_{15}.$
\par We assert that $x_{10}$ is real. In fact
\[
\begin{array}{rll}
\bb_3(r_3\bc_3)&=&b_3\bb_3+x_6\bb_3\\
&=&1+b_8+b_8+x_{10},\\
(\bb_3r_3)\bc_3&=&c_3\bc_3+b_6\bc_3\\
&=&1+b_8+b_8+\bar{x}_{10},
\end{array}
\]
so it follows that $x_{10}$ is real.
\par Since
\[
\begin{array}{rll}
\bb_3(w_3\bc_3)&=&c_3\bb_3+b_6\bb_3\\
&=&b_3+x_6+b_3+x_{15},\\
(\bb_3\bc_3)w_3&=&\bar{r}_3w_3+\bar{s}_6w_3,
\end{array}
\]
we have that
\[
\begin{array}{rll}
\bar{r}_3w_3=b_3+x_6,&w_3\bar{s}_6=b_3+x_{15}.\end{array}
\]
Consequently $\bb_3w_3=r_3+s_6.$

\par Now we have that $(b_3\bb_3)w_3=(\bb_3w_3)b_3$ and
$$w_3b_8=w_3+z_6+z_{15}.$$ Then $w_3\bar{w}_3=1+b_8.$
\par Since
\[
\begin{array}{rll}
(b_3\bb_3)^2&=&1+2b_8+b^2_8\\
b^2_3\bb^2_3&=&(c_3+b_6)(\bc_3+\bb_6)\\
&=&1+b_8+c_3\bb_6+\bc_3b_6+b_6\bb_6\\
&=&1+3b_8+2x_{10}+b_6\bb_6,
\end{array}
\]
 so
$$b^2_8=2x_{10}+b_8+b_6\bb_6.$$
\par By $b^2_3\bb_6=(b_3\bb_6)b_3$, we have that $x_{10}+b_6\bb_6=1+b_3\bar{x}_{15}.$ Then $(b_3x_{10},\ x_{15})=1$. But $(b_3x_{10},\ x_6)=1$. There exits $b_9\in B$ such that
$$b_3x_{10}=x_6+x_{15}+b_9,\ b_9\in \mathit{N}^*B.$$
\par It is easy to see that $b_3x_{10}$ cannot have constituents of degree 3 and 4 by $L_1(B)=\{1\}$ and $L_2(B)=\emptyset$. Thus $b_9\in B$.
\par Since
\[
\begin{array}{rll}
(\bb_3x_6)b_3&=&b_3b_8+b_3x_{10}\\
&=&b_3+x_6+x_{15}+x_6+x_{15}+b_9,\\
(b_3\bb_3)x_6&=&x_6+x_6b_8,
\end{array}
\]
we have $$x_6b_8=b_3+x_6+2x_{15}+b_9.$$
\par Multiplying both sides of the equation $x_{10}+b_6\bb_6=1+b_3\bar{x}_{15}$ by $b_3$, we have that
$$b_3x_{10}+(b_3\bb_6)b_6=b_3+b^2_3\bar{x}_{15}.$$
so
$$x_6+x_{15}+b_9+b_6\bb_3+b_6\bar{x}_{15}=c_3\bar{x}_{15}+b_6\bar{x}_{15}.$$
Hence $$c_3\bar{x}_{15}=2x_{15}+x_6+b_9.$$
\par Now checking associative law of $(b_3\bc_3)b_6=(\bc_3b_6)b_3$, $(r_3\bar{r}_3)x_6=(\bar{r}_3x_6)r_3$ and $(r_3\bb_3)\bar{x}_6=(r_3\bar{x}_6)\bb_3$, we come to
\[
\begin{array}{rll}
b_6\bar{x}_6&=&2x_6+b_9+x_{15},\\
r_3\bar{y}_{15}&=&x_6+2x_{15}+b_9,\\
\bb_3y_{15}&=&x_6+2x_{15}+b_9.
\end{array}
\]
\par Hence
$$c_3\bar{y}_{15}+b_6\bar{y}_{15}=b^2_3\bar{y}_{15}=b_3(b_3\bar{y}_{15})=b_3\bar{x}_6+2b_3\bar{x}_{15}+b_3\bb_9.$$

\par Since  $(b_3\bar{x}_6,\ b_8)=(b_3\bar{x}_6,\ x_{10})=(b_3\bar{x}_{15},\ b_8)=(b_3\bar{x}_{15},\ x_{10})=(b_3\bb_9,\ x_{10})=1$ and $(b_3\bb_9,\ b_8)=0$, we have that
$$(c_3\bar{y}_{15}+b_6\bar{y}_{15},\ b_8)=3\ \mbox{and}\ (c_3\bar{y}_{15}+b_6\bar{y}_{15},\ x_{10})=4.$$
\par Observe the degrees, we know that $(c_3\bar{y}_{15},\ x_{10})\leq 1$. If $(c_3\bar{y}_{15},\ x_{10})=0$, then $(b_6\bar{y}_{15},\ x_{10})=4$, so
$b_6x_{10}=4y_{15}.$ But $(b_6x_{10},\ c_3)=1$, a contradiction.
Hence $(c_3\bar{y}_{15},\ x_{10})=1$ and $(b_6x_{10},\
y_{15})=(b_6\bar{y}_{15},\ x_{10})=3.$
\par Since $(c_3x_{10},\ c_3x_{10})=(b_3x_{10},\ b_3x_{10})=3$ and $(c_3x_{10}, b_6)=1$, there exists $c_9\in B$ such that
$$c_3x_{10}=b_6+y_{15}+c_9.$$
Then
\[
\begin{array}{rll}
(\bb_3c_3)x_6&=&b_3x_6+x^2_6\\
&=&c_3+y_{15}+x^2_6,\\
(\bb_3x_6)c_3&=&c_3b_8+c_3x_{10}\\
&=&c_3+b_6+y_{15}+b_6+y_{15}+c_9,\\
(c_3\bc_3)b_6&=&b_6+b_6b_8,\\
(\bc_3b_6)c_3&=&b_8c_3+x_{10}c_3\\
&=&c_3+b_6+y_{15}+y_{15}+b_6+c_9,\\
(c_3b_6)\bc_3&=&c_3+b_6+\bc_3z_{15},
\end{array}
\]
which implies that
\[
\begin{array}{rll}
x^2_6&=&2b_6+y_{15}+c_9,\\
b_6b_8&=&c_3+b_6+c_9+2y_{15},\\
c_3\bar{z}_{15}&=&\bb_6+\bc_9+2\bar{y}_{15}.
\end{array}
\]
\par Since
\[
\begin{array}{rll}
(b_3\bc_3)^2&=&(\bb_3+\bar{x}_6)^2\\
&=&\bb^2_3+2\bb_3\bar{x}_6+\bar{x}^2_6\\
&=&2\bc_3+3\bb_6+3\bar{y}_{15}+\bc_9,\\
b^2_3\bc^2_3&=&c_3\bar{w}_3+c_3\bar{z}_6+\bar{w}_3b_6+b_6\bar{z}_6\\
&=&\bc_3+\bb_6+\bc_3+\bar{y}_{15}+b_6\bar{w}_3+b_6\bar{z}_6,
\end{array}
\]
but $(\bar{w}_3b_6,\ \bc_3)=1$,  we have that
\[
\begin{array}{rll}
w_3\bb_6&=&c_3+y_{15},\\
b_6\bar{z}_6&=&2\bb_6+\bar{y}_{15}+\bc_9.
\end{array}
\]
 \par Since
\[
\begin{array}{rll}
(b_3\bar{r}_3)\bb_6&=&\bc_3\bb_6+\bb_6^2\\
&=&\bar{w}_3+\bar{z}_{15}+\bb_6^2,\\
(b_3\bb_6)\bar{r}_3&=&\bb_3\bar{r}_3+\bar{x}_{15}\bar{r}_3\\
&=&\bar{w}_3+\bar{z}_6+\bar{r}_3\bar{x}_{15},
\end{array}
\]
and $(\bb^2_6,\ \bar{z}_6)=2$, we have that
$(\bar{r}_3\bar{x}_{15},\ \bar{z}_6)=1$. Hence $(r_3\bar{z}_6,\
\bar{x}_{15})=1$. But  $(r_3\bar{z}_6,\ \bb_3)=1$ by
$b_3r_3=w_3+z_6$. Thus
$$\bar{r}_3z_6=b_3+x_{15}.$$
\par From the following equations
\[\begin{array}{rlll}
\bb_3(w_3\bb_6)&=&\bb_3(c_3+y_{15})\\
&=&b_3+x_6+x_6+2x_{15}+b_9,\\
(\bb_3w_3)\bb_6&=&r_3\bb_6+s_6\bb_6,
\end{array}
\]
and $(r_3\bb_6, r_3\bb_6)=(b_3\bb_6,\ b_3\bb_6)=2$, $(r_3\bb_6,\
b_3)=1$, we come to that
$$r_3\bb_6=b_3+x_{15},\ s_6\bb_6=2x_6+x_{15}+b_9.$$
\par By  $(b_3c_3)\bar{r}_3=(c_3\bar{r}_3)b_3$ and $\bc_3(w_3b_8)=(\bc_3b_8)w_3$, we have that
$$\bar{r}_3s_6=b_8+x_{10},\ \ w_3\bar{y}_{15}=b_6+c_9+2y_{15}.$$
\par Since $(\bar{r}_3b_8,\ \bar{r}_3b_8)=(r_3\bar{r}_3,\ b_8^2)=(b_3\bb_3,\ b_8^2)=(b_3b_8,\ b_3b_8)=3$ by {\bf Lemma \ref{s1}} and $(\bar{r}_3b_8,\ \bar{s}_6)=(\bar{r}_3b_8,\ \bar{r}_3)=1$, there exists $t_{15}\in B$ such that
$$\bar{r}_3b_8=\bar{s}_6+\bar{r}_3+t_{15}.$$
\par  Since
\[
\begin{array}{rll}
(w_3b_8)\bar{r}_3&=&w_3\bar{r}_3+z_6\bar{r}_3+z_{15}\bar{r}_3\\
&=&b_3+x_6+b_3+x_{15}+\bar{r}_3z_{15},\\
(\bar{r}_3w_3)b_8&=&b_3b_8+x_6b_8\\
&=&b_3+x_6+x_{15}+b_3+x_6+2x_{15}+b_9,\\
(\bar{r}_3b_8)w_3&=&\bar{s}_6w_3+\bar{r}_3w_3+t_{15}w_3\\
&=&b_3+x_{15}+b_3+x_6+w_3t_{15},
\end{array}
\]
we have that
\[
\begin{array}{rll}
r_3\bar{z}_{15}&=&\bar{x}_6+2\bar{x}_{15}+\bb_9,\\
w_3t_{15}&=&x_6+b_9+2x_{15}.
\end{array}
\]
\par Since
\[
\begin{array}{rll}
(\bb_3r_3)b_8&=&c_3b_8+b_8b_6\\
&=&c_3+b_6+y_{15}+2y_{15}+b_6+c_3+c_9,\\
(\bb_3b_8)r_3&=&r_3\bb_3+r_3\bar{x}_6+r_3\bar{x}_{15}\\
&=&c_3+b_6+c_3+y_{15}+r_3\bar{x}_{15},\\
(r_3b_8)\bb_3&=&r_3\bb_3+s_6\bb_3+\bar{t}_{15}\bb_3\\
&=&c_3+b_6+c_3+y_{15}+\bb_3\bar{t}_{15},
\end{array}
\]
we have that $w_3=\bc_3,\ z_6=\bb_6,\ z_{15}=\bar{y}_{15}$ and
\[
\begin{array}{rll}
r_3\bar{x}_{15}&=&b_6+c_9+2y_{15},\\
b_3t_{15}&=&\bb_6+\bc_9+2\bar{y}_{15}.
\end{array}
\]
\par Further $b_3r_3=\bc_3+\bb_6$, we have that $(b_3c_3,\ \bar{r}_3)=1$.
But $(b_3c_3,\ r_3)=1$. Hence $r_3=\bar{r}_3$. Moreover $\bar{r}_3+\bar{s}_6=b_3\bar{w}_3=b_3c_3=r_3+s_6$,
which implies that $s_6=\bar{s}_6.$
Therefore $t_{15}$ is real by the expression of $r_3b_8$, and
$$c_3t_{15}=\bar{x}_6+\bc_9+2\bar{x}_{15}.$$
 Now we have that $t_{15}\in Supp\{c_3x_6\}$. Thus
$$c_3x_6=r_3+t_{15}.$$
\par Since
\[
\begin{array}{rll}
(c_3r_3)\bb_6&=&\bb_3\bb_6+\bar{x}_6\bb_6,\\
(r_3\bb_6)c_3&=&b_3c_3+c_3x_{15}\\
&=&r_3+s_6+c_3x_{15},\\
(c_3\bb_6)r_3&=&r_3b_8+r_3x_{10}\\
&=&r_3+s_6+t_{15}+r_3x_{10}
\end{array}
\]
and $(r_3x_{10},\ s_6)=1$, we have $s_6\in Supp\{c_3x_{15}\}$. But
$(c_3x_{15},\ t_{15})=(c_3t_{15},\ \bar{x}_{15})=2$.  Let
$$c_3x_{15}=2t_{15}+s_6+d_9,\ d_9\in\mathit{N}^*B.$$
Then $r_3x_{10}=s_6+t_{15}+d_9$ and by
$(c_3r_3)\bb_6=(c_3\bb_6)r_3$ we obtain
$$\bb_3\bb_6+\bar{x}_6\bb_6=r_3+s_6+t_{15}+s_6+t_{15}+d_9.$$
Further $d_9$ is real for $r_3$ and $x_{10}$ are real, and
\[
\begin{array}{rll}
b_3b_6&=&r_3+t_{15},\\
x_6b_6&=&2s_6+t_{15}+d_9.
\end{array}
\]
\par Since
\[
\begin{array}{rll}
(b_3\bb_3)b_6&=&b_6+b_6b_8\\
&=&b_6+c_3+b_6+c_9+2y_{15},\\
b_3(\bb_3b_6)&=&b^2_3+b_3x_{15}\\
&=&c_3+b_6+b_3x_{15},\\
(b_3b_6)\bb_3&=&r_3\bb_3+t_{15}\bb_3\\
&=&c_3+b_6+\bb_3t_{15},
\end{array}
\]
we have that
\[
\begin{array}{rll}
b_3x_{15}&=&b_6+2y_{15}+c_9,\\
b_3t_{15}&=&\bb_6+\bc_9+2\bar{y}_{15}.
\end{array}
\]
\par Since
\[
\begin{array}{rll}
b^2_3x_6&=&c_3x_6+b_6x_6\\
&=&r_3+t_{15}+2s_6+t_{15}+d_9,\\
(b_3x_6)b_3&=&b_3c_3+b_3y_{15}\\
&=&r_3+s_6+b_3y_{15},
\end{array}
\]
we have that $$b_3y_{15}=2t_{15}+s_6+d_9.$$
\par Because $(b_3y_{15},\ b_3y_{15})=(b_3\bar{y}_{15},\ b_3\bar{y}_{15})=6$, we have that $d_9\in B$.
\par Checking the associative law of $(b_3\bb_3)s_6=(b_3s_6)\bb_3$, we have
$$s_6b_8=r_3+s_6+2t_{15}+d_9.$$
\par Since
\[
\begin{array}{rll}
(b_3\bb_3)y_{15}&=&y_{15}+b_8y_{15},\\
(b_3y_{15})\bb_3&=&\bb_3s_6+2\bb_3t_{15}+\bb_3d_9\\
&=&c_3+y_{15}+2b_6+2c_9+4y_{15}+\bb_3d_9,
\end{array}
\]
it follows that $b_8y_{15}=c_3+2b_6+2c_9+4y_{15}+\bb_3d_9.$
\par This is the end of the proof of the lemma.\\
\par {\bf Remark.} It always follows that $(c_3c_9,\ \bar{y}_{15})=1$ by {\bf Lemma \ref{s2}}. So $(c_3c_9,\ c_3c_9)\geq 2$. At the following, we shall investigate  the expression of $c_3c_9$.\\

\begin{lemma}\label{dir1}
There exist no $x_4$,\ $y_8\in B$ such that
$c_3c_9=x_4+y_8+\bar{y}_{15}$.
\end{lemma}
\par {\bf Proof.} If there exist no $x_4$, $y_8\in B$ such that $$c_3c_9=x_4+y_8+\bar{y}_{15}.$$
Then there exist $z_3\in B$ such that $c_3\bar{x}_4=\bc_9+z_3$.
Further there is $y_5\in B$ such that $c_3\bar{z}_3=x_4+y_5$. Thus
$z_3\bar{z}_3=1+b_8$. By {\bf Lemma \ref{s1}} we have that
\[
\begin{array}{rll}
c_3(z_3\bar{z}_3)&=&c_3+c_3b_8\\
&=&2c_3+b_6+y_{15},\\
z_3(c_3\bar{z}_3)&=&z_3x_4+z_3y_5.
\end{array}
\]
It must follow that $z_3x_4=2c_3+b_6$, which implies that
$(z_3\bc_3,\ x_4)=2$, a contradiction.

\begin{lemma}\label{dir3}
There exist no $x_5$,\ $y_7\in B$ such that
$c_3c_9=x_5+y_7+\bar{y}_{15}$.
\end{lemma}
\par {\bf Proof.} If some $x_5$,\ $y_7\in B$ such that $c_3c_9=x_5+y_7+\bar{y}_{15}$.
Then $(c_3\bar{x}_5,\ \bc_9)=1$. There are three possiblities:
\[
\begin{array}{cl}
(I)&c_3\bar{x}_5=\bc_9+2x_3,\\
(II)&c_3\bar{x}_5=\bc_9+x_3+y_3,\\
(III)&c_3\bar{x}_5=\bc_9+y_6.
\end{array}
\]
It is easy to see that (I) will leads to a contradiction.
\par If (II) follows, then $c_3\bar{x}_5=\bc_9+x_3+y_3$. Set $c_3\bar{x}_3=x_5+z_4$, $c_3\bar{y}_3=x_5+t_4$, where $z_4,\ t_4\in B$. Hence
\[
\begin{array}{rll}
\bc_3(c_3\bar{x}_3)&=&\bc_3x_5+\bc_3z_4\\
&=&c_9+\bar{x}_3+\bar{y}_3+\bc_3z_4,\\
(c_3\bc_3)\bar{x}_3&=&\bar{x}_3+\bar{x}_3b_8.
\end{array}
\]
Therefore $(x_3\bar{y}_3,\ b_8)=1$, which implies that
$\bar{y}_3x_3=1+b_8$. So $y_3=x_3$, a contradiction.
\par If (III) follows, then
\[
\begin{array}{rll}
\bc_3(c_3\bar{x}_5)&=&\bc_3\bc_9+\bc_3y_6\\
&=&y_{15}+\bar{x}_5+\bar{y}_7+\bc_3y_6\\
(c_3\bc_3)\bar{x}_5&=&\bar{x}_5+\bar{x}_5b_8.
\end{array}
\]
Hence $\bar{x}_5b_8=y_{15}+\bar{y}_7+\bc_3y_6.$ Therefore
$(b_8y_{15},\ \bar{x}_5)=1$. By {\bf Lemma \ref{s2}}, we have that
$(\bar{x}_5,\ \bb_3d_9)=1$. Since $(\bb_3d_9,\ y_{15})=1$, we may
set that
$$\bb_3d_9=\bar{x}_5+y_{15}+e_7,\ e_7\in\mathit{N}^*B.$$
\par If $e_7=m_3+n_4$, then $b_3n_4=d_9+p_3$, some $p_3\in B$. Let $\bb_3p_3=n_4+z_5$,\ some $z_5\in B$. Then $p_3\bar{p}_3=1+b_8$. Further
\[
\begin{array}{rll}
b_3(p_3\bar{p}_3)&=&b_3+b_3b_8\\
&=&2b_3+x_6+x_{15},\\
(b_3\bar{p}_3)p_3&=&\bar{n}_4p_3+\bar{z}_5p_3.
\end{array}
\]
Notice that $b_3\in Supp\{\bar{n}_4p_3\}$ and $b_3\in
Supp\{\bar{z}_5p_3\}$, we come to a contradiction.
\par If $e_7\in B$, then  $\bb_3d_9=\bar{x}_5+y_{15}+e_7$. Let $b_3\bar{x}_5=d_9+u_6$.
\par We assert that $u_6\in B$. Otherwise, let $u_6=r_3+s_3$, where $r_3,\ s_3\in B$, then
$$b_3\bar{r}_3=x_5+u_4.$$
Hence $r_3\bar{r}_3=1+b_8$. Moreover
\[
\begin{array}{rll}
b_3(r_3\bar{r}_3)&=&b_3+b_3b_8\\
&=&2b_3+x_6+x_{15},\\
r_3(b_3\bar{r}_3)&=&r_3x_5+r_3u_4.
\end{array}
\]
It is impossible for $b_3\in Supp\{r_3x_5\}$ and $b_3\in
Supp\{r_3u_4\}$.
\par Now we have that
\[
\begin{array}{rll}
\bb_3(c_3\bar{x}_5)&=&\bb_3\bc_9+\bb_3y_6,\\
(\bb_3c_3)\bar{x}_5&=&b_3\bar{x}_5+x_6\bar{x}_5\\
&=&d_9+u_6+x_6\bar{x}_5.
\end{array}
\]
\par If $d_9,\ u_6\in Supp\{\bb_3y_6\}$, then $\bb_3y_6=d_9+u_6+t_3$. Let $b_3t_3=y_6+p_3$. Hence $t_3\bar{t}_3=1+b_8$. Moreover
\[
\begin{array}{rll}
b_3(t_3\bar{t}_3)&=&2b_3+x_6+x_{15},\\
\bar{t}_3(b_3t_3)&=&\bar{t}_3y_6+\bar{t}_3p_3.
\end{array}
\]
So $\bar{t}_3y_6=b_3+x_{15}$. But $3=(\bb_3y_6,\
\bb_3y_6)=(\bar{t}_3y_6,\ \bar{t}_3y_6)=2$ by {\bf Lemma
\ref{la1}}, a contradiction.
\par If  $u_6\in Supp\{\bb_3y_6\}$, $d_9\in Supp\{\bb_3\bc_9\}$ or $u_6\in Supp\{\bb_3\bc_9\}$, $d_9\in Supp\{\bb_3y_6\}$.
\par If  the former one follows, then
there exists a $q_3$ such that
$$b_3c_9=t_{15}+d_9+q_3.$$
Hence $c_9\in Supp\{\bb_3d_9\}$, a contradiction.
\par If $u_6\in Supp\{\bb_3\bc_9\}$, $d_9\in Supp\{\bb_3y_6\}$, then
$$b_3d_9=\bar{y}_{15}+y_6+r_6,$$
where $r_6\in\mathit{N}^*B$. But $d_9$ is real,
$\bb_3d_9=\bar{x}_5+y_{15}+e_7$, a contradiction. This is the end
of the proof.
\begin{lemma}\label{dir4}
There exist no $m_6$ and $n_6$ such that
$c_3c_9=\bar{y}_{15}+m_6+n_6.$
\end{lemma}
\par {\bf Proof.} If there exist $m_6$ and $n_6$ such that $c_3c_9=\bar{y}_{15}+m_6+n_6.$
Then $m_6\neq b_6,\ \bb_6$, $n_6\neq b_6,\ \bb_6$ by {\bf Lemma
\ref{s2}} and $c_3\bar{m}_6$ has three possibilities
\[
\begin{array}{cl}
(I)&c_3\bar{m}_6=\bc_9+x_3+y_3+z_3,,\ x_3,\ y_3,\ z_3\ \mbox{distinct}\\
(II)&c_3\bar{m}_6=\bc_9+x_3+y_6,\ x_3\ \mbox{and}\ y_6\in B, \\
(III)&c_3\bar{m}_6=\bc_9+x_4+y_5 ,\ x_4\ \mbox{and}\ y_5\in B, .
\end{array}
\]
\par Suppose (I) follows. We may assume that $c_3\bar{x}_3=m_6+r_3$, where $r_3\in B$. Then
$$\bar{x}_3+\bar{x}_3b_8=\bar{x}_3(c_3\bc_3)=(\bar{x}_3c_3)\bc_3=m_6\bc_3+r_3\bc _3=c_9+\bar{x} _3+\bar{y}_3+\bar{z} _3+r_3\bc_3.$$
Hence $\bar{y}_3\in Supp\{\bar{x}_3b_8\}$, which implies that
$x_3=y_3$, a contradiction.

\par If (II) follows, then $(c_3\bar{x}_3,\ c_3\bar{x}_3)=2$, from which $x_3\bar{x}_3=1+b_8$ follows. Also we may set that $c_3\bar{m}_3=x_3+t_6$, where $t_6\in B$. Hence
\[
\begin{array}{rll}
(c_3\bc_3)\bar{x}_3&=&\bar{x}_3+\bar{x}_3b_8,\\
\bc_3(c_3\bar{x}_3)&=&\bc_3m_6+\bc_3m_3\\
&=&\bar{x}_4+\bar{y}_5+c_9+\bar{x}_3+\bar{t}_6.
\end{array}
\]
Therefore $x_3b_8=x_4+y_5+\bc_9+t_6$. We have that $(x_3b_8,\
x_3b_8)=4$. But by Lemma \ref{la3}, $(x_3b_8,\ x_3b_8)=(b_3b_8,\
b_3b_8)=3$, a contradiction.

\par If (III) follows, set $c_3\bar{x}_4=m_6+p_6$, where $p_6\in B$. By $\bc_3(c_3\bar{m}_6)=(c_3\bc_3)\bar{m}_6$,
one has that $$\bar{m}_6b_8=y_{15}+\bar{m}_6+\bar{n}_6+\bar{p}_6+\bc_3y_5.$$
Hence $(b_8y_{15},\ \bar{m}_6)=1$. By the expression of
$b_8y_{15}$ in {\bf Lemma \ref{s2}}, $(\bar{m}_6,\ \bb_3d_9)=1$.
But $(y_{15},\ \bb_3d_9)=1$. There exists $q_6\in \mathit{N}^* B$
such that
\begin{eqnarray}\label{cgy1}\bb_3d_9=y_{15}+\bar{m}_6+q_6.\end{eqnarray}
\par If $q_6\not\in B$, then $q=2z_3$ or $q_6=w_3+u_3$. If the first one holds, then $(\bb_3d_9,\ z_3)=2$ and $(b_3z_3,\ d_9)=2$, which is impossible. So the second one follows. Consequently  $b_3w_3=d_9$. Hence
$$w_3+w_3b_8=w_3(b_3\bb_3)=\bb_3(b_3w_3)=\bb_3d_9=y_{15}+\bar{m}_6+w_3+u_3.$$
We can see that $w_3\bar{u}_3=1+b_8$. Then $w_3=u_3$, a
contradiction. Hence $p_6\in B$. It follows by {\bf Lemma
\ref{s2}} and (\ref{cgy1}) that
$$b_8y_{15}=5y_{15}+2c_9+2b_6+c_3+\bar{m}_6+q_6.$$
\par Since
\[
\begin{array}{rll}
c_3(\bb_3d_9)&=&c_3y_{15}+c_3\bar{m}_6+c_3q_6\\
&=&\bb_6+\bc_9+2\bar{y}_{15}+x_4+y_5+\bc_9+c_3q_6,\\
(c_3\bb_3)d_9&=&b_3d_9+x_6d_9\\
&=&\bar{y}_{15}+m_6+\bar{q}_6+x_6d_9,
\end{array}
\]
we have that $(c_3q_6,\ m_6)\geq 1$. Hence $(c_3\bar{m}_6,\
\bar{q}_6)\geq 1$, a contradiction to our assumption. The is the
end of the proof.
\begin{lemma}
No NITA such that $(c_3c_9,\ c_3c_9)\geq 4$.
\end{lemma}
\par {\bf Proof.} The proof is given in three steps.
\par {\bf Step 1.} The following equations hold:
\[
\begin{array}{cl}
(i)&c_3c_9=\bar{y}_{15}+y_3+z_3+w_3+u_3,\\
(ii)&c_3c_9=\bar{y}_{15}+y_3+y_3+y_6,\\
(iii)&c_3c_9=\bar{y}_{15}+y_3+y_4+y_5,\\
(iv)&c_3c_9=\bar{y}_{15}+y_4+z_4+w_4.
\end{array}
\]
\par For any $z\in Supp\{c_3c_9\}\setminus\{\bar{y}_{15}\}$, it is enough to prove that
$(c_3c_9,\ z)=1$. Otherwise $(c_3c_9,\ z)\geq 2$. Since
$c_3c_9-\bar{y}_{15}$ is of degree 12, one has that $(c_3c_9,\
z)=2$ and $z$ is of degree 6. So $c_3c_9=\bar{y}_{15}+2z_6,\
c_3\bar{z}_6=2\bc_9$. Hence
$$\bar{z}_6+\bar{z}_6b_8=(c_3\bc_3)\bar{z}_6=\bc_3(c_3\bar{z}_6)=2\bc_3\bc_9=2y_{15}+4\bar{z}_6.$$
Thus $\bar{z}_6b_8=2y_{15}+3\bar{z}_6$. Consequently $(b_8y_{15},\
\bar{z}_6)=2$. By the expression of $b_6b_8$ in {\bf Lemma
\ref{s2}}, $z_6\neq b_6,\ \bb_6$. the expression of $b_8y_{15}$ in
{\bf Lemma \ref{s2}}, we have that
\[
\begin{array}{rll}
\bb_3d_9&=&y_{15}+2\bar{z}_6,\\
b_8y_{15}&=&5y_{15}+2c_9+2b_6+c_3+2\bar{z}_6,
\end{array}
\]
Therefore $b_3\bar{z}_6=2d_9.$
\par Since
\[
\begin{array}{rll}
(b_3\bb_3)d_9&=&d_9+b_8d_9,\\
b_3(\bb_3d_9)&=&b_3y_{15}+2b_3\bar{z}_6\\
&=&s_6+2t_{15}+5d_9,
\end{array}
\]
it follows that $b_8d_9=s_6+2t_{15}+4d_9.$ Consequently
\[
\begin{array}{rll}
b_8(b_3\bar{z}_6)&=&2b_8b_9\\
&=&2s_6+4t_{15}+8d_9,\\
(b_3b_8)\bar{z}_6&=&b_3\bar{z}_6+x_6\bar{z}_6+x_{15}\bar{z}_6\\
&=&2d_9+x_6\bar{z}_6+x_{15}\bar{z}_6.
\end{array}
\]
Then $x_6\bar{z}_6+x_{15}\bar{z}_6=2s_6+4t_{15}+6d_9.$ We have
exactly three possibilities:
\[
\begin{array}{cl}
\mbox{(I)}&x_6\bar{z}_6=4d_9,\ x_{15}\bar{z}_6=2s_6+4t_{15}+2d_9,\\
\mbox{(II)}&x_6\bar{z}_6=s_6+2t_{15},\ x_{15}\bar{z}_6=s_6+2t_{15}+6d_9,\\
\mbox{(III)}&x_6\bar{z}_6=2s_6+t_{15}+d_9,\
x_{15}\bar{z}_6=3t_{15}+5d_9.
\end{array}
\]
\par The first case implies that
\[
\begin{array}{rll}
\bb_3(x_6\bar{z}_6)&=&4\bb_3d_9\\
&=&4y_{15}+8\bar{z}_6,\\
(\bb_3x_6)\bar{z}_6&=&b_8\bar{z}_6+x_{10}\bar{z}_6\\
&=&2y_{15}+3\bar{z}_6+x_{10}\bar{z}_6.
\end{array}
\]
So $x_{10}\bar{z}_6=2y_{15}+5\bar{z}_6.$ Thus $(z_6\bar{z}_6,\
x_{10})=5$, a contradiction.
\par The second and the third cases mean $(z_6d_9,\ x_{15})=6$ and $(z_6d_9,\ x_{15})=5$ respectively. It is impossible. \\
\par {\bf Step 2.} No NITA such that (i), (ii) and (iii). \\
\par It is enough to show that $y_3\in Supp\{c_3c_9\}$. If $y_3\in Supp\{c_3c_9\}$, then $\bc_3y_3=c_9.$
Hence $(\bb_3y_3,\ \bb_3y_3)=(\bc_3y_3,\ \bc_3y_3)=1$ by {\bf Lemma
\ref{la1}}. Let $\bb_3y_3=g_9$, $g_9\in B$. Then
$$b_3c_9=b_3(\bc_3y_3)=(b_3\bc_3)y_3=\bb_3y_3+\bar{x}_6y_3=g_9+\bar{x}_6y_3.$$
Since $(t_{15},\ b_3c_9)=1$, we have that
$b_3c_9=g_9+t_{15}+\alpha_3$,
some $\alpha_3\in B$. Therefore $(b_3c_9,\ b_3c_9)=3$. But $(b_3c_9,\ b_3c_9)=(c_3c_9,\ c_3c_9)$ by {\bf Lemma \ref{la1}}. It is impossible for $(c_3c_9,\ c_3c_9)\geq 4$.\\
\par {\bf Step 3.} No NITA satisfies (iv).\\
\par If (iv) holds, then there exist $m_3$ such that $\bc_3y_4=c_9+m_3.$
Thus $c_3m_3=y_4+k_5$, some $k_5\in B$. So $(c_3m_3,\ c_3m_3)=2$,
which implies that $m_3\bar{m}_3=1+b_8$. Hence
$$\bar{m}_3y_4+\bar{m}_3k_5=c_3(m_3\bar{m}_3)=c_3+c_3b_8=(c_3\bc_3)\bar{m}_3=2c_3+b_6+y_{15}.$$
It is impossible for $c_3\in Supp\{\bar{m}_3y_4\}$ and
$Supp\{\bar{m}_3k_5\}.$ The lemma follows.
\begin{lemma}
No NITA such that $(c_3c_9,\ c_3c_9)=2$.
\end{lemma}
\par {\bf Proof.} If $(c_3c_9,\ c_3c_9)=2$, then $c_3c_9=\bar{y}_{15}+b_{12},\ b_{12}\in B$ and $c_3\bc_9=x_{10}+x_{17},\ x_{17}\in B$.
\par Since $\bc_3(c_3c_9)=\bc_3\bar{y}_{15}+\bc_3b_{12}=b_6+c_9+2y_{15}+\bc_3b_{12}$, we have that $b_8c_9=b_6+2y_{15}+\bc_3b_{12}.$ On the other hand, since $$c_9+b_8c_9=(c_3\bc_3)c_9=c_3(\bc_3c_9)=c_3(\bar{x}_{10}+\bar{x}_{17})=b_6+y_{15}+c_9+c_3\bar{x}_{17},$$
we have that $b_8c_9=b_6+y_{15}+c_3\bar{x}_{17}.$ By the expression
of $b_8y_{15}$ in {\bf Lemma \ref{s2}}, $(b_8c_9,\ y_{15})=2$. Hence
$y_{15}\in Supp\{c_3\bar{x}_{17}\}$, which implies that $x_{17}\in
Supp\{c_3\bar{y}_{15}\}$. By the expressions of $c_3b_8$ and
$c_3x_{10}$ in {\bf Lemma \ref{s2}}, $(c_3\bar{y}_{15},\
b_8)=(c_3\bar{y}_{15},\ x_{10})=1$. There exists
$t_{10}\in\mathit{N}^*B$ such that
$c_3\bar{y}_{15}=x_{17}+b_8+x_{10}+t_{10}$, some $t_{10}\in
\mathit{N}^*B$. The expression of $c_3y_{15}$ in {\bf  Lemma
\ref{s2}} means that $(c_3y_{15},\ c_3y_{15})=6$. Hence
$(c_3\bar{y}_{15},\ c_3\bar{y}_{15})=6$. So $t_{10}=x_3+y_3+z_4$. It
leads to $y_{15}\in Supp\{c_3\bar{y}_3\}$, a contradiction. The
lemma follows.
\begin{lemma}\label{dir2}$c_3c_9=d_3+\bc_9+\bar{y}_{15}$.
\end{lemma}
\par {\bf Proof.} At first by previous lemma in this section, we know that $(c_3c_9,\ c_3c_9)=3$ and
$c_3c_9$ is a sum of irreducible base elements of degree 3, 9 and
15. So we may assume that $c_3c_9=d_3+y_9+\bar{y}_{15}$, where
$d_3,\ y_9\in B$. Then $c_3\bar{d}_3=\bc_9$. Checking associative
law for $\bc_3(c_3\bar{d}_3)=(c_3\bc_3)\bar{d}_3$, it follows that
$$\bar{d}_3b_8=\bar{y}_9+y_{15}.$$
Since
\[
\begin{array}{rll}
b_6(b_3\bar{x}_6)&=&b_6b_8+b_6x_{10}\\
&=&c_3+b_6+c_9+2y_{15}+b_6x_{10},\\
(b_6\bar{x}_6)b_3&=&2b_3x_6+b_3b_9+b_3x_{15}\\
&=&2c_3+2y_{15}+b_3b_9+b_6+2y_{15}+c_9,
\end{array}
\]
we have that
$$b_6x_{10}=b_3b_9+c_3+2y_{15}.$$
\par Since $(c^2_3)\bar{d}_3=\bc_3\bar{d}_3+\bb_6\bar{d}_3$, $c_3(c_3\bar{d}_3)=c_3\bc_9$ and $x_{10}\in Supp\{c_3\bc_9\}$ by {\bf Lemma \ref{s2}}, we come to  $x_{10}\in Supp\{\bb_6\bar{d}_3\}$. So $(b_6x_{10},\ \bar{d}_3)=1$. Further $(b_3b_9, \bar{d}_3)=1.$ Thus
$$b_3d_3=\bb_9.$$ Therefore
$$\bb_3\bb_9=\bb_3(b_3d_3)=(\bb_3b_3)d_3=d_3+d_3b_8=d_3+y_9+\bar{y}_{15},$$
which implies that
\[
\begin{array}{rll}
b_3b_9&=&\bar{d}_3+\bar{y}_9+y_{15},\\
b_6x_{10}&=&\bar{d}_3+c_3+\bar{y}_9+3y_{15}.
\end{array}
\]
\par Since $(r_3b_3)d_3=c_9+\bb_6d_3$ and $(b_3d_3)r_3=r_3\bb_9$, we have that $c_9\in Supp\{r_3\bb_9\}$. But $y_{15}\in Supp\{r_3\bb_9\}$ by {\bf Lemma \ref{s2}}. Thus there exits $m_3\in B$ such that
$$r_3\bb_9=m_3+c_9+y_{15}.$$
Hence $b_9\in Supp\{r_3\bc_9\}$. But $(r_3\bc_9,\ x_{15})=1$ by {\bf
Lemma \ref{s2}}. Therefore
$$r_3\bc_9=n_3+b_9+x_{15},\ n_3\in B.$$
By $r_3$ real, it follows that
$$r_3n_3=\bc_9.$$  Further
\[
\begin{array}{rll}
r_3(r_3\bc_9)&=&r_3x_{15}+r_3b_9+r_3n_3\\
&=&\bb_6+\bc_9+2\bar{y}_{15}+\bc_9+\bar{y}_{15}+\bar{m}_3+\bc_9,\\
r^2_3\bc_9&=&\bc_9+b_8\bc_9,
\end{array}
\]
so that
$$b_8\bc_9=\bb_6+2\bc_9+3\bar{y}_{15}+\bar{m}_3.$$
 Hence
$(b_8y_{15},\ c_9)=3$.
\par By the expression of $b_8y_{15}$ in {\bf Lemma \ref{s2}}, we have $(\bb_3d_9,\ c_9)=1$.
From the expression of $b_3t_{15}$ in {\bf Lemma \ref{s2}}, one has
that $(t_{15},\ b_3c_9)=1$. Hence there exists $u_3\in B$ such that
$$b_3c_9=u_3+d_9+t_{15}.$$
Thus $\bb_3u_3=c_9.$ Moreover
$$c_3b_9=(\bb_3\bar{d}_3)c_3=\bb_3(c_3\bar{d}_3)=\bb_3\bc_9=\bar{u}_3+d_9+t_{15}.$$
Thus $u_3c_3=\bb_9$.
\par Since
\[
\begin{array}{rll}
\bc_3(b_3\bar{u}_3)&=&\bc_3\bc_9\\
&=&\bar{d}_3+\bar{y}_9+y_{15},\\
(b_3\bc_3)\bar{u}_3&=&\bb_3\bar{u}_3+\bar{x}_6\bar{u}_3,
\end{array}
\]
we have that
$$b_3u_3=y_9\  \mbox{and } u_3x_6=\bar{y}_{15}+\bar{d}_3.$$
Hence
$$c_3y_9=c_3(b_3u_3)=b_3(u_3c_3)=b_3\bb_9.$$
Since $(x_{10},\ b_3\bb_9)=1$ by {\bf Lemma \ref{s2}}, one has that
$(x_{10},\ c_3y_9)=1$, so $(\bar{y}_9,\ c_3x_{10})=1$. Therefore
$\bar{y}_9=c_9$ by the expression of $c_3x_{10}$ (see {\bf Lemma
\ref{s2}}. The lemma follows.

\begin{lemma}\label{subaofc3}NITA generated by $c_3$ is isomorphic to the algebra of characters
of $3\cdot A_6$: $(Ch(3\cdot A_6),\ Irr(3\cdot A_6))$. Further we
have all equations of products of base elements in table subset $D$
listed in \textbf{Section 2}.
\end{lemma}

\par {\bf Proof.}  By {\bf Lemma \ref{s2}} and \textbf{Lemma \ref{dir2}} the
 sub-algebra generated by $c_3$ satisfies the following equations:
\[
\begin{array}{rll}
c_3\bc_3=1+b_8,& c_3^2=\bc_3+\bb_6,& c_3b_6=\bc_3+\bar{y}_{15},\\
c_3\bb_6=b_8+x_{10},\ x_{10}\mbox{ real}, &
c_3x_{10}=b_6+c_9+y_{15},&c_3y_{15}=\bb_6+\bc_9+2\bar{y}_{15}.
\end{array}
\]
Changing $c_3$ into $b_3$, $b_6$ into $\bb_6$, $y_{15}$ into
$b_{15}$, $x_{10}$ into $b_{10}$, we can see above equations are the
same as those in {\bf Hypothesis 5.1} and {\bf Lemma
\ref{lemma4.1}}. Hence NITA generated by $c_3$ is exactly the table
algebra of characters of $3\cdot A_6$ by {\bf Theorem
\ref{theorem4.1}}. Of course we have equations of products of base
elements in  table subset $D$ listed in \textbf{Section 2}.
 \\

\par Now let's continue our calculation, we shall come to a new NITA which is not derived from groups.\\
\begin{theorem}
If $c_3$ is non-real, $(b_3b_8,\ b_3b_8)=3$. Then  NITA generated by
$b_3$ with $b^2_3=c_3+b_6$ is isomorphic to the NITA in $(A(3\cdot
A_6\cdot 2),\ B_{32})$, the NITA of dimension 32 in \textbf{Section
2}.
\end{theorem}

\par {\bf Proof.}  By {\bf Lemma \ref{subaofc3}}, we may say that we have
found a table subset: $ D=\{1,\ b_8,\ x_{10},\ b_5,\ c_5,\ c_8,\
x_9,\ c_3,\ \bc_3,\ d_3,\ \bar{d}_3,\ c_9,\ \bc_9,\ b_6,\ \bb_6,\
y_{5},\ \bar{y}_{15}\}.$ To prove NITA is isomorphic to NITA of
dimension 32 in  {\bf Section 2}, it is enough to find out all basis
elements outside $D$ and remaining expressions of
 all products of basis elements. Of course at this moment we have found a lot of them in
 \textbf{ Lemma \ref{s1}, \ref{s2} and \ref{dir2}}. So what we do is
 to find out others. In order to make the proof easily read, we
 divide the proof into several steps:\\
\par {\bf Step 1.} $b_3b_9=y_{15}+c_9+\bar{d}_3$,  $b_3d_3=\bb_9$,
$b_3\bb_9=x_9+c_8+x_{10}$, $b_3d_9=\bc_9+d_3+\bar{y}_{15}$,
$b_3\bar{d}_3=d_9$, $b_3\bc_9=\bb_9+\bar{x}_{15}+z_3$,
$b_3\bar{z}_3=c_9$, $b_3\bar{x}_{15}=b_5+c_5+b_8+c_8+x_9+x_{10}$, $b_3z_3=x_9$, $b_3y_3=\bc_9$.\\

\par Since by Lemma \ref{s2}
\[
\begin{array}{rll}
b^2_3x_{10}&=&c_3x_{10}+b_6x_{10}\\
&=&b_6+c_9+y_{15}+3y_{15}+c_3+\bar{d}_3+c_9,\\
(b_3x_{10})b_3&=&x_6b_3+x_{15}b_3+b_9b_3\\
&=&c_3+y_{15}+b_6+2y_{15}+c_9+b_3b_9,
\end{array}
\]
we have that $b_3b_9=y_{15}+c_9+\bar{d}_3$. Further
$b_3d_3=\bb_9.$ So
$$b_3\bb_9=(b_3d_3)b_3=b^2_3d_3=c_3d_3+b_6d_3=x_9+c_8+x_{10}.$$
\par Set $b_3\bar{d}_3=y_9$. Since $b_3d_3=\bb_9,$ then $y_9\in B$. And
$$b_3y_9=(b_3\bar{d}_3)b_3=b^2_3\bar{d}_3=c_3\bar{d}_3+b_6\bar{d}_3=\bc_9+d_3+\bar{y}_{15}.$$
We have that $(b_3y_{15},\ \bar{y}_9)=1.$ But
$b_3y_{15}=s_6+2t_{15}+d_9$, thus $d_9=\bar{y}_9$ and
$\bb_3{d}_3=d_9$.
\par Since $(b_3x_{15},\ c_9)=(b_3b_9,\ c_9)=1$, there exists $z_3\in B$ such that
$$b_3\bc_9=\bb_9+\bar{x}_{15}+z_3.$$
Hence $b_3\bar{z}_3=c_9.$
\par Since
\[
\begin{array}{rll}
b^2_3\bb_6&=&c_3\bb_6+b_6\bb_6\\
&=&b_8+x_{10}+1+b_8+b_5+c_5+c_8+x_9,\\
(b_3\bb_6)b_3&=&b_3\bb_3+b_3\bar{x}_{15}\\
&=&1+b_8+b_3\bar{x}_{15},
\end{array}
\]
we have that
$$b_3\bar{x}_{15}=b_5+c_5+b_8+c_8+x_9+x_{10}.$$ Now we have
\[
\begin{array}{rll}
b^2_3\bc_9&=&c_3\bc_9+b_6\bc_9\\
&=&c_8+x_9+x_{10}+2x_9+b_8+x_{10}+c_8+b_5+c_5,\\
(b_3\bc_9)b_3&=&b_3\bb_9+b_3\bar{x}_{15}+b_3z_3\\
&=&c_8+x_9+x_{10}+b_5+c_5+b_8+c_8+x_9+x_{10}+b_3z_3,
\end{array}
\]
which imply $b_3z_3=x_9.$ \\

\par {\bf Step 2.} $b_3c_9=t_{15}+d_9+y_3$, $y_3$ real, $d_3\bar{x}_6=t_{15}+\bar{y}_3,$  $c_3b_9=t_{15}+d_9+y_3$, $c_3y_3=\bb_9$,$b_3y_3=\bc_9$ and
$b_8t_{15}=r_3+\bar{y}_3+5t_{15}+2s_6+3d_9$.

\par We have that
$$b_3c_9=b_3(d_3\bc_3)=(b_3\bc_3)d_3=\bb_3d_3+\bar{x}_3d_3=d_9+\bar{x}_6d_3.$$
Hence $(b_3c_9,\ d_9)=1$. But by Lemma \ref{s2} $(b_3t_{15},\
\bc_9)=1$. Then there exists $y_3\in B$ such that
$$b_3c_9=t_{15}+d_9+y_3,\ d_3\bar{x}_6=t_{15}+\bar{y}_3.$$
Consequently $b_3\bar{y}_3=\bc_9$. Because
$b_3c_9=(d_3\bc_3)b_3=(b_3d_3)\bc_3=\bc_3\bb_9$, we have that
$$c_3b_9=t_{15}+d_9+y_3.$$
Consequently $c_3\bar{y}_3=\bb_9$. Since
\[
\begin{array}{rll}
b^2_3c_9&=&c_3c_9+b_6c_9\\
&=&d_3+\bc_9+\bar{y}_{15}+\bb_6+2\bc_9+2\bar{y}_{15},\\
(b_3c_9)b_3&=&b_3t_{15}+b_3d_9+b_3y_3\\
&=&\bb_6+\bc_9+2\bar{y}_{15}+d_3+\bc_9+\bar{y}_{15}+b_3y_3,
\end{array}
\]
then $b_3y_3=\bc_9.$
\par Because we have that
\[
\begin{array}{rll}
(b_3\bb_3)t_{15}&=&t_{15}+b_8t_{15}\\
(b_3t_{15})\bb_3&=&\bb_3\bb_6+\bb_3\bc_9+2\bb_3\bar{y}_{15}\\
&=&r_3+t_{15}+d_9+\bar{y}_3+t_{15}+2s_6+4t_{15}+2d_9,
\end{array}
\]
so $b_8t_{15}=r_3+\bar{y}_3+5t_{15}+2s_6+3d_9$, which implies that $y_3$ is real.\\
\par {\bf Step 3.}  $b_3d_9=\bc_9+\bar{y}_{15}+d_3$, $r_3t_{15}=b_5+c_5+b_8+c_8+x_9+x_{10},
r_3d_9=c_8+x_9+x_{10}$.

\par Since
\[
\begin{array}{rll}
b^2_3y_{15}&=&c_3y_{15}+b_6y_{15}\\
&=&\bb_6+\bc_9+2\bar{y}_{15}+\bc_3+d_3+\bb_6+2\bc_9+4\bar{y}_{15},\\
(b_3y_{15})b_3&=&b_3s_6+2b_3t_{15}+b_3d_9\\
&=&\bc_3+\bar{y}_{15}+b_3d_9+2\bb_6+2\bc_9+4\bar{y}_{15},
\end{array}
\]
then $b_3d_9=\bc_9+\bar{y}_{15}+d_3$.
\par Since $r^2_3b_8=(r_3b_8)r_3$ and $r^2_3x_{10}=(r_3x_{10})r_3$, we have that
\[
\begin{array}{rll}
r_3t_{15}&=&b_5+c_5+b_8+c_8+x_9+x_{10},\\
r_3d_9&=&c_8+x_9+x_{10}.
\end{array}
\]
\par {\bf Step 4.} $r_3d_3=b_9$ and $r_3b_9=d_3+\bc_9+\bar{y}_{15}$.\\
\par There are three possibilities: $r_3d_3=t_9,\ m_3+n_6,\ m_4+n_5.$ Thus
\[
\begin{array}{rll}
r^2_3d_3&=&d_3+d_3b_8\\
&=&d_3+\bc_9+\bar{y}_{15},\\
(r_3d_3)r_3&=&r_3t_9\\
&\mbox{or}&m_3r_3+n_6r_3\\
&\mbox{or}&m_4r_3+n_5r_3,
\end{array}
\]
\par If $r_3d_3=m_3+n_6$, then $m_3r_3=\bc_9$ and $r_3n_6=d_3+\bar{y}_{15}.$ But
$r_3\bar{y}_{15}=x_6+2x_{15}+b_9.$ Hence $n_6=x_6$, so that
$(r_3x_6,\ d_3)=1$, a contradiction.
\par If $r_3d_3=m_4+n_5$, then $m_4r_3=d_3+\bc_9$ and $r_3n_5=\bar{y}_{15}$. But
$(r_3\bar{y}_{15},\ n_5)=0$,  a contradiction.
Hence $r_3d_3=t_9$ and $r_3t_9=d_3+\bc_9+\bar{y}_{15}$, which implies that $(r_3\bar{y}_{15},\ t_9)=1$. Then $t_9=b_9$.\\
\par {\bf Step 5.} $x_6y_{15}=r_3+4t_{15}+s_6+2d_9+y_3,$
$b_8x_{15}=b_3+\bar{z}_3+2x_6+5x_{15}+3b_9,$
$b_8b_9=\bar{z}_3+x_6+2b_9+3x_{15},$
$b_8d_9=y_3+s_6+2d_9+3t_{15}$, $y_3b_8=d_9+t_{15}$ and
$b_6b_9=s_6+2d_9+2t_{15}$
and $b_6x_{15}=r_3+s_6+y_3+2d_9+4t_{15}$.\\
\par The above equations follow from $b_3x^2_6=(b_3x_6)x_6$, $(b_3\bb_3)x_{15}=(b_3x_{15})\bb_3$, $(b_3\bb_3)b_9=b_3b_9)\bb_3$, $(b_3\bb_3)d_9=(b_3d_9)\bb_3$, $(b_3\bb_3)y_3=(b_3y_3)\bb_3$, $b^2_3b_9=(b_3b_9)b_3$
and $b^2_3x_{15}=(b_3x_{15})b_3$.
\par \par {\bf Step 6.} $c_3d_9=\bb_9+z_3+\bar{x}_{15},$ $b_6d_9=\bar{x}_6+3\bb_9+2\bar{x}_{15}$, $z_3\bar{c}_3=d_9$.
\par Since
\[
\begin{array}{rll}
b^2_3d_9&=&c_3d_9+b_6d_9,\\
(b_3d_9)b_3&=&b_3d_3+b_3\bc_9+b_3\bar{y}_{15}\\
&=&2\bb_9+3\bar{x}_{15}+z_3+\bar{x}_6,
\end{array}
\]
and $(c_3x_{15},\ d_9)=1$, which the {\bf Step 6} follows from.\\
\par {\bf Step 7.} $b_8x_{15}=b_3+\bar{z}_3+2x_6+5x_{15}+3b_9$, $b_3x_9=\bar{z}_3+b_9+x_{15}$, $b_3c_8=x_{15}+b_9$, $b_3b_5=x_{15}$ and $b_3c_5=x_{15}.$\\
\par By \textbf{Lemma 7.2} and \textbf{Step 1}, we have
\begin{center}$(b_3x_{15})\bb_3=b_6\bb_3+2y_{15}\bb_3+c_9\bb_3=b_3+x_{15}+x_6+2x_{15}+b_9+b_9+x_{15}+\bar{z}_3$
$(b_3\bb_3)x_{15}=x_{15}+b_8x_{15}$,\end{center} hence
$$b_8x_{15}=b_3+\bar{z}_3+2x_6+5x_{15}+3b_9$$
\par Consequently we have the following equations:
\[
\begin{array}{rll}
b^2_8b_3&=&b_3+2b_3b_8+2b_3x_{10}+b_3b_5+b_3c_5+b_3c_8+b_3x_9\\
&=&3b_3+2x_6+2x_{15}+2x_6+2x_{15}+2b_9+b_3b_5+b_3c_5+b_3c_8+b_3x_9,\\
(b_3b_8)b_8&=&b_3b_8+x_6b_8+x_{15}b_8\\
&=&b_3+x_6+x_{15}+b_3+x_6+b_9+2x_{15}+b_3+\bar{z}_3+2x_6+5x_{15}+3b_9,
\end{array}
\]
 which imply that $$b_3b_5+b_3c_5+b_3c_8+b_3x_9=4x_{15}+2b_9+\bar{z}_3.$$
By the expression of $b_3\bar{x}_{15}$ in {\bf Step 3}, it can be easily show that {\bf Step 7} follows.\\
\par {\bf Step 8.} $r_3b_5=t_{15}$, $r_3c_5=t_{15}$, $r_3c_8=t_{15}+d_9$, $r_3x_9=t_{15}+d_9+y_3$ and $r_3y_3=x_9$.\\
\par Since
\[
\begin{array}{rll}
b^2_8r_3&=&r_3+2r_3b_8+2r_3x_{10}+r_3b_5+r_3c_5+r_3c_8+r_3x_9\\
&=&r_3+2r_3+2s_6+2t_{15}+2s_6+2t_{15}+2d_9+r_3b_5+r_3c_5+r_3c_8+r_3x_9,\\
(r_3b_8)b_8&=&r_3b_8+b_8s_6+b_8t_{15}\\
&=&r_3+s_6+t_{15}+r_3+s_6+2t_{15}+d_9+r_3+y_3+5t_{15}+2s_6+3d_9,
\end{array}
\]
it holds that
$$r_3b_5+r_3c_5=r_3c_8+r_3x_9=4t_{15}+2d_9+y_3.$$
But $(r_3t_{15},\ x_9)=1$. The first four equations follow from above equation. Consequently, $r_3y_3=x_9$.\\
\par {\bf Step 9.} $r_3c_9=\bb_9+\bar{x}_{15}+z_3$, $r_3z_3=c_9$, $x_6\bar{y}_{15}=\bb_3+z_3+\bar{x}_6+4\bar{x}_{15}+2\bb_9$,
$x_6c_9=2t_{15}+2d_9+s_6$, $x_6x_{10}=b_9+b_3+\bar{z}_3+3x_{15}$, $x_6\bar{x}_6=1+b_8+b_5+c_5+c_8+x_9$, $x_6s_6=2\bb_6+\bc_9+\bar{y}_{15}$ and $x_6\bar{c}_9=2\bar{x}_{15}+\bar{x}_6+2\bb_9.$\\
\par Since $r^2_3b_9=(r_3b_9)r_3$, we have the first equation, which the second equation follows from. The rest of equations follow from associativity of $(c_3b_8)\bar{x}_6=b_8(c_3\bar{x}_6)$, $(d_3b_8)\bar{x}_6=b_8(d_3\bar{x}_6)$, $(b_3\bar{x}_6)x_6=(b_3x_6)\bar{x}_6$, $(b_3\bc_3)x_6=(b_3x_6)\bc_3$, $(b_3c_3)x_6=b_3(c_3x_6)$ and $(r_3c_3)c_9=r_3(c_3c_9)$ respectively.\\
\par {\bf Step 10.} $x_6y_3=d_3+\bar{y}_{15}.$\\
\par By known equations, we have that $(x_6y_3,\ \bar{y}_{15})=(x_6y_{15},\ y_3)=1$ and $(x_6y_3,\ d_3)=(x_6\bar{d}_3,\ y_3)=1$, which {\bf Step 10} follows from.\\
\par {\bf Step 11.} $x_6b_5=x_6+b_9+x_{15}$, $x_6c_5=x_6+b_9+x_{15}$, $x_6c_8=x_6+b_9+2x_{15}+z_3$ and $z_3x_6=x_{10}+c_8$.\\
\par Since
\[
\begin{array}{rll}
b^2_8x_6&=&x_6+2x_6b_8+2x_6x_{10}+x_6b_5+x_6c_5+x_6c_8+x_6x_9\\
&=&x_6+2b_3+2x_6+2b_9+4x_{15}+2b_9+2b_3+2\bar{z}_3+6x_{15}+\\
&&x_6b_5+x_6c_5+x_6c_8+x_6+2b_9+2x_{15},\\
(b_8x_6)b_8&=&b_3b_8+x_6b_8+b_9b_8+2b_8x_{15}\\
&=&b_3+x_6+x_{15}+b_3+x_6+b_9+2x_{15}+\bar{z}_3+x_6+2b_9+\\
&&+3x_{15}+2b_3+2\bar{z}_3+4x_6+10x_{15}+6b_9,
\end{array}
\]
we have that
$$x_6b_5+x_6c_5+x_6c_8=3x_6+3b_9+\bar{z}_3+4x_{15}.$$
By the expression of $x_6\bar{x}_6$ in {\bf Step 9}, we have that $(x_6b_5,\ x_6)=(x_6c_5,\ x_6)=(x_6c_8,\ x_6)=1$. The first three equations of {\bf Step 11} follows by comparing the degrees. Hence $(x_6z_3,\ c_8)=1$. In addition to $(z_3x_6,\ x_{10})=(x_6x_{10},\ \bar{z}_3)=1$, which the last equation follows from.\\
\par {Step 11.} $x_6x_9=x_6+2b_9+2x_{15}$, $x_6d_9=\bb_6+2\bc_9+2\bar{y}_{15}$, $x_6d_3=z_3+\bar{x}_{15}$ and $x_6b_9=2\bc_9+\bb_6+2\bar{y}_{15}$.\\
\par It follows from $(b_3z_3)x_6=b_3(z_3x_6)$, $(z_3\bc_3)x_6=\bc_3(z_3x_6)$, $(r_3c_3)\bar{d}_3=r_3(c_3\bar{d}_3)$ and $(r_3c_3)\bb_9=(r_3\bb_9)c_3$ respectively.\\
\par {\bf Step 12.} Based on the equations known by us, we can check the inner products of some products of basis elements and come to the following table:
\[
\begin{array}{|l|l|}
\hline
\mbox{values of inner products}& \mbox{new equations}\\
\hline &\\
(z_3b_9,\ b_8)=(b_8b_9,\ \bar{z}_3)=1&\\
(z_3x_{15},\ b_8)=(b_8x_{15},\ \bar{z}_3)=1&z_3b_8=\bb_9+\bar{x}_{15}\\
\hline
& \\
(x_6x_{10},\ b_9)=(x_6c_5,\ b_9)=(x_6b_5,\ b_9)=1&\\
(x_6b_8,\ b_9)=(x_6c_8,\ b_9)=1,\ (x_6x_9,\ b_9)=2&x_6\bb_9=b_5+c_5+x_{10}+b_8+c_8+2x_9\\
\hline &\\
(x_6b_8,\ x_{15})=2,\ (x_6x_{10},\ x_{15})=3,&\\
(x_6b_5,\ x_{15})=(x_6c_5,\ x_{15})=1&x_6\bar{x}_{15}=b_5+c_5+2b_8+3x_{10}+2c_8+2x_9\\
(x_6c_8,\ x_{15})=(x_6x_9,\ x_{15})=2&\\
\hline
 &\\
(x_6c_3,\ t_{15})=(x_6\bar{d}_3,\ t_{15})=(x_6b_6,\ t_{15})=1&\\
(x_6c_9,\ t_{15})=2,\ (x_6y_{15},\ t_{15})=4&x_6t_
{15}=\bc_3+d_3+2\bc_9+4\bar{y}_{15}+\bb_6\\
(x_6\bar{z}_3,\ \bar{d}_3)=(x_6\bar{y}_{15},\ z_3)=1&x_6\bar{z}_3=\bar{d}_3+y_{15}\\
(x_6d_3,\ \bar{x}_{15})=(x_6\bb_6,\ \bar{x}_{15})=(x_6\bc_3,\ \bar{x}_{15})=1&\\
(x_6\bc_9,\ \bar{x}_{15})=2,\ (x_6\bar{y}_{15},\ \bar{x}_{15})=4&x_6x_{15}=4y_{15}+b_6+2c_9+d_3+c_3\\
\hline
\end{array}
\]

\par {\bf Step 13.} Repeatedly considering the associativity of basis elements, we can come to new equations. Since the process is trivial and repeated, the concrete process are omitted and corresponce between associativity of basis elements and new equations are listed in following table. Occasionally we need check the inner products based on the old and new equations to obtain more new equation, at this moment, an explanation is given:

\[\begin{array}{|l|l|}  \hline
\mbox{associativity checked} &\mbox{new equations}\\
\mbox{or inner products}&\\
\hline &\\
(b_3c_3)\bd_3=b_3(c_3\bd_3) &d_3s_6=\bar{z}_3+x_{15}\\

(b_3c_3)y_3=b_3(c_3y_3)&y_3s_6=c_8+x_{10}\\

(b_3c_3)z_3=(b_3z_3)c_3&z_3s_6=\bar{d}_3+y_{15}\\

(b_3c_3)\bar{z}_3=b_3(c_3\bar{z}_3)&\bar{z}_3s_6=d_3+\bar{y}_{15}\\

(b_3c_3)s_6=b_3(c_3s_6)&s^2_6=1+b_5+c_5+b_8+c_8+x_9\\

(b_3c_3)b_9=b_3(c_3b_9)&b_9s_6=\bb_6+2\bc_9+2\bar{y}_{15}\\

(b_3c_3)x_{15}=b_3(c_3x_{15})&s_6x_{15}=\bc_3+d_3+\bb_6+2\bc_9+4\bar{y}_{15}\\

(b_3c_3)t_{15}=b_3(c_3t_{15})&s_6t_{15}=b_5+c_5+2b_8+2c_8+2x_9+3x_{10}\\

(b_3c_3)d_9=b_3(c_3d_9)&s_6d_9=b_5+c_5+2x_9+b_8+c_8+x_{10}\\

(b_3c_3)x_{10}=b_3(c_3x_{10})&s_6x_{10}=r_3+y_3+d_9+3t_{15}\\

(b_3c_3)b_5=b_3(c_3b_5)&s_6b_5=s_6+d_9+t_{15}\\

(b_3c_3)c_5=b_3(c_3c_5)&s_6c_5=s_6+t_{15}+d_9\\

(b_3c_3)c_8=b_3(c_3c_8)&s_6c_8=y_3+s_6+d_9+2t_{15}\\

(b_3c_3)x_9=b_3(c_3x_9)&s_6x_9=s_6+2d_9+2t_{15}\\

(b_3c_3)c_9=b_3(c_3c_9)&s_6c_9=\bar{x}_6+2\bb_9+2\bar{x}_{15}\\
\hline & \\
(s_6b_3,\ \bar{y}_{15})=1\\
(s_6b_9,\ \bar{y}_{15})=2,\ (s_6x_{15},\ \bar{y}_{15})=4&\\
(s_6\bar{z}_3,\ \bar{y}_{15})=(s_6x_6,\ \bar{y}_{15})=1&s_6y_{15}=2\bb_9+4\bar{x}_{15}+z_3+\bar{x}_6+\bar{b}_3\\
\hline
\end{array}
\]
\[
\begin{array}{|l|l|}

\hline
\mbox{associativity checked} &\mbox{new equations}\\
\mbox{or inner products}& \\
\hline &\\
(b_3d_3)b_9=d_3(b_3b_9)&b_9\bb_9=1+b_5+c_5+2b_8+2c_8+2x_9+2x_{10}\\
(b_3d_3)\bb_9=d_3(b_3\bb_9)&b^2_9=\bar{d}_3+c_3+2b_6+2c_9+3y_{15}\\
(b_3d_3)\bar{x}_{15}=d_3(b_3\bar{x}_{15})&b_9x_{15}=6y_{15}+3c_9+\bar{d}_3+2b_6+c_3\\
(b_3d_3)x_{15}=d_3(b_3x_{15})&b_9\bar{x}_{15}=2b_5+2c_5+3c_8+3b_8+4x_{10}+3x_9\\
(b_3d_3)t_{15}=d_3(b_3t_{15})&b_9t_{15}=d_3+\bc_3+6\bar{y}_{15}+3\bc_9+2\bb_6\\
(b_3d_3)d_9=d_3(b_3d_9)&b_9d_9=\bc_3+d_3+2\bb_6+2\bc_9+3\bar{y}_{15}\\
(b_3d_3)y_3=d_3(b_3y_3)&b_9y_3=\bc_3+\bar{y}_{15}+\bc_9\\
(b_3d_3)z_3=d_3(b_3z_3)&b_9\bar{z}_3=c_9+c_3+y_{15}\\
(b_3d_3)\bar{z}_3=d_3(b_3\bar{z}_3)&b_9z_3=b_8+x_9+x_{10}\\
(b_3d_3)x_{10}=b_3(d_3x_{10})&b_9x_{10}=\bar{z}_3+b_3+2b_9+x_6+4x_{15}\\
(b_3d_3)b_5=b_3(d_3b_5)&\bb_9c_9=b_3+\bar{z}_3+2x_6+3x_{15}+2b_9\\
(b_3d_3)c_5=b_3(d_3c_5)&b_9c_9=r_3+y_3+2s_6+2d_9+3t_{15}\\
(b_3d_3)b_6=b_3(d_3b_6)&\bb_9b_6=2b_9+2x_{15}+x_6\\
(b_3d_3)c_8=b_3(d_3c_8)&\bb_9y_{15}=6x_{15}+b_3+\bar{z}_3+2x_6+3b_9\\
(b_3d_3)x_9=b_3(d_3x_9)&\bb_9x_9=2\bb_9+3\bar{x}_{15}+z_3+\bb_3+2\bar{x}_6\\
(b_3d_3)\bar{y}_{15}=b_3(d_3\bar{y}_{15})&\bb_9\bar{y}_{15}=r_3+y_3+2s_6+6t_{15}+3d_9\\
(b_3d_3)d_3=b_3(d^2_3)&\bb_9d_3=r_3+t_{15}+d_9\\
(b_3d_3)c_3=b_3(d_3c_3)&\bb_9c_3=\bar{z}_3+b_9+x_{15}\\
(b_3d_3)\bar{d}_3=b_3(d_3\bar{d}_3)&b_9d_3=\bb_3+\bb_9+\bar{x}_{15}\\
(b_3d_3)c_3=b_3(d_3c_3)&b_9\bc_3=z_3+\bb_9+\bar{x}_{15}\\
\hline
\end{array}
\]
\[
\begin{array}{|l|l|}
\hline
\mbox{associativity checked}&\mbox{new equations}\\
\mbox{or inner products}
& \\
\hline &\\
(b_3b_5)x_{15}=(b_3x_{15})b_5&x^2_{15}=2c_3+2\bar{d}_3+4b_6+6c_9+9y_{15}\\
(b_3b_5)\bar{x}_{15}=(b_3\bar{x}_{15})b_5&x_{15}\bar{x}_{15}=1+6x_9+3b_5+3c_5+5b_8+\\
&+5c_8+6x_{10}\\
(b_3b_5)t_{15}=(b_3t_{15})b_5&x_{15}t_{15}=4\bb_6+6\bc_9+9\bar{y}_{15}+\\
&+2\bc_3+2d_3\\
(b_3b_5)d_9=(b_3d_9)b_5&x_{15}d_9=\bc_3+d_3+2\bb_6+3\bc_9+6\bar{y}_{15}\\
(b_3b_5)z_3=(b_3z_3)b_5&x_{15}z_3=b_5+c_5+b_8+c_8+x_9+x_{10}\\
(b_3b_5)b_5=b_3b_5^2&x_{15}b_5=\bar{z}_3+b_3+x_6+2b_9+3x_{15}\\
(b_3b_5)x_{10}=(b_3x_{10})b_5&x_{15}x_{10}=b_3+\bar{z}_3+3x_6+4b_9+6x_{15}\\
(b_3b_5)c_5=(b_5c_5)b_3&x_{15}c_5=b_3+\bar{z}_3+x_6+2b_9+3x_{15}\\
(b_3b_5)c_8=(b_5c_8)b_3&x_{15}c_8=b_3+\bar{z}_3+2x_6+3b_9+5x_{15}\\
(b_3b_5)x_9=(b_5x_9)b_3&x_{15}x_9=b_3+\bar{z}_3+2x_6+3b_9+6x_{15}\\
(b_3b_5)d_3=(b_5d_3)b_3&x_{15}d_3=\bar{x}_6+2\bar{x}_{15}+\bb_9\\
(b_3b_5)\bar{d}_3=(b_5\bar{d}_3)b_3&x_{15}d_3=s_6+d_9+2t_{15}\\
(b_3b_5)c_9=(b_5c_9)b_3&x_{15}c_9=r_3+y_3+2s_6+6t_{15}+3d_9\\
(b_3b_5)\bc_9=(b_5\bc_9)b_3&x_{15}\bb_9=\bb_3+z_3+2\bar{x}_6+3\bb_9+6\bar{x}_{15}\\
(b_3b_5)\bb_6=(b_5\bb_6)b_3&x_{15}\bb_6=\bb_3+z_3+2\bb_9+4\bar{x}_{15}+\bar{x}_6\\
(b_3b_5)y_{15}=(b_5y_{15})b_3&x_{15}y_{15}=2y_3+2r_3+4s_6+6d_9+9t_{15}\\
\hline
\end{array}
\]

\[
\begin{array}{|l|l|}
\hline
\mbox{associativity checked} &\mbox{new equations}\\
\mbox{or inner products}& \\
\hline &\\
(b_3b_6)t_{15}=(b_3t_{15})b_6&t^2_{15}=1+5b_8+5c_8+3b_5+3c_5+6x_{10}+6x_9\\
(b_3b_6)d_9=(b_3d_9)b_6&d_9t_{15}=2b_5+2c_5+3b_8+3c_8+3x_9+4x_{10}\\
(b_3b_6)y_3=(b_3y_3)b_6&y_3t_{15}=b_5+c_5+b_8+c_8+x_9+x_{10}\\
(b_3b_6)z_3=(b_3z_3)b_6&z_3t_{15}=b_6+c_9+2y_{15}\\
(b_3b_6)x_{10}=(b_6x_{10})b_3&t_{15}x_{10}=r_3+y_3+4d_9+3s_6+6t_{15}\\
(b_3b_6)b_5=(b_6b_5)b_3&t_{15}b_5=r_3+y_3+s_6+2d_9+3t_{15}\\
(b_3b_6)c_5=(b_6c_5)b_3&t_{15}c_5=r_3+y_3+s_6+2d_9+3t_{15}\\
(b_3b_6)c_8=(b_6c_8)b_3&t_{15}c_8=r_3+y_3+2s_6+3d_9+5t_{15}\\
(b_3b_6)x_9=(b_6x_9)b_3&t_{15}x_9=r_3+y_3+3d_9+2s_6+6t_{15}\\
(b_3b_6)d_3=(b_6d_3)b_3&t_{15}d_3=x_6+b_9+2x_{15}\\
(b_3b_6)c_9=(b_6c_9)b_3&t_{15}c_9=\bb_3+z_3+3\bb_9+2\bar{x}_6+6\bar{x}_{15}\\
(b_3b_6)b_6=b_3b^2_6&t_{15}b_6=\bb_3+z_3+\bar{x}_6+2\bb_9+4\bar{x}_{15}\\
(b_3b_6)y_{15}=(b_6y_{15})b_3&t_{15}y_{15}=2\bb_3+2z_3+4\bar{x}_6+6\bb_9+9\bar{x}_{15}\\
\hline
\end{array}
\]

\[
\begin{array}{|l|l|}
\hline
\mbox{associativity checked} &\mbox{new equations}\\
\mbox{or inner products}&\\
\hline &\\
(b_3\bar{d}_3)d_9=(b_3d_9)\bar{d}_3&d^2_9=1+b_5+c_5+2b_8+2c_8+2x_9+2x_{10}\\
(b_3\bar{d}_3)x_{10}=(b_3x_{10})\bar{d}_3&d_9x_{10}=y_3+r_3+s_6+2d_9+4t_{15}\\
(b_3\bar{d}_3)b_5=(b_3b_5)\bar{d}_3&d_9b_5=s_6+d_9+2t_{15}\\
(b_3\bar{d}_3)c_5=(b_3c_5)\bar{d}_3&d_9c_5=s_6+d_9+2t_{15}\\
(b_3\bar{d}_3)c_8=(b_3c_8)\bar{d}_3&d_9c_8=r_3+s_6+2d_9+3t_{15}\\
(b_3\bar{d}_3)d_3=(b_3d_3)\bar{d}_3&d_9d_3=b_3+b_9+x_{15}\\
(b_3\bar{d}_3)c_9=(b_3c_9)\bar{d}_3&d_9c_9=b_3+z_3+2\bar{x}_6+2\bb_9+3\bar{x}_{15}\\
(b_3\bar{d}_3)y_{15}=(b_3y_{15})\bar{d}_3&d_9y_{15}=\bb_3+z_3+2\bar{x}_6+3\bb_9+6\bar{x}_{15}\\
(b_3\bar{d}_3)x_9=(\bar{d}_3x_9)b_3&d_9x_9=r_3+y_3+2s_6+2d_9+3t_{15}\\
\hline
\end{array}
\]

\[
\begin{array}{|l|l|}
\hline
\mbox{associativity checked} &\mbox{new equations}\\
\mbox{or inner products}&\\
\hline &\\
(x_6d_3)z_3=(z_3x_6)d_3&z^2_3=d_3+\bb_6\\
(x_6d_3)\bar{z}_3=(\bar{z}_3x_6)d_3&z_3\bar{z}_3=1+b_8\\
(x_6\bar{d}_3)z_3=(y_3x_6)\bar{d}_3&y^2_3=1+c_8\\
(x_6d_3)y_3=(y_3x_6)d_3&y_3z_3=\bar{d}_3+b_6\\
(x_6d_3)\bar{d}_3=(d_3\bar{d}_3)x_6&z_3\bar{d}_3=\bar{z}_3+x_6\\
(x_6d_3)\bc_9=(d_3\bc_9)x_6&z_3\bc_9=r_3+d_9+t_{15}\\
(x_6d_3)\bb_6=(d_3\bb_6)x_6&z_3\bb_6=y_3+t_{15}\\
(x_6d_3)b_6=(d_3b_6)x_6&z_3b_6=\bar{z}_3+x_{15}\\
\hline
\end{array}
\]

\[
\begin{array}{|l|l|}
\hline
\mbox{associativity checked} &\mbox{new equations}\\
\mbox{or inner products}&\\
\hline &\\
(r_3x_{10},\ s_6)=(s_6b_8,\ r_3)=1 & r_3s_6=x_{10}+b_8\\
(d_9y_{15},\ \bar{x}_6)=2,\ (d_9b_6,\ \bar{x}_6)=1,\ (d_9c_9,\bar{x}_6)=2 & x_6d_9=\bb_6+2\bar{y}_{15}+2\bar{c}_9\\
(x_6y_3,\ d_3)=(x_6t_{15},\ d_3)=1&x_6\bar{d}_3=y_3+t_{15}\\
(r_3c_3)\bar{d}_3=r_3(c_3\bar{d}_3)&x_6d_3=z_3+\bar{x}_{15}\\
(\bar{x}_6b_6,\ b_9)=1,\ (\bb_9y_{15},\ x_6)=(\bb_9c_9,\ x_6)=2&x_6b_9=b_6+2y_{15}+2c_9\\
(x_6\bb_9,\ x_9)=(x_6\bar{x}_{15},\ x_9)=2,\ (x_6\bar{x}_6,\ x_9)=1&x_6x_9=2b_9+x_6+2x_{15}\\
(b_9\bb_9,\ c_8)=2,\ (b_9\bar{x}_{15},\ c_8)=1,\ (x_6c_8,\ b_9)=1&b_9c_8=2b_9+3x_{15}+x_6+b_3\\
(b_9\bb_9,\ b_8)=2,\ (b_9\bar{x}_{15},\ b_8)=3,&\\
(b_9z_3,\ b_8)=(x_6b_8,\ b_9)=1& b_9b_8=2b_9+3x_{15}+\bar{z}_3+x_6\\
(b_9\bb_9,\ c_5)=(x_6c_5,\ b_9)=1,\ (b_9\bar{x}_{15},c_5)=2&c_5b_9=b_9+2x_{15}+x_6\\
(x_{15}c_9,\ y_3)=(x_{15}b_6,\ y_3)=1,\ (x_{15}y_{15},\ y_3)=2& x_{15}y_3=\bb_6+\bc_9+2\bar{y}_{15}\\
(b_3b_5)\bar{z}_3=(b_3\bar{z}_3)b_5&x_{15}\bar{z}_3=b_6+c_9+2y_{15}\\
(b_3b_5)\bar{d}_3=(b_5\bar{d}_3)b_3&x_{15}\bar{d}_3=s_6+2t_{15}+d_9\\
(b_3b_5)\bc_9=(b_5\bc_9)b_3& x_{15}\bc_9=\bb_3+z_3+3\bb_9+2\bar{x}_6+6\bar{x}_{15}\\
(b_3b_5)y_{15}=(b_5y_{15})b_3&x_{15}\bar{y}_{15}=9\bar{x}_{15}+2z_3+2\bb_3+4\bar{x}_6+6\bb_9\\
(x_6t_{15},\ c_3)=(b_9t_{15}, \bc_3)=1,\ (x_{15}t_{15},\ \bc_3)=2&t_{15}c_3=\bar{x}_6+\bb_9+2\bar{x}_{15}\\
(d_9b_8,\ y_3)=(d_9x_{10},\ y_3)=(d_9x_9,\ y_3)=1&d_9y_3=x_9+b_8+x_{10}\\
(d_9y_{15},\ z_3)=(d_9c_3,\ z_3)=(d_9c_9,\ z_3)=1&d_9z_3=c_3+c_9+y_{15}\\
\hline
\end{array}
\]

\[
\begin{array}{|l|l|}
\hline
\mbox{associativity checked} &\mbox{new equations}\\
\mbox{or inner products}&\\
\hline &\\
(b_3\bb_3)d_9=(b_3d_9)\bb_3&d_9b_8=y_3+s_6+2d_9+3t_{15}\\
(b_9d_9,\ \bc_3)=(d_9z_3,\ c_3)=(x_{15}d_9,\ \bc_3)=1&d_9c_3=\bb_9+z_3+\bar{x}_{15}\\
(d_9b_8,\ y_3)=(t_{15}b_8,\ y_3)=1&y_3b_8=d_9+t_{15}\\
(s_6x_{10},\ y_3)=(t_{15}x_{10},\ y_3)=(d_9x_{10},\ y_3)=1&y_3x_{10}=s_6+t_{15}+d_9\\
(t_{15}y_3,\ b_5)=1&y_3b_5=t_{15}\\
(t_{15}y_3,\ c_5)=1&y_3c_5=t_{15}\\
(t_{15}y_3,\ c_8)=(s_6y_3,\ c_8)=(y_3c_8,\ y_3)=1&y_3c_8=s_6+t_{15}+y_3\\
(t_{15}y_3,\ x_9)=(d_9y_3,\ x_9)=(r_3y_3,\ x_9)=1&y_3x_9=t_{15}+d_9+r_3\\
(b_9y_3,\ \bc_3)=1&y_3c_3=\bb_9\\
(b_9y_3,\ \bar{d}_3)=1&y_3d_3=\bb_9\\
(b_3y_3,\ c_9)=(b_9c_9,\ y_3)=(x_{15}c_9,\ y_3)=1&y_3c_9=\bb_3+\bb_9+\bar{x}_{15}\\
(x_{15}y_3,\ b_6)=(y_3z_3,\ b_6)=1&y_3b_6=\bar{x}_{15}+z_3\\
(x_6y_3,\ \bar{y}_{15})=(b_9y_3,\ \bar{y}_{15})=1,\ (x_{15}y_3,\ \bar{y}_{15})=2&y_3y_{15}=\bar{x}_6+\bb_9+2\bar{x}_{15}\\
(x_6x_{10},\ z_3)=(b_9x_{10},\ z_3)=(x_{15}x_{10},\ z_3)=1&z_3x_{10}=\bar{x}_6+\bb_9+\bar{x}_{15}\\
(x_{15}b_5,\ z_3)=1&z_3b_5=\bar{x}_{15}\\
(x_{15}c_5,\ z_3)=1&z_3c_5=\bar{x}_{15}\\
(x_6c_8,\ z_3)=(x_{15}c_8,\ z_3)=(z_3\bar{z}_3,\ c_8)=1&z_3c_8=\bar{x}_6+\bar{x}_{15}+z_3\\
(b_3x_9,\ z_3)=(b_9x_9,\ z_3)=(x_{15}x_9,\ z_3)=1&z_3x_9=\bb_3+\bb_9+\bar{x}_{15}\\
(b_9c_3,\ z_3)=1&z_3c_3=b_9\\
(s_6\bar{d}_3,\ z_3)=(d_3y_3,\ \bar{z}_3)=1&z_3d_3=s_6+y_3\\
(\bb_3c_9,\ \bar{z}_3)=(\bb_9c_9,\ \bar{z}_3)=(x_{15}\bc_9,\ z_3)=1&z_3y_{15}=b_3+b_9+x_{15}\\
(x_6\bar{y}_{15},\ z_3)=(b_9\bar{y}_{15},\ z_3)=1,\ (x_{15}\bar{y}_{15},\ z_3)=2&z_3y_{15}=x_6+b_9+2x_{15}\\
(s_6y_{15},\ z_3)=(d_9y_{15},\ z_3)=1,\ (t_{15}y_{15},\ z_3)=2&z_3\bar{y}_{15}=s_6+2t_{15}+d_9\\
\hline
\end{array}
\]
\newpage
\section{Structure of NITA generated by $b_3$ and satisfying
$b^2_3=c_3+b_6$, $c_3\neq b_3,\ \bb_3$, $(b_3b_8,\ b_3b_8)=3$ and
$c_3$ real}

In this section we discuss NITA satisfying the following
hypothesis.
\begin{hyp}\label{hyp2}
Let $(A,\ B)$ be NITA generated by a non-real element $b_3\in B$
of degree 3 and without non-identity basis element of degree 1 or
2,\ such that :
 \begin{equation*}
 \overline{b}_3b_3=1+b_8
 \end{equation*}
 and
 \begin{equation*}
 b_3^2=c_3+b_6,\  c_3=\overline{c}_3
 \end{equation*}
 where $b_i,\  c_i ,\ d_i\in B$ are of degree i.\end{hyp}

\begin{lemma} \label{lem5.1}Let $(A,\  B)$ satisfy Hypothesis \ref{hyp2},\  then

1) $c_3^2 = 1+b_8$,\

2) $\overline{b}_3c_3 = b_3+x_6,\  x_6\neq\overline{x}_6\in B,\ $

3) $b_8b_3 = x_6+b_6\overline{b}_3,\ $

4) $ b_8+b_8^2 =c_3\overline{b}_6+b_6c_3+b_6\overline{b}_6,\  $

5) $\overline{b}_6c_3 = x_6\overline{b}_3,\ $

6) $c_3b_8 =b_6+x_6b_3$,\

7) $( b_8b_3,\  b_8b_3)=3,\ 4.$
\end{lemma}

\begin{proof}: By {\bf Hypothesis \ref{hyp2}},\  we have $(\overline{b}_3\overline{c}_3,\  b_3)=(c_3,\  b_3^2)=1$. So
$(\overline{b}_3\overline{b}_3,\ c_3^2)=(\overline{b}_3c_3,\
\overline{b}_3c_3) \geq2,\ $ hence
$$c_3^2=1+b_8.  \eqno (1)$$
Further we have
$$\overline{b}_3c_3=b_3+x_6. \eqno (2)$$
 If $x_6=\bar{x}_6$, then $1=(\overline{b}_3c_3,\ b_3c_3)=( c_3^2,\  b_3^2),\ $
 a contradiction to {\bf Hypothesis \ref{hyp2}}.
 By the associative law $(\overline{b}_3b_3)b_3=\overline{b}_3b_3^2$ and {\bf Hypothesis
 \ref{hyp2}}, one concludes that
 $b_3+b_8b_3=
 c_3\overline{b}_3+b_6\overline{b}_3.$ So
 $b_8b_3=x_6+b_6\overline{b}_3$ holds by (2).
By the associative law
$(\overline{b}_3b_3)^2=\overline{b}_3^2b_3^2$ and {\bf Hypothesis
\ref{hyp2}} and (1), we see that
\begin{eqnarray*}
1+2b_8+b_8^2&=&(c_3+\overline{b}_6)(c_3+b_6),\  \\
1+2b_8+b_8^2&=&c_3^2+c_3b_6+\overline{b}_6c_3+\overline{b}_6b_6,\ \\
b_8^2+b_8&=&c_3b_6+\overline{b}_6c_3+\overline{b}_6b_6.
\end{eqnarray*}

By the association law $(\overline{b}_6b_3)b_3 =
\overline{b}_3b_3^2$ and {\bf Hypothesis \ref{hyp2}} and (2), it
holds that $ b_3+b_8b_3 = c_3\overline{b}_3+b_6\overline{b}_3$,
hence
$$
b_8b_3 = x_6+b_6\overline{b}_3. \eqno (3)
$$
By {\bf Hypothesis \ref{hyp2}} ,\ $(b_6\overline{b}_3,\
b_3)=(b_6,\ b_3^2)=1$ ,\  so by (3) $$(b_8b_3,\  b_8b_3) \geq3,\
$$ we shall show that $(b_8b_3,\  b_8b_3)=3,\  4.$ By the
associative law and {\bf Hypothesis \ref{hyp2}} and equations (1)
 and (2), we have that $
(\overline{b}_3^2)c_3 = \overline{b}_3(\overline{b}_3c_3),$ $
c_3^2+\overline{b}_6c_3 = b_3\overline{b}_3+x_6\overline{b}_3,$
thus
\begin{equation*}
\overline{b}_6c_3 = x_6\overline{b}_3.\eqno (4)
\end{equation*}
By (3),\ $ (x_6\overline{b}_3,\  b_8)=(x_6,\  b_3b_8)=1$ holds. So
$( x_6\overline{b}_3,\  x_6\overline{b}_3)\geq 2.$ We shall show
that $(x_6\overline{b}_3,\  x_6\overline{b}_3)=2,\  3.$ If
$(x_6\overline{b}_3,\  x_6\overline{b}_3)=4$,\  then
 $x_6\overline{b}_3=b_8+g_3+w_3+v_4,\ $  hence $$b_3g_3=x_6+w_3,\  b_3w_3=x_6+k_3,\
 b_3v_4=x_6+z_6.$$
Therefore $1\leq (b_3g_3,\  b_3w_3)=( b_3\overline{b}_3,\
\bar{g}_3w_3).$ We come to $w_3=g_3$ by {\bf Hypothesis
\ref{hyp2}},\ hence $ 2=(x_6\overline{b}_3,\ w_3)=(x_6,\ b_3w_3),\
$ so that $b_3w_3=2x_6$, a contradiction by (4). Now we can
announce that $ (\overline{b}_6c_3,\
\overline{b}_6c_3)=(\overline{b}_6b_6,\ c_3^2)=2,\  3$. By (1), it
follows that $(\overline{b}_6b_6,\ b_3\overline{b}_3)=(
b_6\overline{b}_3,\ b_6\overline{b}_3)=$ 2 or 3. Therefore $(
b_8b_3,\ b_8b_3)=3,\  4$ by (3).

By the associative law $( \bar{b_3}c_3)b_3 =
c_3(\overline{b}_3b_3)$ and {\bf Hypothesis \ref{hyp2}}, we get
 $b_3^2+x_6b_3 = c_3+c_3b_8$ and $ c_3b_8=b_6+x_6b_3.$
\end{proof}

\begin{lemma}\label{lem5.2}Let $(A,\  B)$ satisfy {\bf Hypothesis \ref{hyp2}} and $( b_8b_3,\
b_8b_3)=3$,\  then we get:

1) $b_6\overline{b}_3=b_3+x_{15}$,\

2) $ b_8b_3=x_6+b_3+x_{15},\  x_6\neq\overline{x}_6$,\

3) $\overline{x}_6b_3=b_8+x_{10}=b_6c_3,\  x_6\neq b_6,\
x_{10}=\overline{x}_{10}$,\

4) $x_6c_3=\overline{b}_3+\overline{x}_{15}$,\

5) $ b_3x_6=c_3+y_{15}$,\

6) $ b_8c_3=b_6+c_3+y_{15},\  b_6=\overline{b}_6,\
y_{15}=\overline{y}_{15}.$
\end{lemma}

\begin{proof} By our assumption in this lemma and part 3) in {\bf Lemma \ref{lem5.1}},\  and
{\bf Hypothesis \ref{hyp2}}: $$b_6\overline{b}_3=b_3+x_{15},\
\eqno(1)
$$then
$$b_8b_3=x_6+b_3+x_{15}$$ by (1),\ $$( b_6\overline{b}_6,\
b_3\overline{b}_3)=( b_6\overline{b}_3,\  b_6\overline{b}_3)=2,\
$$ then by 1) in {\bf Lemma \ref{lem5.1}},\  $$2=(
b_6\overline{b}_6,\ c_3^2)=(\overline{b}_6c_3,\  \overline{b}_6c_3
). \eqno(2)$$ By 6) in {\bf Lemma \ref{lem5.1}},\
 $$(\overline{b}_6c_3,\  b_8)=(c_3b_8,\  b_6)=1,\  $$ so by (2)
 $$\overline{b}_6c_3=b_8+x_{10}. \eqno(3)$$ so by 5) in {\bf Lemma \ref{lem5.1}}
 $$x_6\overline{b}_3=b_8+x_{10},\   \eqno(4)$$ so by (1) we get that $x_6\neq b_6$.
 By the associative law $(c_3^2)\overline{b}_3 = (c_3\overline{b}_3)c_3$ and 1),\  2) in {\bf Lemma \ref{lem5.1}}
  and (2)
 , it follows that
  \begin{eqnarray*}
\overline{b}_3+\overline{b}_3b_8 &=& b_3c_3+x_6c_3,\ \\
\overline{b}_3+\overline{x}_6+\overline{b}_3+\overline{x}_{15} &=&
\overline{b}_3+\overline{x}_6+x_6c_3,
\end{eqnarray*}
so \[x_6c_3 = \overline{b}_3+\overline{x}_{15}, \eqno (5)\]  By
(4) and 2) in {\bf Lemma \ref{lem5.1}}
$$2=(\overline{x}_6b_3,\  x_6b_3)=( x_6b_3,\  x_6b_3),\  (b_3x_6,\
c_3)=(x_6,\  \overline{b}_3c_3)=1.$$ Therefore $$b_3x_6=c_3+y_{15}
. \eqno(6)$$ By the associative law $(b_3\overline{b}_3)c_3 =
(\overline{b}_3c_3)b_3$ and 2) in {\bf Lemma \ref{lem5.1}} and
{\bf Hypothesis \ref{hyp2}}, one has
\begin{eqnarray*}
c_3+b_8c_3 &=& b_3^2+x_6b_3,\ \\
c_3+b_8c_3 &=& c_3+b_6+c_3+y_{15},\ \\
b_8c_3 &=& b_6+c_3+y_{15}.
\end{eqnarray*}
$ b_8,\  c_3 $ are reals ,\ then $$b_6=\overline{b}_6,\
y_{15}=\overline{y}_{15}. $$ By (3), we have
$x_{10}=\overline{x}_{10}$.
\end{proof}

\begin{lemma} \label{lem5.3}Let $(A,\  B)$ satisfy {\bf Hypothesis \ref{hyp2}} and $( b_8b_3,\
b_8b_3)=3$,\  then {\bf Lemma \ref{lem5.2}} hold and we get :

1) $ x_{10}b_3 =x_{15}+x_6+x_9,\  x_{10}=\overline{x}_{10},\ $

2) $x_6b_8 = 2x_{15}+b_3+x_6+x_9,\ $

3) $ \overline{x}_{15}c_3 = 2x_{15}+x_6+x_9,\ $

4) $ y_{15}\overline{b}_3 = 2x_{15}+x_6+x_9,\ $

5) $ x_6b_6 = 2\overline{x}_6+\overline{x}_{15}+\overline{x}_9,\ $

6) $ c_3x_{10} = b_6+y_{15}+y_9,\ $

7) $x_6^2 = 2b_6+y_{15}+y_9,\ $

8) $ x_{15}b_3 = 2y_{15}+b_6+y_9,\ $

9) $ b_8b_6 = 2y_{15}+b_6+c_3+y_9,\  y_9=\overline{y}_9,\ $

10) $b_8x_{15} = 4x_{15}+2x_6+2x_9+b_3+y_9\overline{b}_3.$
\end{lemma}

\begin{proof} By the associative law $(b_3^2)b_6 = (b_3b_6)b_3$ and {\bf Hypothesis \ref{hyp2}} and 1),\  3) in {\bf Lemma
\ref{lem5.2}}, we get
\begin{eqnarray*}
c_3b_6+b_6^2&=& \overline{b}_3b_3+\overline{x}_{15}b_3,\ \\
b_8+x_{10}+b_6^2&=& 1+b_8+\overline{x}_{15}b_3,\
\end{eqnarray*}
 Hence\[ x_{10}+b_6^2 =
1+\overline{x}_{15}b_3. \eqno(1)\] By 3) in {\bf Lemma
\ref{lem5.2}} ,\ $$(b_3x_{10},\  x_6)=(x_{10},\
\overline{b}_3x_6)=1,\ $$ hence by (1),\
$$(x_{10}b_3,\  x_{15})=(b_3\overline{x}_{15},\  x_{10})=1,\ $$ So
$$x_{10}b_3=x_{15}+x_6+x_9,\  $$ where $x_9$ is linear combination of
B. we shall show that $x_9\in B$. If $x_{10}b_3$ have element of
degree 3 ,\ then $$1=(x_{10}b_3,\  y_3)=(x_{10},\
\overline{b}_3y_3),\ $$
 hence $\overline{b}_3y_3=x_{10}$, a contradiction. So $$x_{10}b_3=x_{15}+x_6+x_9,\  \eqno(2)$$ where $x_9\in B$.
 By the associative law and
 {\bf Hypothesis \ref{hyp2}} and  2),\  3) in {\bf Lemma \ref{lem5.2}},\  we get
\begin{eqnarray*}
 (b_3\overline{x}_6)\overline{b}_3 &=& \overline{x}_6(b_3\overline{b}_3),\ \\
 b_8\overline{b}_3+x_{10}\overline{b}_3 &=& \overline{x}_6+\overline{x}_6b_8,\ \\
\overline{x}_6+\overline{x}_6b_8&=&\overline{x}_6+\overline{b}_3+\overline{x}_{15}+\overline{x}_{15}+\overline{x}_6+\overline{x}_9 ,\ \\
 \end{eqnarray*}
 Hence \[\overline{x}_6b_8 = 2\overline{x}_{15}+\overline{b}_3+\overline{x}_6+\overline{x}_9.\eqno (3)\]

 By the associative law $(\bar{b_3}c_3)b_8 = \bar{b_3}(b_8c_3)$ and 1),\  2),\  6) in {\bf Lemma \ref{lem5.2}} and 2) in {\bf Lemma \ref{lem5.1}},\
 we come to
 \begin{eqnarray*}
 b_3b_8+x_6b_8 &=& b_6\bar{b_3}+c_3\bar{b_3}+y_{15}\bar{b_3},\ \\
 x_6+b_3+x_{15}+2x_{15}+b_3+x_6+x_9 &=& b_3+x_{15}+b_3+x_6+y_{15}\bar{b_3},\
\end{eqnarray*}
 So \[y_{15}\bar{b_3} = 2x_{15}+x_6+x_9. \eqno(4)\]

 By the associative $(\overline{b}_3c_3)b_8 = (\overline{b}_3b_8)c_3$ and
 1),\  2),\  6) in {\bf Lemma \ref{lem5.2}} and in {\bf Lemma \ref{lem5.1}} and
 (3), the following equations hold:
 \begin{eqnarray*}
 b_3b_8+x_6b_8 &=& \overline{x}_6c_3+\overline{b}_3c_3+\overline{x}_{15}c_3,\ \\
 x_6+b_3+x_{15}+2x_{15}+b_3+x_6+x_9 &=& b_3+x_{15}+b_3+x_6+\overline{x}_{15}c_3,
  \end{eqnarray*}
  So\[\overline{x}_{15}c_3 = 2x_{15}+x_6+x_9. \eqno(5)\]

 By the associative law $(\overline{b}_3c_3)b_6=c_3(\overline{b}_3b_6)$ and 2) in {\bf Lemma \ref{lem5.1}} and 1) in {\bf Lemma \ref{lem5.2}} and
 (5), one should have
 \begin{eqnarray*}
 b_3b_6+x_6b_6 &=& b_3c_3+x_{15}c_3,\ \\
 \overline{b}_3+\overline{x}_{15}+x_6b_6 &=& \overline{b}_3+\overline{x}_6+2\overline{x}_{15}+\overline{x}_6+\overline{x}_9,\ \\
 x_6b_6 &=& 2\overline{x}_6+\overline{x}_{15}+\overline{x}_9.
 \end{eqnarray*}
 By the associative law and 5),\  4) in {\bf Lemma \ref{lem5.2}} and 1) in {\bf Lemma
 \ref{lem5.1}}, we get
$ (b_3x_6)c_3 = b_3(x_6c_3),$
 $c_3^2+y_{15}c_3 = b_3\overline{b}_3+b_3\overline{x}_{15}$ and $
 y_{15}c_3 = b_3\overline{x}_{15}.$
  Now by (2) it follows that $(\overline{x}_{15}b_5,\  x_{10})=(b_3x_{10},\  x_{15})=1$.
  So $$1=(y_{15}c_3,\  x_{10})=(y_{15},\  c_3x_{10}).
   \eqno(7)$$
By 3) in {\bf Lemma \ref{lem5.2}},\  $$(c_3x_{10},\ b_6)=(c_3b_6,\
x_{10})=1. \eqno (8)$$ By (2) and 1) in {\bf Lemma \ref{lem5.1}},
we have $(c_3x_{10},\ c_3x_{10})=(c_3^2,\
x_{10}^2)=(b_3\overline{b}_3,\ x_{10}^2)=(b_3x_{10},\
b_3x_{10})=3$. Therefore by (7) and (8), the following equation
holds
$$c_3x_{10}=b_6+y_{15}+y_9, \eqno(9)$$ where $y_9\in B$. By the
associative law and 3), 5), 6) in {\bf Lemma \ref{lem5.2}} and 2)
in {\bf Lemma \ref{lem5.1}} and (9), we have that
$(\overline{x}_6b_3)c_3 = \overline{x}_6(b_3c_3)$  and
\begin{eqnarray*}
b_8c_3+x_{10}c_3&=&\overline{x}_6\overline{b}_3+\overline{x}_6^2,\ \\
b_6+c_3+y_{15}+b_6+y_{15}+y_9 &=&c_3+y_{15}+\overline{x}_6^2,\
\end{eqnarray*}
so \[ x_6^2 = 2b_6+y_{15}+y_9.  \eqno(10)\]

By the associative law and 3),\  4) ,\ 6) in {\bf Lemma
\ref{lem5.2}} and Hypothesis \ref{hyp2} and (9), we have
$(\overline{x}_6b_3)c_3 = (\overline{x}_6c_3)b_3$, hence
\begin{eqnarray*}
b_8c_3+x_{10}c_3&=& b_3^2+x_{15}b_3,\ \\
b_6+c_3+y_{15}+b_6+y_{15}+y_9 &=&c_3+b_6+x_{15}b_3,\
\end{eqnarray*} so \[x_{15}b_3=2y_{15}+b_6+y_9.\eqno (11)\]

By the associative law and 2),\  5) ,\ 6) in {\bf Lemma
\ref{lem5.2}} and (11) and the {\bf Hypothesis \ref{hyp2}},
$b_3(b_8b_3) = b_8(b_3^2)$ holds, which concludes
\begin{eqnarray*}
x_6b_3+b_3^2+x_{15}b_3 &=& b_8c_3+b_8b_6,\ \\
c_3+y_{15}+c_3+b_6+2y_{15}+b_6+y_9 &=&  b_6+c_3+y_{15}+b_8b_6,\
\end{eqnarray*}
 so \[
b_8b_6 = 2y_{12}+b_6+c_3+y_9.\eqno (12)\]

By the associative law and 1),\  2) in {\bf Lemma \ref{lem5.2}}
and (4) and 2) in {\bf Lemma \ref{lem5.1}} :
\begin{center}
$(b_6\overline{b}_3)b_8 = \overline{b}_3(b_6b_8),$\ \
$b_3b_8+x_{15}b_8 =
2y_{12}\overline{b}_3+b_6\overline{b}_3+c_3\overline{b}_3+y_9\overline{b}_3.$
 \end{center}
\end{proof}

\begin{lemma}\label{lem5.4} Let $(A,\  B)$ satisfy {\bf Hypothesis \ref{hyp2}} and $(b_8b_3,\
b_3b_8)=3$,\  then {\bf Lemma \ref{lem5.3}} holds and
$x_9c_3=\bar{x}_{15}+d_3+d_9   $,\  $d_3\neq c_3,\  b_3,\
\bar{b_3}$.
\end{lemma}

\begin{proof} By 3) in {\bf Lemma \ref{lem5.3}} ,\ $$(c_3x_9,\ \overline{x}_{15})=(x_9,\
c_3\overline{x}_{15})=1,\ $$
 so $$c_3x_9=\overline{x}_{15}+\alpha,\ \eqno  (1)  $$where $\alpha$ is linear combination at B.
 We shall show that $\alpha$ does not have degree 4 and also two elements of degree 3. If $x_4\in\alpha$,\
 then $1\leq(c_3x_9,\  x_4)=(x_9,\  c_3x_4)$,\  hence $1=(x_9,\  c_3x_4)$. So $$c_3x_4=x_9+m_3.\eqno (2)$$ So
 $$(c_3m_3,\  x_4)=(m_3,\  c_3x_4)=1,\ $$  hence $$c_3m_3=y_4+k_5.\eqno (3) $$Hence$$2=( c_3m_3,\ c_3m_3)
 =( c_3^2,\  \overline{m}_3m_3),\ $$ therefore $\overline{m}_3m_3=1+b_8$. So by the associative law and (3) and
 6) in {\bf Lemma \ref{lem5.2}},\  we get
 $(\overline{m}_3m_3)c_3 =\overline{m}_3(m_3c_3),\ $
$ c_3+c_3b_8 = y_4\overline{m}_3+k_5\overline{m}_3$. Hence \[
2c_3+b_6+y_{15} = y_4\overline{m}_3+k_5\overline{m}_3\eqno (4)\]

 By (3) $$(y_4\overline{m}_3,\  c_3)=(y_4,\  m_3c_3)=1,\  (k_5\overline{m}_3,\  c_3)=(k_5,\  m_3c_3)=1,\ $$ hence by (4)
 we get a contradiction. If $x_3,\  y_3\in\alpha$ ,\  then by  (1) ,\ $$1=(c_3x_9,\  x_3)=(x_9,\  c_3x_3),\ $$
 $$1=(c_3x_9,\  y_3)=(x_9,\  c_3y_3).$$ So $$c_3x_3=x_9,\  c_3y_3=x_9,\ $$ hence $$1=( c_3x_3,\  c_3y_3)
 =( c_3^2,\  \overline{x}_3y_3),\ $$ therefore by 1) in {\bf Lemma \ref{lem5.1}},\  $\bar{x_3}y_3=b_8+1,\ $ so $x_3=y_3$,\
 then by (1) $2=(c_3x_9,\  x_3)=(x_9,\  c_3x_3)$ a contradiction. Therefore  by
 (1) one of following holds:
\begin{center}
a) $\alpha=d_3+d_9$ or b) $\alpha=d_5+d_7$ or c) $\alpha=d_6+k_6$
or d) $\alpha=2d_6$ or e) $ \alpha=d_{12}$. \end{center} We shall
show
 that a) is the only possibility. Assume b) holds,\  then by (1) we get $c_3x_9=\overline{x}_{15}+d_5+d_7.$ So
 $(c_3d_5,\  x_9)=(d_5,\  c_3x_9)=1,\ $ therefore $c_3d_5=x_9+\beta,\ $ where $\beta=2x_3$ or $\beta=x_3+y_3$ or $\beta=m_6$.
 We shall show that  this possibilities are impossible if $\beta=2x_3$,\  then $c_3d_5=x_9+2x_3$,\  therefore
 $(c_3x_3,\  d_5)=(x_3,\  c_3d_5)=2,\ $ so $c_3x_3=2d_5$ a contradiction.
 If $\beta=x_3+y_3,\ $ then
 $c_3d_5=x_9+x_3+y_3,\ $
 hence $(c_3x_3,\  d_5)=(x_3,\  c_3d_5)=1,\  (c_3y_3,\  d_5)=(y_3,\  c_3d_5)=1,\ $  so $1\leq( c_3x_3,\  c_3y_3)
 =( c_3^2,\  \overline{x}_3y_3),\ $ then by 1) in {\bf Lemma \ref{lem5.1}} $x_3=y_3$ a contradiction. If $\beta=m_6,\ $ then
$c_3d_5=x_9+m_6,\ $ so by (1) and  (6), we have $
(c_3^2)d_5=c_3(c_3d_5),\
 d_5+b_8d_5 = x_9c_3+m_6c_3$ and $ b_8d_5 = \overline{x}_{15}+d_7+m_6c_3.$
 So $(b_8x_{15},\  \overline{d}_{5})=(b_8d_5,\  \overline{x}_{15})=1,\ $ hence by 10) in {\bf Lemma \ref{lem5.3}}:
 $(y_9\overline{b}_3,\  \overline{d}_5)=1$ and by 8) in {\bf Lemma \ref{lem5.3}}: $(y_9\overline{b}_3,\  x_{15})=(y_9,\
 b_3x_{15})=1,\ $
 so $y_9\overline{b}_3=x_{15}+\overline{d}_5+k_7,\ $ where $k_7$ is linear combination.
  When $k_7=e_3+m_4$,\  then $y_9\overline{b}_3=x_{15}+\overline{d}_5+e_3+t_3,\ $ hence $(b_3m_4,\  y_9)=(m_4,\
  \overline{b}_3y_9)=1,\ $
  so $b_3m_4=y_9+t_3$,\  hence $(t_3\overline{b}_3,\  m_4)=(t_3,\  b_3m_4)=1,\ $
  therefore $$t_3\bar{b_3}=m_4+t_5,\  \eqno(5)$$
  so $ 2=( t_3\overline{b}_3,\  t_3\overline{b}_3)=( t_3\overline{t}_3,\  b_3\overline{b}_3),\ $ so by {\bf Hypothesis \ref{hyp2}}
  $\bar{t_3}t_3=1+b_8$,\  hence by the associative law and (5) we come to:
  \begin{equation*}
  (\overline{t}_3t_3)b_3 = (\overline{t}_3b_3)t_3,\ \\
  b_3+b_8b_3 = \overline{m}_4t_3+\overline{t_5}t_3.
  \end{equation*}
  and by 2) in {\bf Lemma \ref{lem5.2}} it follows that $2b_3+x_6+x_{15}=\overline{m}_4t_3+\overline{t_5}t_3,\ $
  so $\bar{m_4}t_3=2b_3+x_6$,\
  hence $5=(\overline{m}_4t_3,\  \overline{m}_4t_3)=(\overline{m}_4m_4,\  \overline{t}_5t_3),\ $ then $\overline{m}_4m_4=1+4b_8,\ $ a contradiction.
  If $$y_9\overline{b}_3=x_{15}+\overline{d}_5+k_7,\ \eqno (6)$$ where $k_7\in B$. Then $(b_3\overline{d}_5,\  y_9)=(\overline{d}_5,\
  \overline{b}_3y_9)=1,\ $
   hence $b_3\overline{d}_5=y_9+\gamma,\ $ where $\gamma=v_6$ or $\gamma=v_3+w_3$.
   When $b_3\overline{d}_5=y_9+v_3+w_3$,\  then $(b_3\overline{v}_3,\  d_5)=(b_3\overline{d}_5,\  v_3)=1,\ $ also $(b_3\overline{m}_3,\  d_5)=1$,\
   so $1\leq( b_3\overline{v}_3,\  b_3\overline{m}_3)=( b_3\overline{b}_3,\
   v_3\overline{m}_3),\ $
   by {\bf Hypothesis \ref{hyp2}},\  $v_3=w_3$,\  then $$2=(b_3\overline{d}_5,\  v_3)=(\overline{d}_5,\  \overline{b}_3v_3),\ $$ so $$\overline{b}_3v_3=2\overline{d}_5$$ a contradiction.
   When $b_3\bar{d_5}=y_9+v_6$,\  Then by the associative law and (6) and our assumption:
  \begin{eqnarray*}
  (b_3^2)\overline{d}_5 &=& b_3(b_3\overline{d}_5),\ \\
  \overline{d}_5c_3+\overline{d}_5b_6 &=& y_9b_3+v_6b_3,\ \\
  \overline{x}_9+\overline{m}_6+\overline{d}_5b_6 &=& \overline{x}_{15}+d_5+\overline{k}_7+v_6b_3.
  \end{eqnarray*}
  Hence $$v_6b_3=\overline{x}_9+\overline{m}_6+t_3,\ \eqno (7)$$ so
  $(t_3\overline{b}_3,\  v_6)=(t_3,\ b_3v_6)=1,\ $ hence
  $t_3\overline{b}_3=v_6+p_3.$ Therefore $2=( \overline{b}_3t_3,\  \overline{b}_3t_3)=( t_3\overline{t}_3,\  b_3\overline{b}_3),\ $
  Now by the {\bf Hypothesis \ref{hyp2}} we get $$t_3\overline{t}_3=1+b_8.\eqno (8)$$
  Again by the associative law:
   \begin{equation*}
   b_3(t_3\overline{t}_3) = t_3(b_3\overline{t}_3),\ \\
   b_3+b_3b_8 = t_3\overline{v}_6+t_3\overline{p}_3.
  \end{equation*}
and by 2) in {\bf Lemma \ref{lem5.2}} we have
$2b_3+x_6+x_{15}=t_3\overline{v}_6+t_3\overline{p}_3,\ $ but
$(t_3\overline{v}_6,\  b_3)=(t_3\overline{b}_3,\ v_6)=1,\
(t_3\overline{p}_3,\  b_3)=(t_3\overline{b}_3,\  p_3)=1.$ This
means $\overline{v}_6t_3=b_3+x_{15},\  t_3\overline{p}_3=b_3+x_6.$
Then by (7) and (8) it follows that $2=(\overline{v}_6t_3,\
\overline{v}_6t_3)=(\overline{v}_6v_6,\
\overline{t}_3t_3)=(\overline{v}_6v_6,\  \overline{b}_3b_3) =(
v_6b_3,\  v_6b_5)=3,\ $ a contradiction.

Assume c) holds ,\  then by (1) one has that
$c_3x_9=\overline{x}_{15}+d_6+k_6,\ $ so $(c_3d_6,\  x_9)=(d_6,\
c_3x_9)=1,\ $ hence $c_3d_6=x_9+\gamma,\ $
 where $\gamma=m_9$ or $\gamma=m_3+v_3+k_3$ or $\gamma=v_3+y_6$ or $\gamma=v_4+v_5$ or $\gamma=2v_3+m_3$ or $\gamma=3v_3$.
 If $\gamma=m_9$,\  then $c_3d_6=x_9+m_9$,\  so  by the associative law and 1) in {\bf Lemma \ref{lem5.2}} :
 \begin{eqnarray*}
 (c_3^2)d_6 &=& (c_3d_6)c_3,\ \\
 d_6+d_6b_8 &=&x_9c_3+m_9c_3,\ \\
 d_6+d_6b_8 &=& \overline{x}_{15}+d_6+k_6+m_9c_3,\ \\
\end{eqnarray*}
$$ d_6b_8 = \overline{x}_{15}+k_6+m_9c_3.  \eqno (9)$$

 We shall that $d_6\neq x_6,\  \overline{x}_6,\  k_6\neq x_6,\  \overline{x}_6.$ By 4) in {\bf Lemma \ref{lem5.2}} and our assumption if $d_6=x_6$,\  then
 $$
 1=(c_3x_9,\  x_6)=(x_9,\  c_3x_6)=0,\  1=(c_3x_9,\  \overline{x}_6)=(c_3x_6,\  \overline{x}_9)=0 ,\ $$ a contradiction,\  so $d_6\neq x_6,\  \overline{x}_6$.
 in the same way $k_6\neq x_6,\  \overline{x}_6.$ By (9),\  $(b_8x_{15},\  \overline{d}_6)=(b_8d_6,\  \overline{x}_{15})=1,\ $ in the same way
  $(b_8x_{15},\  \overline{k}_6)=(b_8k_6,\  \overline{x}_{15}) =1,\ $ then by 8),\  10) in {\bf Lemma \ref{lem5.3}},\  $(y_9\overline{b}_3,\  x_{15})=(y_9,\
  b_3x_{15})=1,\ $
  $(y_9\overline{b}_3,\  \overline{d}_6)=1,\  (y_9\overline{b}_3,\  \overline{k}_6)=1. $ So
  $$y_9\overline{b}_3=x_{15}+\overline{d}_6+\overline{k}_6.\eqno (10)$$ So by the
  associative law and 2) in {\bf Lemma \ref{lem5.1}} :
  \begin{eqnarray*}
  (\overline{b}_3c_3)y_9 &=& c_3(\overline{b}_3y_9),\ \\
  b_3y_9+y_9x_6 &=&x_{15}c_3+\overline{d}_6c_3+\overline{k}_6c_3,\ \\
  \overline{x}_{15}+d_6+k_6+y_9x_6&=&2\overline{x}_{15}+\overline{x}_6+\overline{x}_9+\overline{x}_9+\overline{m}_9+\overline{k}_6c_3,\ \\
  d_6+k_6+y_9x_6&=&\overline{x}_{15}+\overline{x}_6+2\overline{x}_9+\overline{m}_9+\overline{k}_6c_3.
 \end{eqnarray*}
  And by 3) in {\bf Lemma \ref{lem5.3}} and (10) and that $c_3d_6=x_9+m_9$,\  since $k_6,\  d_6\neq x_6 $,\  then $1\leq (\overline{k}_6c_3,\  d_6)=(\overline{k}_6,\
  c_3d_6),\ $
   a contradiction to the assumption that $c_3d_6=x_9+m_9$. If $\gamma=m_3+v_3+k_3$,\  then $c_3d_3=x_9+m_3+v_3+k_3,\ $ so
   $(c_3m_3,\  d_3)=(m_3,\  c_3d_3)=1,\  (c_3v_3,\  d_3)=(v_3,\  c_3d_3)=1,\ $ hence $1\leq( c_3m_3,\  c_3v_3)=
   ( m_3\overline{v}_3,\  c_3^2),\ $ then by 1) in {\bf Lemma \ref{lem5.1}} $m_3=v_3$ a contradiction.
    If $\gamma=v_3+y_6$,\  then $c_3d_6=x_9+v_3+y_6,\ $ so $(c_3v_3,\  d_6)=(v_3,\  c_3d_6)=1,\ $ hence $c_3v_3=d_6+t_3$,\
    so by  the associative law and 1) in {\bf Lemma \ref{lem5.1}},
    we get equations: $c_3^2v_3 = c_3(c_3v_3)$ and
  \begin{equation*}
   v_3+v_3b_8 = d_6c_3+t_3c_3,\ \\
  v_3+v_3b_8 = x_9+v_3+y_9+t_3c_3.
 \end{equation*}
 Hence $v_3b_8=x_9+y_6+t_3c_3.$  Since $(t_3c_3,\  v_3)=(t_3,\  c_3v_3)=1$,\  then $( v_3b_8,\
 v_3b_8)\geq4$,\
 in other hand $2=( c_3v_3,\  c_3v_3)=( c_3^2,\  \overline{v}_3v_3),\ $ so by 1) in {\bf Lemma \ref{lem5.1}} :
 $\bar{v_3}v_3=1+b_8=\bar{b_3}b_3$,\  therefore $( v_3b_8,\  v_3b_8)=( b_8^2,\  v_3\overline{v}_3)
 =( b_8^2,\  \overline{b}_3b_3)=( b_8b_3,\  b_8b_3),\ $ by 2) in {\bf Lemma \ref{lem5.2}},\  $( b_8b_3,\
 b_8b_3)=3,\ $
  a contradiction.
   If $\gamma=v_4+v_5$,\  then $c_3d_6=x_4+v_4+v_5,\ $ so by the associative law and 1) in {\bf Lemma \ref{lem5.1}} :
   \begin{equation*}
   c_3^2d_6 = c_3(c_3d_6),\ \\
   d_6+d_6b_8 = x_9c_3+v_4c_3+v_5c_3,\
  \end{equation*}
  hence by our assumption we have
  $d_6+d_6b_8 = \overline{x}_{15} +d_6+k_6+v_4c_3+v_5c_3$,\  thus
  $$d_6b_8=\bar{x}_{15}+k_6+v_4c_3+v_5c_3.\eqno (11)$$

 We shall show that $d_6\neq x_6,\  \overline{x}_6,\  k_6\neq x_6,\  \overline{x}_6.$ By 4) in {\bf Lemma \ref{lem5.2}} ,\ if $d_6=x_6$,\  or $d_6=\overline{x}_6$,\
  then $1=(c_3x_9,\  x_6)=(x_9,\  c_3x_6)=0,\  1=(c_3x_9,\  \overline{x}_6)=(c_3x_6,\  \overline{x}_9)=0,\ $ a contradiction .In the same way
  $ k_6\neq x_6,\  \overline{x}_6$. By (11) and 10) in {\bf Lemma \ref{lem5.3}},\  $$(b_8x_{15},\  \overline{d}_6)=(b_8d_6,\
  \overline{x}_{15})=1$$
 in the same way $(b_8x_{15},\  \overline{k}_6)=1$ and by 10),\  8) in {\bf Lemma \ref{lem5.3}} $(y_9\overline{b}_3,\  \overline{k}_6)=1,\
 (y_9\overline{b}_3,\  \overline{d}_6)=1,\  (y_9\overline{b}_3,\  x_{15})=1,\ $ so $y_9\overline{b}_3=x_{15}+\overline{k}_6+\overline{d}_6.$ Hence
 by the associative law an 2) in {\bf Lemma \ref{lem5.1}} :
  \begin{eqnarray*}
 (\overline{b}_3c_3)y_9 &=& (\overline{b}_3y_9)c_3,\ \\
 b_3y_9+x_6y_9 &=& x_{15}c_3+\overline{k}_6c_3+\overline{d}_6c_3,\
\end{eqnarray*}
 so by 3) in {\bf Lemma \ref{lem5.3}} and our assumption:
 \begin{eqnarray*}
 \overline{x}_{15}+k_6+d_6+x_6y_9 &=& 2\overline{x}_{15}+\overline{x}_6+\overline{x}_9+\overline{k}_6c_3+\overline{d}_6c_3,\ \\
 k_6+d_6+x_6y_9 &=& \overline{x}_{15}+\overline{x}_6+\overline{x}_9+\overline{v}_4+\overline{v}_5+\overline{x}_9+\overline{k}_6c_3.
  \end{eqnarray*}
  hence $1=(\overline{k}_6c_3,\  d_6)=(\overline{k}_6,\  c_3d_6),\ $ a contradiction to our assumption .
  If $\gamma=2v_3+m_3$,\  or $\gamma=3v_3$,\  then $c_3d_6=x_9+2v_3+m_3$ or $c_3d_6=x_9+3v_3$. In two
  cases,\  it always follows that  $(c_3v_3,\  d_6)=(v_3,\  c_3d_6)\geq2$,\  a contradiction.
  assumed then by (1) $c_3x_9=\overline{x}_{15}+2d_6$. Hence $(c_3d_6,\  x_9)=(d_6,\
  c_3x_9)=2$
  therefore $c_3d_6=2x_9$. so by the associative law and 1) in {\bf Lemma \ref{lem5.1}}, we have that
  $(c_3^2)d_6 = c_3(c_3d_6)$
  \begin{eqnarray*}
  d_6+d_6b_8 &=& 2x_9c_3\\
 &=& 2\overline{x}_{15}+4d_6,\ \\
 d_6b_8 &=& 2\overline{x}_{15}+3d_6.
 \end{eqnarray*}
  so $(b_8x_{15},\  \overline{d}_6)=(b_8d_6,\  \overline{x}_{15})=2 (12).$ by 4) in {\bf Lemma \ref{lem5.2}},\  if $x_6=d_6$,\  then
  $2=(c_3x_9,\  x_6)=(x_9,\  c_3x_6)=0,\ $ then $d_6\neq x_6,\  $ also $d_6\neq \overline{x}_6$.
  then by (12) and 8),\  10) in {\bf Lemma \ref{lem5.3}}:
  $(y_9\overline{b}_3,\  x_{15})=(y_9,\  b_3x_{15})=1,\  (y_9\overline{b}_3,\
  \overline{d}_6)=2.$
  So $y_9\overline{b}_3=x_{15}+2\overline{d}_6.$ Hence by the associative law and 2) in lemma \ref{lem5.1}:
  \begin{eqnarray*}
  y_9(\overline{b}_3c_3)& =& (y_9\overline{b}_3)c_3,\ \\
  y_9b_3+y_9x_6 &=& x_{15}c_3+2\overline{d}_6c_3,\
 \end{eqnarray*}
  then by 3) in {\bf Lemma \ref{lem5.3}} $\overline{x}_{15}+2d_6+y_9x_6=2\overline{x}_{15}+\overline{x}_6+\overline{x}_9+2\overline{x}_9,\ $ a contradiction.
 Assume e) then by 1) ,\ $c_3x_9=\overline{x}_{15}+d_{12}.$ Then by 1) in {\bf Lemma \ref{lem5.1}} :
  $2=( x_9c_3,\  x_9c_3)=( x_9\overline{x}_9,\  c_3^2)=( x_9\overline{x}_9,\
  b_3\overline{b}_3)=(
  x_9\overline{b}_3,\  x_9\overline{b}_3).$ So by 1) in {\bf Lemma \ref{lem5.3}} $(x_9\overline{b}_3,\  x_{10})=(x_9,\  b_3x_{10})
 =1,\ $ hence $x_9\overline{b}_3=x_{10}+x_{17} (13) .$ By the associative law and 1) in {\bf Lemma \ref{lem5.1}}:
 $(c_3^2)x_9=(b_3\overline{b}_3)x_9=b_3(\overline{b}_3x_9)=(c_3x_9)c_3,\
 x_{10}b_3+b_3x_{17}=\overline{x}_{15}c_3+d_{12}c_3.$
 By 1),\  3) and {\bf Lemma \ref{lem5.3}} ,\ $x_{15}+x_6+x_9+b_3x_{17}=2x_{15}+x_6+x_9+d_{12}c_3,\ $ $ b_3x_{17}=x_{15}+d_{12}c_3.$ Then
 $(x_{15}\overline{b}_3,\  x_{17})=(x_{15},\  b_3x_{17})=1$ and by 2) in {\bf Lemma \ref{lem5.2}} and 1),\  8) in {\bf Lemma \ref{lem5.3}}:
 \begin{equation*}
           \begin{array}{rll}
(\overline{b}_3x_{15},\  b_8)=(x_{15},\  b_3b_8)=1,\ \\
(\overline{b}_3x_{15},\  x_{10})=(x_{15},\  b_3x_{10})=1,\  \\
( \overline{b}_3x_{15},\  \overline{b}_3x_{15})=( x_{15}b_3,\
x_{15}b_3)=6,\
\end{array}
       \end{equation*}
  Therefore the only possibility $\overline{b}_3x_{15}=b_8+x_{10}+x_{17}+m_3+z_3+z_4.$ So $(b_3m_3,\  x_{15})=(m_3,\
  \overline{b}_3x_{15})=1$
  a contradiction.
  Hence the only possibility is that $\alpha=d_3+d_9$,\  hence by 1) $x_9c_3=\overline{x}_{15}+d_3+d_9$. If $d_3=c_3$,\  then
  $(c_3^2,\  x_9)=(c_3,\  c_3x_9)=1$,\  a contradiction to 1) in {\bf Lemma \ref{lem5.1}}. If $d_3=b_3$,\  then $
  (c_3b_3,\  x_9)=(b_3 ,\  c_3x_9)=1$ a contradiction to 2) in {\bf Lemma \ref{lem5.1}} in the same way . If $d_3=\bar{b_3}$,\  a  contradiction
  so $d_3\neq c_3,\  b_3,\  \overline{b}_3$.
\end{proof}

  \begin{lemma} \label{lem5.5}Let$(A,\  B)$ satisfy {\bf Hypothesis \ref{hyp2}} and $( b_8b_3,\
b_8b_3)=3$,\  then {\bf Lemma \ref{lem5.4}} hold and we get the
following :

 1) $ d_3c_3 = x_9,\  x_9 \neq \overline{x}_9,\ $

 2) $ d_3b_8 = \overline{x}_{15}+\overline{x}_9,\ $

 3) $b_3\overline{d}_3 = y_9,\  y_9=\overline{y}_9,\  y_9\neq x_9,\  $

 4) $x_9c_3 = \overline{x}_{15}+d_3+\overline{x}_9,\ $

 5) $ y_9\overline{b}_3 = \overline{d}_3+x_9+x_{15},\ $

 6) $ \overline{d}_3b_6 = d_3+\overline{x}_{15},\ $

 7) $b_8x_{15} =5 x_{15}+3x_9+2x_6+b_3+\overline{d}_3,\ $

 8) $ x_9b_8 = 3x_{15}+2x_9+x_6+\overline{d}_3,\ $

 9) $y_{15}x_6 =4\overline{x}_{15}+2\overline{x}_9+\overline{x}_6+\overline{b}_3+d_3,\ $

 10) $b_6x_{15} =4\overline{x}_{15}+\overline{x}_6+2\overline{x}_9+\overline{b}_3+d_3,\ $

 11) $d_3y_{15} = 2x_{15}+x_6+x_9,\ $

 12) $x_9x_{15} =6\overline{x}_{15}+2\overline{x}_6+3\overline{x}_9+\overline{b}_3+d_3,\ $

 13) $d_3b_8= \overline{x}_{15}+\overline{x}_9,\ $

 14) $\overline{d}_3d_3 = 1+c_8,\  c_8\neq b_8,\ $

 15) $ d_3^2 = b_6+y_3,\  d_3\neq c_3,\  b_3,\  \overline{b}_3,\  y_3\neq b_3,\  c_3,\  \overline{b}_3,\ $

 16) $ b_3d_3=z_9,\  z_9\neq y_9,\  x_9,\  \overline{x}_9,\  z_9=\overline{z}_9,\ $

 17) $y_3^2=1+c_8.$
\end{lemma}

\begin{proof} By {\bf Lemma \ref{lem5.4}},\  $(c_3d_3,\  x_9)=(d_3,\  c_3x_9)=1$,\  hence
$$d_3c_3=x_9. \eqno (1)$$
By the associative law and 1) in {\bf Lemma \ref{lem5.4}},\  we
have
 $ d_3(c_3^2) = x_9c_3$ and $
  d_3+d_3b_8=\overline{x}_{15}+d_3+d_9,\ $ which leads to
  $$d_3b_8 = \overline{x}_{15}+d_9.\eqno (2)$$

So $(b_8x_{15},\  \overline{d}_3)=(b_8d_3,\  \overline{x}_{15})=1$
and by 10) in {\bf Lemma \ref{lem5.3}} ,\ $ (y_9\overline{b}_3,\
\overline{d}_3)=1=(y_9,\  b_3\overline{d}_3).$ So
$$b_3\overline{d}_3=y_9.\eqno (3)$$
 By the
associative law
$(\overline{b}_3b_3)\overline{d}_3=\overline{b}_3(b_3\overline{d}_3)$
and {\bf Hypothesis \ref{hyp2}},\  one has
  $\overline{d}_3+\overline{d}_3b_8 =y_9\overline{b}_3. $ By (2),\
  we get
$$y_9\overline{b}_3=\overline{d}_3+\overline{d}_9+x_{15}.\eqno (4)$$ By the
associative law and {\bf Hypothesis \ref{hyp2}} and (3) we have
$(b_3^2)\overline{d}_3=(b_3\overline{d}_3)b_3=y_9b_3,\
\overline{d}_3c_3+\overline{d}_3b_6=y_9b_3=d_3+d_9+\overline{x}_{15}.$
So
$$\overline{d}_3c_3=d_9 ,\ \eqno(5)$$
$$ \overline{d}_3b_6=d_3+\overline{x}_{15}.\eqno (6)$$ Hence by (1)
we get $\overline{d}_9=x_9$. By {\bf Lemma \ref{lem5.4}},\  we may
assume
$$x_9c_3=\overline{x}_{15}+d_3+\overline{x}_9.$$ By 10) in {\bf Lemma \ref{lem5.3}} and (4)
we meet $$ b_8x_{15}=5x_{15}+3x_9+2x_6+b_3+\overline{d}_3.\eqno
(7)$$ By the associative law and 1) in {\bf Lemma \ref{lem5.1}}
and {\bf Lemma \ref{lem5.4}}:
\begin{equation*}
(c_3^2)x_9 = (c_3x_9)c_3,\ \\
x_9+x_9b_8 = \overline{x}_{15}c_3+d_3c_3+\overline{x}_9c_3.
\end{equation*}
By (1)  and 3) in {\bf Lemma \ref{lem5.3}} and  {\bf Lemma
\ref{lem5.4}},\  we get $$ x_9+x_9b_8 =
2x_{15}+x_6+x_9+x_{15}+\overline{d}_3+x_9,\ $$ so
$$x_9b_8 =3x_{15}+2x_9+x_6+\overline{d}_3.\eqno (8)$$

By the associative law and 4) in {\bf Lemma \ref{lem5.3}} and 2)
in {\bf Lemma \ref{lem5.1}},\  we have that $
y_{15}(\overline{b}_3c_3) = (y_{15}\overline{b}_3)c_3,\ $
$y_{15}b_3+y_{15}x_6 =2x_{15}c_3+x_6c_3+x_9c_3.  $ Moreover by 4)
in {\bf Lemma \ref{lem5.2}} and 3),\  4) in {\bf Lemma
\ref{lem5.3}} and {\bf Lemma \ref{lem5.4}} we get $$
y_{15}b_3+y_{15}x_6 =
4\overline{x}_{15}+2\overline{x}_6+2\overline{x}_9+\overline{b}_3+\overline{x}_{15}+\overline{x}_{15}+d_3+\overline{x}_9,\
$$ then $$y_{15}x_6 =
4\overline{x}_{15}+2\overline{x}_9+\overline{x}_6+\overline{b}_3+d_3.\eqno
(9)$$
 By the associative law and 8) in {\bf Lemma \ref{lem5.3}} and {\bf Hypothesis \ref{hyp2}} we get
$ (b_3^2)x_{15} = (b_3x_{15})b_3,\   x_{15}c_3+b_6x_{15} =
2y_{15}b_3+b_6b_3+y_9b_3,\  $
 By 3),\  4) in {\bf Lemma \ref{lem5.3}}  and 10) in {\bf Lemma \ref{lem5.2}} and equation (4) we have
 $2\overline{x}_{15}+\overline{x}_6+\overline{x}_9+b_6x_{15}=4\overline{x}_{15}+2\overline{x}_6+2\overline{x}_9+\overline{b}_3+\overline{x}_{15}+d_3+\overline{x}_9+\overline{x}_{15},\ $
 so $b_6x_{15}=4\overline{x}_{15}+\overline{x}_6+2\overline{x}_9+\overline{b}_3+d_3.$
 By the associative law and 6) in {\bf Lemma \ref{lem5.2}} and equations (1),\  (6) and (8),\  we have
$(c_3b_8)d_3 = (c_3d_3)b_8= x_9b_8$. Hence
 \begin{eqnarray*}
 b_6d_3+c_3d_3+y_{15}d_3 &=& 3x_{15}+2x_9+x_6+\overline{d}_3,\ \\
\overline{d}_3+x_{15}+x_9+y_{15}d_3 &=& 3x_{15}+2x_9+x_6+\overline{d}_3,\ \\
 y_{15}d_3 &=& 2x_{15}+x_6+x_9.
 \end{eqnarray*}
 So by the associative law and (1),\  it follow  $(d_3c_3)y_{15}=
 c_3(d_3y_{15})$ and $x_9y_{15}= 2x_{15}c_3+x_9c_3+x_6c_3.$
Hence by 3) in {\bf Lemma \ref{lem5.3}} and 4) in {\bf Lemma
\ref{lem5.2}} we get that
$$x_9c_3=\overline{x}_{15}+d_3+\overline{x}_9,\ \ \ x_9x_{15}=6\overline{x}_{15}+2\overline{x}_6+3\overline{x}_9+\overline{b}_3+d_3.$$
By the associative law 1) and  in {\bf Lemma \ref{lem5.1}} and
Lemma \ref{lem5.4} and equation (1) on has that
   $d_3(c_3^2)=(d_3c_3)c_3,\ $
   $d_3+d_3b_8 = c_3x_9 =\overline{x}_{15}+d_3+\overline{x}_9,\ $ and $d_3b_8=\overline{x}_{15}+\overline{x}_9.$
Then we have by (6) that $(d_3^2,\  b_6)=(d_3,\
\overline{d}_3b_6)=1,\ $ so $d_3^2=b_6+y_3(10),\ $ hence
    $( d_3\overline{d}_3,\  d_3\overline{d}_3)  =( d_3^2,\  d_3^2)=2.$
    Now by equation (1) and 1) in {\bf Lemma \ref{lem5.1}} we get
    $( d_3c_3,\  d_3c_3)=( \overline{d}_3d_3,\  c_3^2)=1.$ Hence
    $$d_3\overline{d}_3=1+c_8.\eqno (11)$$
where $c_8\neq b_8$ by main theorem in \cite{cgyarad},\  so
$d_3\neq c_3,\  b_3,\  \bar{b_3} $. Also $y_3\neq
      c_3,\  b_3,\  \overline{b}_3,\ $ since if $y_3=c_3$ ,\ then $1=(d_3^2,\  c_3)=(\overline{d}_3c_3,\  d_3),\ $ therefore
      $
      (\overline{d}_3d_3,\  c_3\overline{d}_3)=(\overline{d}_3c_3,\  \overline{d}_3c_3)\geq2$,\  a contradiction.
      in the same way if $y_3=b_3,\   \overline{b}_3  $,\  If $x_9=\overline{x}_9$,\  then
      by (1) we have $1=( c_3d_3,\  c_3\overline{d}_3) =(\overline{d}_3^2,\
      c_3^2)$,\
      a contradiction to 1) in {\bf Lemma \ref{lem5.1}} and that $d_3
      ^2=b_6+y_3$,\  so $x_9\neq \overline{x}_9.$ By equation(11) and {\bf Hypothesis \ref{hyp2}},\  it holds that
      $1=(\overline{d}_3d_3,\  \overline{b}_3b_3)
      =( d_3b_3,\  d_3b_3),\ $ hence $d_3b_3=z_9$. By (10) and
      hypothesis \ref{hyp2}
     ,\  we get $
      1=( d_3^2,\  \overline{b}_3^2)=( d_3b_3,\  \overline{d}_3\bar{b}_3),$ so $z_9=\bar{z_9},\ $ hence $z_9\neq x_9,\
     \overline{x}_9$. If $z_9=y_9$ ,\ then by (3) we have
     $1=( d_3b_3,\  b_3\overline{d}_3)=( d_3^2,\
     \overline{b}_3b_3),\ $
       a contradiction to (10) and {\bf Hypothesis \ref{hyp2}},\  hence $z_9\neq y_9$.
       By the associative law and equations (11),\  (6) and (10):
       \begin{equation*}
       (d_3b_6)\overline{d}_3 = b_6(d_3\overline{d}_3),\ \\
      \overline{d}_3^2+x_{15}\overline{d}_3 = b_6+c_8b_6,\
       \end{equation*}
       then$$ c_8b_6=\overline{y}_3+x_{15}\overline{d}_3,\  c_8,\  b_6$$ are real s and $$(x_{15}\overline{d}_3,\  y_3)=(x_{15},\
       d_3y_3)=0,\ $$
        so $y_3=\overline{y}_3$,\  by (10) ,\ $(\overline{d}_3y_3,\  d_3)=1$ hence $$2\leq(\overline{d}_3y_3,\  \overline{d}_3y_3)=
        ( d_3\overline{d}_3,\  y_3^2),\ $$ so $y_3^2=1+c_8$.
        \end{proof}

       \begin{lemma}\label{lem5.6} Let$(A,\  B)$ satisfy {\bf Hypothesis \ref{hyp2}} and $( b_8b_3,\  b_8b_3)=3$,\  then {\bf Lemma \ref{lem5.5}} hold and we get the following :

       1) $ b_8^2 = 1+2b_8+2x_{10}+z_4+c_8+b_5+c_5,\ $

       2) $\overline{x}_{15}b_3 = x_{10}+b_8+c_8+z_9+b_5+c_5,\ \mbox{either } b_5,\  c_5 \mbox{are reals
       or }
       b_5=\overline{c}_5,\ $

       3) $ b_6^2 = 1+b_8+z_9+c_8+b_5+c_5,\ $

       4) $\overline{x}_6x_6 = 1+b_8+c_8+z_9+b_5+c_5,\ $

       5) $x_{10}^2 = 1+2b_8+2x_{10}+2c_8+2b_5+2c_5+3z_9,\ $

       6) $ x_{10}\overline{x}_6 =
       \overline{b}_3+3\overline{x}_{15}+d_3+\overline{x}_9,\ $

       7) $ b_8x_{10} =
       2b_8+x_{10}+\overline{x}_9b_3+c_8+z_9+b_5+c_5,\ $

       8) $c_3y_{15} = x_{10}+b_8+c_8+z_9+b_5+c_5,\ $

       9) $ c_3b_5 = y_{15},\ $

       10) $c_3c_5 = y_{15}.$
\end{lemma}

     \begin{proof} At first we notice that
      either $b_5$ and $c_5$ are reals or $b_5=\overline{c_5}$ holds. But in these two cases,\
      it always follows that $\overline{b_5+c_5}=\overline{b}_5+\overline{c_5}=b_5+c_5$.
     By 2) in {\bf Lemma \ref{lem5.2}} and 1),\  8) in {\bf Lemma \ref{lem5.3}},\  we have
     $(\overline{x}_{15}b_3,\  \overline{x}_{15}b_3)=6,\ \ (\overline{x}_{15}b_3,\  x_{10})=(b_3x_{10},\ \ x_{15})=1,\
     (\overline{x}_{15}b_3,\  b_8)=(b_3b_8,\  x_{15})$,\
     then $$\overline{x}_{15}b_3=x_{10}+b_8+x+y+z+w,\  \eqno  (1) $$
     where $|x+y+z+w|=27$. Here the case $\bar{x}_{15}b_3=x_{10}+b_8+2x$ is impossible
      since $|2x|=27$ will implies that $|x|=13.5$.
     But $x$ is an integer, a contradiction. By the associative law $(\overline{b}_3b_8)b_3=b_8(\overline{b}_3b_3)$
      and \textbf{Hypothesis 8.1} and 2),\  3) in {\bf Lemma \ref{lem5.2}},\  we get
     \begin{equation*}
     \overline{x}_6b_3+\overline{b}_3b_3+\overline{x}_{15}b_3 = b_8+b_8^2,\ \\
     b_8+x_{10}+1+b_8+\overline{x}_{15}b_3 = b_8+b_8^2,\ \\
     1+x_{10}+b_8+\overline{x}_{15}b_3 = b_8^2.
     \end{equation*}
     Thus
     $$ b_8^2=1+2b_8+2x_{10}+x+y+z+w .\eqno  (2)$$
     So by 4) in {\bf Lemma \ref{lem5.1}} and 3) in {\bf Lemma \ref{lem5.2}},\  $$b_6^2=1+b_8+x+y+z+w.\eqno(3)$$
      By the associative law $(\overline{b}_3c_3)(b_3c_3)=(\overline{b}_3b_3)c_3^2$ and 1),\  2) in {\bf Lemma \ref{lem5.1}}
      and {\bf Hypothesis \ref{hyp2}} and 3) in {\bf Lemma \ref{lem5.2}} :
     \begin{eqnarray*}
     (b_3+x_6)(\overline{b}_3+\overline{x}_6) = (1+b_8)^2 &=& 1+2b_8+b_8^2,\ \\
     b_3\overline{b}_3+b_3\overline{x}_6+\overline{b}_3x_6+x_6\overline{x}_6 &=& 1+2b_8+b_8^2.
     \end{eqnarray*}
      So by (2) and (1) and 3) in\textbf{ Lemma 8.2}:  $$x_6\overline{x}_6=1+b_8+x+y+z+w. \eqno (4)$$
      By the associative law $(\overline{x}_6b_8)b_3=(\overline{x}_6b_3)b_8$ and 3) in {\bf Lemma \ref{lem5.2}} and 2) in {\bf Lemma \ref{lem5.3}} :
      \begin{equation*}
 2\overline{x}_{15}b_3+\overline{b}_3b_3+\overline{x}_6b_3+\overline{x}_9b_3 = b_8^2+x_{10}b_8,\
     \end{equation*}
     then by  3) in {\bf Lemma \ref{lem5.2}} and equations (1) and (2) we have\\
\begin{center}$2x_{10}+2b_8+2x+2y+2z+2w+1+b_8+b_8+x_{10}+\overline{x}_9b_3$\\
     $= 1+2b_8+2x_{10}+x+y+z+w+x_{10}b_8,\ $\end{center}
 so    $$2b_8+x_{10}+\overline{x}_9b_3+x+y+z+w= x_{10}b_8. \eqno
 (5)$$
     By the associative law $(\overline{x}_6^2)b_3 = \overline{x}_6(\overline{x}_6b_3)$
     and 7) in {\bf Lemma \ref{lem5.3}} and 3) in {\bf Lemma \ref{lem5.2}},\  we get
     $2b_6b_3+y_{15}b_3+y_9b_3=
     \overline{x}_6b_8+\overline{x}_6x_{10}$. Then
      So by 1) in {\bf Lemma \ref{lem5.2}} and 2),\  4) in {\bf Lemma \ref{lem5.3}} and 5) in {\bf Lemma \ref{lem5.5}},\
      we have equation:
     \begin{equation*}
  2\overline{b}_3+2\overline{x}_{15}+2\overline{x}_{15}+\overline{x}_6+\overline{x}_9+d_3+\overline{x}_{15}+\overline{x}_9 =
  2\overline{x}_{15}+\overline{b}_3+\overline{x}_6+\overline{x}_9+x_{10}\overline{x}_6,\
  \end{equation*}
  so $$\overline{b}_3+3\overline{x}_{15}+d_3+\overline{x}_9 =x_{10}\overline{x}_6.\eqno
  (6)$$

  By the associative $(\overline{x}_6b_3)x_{10}= b_3(\overline{x}_6x_{10})$ and 3) in {\bf Lemma \ref{lem5.2}} :
  we get $b_8x_{10}+x_{10}^2
  =\overline{b}_3b_3+3\overline{x}_{15}b_3+d_3b_3+\overline{x}_9b_3.$
   Then by (3) and {\bf Hypothesis \ref{hyp2}} and (1) and 5) in {\bf Lemma \ref{lem5.5}} it
   follows that
   \begin{center}
   $2b_8+x_{10}+\overline{x}_9b_3+x+y+z+w+x_{10}^2$\\
   $=1+b_8+3x_{10}+3b_8+3x+3y+3z+3w+z_9+\overline{x}_9b_3.$
   \end{center}
   Hence $$ x_{10}^2 =1+ 2b_8+2x_{10}+z_9+2x+2y+2z+2w.\eqno
   (7)$$

   By 13) and 14)in {\bf Lemma \ref{lem5.5}},\  one has that $ ( d_3\overline{d}_3,\  b_8^2)=( d_3b_8,\  d_3b_8)
   =2$,\  hence we get by (2) that $(b_8^2,\  c_8)=1 $. Without loss of generality,\  say  $x=c_8$.
   By (4) and (7),\  it holds that $11\leq( x_{10}^2,\  x_6\overline{x}_6)=
   ( x_{10}\overline{x}_6,\  x_{10}\overline{x}_6).$
    But by (6),\  $( x_{10}\overline{x}_6,\  x_{10}\overline{x}_6)=12$ holds,\
    thus we may assume $y=z_9$. By (1)
    ,\ $(\overline{b}_3\alpha,\  \overline{x}_{15})=(\alpha,\  b_3\overline{x}_{15})=1$ holds,\
     where $\alpha\in\{x,\  y,\  z,\  w\}$.
    So $|\alpha|\geq5$. Hence $z=b_5$ and $w=c_5$. Hence by equations
    (1),\  (3),\  (4),\  (5),\  (6) and (7),\  we come to
  \begin{eqnarray*}
  \overline{x}_{15}\bar{b}_3 &=& x_{10}+b_8+c_8+z_9+b_5+c_5,\ \\
  b_8^2 &=& 1+2b_8+2x_{10}+c_8+z_9+b_5+c_5,\ \\
  b_6^2 &=& 1+b_8+c_8+z_9+b_5+c_5,\ \\
  x_6\overline{x}_6 &=& 1+b_8+c_8+z_9+b_5+c_5,\ \\
  x_{10}^2 &=& 1+2b_8+2x_{10}+3z_9+2c_8+2b_5+2c_5,\ \\
  b_8x_{10} &=& 2b_8+x_{10}+\overline{x}_9b_3+c_8+z_9+b_5+c_5.
  \end{eqnarray*}
  So by the associative law $(\overline{b}_3\overline{x}_6)c_3=\overline{x}_6(\overline{b}_3c_3)$ and 1),\  2) in {\bf Lemma \ref{lem5.1}} and 3),\  5) in {\bf Lemma \ref{lem5.2}} :
  \begin{eqnarray*}
   c_3^2+c_3y_{15}&=& b_3\overline{x}_6+x_6\overline{x}_6,\ \\
  1+b_8+c_3y_{15} &=& b_8+x_{10}+1+b_8+c_8+z_9+b_5+c_5,\ \\
  c_3y_{15} &=& x_{10}+b_8+c_8+z_9+b_5+c_5.
   \end{eqnarray*}
   So $(c_3b_5,\  y_{15})=(b_5,\  c_3y_{15})=1,\ \ (c_3c_5,\  y_{15})=(y_{15}c_3,\  c_5)=1,\ $ hence
   $c_3b_5=y_{15},\ \ c_3c_5=y_{15}.$
\end{proof}
   \begin{lemma}\label{lem5.7} Let $(A,\  B) $ satisfy {\bf Hypothesis \ref{hyp2}} and $( b_8b_3,\  b_8b_3)=3$ ,\  then {\bf Lemma \ref{lem5.6}} hold and we get the following :

   1) $ b_3z_9 = x_9+\overline{d}_3+x_{15},\ $

   2) $ \overline{x}_9b_3 = z_9+x_{10}+c_8,\ $

   3) $ b_8x_{10} = 2b_8+2x_{10}+2z_9+2c_8+b_5+c_5,\ $

   4) $b_3c_8 = x_9+x_{15},\ $

   5) $d_3x_6 = x_{10}+c_8,\ $

   6) $ d_3b_6 =\overline{d}_3+x_{15},\ $

   7) $d_3c_8 = \overline{x}_{15}+\overline{x}_6 +d_3,\ $

   8) $y_3\overline{d}_3 = \overline{x}_6+d_3,\ $

   9) $d_3\overline{x}_6 = y_{15}+y_3.$
\end{lemma}

  \begin{proof}By the associative law $(b_3^2)d_3 = b_3(b_3d_3)$ and 6),\  4) in {\bf Lemma \ref{lem5.5}} and {\bf Hypothesis \ref{hyp2}},\  we
  get $ c_3d_3+b_6d_3= b_3z_9.$
   Then by 1),\  6) in {\bf Lemma \ref{lem5.5}} it follows that $x_9+\overline{d}_3+x_{15}=b_3z_9(1). $
   Then $(\overline{b}_3x_9,\  z_9)=(b_3z_9,\
   x_9)=1$.
  And by 1) in {\bf Lemma \ref{lem5.3}},\  $(\overline{x}_9b_3,\  x_{10})=(b_3x_{10},\  x_9)=1$ holds.
    It follows by 4) in {\bf Lemma \ref{lem5.5}} and 1) in {\bf Lemma \ref{lem5.1}} and {\bf Hypothesis \ref{hyp2}}
    that
    $(\overline{x}_9b_3,\  \overline{x}_9b_3)=(\overline{x}_9x_9,\  \overline{b}_3b_3)=
    (\overline{x}_9x_9,\  c_3^2)=( x_9c_3,\  x_9c_3)=3.$
    Hence $$\overline{x}_9b_3=z_9+x_{10}+d_8.\eqno (2)$$
     By 1) and 5) in {\bf Lemma \ref{lem5.6}},\  we get
      $$( b_8^2,\  x_{10}^2)=18=( b_8x_{10},\  b_8x_{10}).\eqno     (3)$$
     By 7) in {\bf Lemma \ref{lem5.6}} and equation (2),\  we have
     $b_8x_{10}=2b_8+2x_{10}+2z_9+c_8+b_5+c_5+d_8,\ $
      then by (3) we get that $c_8=d_8$,\  so
      $$b_8x_{10}=2b_8+2x_{10}+2z_9+2c_8+b_5+c_5. \eqno
      (4)$$
Now we have $(b_3c_8,\  x_9)=(b_3\overline{x}_9,\  c_8)=1 $ by (2)
and $(b_3c_8,\  x_{15})=(b_3\overline{x}_{15},\  c_8)=1$ by 2) in
{\bf Lemma \ref{lem5.6}},\  so $$b_3c_8=x_9+x_{15}.\eqno (5)$$ By
the associative law
$d_3(c_3\overline{b}_3)=(d_3c_3)\overline{b}_3$ and 1) in {\bf
Lemma \ref{lem5.5}} and 2) in {\bf Lemma \ref{lem5.1}} and
equation (2) we see that $ d_3(c_3\overline{b}_3) =
(d_3c_3)\overline{b}_3,\ \
     d_3b_3+d_3x_6= x_9\overline{b}_3=z_9+x_{10}+c_8.$
Hence by 14) in {\bf Lemma \ref{lem5.5}} we have
$$d_3x_6=x_{10} +c_8. \eqno(6)$$
       By the associative law $(b_3^2)d_3=(b_3d_3)b_3$ and 14) in {\bf Lemma \ref{lem5.5}}
       and equation (1) and {\bf Hypothesis \ref{hyp2}},\  it follows
    $ d_3c_3+d_3b_6= z_9b_3=x_9+\overline{d}_3+x_{15}$.
  Hence by 1) in {\bf Lemma \ref{lem5.5}} it holds:
  $$d_3b_6=\bar{d_3}+x_{15}.\eqno (7) $$
  By the associative law $(\overline{d}_3d_3)d_3 = \overline{d}_3(d_3^2)$
  and 12) and 13) in {\bf Lemma \ref{lem5.5}} we have
   $d_3+d_3c_8 = b_6\overline{d}_3+y_3\overline{d}_3,\ $ and get
   $d_3c_8=\overline{x}_{15}+y_3\overline{d}_3$ by (7) and
     $(d_3c_8,\  \overline{x}_6)=(d_3x_6,\
   c_8)=1$ by (6).
   By 12) in {\bf Lemma \ref{lem5.5}},\  $(d_3c_8,\  d_3)=(c_8,\  \overline{d}_3d_3)=1$ holds,\
    hence $$d_3c_8=\overline{x}_{15}+\overline{x}_6+d_3. \eqno (8)
   $$
    Therefore $$y_3\overline{d}_3=\overline{x}_6+d_3.\eqno (9)$$
    Moreover $(d_3\overline{x}_6,\  y_3)=(\overline{x}_6,\  \overline{d}_3y_3)=1,\ $
      and  by 11) in {\bf Lemma \ref{lem5.5}} it follows that $(d_3\overline{x}_6,\  y_{15})=1$,\  so $d_3\overline{x}_6=y_{15}+y_3$.
\end{proof}

\begin{lemma} \label{lem5.8}Let $(A,\  B)$ satisfy {\bf Hypothesis \ref{hyp2}} and $( b_8b_3,\  b_8b_3)=3$ ,\  then {\bf Lemma \ref{lem5.7}} hold and
    we get the following :

   1) $ x_9b_3 = y_9+y_{15}+y_3,\ $

   2) $ x_9\overline{d}_3 = y_{15}+c_3+y_9,\ $

   3) $ x_{15}\overline{d}_3 = y_{15}+y_6+y_9,\ $

    4) $c_8b_6 =2y_{15}+y_3+y_9+b_6,\ $

    5) $z_9b_6 = 2y_{15}+2y_9+b_6,\ $

    6) $ x_{10}b_6 = 3y_{15}+y_3+c_3+y_9,\ $

    7) $ y_3b_6 = c_8+x_{10},\ $

    8) $ d_3x_{15} = b_8+z_9+b_5+c_5+x_{10}+c_8,\ $

    9) $y_{15}b_6 =2b_8+3x_{10}+2z_9+2c_8+b_5+c_5,\ $

    10) $ y_9b_6 = b_8+2z_9+b_5+c_5+x_{10}+c_8,\ $

    11) $\overline{x}_{15}x_{15} = 5
    b_8+6x_{10}+6z_9+5c_8+3b_5+3c_5+1,\ $

    12) $b_6b_5 = y_{15}+y_9+b_6,\ $

    13) $ b_3y_3 = \overline{x}_9$
\end{lemma}

\begin{proof}By the associative law $d_3(c_3b_3)=(d_3c_3)b_3=x_9b_3$ and 1) in {\bf Lemma \ref{lem5.5}} and 2) in {\bf Lemma \ref{lem5.1}},\  we have that $ d_3\overline{b}_3+d_3\overline{x}_6=x_9b_3$ and
by 3) in {\bf Lemma \ref{lem5.5}}  that $y_9+d_3\overline{x}_6
=x_9b_3$. At last,\ we come to  $$x_9b_3=y_9+y_{15}+y_3,\
\eqno(1)$$ by 9) in {\bf Lemma \ref{lem5.7}}.
\par By the associative law $b_3(c_8\overline{d}_3)=\overline{d}_3(b_3c_8)$ and 4),\  7) in {\bf Lemma \ref{lem5.7}},\  we get
    $x_{15}b_3+b_3x_6+b_3\overline{d}_3=
    x_9\overline{d}_3+x_{15}\overline{d}_3,\  $
     and by 5) in {\bf Lemma \ref{lem5.2}} and 3) in {\bf Lemma \ref{lem5.5}} and 8) in {\bf Lemma \ref{lem5.3}},\  we
     get
  $2y_{15}+b_6+y_9+c_3+y_{15}= x_9\overline{d}_3+x_{15}\overline{d}_3 ,\ $
  $ x_9\overline{d}_3+x_{15}\overline{d}_3=3y_{15}+2y_9+b_6+c_3.$
   Now from 10,\  11),\  1) in {\bf Lemma \ref{lem5.5}},\  we see that
    $(x_9\overline{d}_9,\  c_3)=(x_9,\  d_3c_3)=1$ and $(x_9\overline{d}_3,\  y_{15})=(x_9,\
    d_3y_{15})=1$ and $
    (x_9\overline{d}_3,\  b_6)=(x_9,\  d_3b_6)=0$,\  hence $x_9\overline{d}_3=y_{15}+c_3+y_9 (2) $,\
   which gives $$x_{15} \overline{d}_3=2y_{15}+y_9+b_6.\eqno (3) $$
    By the associative law $(\overline{d}_3d_3)b_6= \overline{d}_3(d_3b_6) $ and 6),\  12) in {\bf Lemma \ref{lem5.5}},\  one has that
  $
   b_6+c_8b_6 =\overline{d}_3^2+\overline{d}_3x_{15}$. And by (3) and 13) in
   {\bf Lemma \ref{lem5.5}},\  we get $b_6+c_8b_6= b_6+y_3+2y_{15}+y_9+b_6,\ $ hence
    $$c_8b_6=2y_{15}+y_3+y_9+b_6.\eqno (4)$$
    By the associative law $(d_3b_3)b_6 = d_3(b_3b_6)$ and 1) in {\bf Lemma \ref{lem5.2}} and
    3)  and 16) in {\bf Lemma \ref{lem5.5}},\  we obtain
    $z_9b_6=d_3\overline{b}_3+d_3\overline{x}_{15},\ $ and hence by
     (3) and 3) in {\bf Lemma \ref{lem5.5}},\  we have
   $$z_9b_6=2y_{15}+2y_9+b_6.\eqno (5)$$
   By the associative law $(x_6b_6)d_3= b_6(x_6d_3)$ and 5) in {\bf Lemma \ref{lem5.3}} and 5) in \textbf{Lemma
   8.7},\  we have
   $2\overline{x}_6d_3+\overline{x}_{15}d_3+\overline{x}_9d_3=
   x_{10}b_6+c_8b_6$. Then by 9) in {\bf Lemma \ref{lem5.7}} and equations (2),\  (3) and
   (4) it follows that
     $$2y_{15}+2y_3+2y_{15}+y_9+b_6+y_{15}+c_3+y_9=x_{10}b_6+2y_{15}+y_3+y_9+b_6.$$
     Hence $$x_{10}b_6=3y_{15}+y_3+c_3+y_8.\eqno (6) $$
    Now we have by (4) $(y_3b_6,\  c_8)=(y_3,\  b_6c_8)=1,\  (y_3b_6,\  x_{10})=(y_3,\  b_6x_{10})=1$ by (4) and (6),\
    from which we get  $$y_3b_6=c_8+x_{10}.\eqno (7)$$
     By the associative law $(d_3^2)b_6=d_3(d_3b_6)$
     and 15) in {\bf Lemma \ref{lem5.5}} and 6) in {\bf Lemma \ref{lem5.7}},\  we get
     $b_6^2+y_3b_6=d_3\overline{d}_3+x_{15}d_3$.
Moreover by 14) in {\bf Lemma \ref{lem5.5}} and 3) in {\bf Lemma
\ref{lem5.6}} and equation (7) we have
   $$
    1+b_8+z_9+c_8+b_5+c_5+x_{10}+c_8 = 1+c_8+x_{15}d_3,$$
    $$
    b_8+z_9+b_5+c_5+x_{10}+c_8 = x_{15}d_3. \eqno (8)
   $$
     By the associative law and 5) in {\bf Lemma \ref{lem5.2}} and 5) in {\bf Lemma \ref{lem5.3}} :
      \begin{equation*}
      b_3(x_6b_6)=b_6(b_3x_6),\ \\
      2\overline{x}_6b_3+\overline{x}_{15}b_3+\overline{x}_9b_3 = b_6c_3+y_{15}b_6 ,\
      \end{equation*}
      By 2) in {\bf Lemma \ref{lem5.7}} and 2) in {\bf Lemma \ref{lem5.6}} and  3) in {\bf Lemma \ref{lem5.2}} :
     \[
      2b_8+2x_{10}+x_{10}+b_8+z_9+c_8+b_5+c_5+x_{10}+z_9+c_8 =
      b_8+x_{10}+y_{15}b_6,\]
     \[ y_{15}b_6=2b_8+3x_{10}+2z_9+2c_8+b_5+c_5. \eqno(9) \]
     By 3) in {\bf Lemma \ref{lem5.5}} and 1) in {\bf Lemma \ref{lem5.2}} :
     \begin{equation*}
     (b_3b_6)\overline{d}_3 = (b_3\bar{d}_3)b_6=y_9b_6,\ \\
     \overline{b}_3\overline{d}_3+\overline{x}_{15}\overline{d}_3 = y_9b_6,\
     \end{equation*}
     by 16) in {\bf Lemma \ref{lem5.5}} and (8) $$b_8+2z_9+b_5+c_5+x_{10}+c_8=y_9b_6.\eqno
     (10)$$
      By the associative law and 8) in {\bf Lemma \ref{lem5.3}} and 1) in {\bf Lemma \ref{lem5.2}} :
       \begin{equation*}
       (b_3b_6)x_{15}=b_6(b_3x_{15}),\ \\
       \overline{b}_3x_{15}+\overline{x}_{15}x_{15} = 2y_{15}b_6+b_6^2+y_9b_6,\
      \end{equation*}
     so by (10),\  (9) and 2),\  3) in {\bf Lemma \ref{lem5.6}}
       \begin{eqnarray*}
       x_{10}+b_8+z_9+c_8+b_5+c_5+\overline{x}_{15}x_{15}&=&4b_8+6x_{10}+4z_9+4c_8+2c_5+1+b_8+z_9\\
       & &+c_8+b_5+c_5+b_8+2z_9+c_8+x_{10}+b_5+c_5,\   \\
       \overline{x}_{15}x_{15}&=&1+5b_8+6x_{10}+6z_9+5c_8+3b_5+4c_5,\
      \end{eqnarray*}
      By (9) and (10), we get $$ (b_6b_5,\  y_9)=(b_6y_9,\  b_5)=1,\  (b_6b_5,\  y_{15})=(b_5,\  y_{15}b_6)=1
      $$
      and by (3) $$(b_5b_6,\  b_6)=(b_5,\  b_6^2)=1,\ $$ so $$b_5b_6=y_9+y_{15}+b_6.$$
      By (1) $$(b_3y_3,\  \overline{x}_9)=(b_3x_9,\  y_3)
      =1,\ $$
      then $b_3y_3=\bar{x_9}$.
      \end{proof}

     \begin{lemma}\label{lem5.9}  Let $(A,\  B) $ satisfy {\bf Hypothesis \ref{hyp2}} and $( b_8b_3,\  b_8b_3)=3$ ,\  then lemma 8 hold and we get the following:

       1) $y_9d_3 = b_3+x_{15}+x_9,\ $

       2) $ z_9d_3 = \overline{b}_3+\overline{x}_{15}+\overline{x}_9,\ $

       3) $z_9y_9 =3y_{15}+2b_6+2y_9+c_3+y_3,\ $

       4) $ c_8b_8 = 2x_{10}+2z_9+b_8+c_8+b_5+c_5,\ $

       5) $ y_9c_3 = z_9+x_{10}+c_8,\ $

       6) $\overline{x}_{15}x_6 =2b_8+3x_{10}+2z_9+2c_8+b_5+c_5,\ $

       7) $x_{10}d_3 = \overline{x}_{15}+\overline{x}_6+\overline{x}_9,\ $

       8) $ d_3x_9 = x_{10}+b_8+z_9,\ $

       9) $y_3y_9 = z_9+b_8+x_{10},\ $

       10) $ b_6y_9 = 2z_9+b_8+x_{10}+b_5+c_5+c_8,\ $

       11) $x_{15}^2 =9y_{15}+4b_6+6y_9+2c_3+2y_3,\ $

       12) $x_{15}x_6 = c_3 +4y_{15}+b_6+2y_9+y_3$
\end{lemma}

     \begin{proof}:  By the associative law and 13),\  3) in {\bf Lemma \ref{lem5.5}}:
       $$(b_3\overline{d}_3)d_3=(\overline{b}_3d_3)d_3=\overline{b}_3(d_3^2),\  y_9d_3=b_6\overline{b}_3+y_3\overline{b}_3,\ $$ then by 1) in {\bf Lemma \ref{lem5.2}}
       and 13) in {\bf Lemma \ref{lem5.8}} ,\ $$y_9d_3=b_3+x_{15}+x_9.\eqno (1)$$
       By the associative law and 16) ,\  15) in {\bf Lemma \ref{lem5.5}} :
       \begin{equation*}
       (d_3^2)b_3 = (d_3b_3)d_3,\ \\
       b_6b_3+y_3b_3 = z_9d_3,\
       \end{equation*}
       then by 13) in {\bf Lemma \ref{lem5.8}} and 1) in {\bf Lemma \ref{lem5.2}},\
        $$z_9d_3=\bar{b_3}+\bar{x}_{15}+\bar{x_9}\eqno (2)$$ By 16) in {\bf Lemma \ref{lem5.5}} and (1) :
       \begin{equation*}
       (b_3d_3)y_9=z_9y_9=(d_3y_9)b_3 =b_3^2+x_{15}b_3+x_9b_3,\
       \end{equation*}
       by 8) in {\bf Lemma \ref{lem5.3}} and {\bf Hypothesis \ref{hyp2}} and 1) in {\bf Lemma \ref{lem5.8}} :
       \begin{equation*}
       c_3+b_6+2y_{15}+b_6+y_9+y_9+y_{15}+y_3=z_9y_9,\end{equation*}
      \begin{equation*} 3y_{15}+2b_6+2y_9+c_3+y_3=z_9y_9 .\eqno(3)
      \end{equation*}
       By the associative law and {\bf Hypothesis \ref{hyp2}} and 4) in {\bf Lemma \ref{lem5.7}}:
        \begin{equation*}
        (b_3\overline{b}_3)c_8 = \overline{b}_3(b_3c_8),\ \\
        c_8+c_8b_8 = x_9\overline{b}_3+x_{15} \overline{b}_3,\
        \end{equation*}
        then by 2) in {\bf Lemma \ref{lem5.7}} and 2) in {\bf Lemma \ref{lem5.6}}:
        $$c_8b_8=2x_{10}+2z_9+b_8+c_8+b_5+c_5 . \eqno (4)$$  By the associative law and 2) in {\bf Lemma \ref{lem5.1}} and 3) in \textbf{Lemma 8.5}:
       \begin{equation*}
       (b_3c_3)\overline{d}_3 = (b_3\overline{d}_3 )c_3,\ \\
      \overline{b}_3\bar{d}_3+\overline{x}_6\overline{d}_3 = y_9c_3,\
    \end{equation*}
    then by 16) in {\bf Lemma \ref{lem5.5}} and 5) in {\bf Lemma \ref{lem5.7}} : $$y_9c_3=z_9+x_{10}+c_8.\eqno (5)$$ By the associative law and 7) in {\bf Lemma \ref{lem5.3}} and 4) in {\bf Lemma \ref{lem5.2}} :
     \begin{equation*}
     (x_6^2)c_3 = x_6(x_6c_3),\ \\
     2b_6c_3+y_{15}c_3+y_9c_3 = \overline{b}_3x_6+\overline{x}_{15}x_6,\
     \end{equation*}
     then by (5) and 3) in {\bf Lemma \ref{lem5.2}} and 8)in {\bf Lemma \ref{lem5.6}} :
    \[
     2b_8+2x_{10}+x_{10}+b_8+z_9+c_8+b_5+c_5+z_9+x_{10}+c_8=     b_8+x_{10}+\overline{x}_{15}x_6,\]
    \[\overline{x}_{15}x_6=2b_8+3x_{10}+2z_9+2c_8+b_5+c_5. \eqno(6)\]

  By associative law ,\  9) in {\bf Lemma \ref{lem5.7}} and 3) in {\bf Lemma \ref{lem5.2}} :
  \begin{equation*}
  (\overline{x}_6b_3)d_3 = (\overline{x}_6d_3)b_3,\ \\
  b_8d_3+x_{10}d_3 = y_{15}b_3 +y_3b_3,\
  \end{equation*}
  then by 4) in {\bf Lemma \ref{lem5.3}} and 13) in {\bf Lemma \ref{lem5.8}} $$b_8d_3+x_{10}d_3=2\overline{x}_{15}+\overline{x}_6+\overline{x}_9+\overline{x}_9,\ $$
  by 13) in
  {\bf Lemma \ref{lem5.5}} $$x_{10}d_3=\overline{x}_{15}+\overline{x}_6+\overline{x}_9.\eqno (7)$$
  So$$(d_3x_9,\  x_{10})=(d_3x_{10},\  \overline{x}_9)=1,\ $$ and
  by 13) in {\bf Lemma \ref{lem5.5}} ,\  $$ (d_3x_9,\  b_8)=(d_3b_8,\  \overline{x}_9)=1$$ and by (2) $$(d_3x_9,\  z_9)=(d_3z_9,\
  \overline{x}_9)=1,\ $$
  so $$d_3x_9=x_{10}+b_8+z_9.\eqno (8)$$By the associative law and 15) in {\bf Lemma \ref{lem5.5}} and (1) :
 \begin{equation*}
 (d_3^2)y_9 = d_3(d_3y_9),\ \\
 b_6y_9+y_3y_9 = b_3d_3+x_{15}d_3+x_9d_3,\
\end{equation*}
then by 16) in {\bf Lemma \ref{lem5.5}} and (8) and 8) in {\bf
Lemma \ref{lem5.8}} :
  \begin{eqnarray*}
  b_6y_9+y_3y_9 &=& z_9+b_8+z_9+b_5+c_5+x_{10}+c_8+x_{10}+b_8+z_9,\ \\
  b_6y_9+y_3y_9&=&3z_9+2b_8+b_5+c_5+c_8+2x_{10}.
\end{eqnarray*}
by 5),\  6) in {\bf Lemma \ref{lem5.8}} and 9) in {\bf Lemma
\ref{lem5.3}} :
$$(b_6y_9,\  b_8)=(y_9,\  b_6b_8)=1,\  (b_6y_9,\  x_{10})=(y_9,\
b_6x_{10})=1,\  (b_6y_9,\  z_9)=(y_9,\  z_9b_6) =2,\ $$ hence
$$y_3y_9=z_9+b_8+x_{10}.\eqno (9) $$
 therefore $$b_6y_9=2z_9+b_8+x_{10}+b_5+c_5+c_8.\eqno  (10) $$By the associative law and 10) in {\bf Lemma \ref{lem5.5}} and 1) in {\bf Lemma \ref{lem5.2}}:
 \begin{equation*}
 (b_6\overline{b}_3)x_{15} = \overline{b}_3(b_6x_{15}),\ \\
 b_3x_{15}+x_{15}^2  = 4\overline{x}_{15}\overline{b}_3+\overline{x}_6\overline{b}_3+2\overline{x}_9\overline{b}_3+\overline{b}_3^2+d_3\overline{b}_3.
\end{equation*}
by 8) in {\bf Lemma \ref{lem5.3}} and 5) in {\bf Lemma
\ref{lem5.2}} and 1) in {\bf Lemma \ref{lem5.8}} and {\bf
Hypothesis \ref{hyp2}} and 3) in {\bf Lemma \ref{lem5.5}}
\begin{equation*}
x_{15}^2 = 6y_{15}+3b_6+3y_9+c_3+y_{15}+2y_9+2y_{15}+2y_3+c_3+b_6+y_9,\ \\
x_{15}^2 = 9y_{15}+4b_6+6y_9+2c_3+2y_3.
\end{equation*}
By the associative law and 5) in {\bf Lemma \ref{lem5.3}} and 1)
in {\bf Lemma \ref{lem5.2}} :
\begin{equation*}
(\overline{x}_6b_6)b_3 = \overline{x}_6(b_6b_3),\ \\
2x_6b_3+x_{15}b_3+x_9b_3 =
\overline{b}_3\overline{x}_6+\overline{x}_{15}\overline{x}_6,\
\end{equation*}
then by 5) in {\bf Lemma \ref{lem5.2}} and 8) in {\bf Lemma
\ref{lem5.3}} and 1) in {\bf Lemma \ref{lem5.8}} :
\begin{equation*}
2c_3+2y_{15}+2y_{15}+b_6+y_9+y_{15}+y_3 = c_3+y_{15}+\overline{x}_{15}\overline{x}_6,\ \\
x_{15}x_6 = c_3 +4y_{15}+b_6+2y_{9}+y_3.
\end{equation*}
\end{proof}

\begin{lemma}\label{lem5.10}Let$ (A,\  B) $ satisfy {\bf Hypothesis \ref{hyp2}} and $( b_8b_3,\
b_8b_3)=3$,\  then {\bf Lemma \ref{lem5.9}} hold and we get the
following :

 1) $y_3c_3 =z_9,\ $

 2) $z_9c_3 = y_9+y_{15}+y_3,\ $

 3) $y_3b_8 = y_9+y_{15},\ $

 4) $ z_9b_8 = 2z_9+2x_{10}+2c_8+b_8+b_5+c_5,\ $

 5) $ z_9^2 = 1+2x_{10}+2z_9+2b_8+2c_8+b_5+c_5,\ $

 6) $c_8^2 = 1+2c_8+2x_{10}+b_8+z_9+b_5=c_5,\ $

 7) $ c_8c_3 = y_{15}+y_9,\ $

 8) $ x_9^2 = 2b_6+3y_{15}+2y_9+y_3+c_3,\ $

 9) $c_8y_3 = b_6+y_3+y_{15},\ $

 10) $ x_{10}y_3 = b_6+y_3+y_{15},\ $

 11) $ d_3y_3 = \overline{d}_3+x_6,\ $

 12) $\overline{x}_6y_3 = \overline{d}_3+x_{15},\ $

 13) $x_{15}y_3 = 2\overline{x}_{15}+\overline{x}_6+\overline{x}_9.$
\end{lemma}

\begin{proof}:By the associative law and {\bf Hypothesis \ref{hyp2}} and 13) in {\bf Lemma \ref{lem5.8}}:
\begin{equation*}
y_3 (b_3^2) = \overline{x}_9b_3,\ \\
y_3c_3+y_3b_6 = \overline{x}_9b_3,\
\end{equation*}
then by 2) in {\bf Lemma \ref{lem5.7}}
$$y_3c_3+y_3b_6=z_9+x_{10}+c_8,\ $$ then by 7) in {\bf Lemma \ref{lem5.8}}
$$y_3c_3=z_9.\eqno
(1)$$ By 1),\  16) in {\bf Lemma \ref{lem5.5}} and the associative
law and 1) in {\bf Lemma \ref{lem5.8}}:
\begin{equation*}
b_3(c_3d_3) = (b_3c_3)d_3=(b_3d_3)c_3= z_9c_3,\ \ \\
 b_3x_9 = z_9c_3,\
\end{equation*}
then $$z_9c_3=y_9+y_{15}+y_3.\eqno (2)$$  So by (1) and 1) in {\bf
Lemma \ref{lem5.1}} and the associative law :
\begin{equation*}
y_3(c_3^2) = z_9c_3,\ \\
y_3+y_3b_8 = y_9+y_{15}+y_3,\ \end{equation*} Thus
$$y_3b_8=y_9+y_{15}. (3)$$

By (2) and the associative law and 1) in {\bf Lemma \ref{lem5.1}}:
\begin{equation*}
z_9(c_3^2) = y_9c_3+y_{15}c_3+y_3c_3,\ \\
z_9+z_9b_8 = y_9c_3+y_{15}c_3+y_3c_3,\
\end{equation*}
so by (1) and 8) in {\bf Lemma \ref{lem5.6}} and 5) in {\bf Lemma
\ref{lem5.9}} :
\[
z_9+z_9b_8=z_9+x_{10}+c_8+x_{10}+b_8+c_8+z_9+b_5+c_5+z_9,\]
 \[
z_9b_8=2z_9+2x_{10}+2c_8+b_8+b_5+c_5.\eqno (4) \]

By 17) in {\bf Lemma \ref{lem5.5}} and 1) in {\bf Lemma
\ref{lem5.1}} and (1) :
\begin{equation*}
y_3^2c_3^2 = z_9^2,\ \\
(1+c_8)(1+b_8) = z_9^2,\ \\
1+b_8+c_8+c_8b_8 = z_9^2,\
\end{equation*}
by 4) in {\bf Lemma \ref{lem5.9}}
$$z_9^2=1+2x_{10}+2z_9+2b_8+2c_8+b_5+c_5.  \eqno  (5)$$
By the associative law and 14) in {\bf Lemma \ref{lem5.5}} and 7)
in {\bf Lemma \ref{lem5.7}}:
\begin{equation*}
(\overline{d}_3d_3)c_8 = \overline{d}_3(d_3c_8),\ \\
c_8+c_8^2 =
\overline{d}_3\overline{x}_{15}+\overline{x}_6\overline{d}_3+d_3\overline{d}_3,\
\end{equation*}
then by 5) in {\bf Lemma \ref{lem5.7}} and 8) in {\bf Lemma
\ref{lem5.8}} :
\[
c_8+c_8^2=b_8+z_9+b_5+c_5+x_{10}+c_8+x_{10}+c_8+1+c_8,\] \[ c_8^2
=1+2c_8+2x_{10}+b_8+z_9+b_5+c_5.\eqno   (6)\]

By 1),\  12) in {\bf Lemma \ref{lem5.5}} and the associative law :
\begin{equation*}
\overline{d}_3(d_3c_3) = c_3+c_8c_3,\ \\
\overline{d}_3x_9 = c_3+c_8c_3,\ \\
\end{equation*}
then by 2) in {\bf Lemma \ref{lem5.8}} $$c_8c_3=y_{15}+y_9.\eqno
(7)$$ By the associative law and 13) in {\bf Lemma \ref{lem5.8}}
and 17) in {\bf Lemma \ref{lem5.5}} and {\bf Hypothesis
\ref{hyp2}}:
\begin{equation*}
b_3^2y_3^2 = \overline{x}_9^2,\ \\
(c_3+b_6)(1+c_8) = \overline{x}_9^2,\ \\
c_3+b_6+c_3c_8+b_6c_8 = \overline{x}_9^2,\
\end{equation*}
then by (7) and 4) in {\bf Lemma \ref{lem5.8}} :
\begin{equation*}
 c_3+b_6+y_{15}+y_9+2y_{15}+y_3+y_9+b_6 = \overline{x}_9^2,\ \\
 x_9^2 = 2b_6+3y_{15}+2y_9+y_3+c_3.
\end{equation*}
By the associative law and 4), 7) in {\bf Lemma \ref{lem5.8}} and
17) in {\bf Lemma \ref{lem5.5}}:
\begin{equation*}
(y_3^2)b_6 = y_3(y_3b_6),\ \\
b_6+c_8b_6 = c_8y_3+x_{10}y_3,\ \\
2y_{15}+y_3+y_9+2b_6 = c_8y_3+x_{10}y_3,\
\end{equation*}
by 9) in {\bf Lemma \ref{lem5.9}},\  4) in {\bf Lemma
\ref{lem5.8}} and 15) in {\bf Lemma \ref{lem5.5}}: $(c_8y_3,\
b_6)=(y_3,\ c_8b_6)=1,\ (c_8y_3,\  y_3)=(c_8,\  y_3^2)=1,\
(c_8y_3,\ y_9)=(c_8,\ y_3y_9)=0.$ So $$c_8y_3=b_6+y_3+y_{15},
\eqno (9)$$ then
$$x_{10}y_3=b_6+y_9+y_{15}.\eqno (10)$$
By the associative law and 8) in {\bf Lemma \ref{lem5.7}} and 17)
in {\bf Lemma \ref{lem5.5}}:
\begin{equation*}
(y_3^2)\overline{d}_3 = y_3(y_3\overline{d}_3),\ \\
\overline{d}_3+\overline{d}_3c_8 = \overline{x}_6y_3+d_3y_3,\
\end{equation*}
by 7) in  {\bf Lemma \ref{lem5.7}},\
$2\overline{d}_3+x_{15}+x_6=\overline{x}_6y_3+d_3y_3$,\  by 15) in
{\bf Lemma \ref{lem5.5}} ,\  $(d_3y_3,\  \overline{d}_3) =(d_3^2,\
y_3)=1$. So $$d_3y_3=\overline{d}_3+x_6,\eqno (11)$$ then
$$\overline{x}_6y_3=\overline{d}_3+x_{15}.\eqno (12)$$
By the associative law and 6) in {\bf Lemma \ref{lem5.7}} and 7)
in {\bf Lemma \ref{lem5.8}}:
\begin{equation*}
(d_3b_6)y_3 = d_3(y_3b_6),\ \\
\overline{d}_3y_3+x_{15}y_3 = c_8d_3+x_{10}d_3,\
\end{equation*}
then by 7) in {\bf Lemma \ref{lem5.9}} and 7) in {\bf Lemma
\ref{lem5.7}} and (11):
\begin{equation*}
d_3+\overline{x}_6+x_{15}y_3 =
\overline{x}_{15}+\overline{x}_6+d_3+\overline{x}_{15}+\overline{x}_6+\overline{x}_9,\ \\
y_{15}y_3 = 2\overline{x}_{15}+\overline{x}_6+\overline{x}_9.
\end{equation*}
\end{proof}

\begin{lemma}\label{lem5.11}Let $(A,\  B)$ satisfy {\bf Hypothesis \ref{hyp2}} and $( b_8b_3,\
b_8b_3 )=3$,\  then {\bf Lemma \ref{lem5.10}} hold and we get the
following :

(1) $z_9x_{10}=3x_{10}+2b_8+2z_9+2c_8+b_5+c_5$,\

(2) $z_9c_8=2x_{10}+2b_8+2z_9+b_5+c_5+c_8,\ $

(3) $c_8x_{10}=2x_{10}+2c_8+2b_8+2z_9+b_5+c_5$,\

 (4) $ y_9^2=1+2z_9+2b_8+2c_8+2x_{10}+b_5+c_5$,\

 (5) $ x_9c_8=3x_{15}+x_6+2x_9+b_3$,\

 (6) $\overline{x}_9y_3=b_3+x_9+x_{15}$,\

 (7) $\overline{x}_9x_9=1+2b_8+2x_{10}+2c_8+2z_9+b_5+c_5$,\

 (8) $x_9y_9=d_3+3\overline{x}_{15}+2\overline{x}_9+\overline{b}_3+2\overline{x}_6$,\

 (9) $z_9y_3=y_9+c_3+y_{15}$,\

 (10) $ x_9z_9=3x_{15}+2x_9+2x_6+b_3+\overline{d}_3$.
\end{lemma}

\begin{proof}: By the associative law and 7) in {\bf Lemma \ref{lem5.9}} and 16) in {\bf Lemma \ref{lem5.5}} :
\begin{equation*}
(b_3d_3)x_{10} = b_3(d_3x_{10}),\ \\
z_9x_{10} =
\overline{x}_{15}b_3+\overline{x}_6b_3+\overline{x}_9b_3.
\end{equation*}
Then by 2) in {\bf Lemma \ref{lem5.6}} and 3) in {\bf Lemma
\ref{lem5.2}} and 2) in {\bf Lemma \ref{lem5.7}}:
\[
z_9x_{10}=x_{10}+b_8+c_8+z_9+b_5+c_5+b_8+x_{10}+z_9+x_{10}+c_8,\]
\[
z_9x_{10}=3x_{10}+2b_8+2z_9+2c_8+b_5+c_5.\eqno  (1)\]
 By the associative law
and 16) in {\bf Lemma \ref{lem5.5}} and 7) in {\bf Lemma
\ref{lem5.7}}:
\begin{equation*}
(b_3d_3)c_8 = b_3(d_3c_8),\ \\
z_9c_8 = \overline{x}_{15}b_3+\overline{x}_6b_3+d_3b_3,\
\end{equation*}
then by 2) in {\bf Lemma \ref{lem5.6}} and 3) in {\bf Lemma
\ref{lem5.2}} and 16) in {\bf Lemma \ref{lem5.5}}, we get $z_9c_8=
x_{10}+b_8+c_8+z_9+b_5+c_5+b_8+x_{10}+z_9$ and $$ z_9c_8 =
2x_{10}+2b_8+2z_9+b_5+c_5+c_8 .\eqno (2)$$

By the associative law and 7) in {\bf Lemma \ref{lem5.9}} and 14)
in {\bf Lemma \ref{lem5.5}} :
\begin{equation*}
(\overline{d}_3d_3)x_{10} = \overline{d}_3(d_3x_{10}),\ \\
x_{10}+c_8x_{10} =
\overline{d}_3\overline{x}_{15}+\overline{d}_3\overline{x}_6+\overline{d}_3\overline{x}_9,\
\end{equation*}
then by 5) in {\bf Lemma \ref{lem5.7}} and 8) in {\bf Lemma
\ref{lem5.8}} and 8) in {\bf Lemma \ref{lem5.9}},it follows that
$x_{10}+c_8x_{10}=
b_8+z_9+b_5+c_5+x_{10}+c_8+x_{10}+c_8+x_{10}+b_8+z_9$ and
 $$c_8x_{10}=2x_{10}+2c_8+2b_8+2z_9+b_5+c_5.\eqno   (3)$$

By the associative law and 3),\  15) in {\bf Lemma \ref{lem5.5}}
and {\bf Hypothesis \ref{hyp2}},one has
\begin{equation*}
b_3^2\overline{d}_3^2 = y_9^2 ,\ \\
(c_3+b_6)(b_6+y_3) = y_9^2,\ \\
c_3b_6+c_3y_3+b_6^2+b_6y_3 = y_9^2,\
\end{equation*}
then by 3) in {\bf Lemma \ref{lem5.2}} and 1) in {\bf Lemma
\ref{lem5.10}} and 3) in {\bf Lemma \ref{lem5.6}} and 7) in {\bf
Lemma \ref{lem5.8}}we have $ y_9^2=
b_8+x_{10}+z_9+1+b_8+z_9+c_8+b_5+c_5+c_8+x_{10}$, so $$ y_9^2=
1+2z_9+2b_8+2c_8+2x_{10}+b_5+c_5.\eqno  (4)$$

By the associative law and 13) in {\bf Lemma \ref{lem5.8}} and 9)
in {\bf Lemma \ref{lem5.10}} :
\begin{equation*}
(b_3y_3)c_8 = b_3(y_3c_8),\ \\
\overline{x}_9c_8 = b_6b_3+y_3b_3+y_{15}b_3,\
\end{equation*}
then by 13) in {\bf Lemma \ref{lem5.8}} and 1) in {\bf Lemma
\ref{lem5.2}} and 4) in {\bf Lemma \ref{lem5.3}},$
\overline{x}_9c_8=
\overline{b}_3+\overline{x}_{15}+\overline{x}_9+2\overline{x}_{15}+\overline{x}_6+\overline{x}_9$
and $$ \overline{x}_9c_8=
\overline{b}_3+3\overline{x}_{15}+2\overline{x}_9+\overline{x}_6.\eqno
(5)$$

By the associative law and 13) in {\bf Lemma \ref{lem5.8}} and 17)
in {\bf Lemma \ref{lem5.5}} :
\begin{equation*}
(y_3^2)b_3 = y_3(y_3b_3),\ \\
 b_3+c_8b_3 = \overline{x}_9y_3,\
\end{equation*}
then 4) in {\bf Lemma \ref{lem5.7}}
$$\overline{x}_9y_{13}=b_3+x_9+x_{15}.\eqno (6)$$ By the
associative law and 13) in {\bf Lemma \ref{lem5.8}} and (6) :
\begin{equation*}
(b_3y_3)x_9 = b_3(y_3x_9),\ \\
 \overline{x}_9x_9 =
 b_3\overline{b}_3+b_3\overline{x}_9+b_3\overline{x}_{15},\
\end{equation*}
then by {\bf Hypothesis \ref{hyp2}} and 2) in {\bf Lemma
\ref{lem5.7}} and 2) in {\bf Lemma \ref{lem5.6}} :
\begin{equation*}
\overline{x}_9x_9=
1+b_8+z_9+x_{10}+c_8+x_{10}+b_8+c_8+z_9+b_5+c_5,\end{equation*}
\begin{equation*}
\overline{x}_9x_9=1+2b_8+2x_{10}+2x_{10}+2c_8+2z_9+b_5+c_5.\eqno
(7)
\end{equation*}
By the associative law and 13) in {\bf Lemma \ref{lem5.8}} and 9)
in {\bf Lemma \ref{lem5.9}} :
\begin{equation*}
(b_3y_3)y_9 = b_3 (y_3y_9),\ \\
\overline{x}_9y_9=z_9b_3+b_8b_3+x_{10}b_3,\
\end{equation*}
then by 1) in {\bf Lemma \ref{lem5.7}} and 2) in {\bf Lemma
\ref{lem5.2}} and 1) in {\bf Lemma \ref{lem5.3}}:
\[
\overline{x}_9y_9 =
x_9+\overline{d}_3+x_{15}+x_6+b_3+x_{15}+x_{15}+x_6+x_9=
2x_9+3x_{15}+2x_6+\overline{d}_3+b_3. \eqno (8)\] By the
associative law and 16) in {\bf Lemma \ref{lem5.5}} and 11) in
{\bf Lemma \ref{lem5.1}}0:
\begin{equation*}
(b_3d_3)y_3 = (d_3y_3)b_3 ,\ \\
z_9y_3 = \overline{d}_3b_3+x_6b_3,\
\end{equation*}
then by 3) in {\bf Lemma \ref{lem5.5}} and 5) in {\bf Lemma
\ref{lem5.2}}
$$z_9y_3=y_9+c_3+y_{15}.\eqno (9)$$ By the associative law and 13) in
{\bf Lemma \ref{lem5.8}} and (9) :
\begin{equation*}
(z_9y_3)b_3 = z_9(y_3b_3),\ \\
y_9b_3+c_3b_3+y_{15}b_3 = z_9\overline{x}_9,\
\end{equation*}
then 5) in {\bf Lemma \ref{lem5.5}} and 2) in {\bf Lemma
\ref{lem5.11}} and 4) in {\bf Lemma \ref{lem5.3}} :
\begin{equation*}
z_9\overline{x}_9=d_3+\overline{x}_9+\overline{x}_{15}+\overline{b}_3+\overline{x}_6+2\overline{x}_{15}+\overline{x}_6+\overline{x}_9
=3\overline{x}_{15}+2\overline{x}_9+2\overline{x}_6+\overline{b}_3+d_3.
\eqno(10)
\end{equation*}
\end{proof}

\begin{lemma}\label{lem5.12} Let $(A,\  B)$ satisfy {\bf Hypothesis \ref{hyp2}} and $( b_8b_3,\
b_8b_3 )=3$,\  then {\bf Lemma \ref{lem5.1}} hold and we get the
following :

 (1) $x_{15}c_8=5x_{15}+2x_6+3x_9+b_3+\overline{d}_3$,\

 (2) $x_{15}x_9=6y_{15}+3y_9+2b_6+c_3+y-3,\ $

 (3) $ \overline{x}_9x_{15}=3b_8+4x_{10}+3z_9+3c_8+2b_5+2c_5$,\

 (4) $x_{15}y_{15}=4\overline{x}_6+9\overline{x}_{15}+6\overline{x}_9+2\overline{b}_3+2d_3$,\

 (5) $x_{15}y_9=6\overline{x}_{15}+2\overline{x}_6+3\overline{x}_9+\overline{b}_3+d_3$,\

 (6) $\overline{x}_9x_6=2z_9+x_{10}+b_8+c_8+b_5+c_5,\ $

 (7) $ x_6x_9=2y_9+2z_{15}+b_6$,\

 (8) $ x_6y_9=2\overline{x}_{15}+2\overline{x}_9+\overline{x}_6$,\

 (9) $ x_6z_9=2x_{15}+2x_9+x_6$,\

 (10) $ x_6c_8=2x_{15}+x_6+x_9+\overline{d}_3$,\

 (11) $ y_{15}y_9=4x_{10}+3b_8+3z_9+3c_8+2c_5+2b_5$.
 \end{lemma}

\begin{proof}: By the associative law and 6),\  7) in {\bf Lemma \ref{lem5.7}} :
\begin{eqnarray*}
(d_3b_6)c_8 &=& b_6(d_3c_8),\ \\
\overline{d}_3c_8+x_{15}c_8 &=&
\overline{x}_{15}b_6+\overline{x}_6b_6+d_3b_6,\
\end{eqnarray*}
then by 7) in {\bf Lemma \ref{lem5.7}} and 6),\  10) in {\bf Lemma
\ref{lem5.5}} ,\  5) in {\bf Lemma \ref{lem5.3}} :
\begin{equation*}
x_{15}+x_6\overline{d}_3+x_{15}c_8 =
4x_{15}+x_6+2x_9+b_3+\overline{d}_3+2x_6+x_{15}+x_9+\overline{d}_3+x_{15}
,\end{equation*}
\[ x_{15}c_8 = 5x_{15}+2x_6+3x_9+b_3+\overline{d}_3.\eqno
(1) \]

By the associative law and 8) in {\bf Lemma \ref{lem5.9}} and 6)
in {\bf Lemma \ref{lem5.7}} :
\begin{equation*}
(d_3b_6)x_9 = b_6(d_3x_9),\ \\
\overline{d}_3x_9+x_{15}x_9 = x_{10}b_6+b_8b_6+z_9b_6,\
\end{equation*}
then by and 9) in {\bf Lemma \ref{lem5.3}} and 2),\  5),\  6) in
{\bf Lemma \ref{lem5.8}} :
\begin{equation*}
y_{15}+c_3+y_9+x_{15}x_9=
3y_{15}+y_3+c_3+y_9+2y_{15}+b_6+c_3+y_9+2y_{15}+2y_9+b_6,\end{equation*}
\begin{equation*}x_{15}x_9=6y_{15}+3y_9+2b_6+c_3+y_3.\eqno  (2)
\end{equation*}
By 13) in {\bf Lemma \ref{lem5.8}} and 13) in {\bf Lemma
\ref{lem5.10}} :
\begin{equation*}
b_3(y_3x_{15}) = (b_3y_3)x_{15}\\
= \overline{x}_9x_{15},\ \\
2\overline{x}_{15}b_3+\overline{x}_6b_3+\overline{x}_9b_3 =
\overline{x}_9x_{15}.
\end{equation*}
and 2) in {\bf Lemma \ref{lem5.7}} and 3) in {\bf Lemma
\ref{lem5.2} and 2) in {\bf Lemma \ref{lem5.6}}} :
$$ 2x_{10}+2b_8+2c_8+2z_9+2b_5+2c_5+b_8+x_{10}+z_9+x_{10}+c_8=
\overline{x}_9x_{15},\ $$
$$ \overline{x}_9x_{15}=3b_8+4x_{10}+3z_9+3c_8+2b_5+2c_5.\eqno   (3)
$$
By 11),\  12) in {\bf Lemma \ref{lem5.9}} and by (2) and 8) in
{\bf Lemma \ref{lem5.3}} and 3) in {\bf Lemma \ref{lem5.8}} :
$$(x_{15}y_{15},\  \overline{x}_6)=(x_{15}x_6,\  y_{15} ) =4,\
(x_{15}y_{15},\  \overline{x}_{15})=(x_{15}^2,\  y_{15})=9 ,\ $$
 $$(x_{15}y_{15},\  \overline{x}_9)=(x_{15}x_9,\  y_{15})=6,\
 (x_{15}y_{15},\  \overline{b}_3)=(x_{15}b_3,\  y_{15})=2,\  $$
 $$(x_{15}y_{15},\  d_3)=(x_{15}\overline{d}_3,\  y_{15})=2.$$
Then
$$x_{15}y_{15}=4\overline{x}_6+9\overline{x}_{15}+6\overline{x}_9+2\overline{b}_3+2d_3
.\eqno(4)$$  By (2) and 11),\  12) in {\bf Lemma \ref{lem5.9}} and
8) in {\bf Lemma \ref{lem5.3}} and 3) in {\bf Lemma \ref{lem5.8}}:
$$(x_{15}y_9,\  \overline{x}_{15})=(y_9,\  x_{15}^2) =6,\  (x_{15}y_9,\
\overline{x}_6)=(x_{15}x_6,\  y_9)=2,\  $$
$$ (x_{15}y_9,\  \overline{x}_9)=(x_{15}x_9,\  y_9)=3,\  (x_{15}y_9,\
\overline{b}_3)=(x_{15}b_3,\  y_9)=1,\  $$
$$ (x_{15}y_9,\  d_3)=(x_{15}\overline{d}_3,\  y_9)=1.$$
Then
$$x_{15}y_9=6\overline{x}_{15}+2\overline{x}_6+3\overline{x}_9+\overline{b}_3+d_3
.\eqno (5)$$
 By the associative law and 13) in {\bf Lemma \ref{lem5.8}} and 12) in {\bf Lemma \ref{lem5.10}}:
$$b_3(y_3x_6)=(b_3y_3)x_6=\overline{x}_9x_6,\ $$
$$ b_3d_3+b_3\overline{x}_{15}=\overline{x}_9x_6,\ $$
then by 16) in {\bf Lemma \ref{lem5.5}} and 2) in {\bf Lemma
\ref{lem5.6}}
$$2z_9+x_{10}+b_8+c_8+b_5+c_5=\overline{x}_9x_6.\eqno  (6)$$
By the associative law and 13) in {\bf Lemma \ref{lem5.8}} and 12)
in {\bf Lemma \ref{lem5.10}}:
\begin{equation*}
b_3(y_3\overline{x}_6) = \overline{x}_9\overline{x}_6,\ \\
b_3\overline{d}_3+x_{15}b_3 = \overline{x}_9\overline{x}_6,\
\end{equation*}
so by 3) in {\bf Lemma \ref{lem5.5}} and 8) in {\bf Lemma
\ref{lem5.3}}
$$\overline{x}_6\overline{x}_9=2y_9+2y_{15}+b_6.\eqno  (7) $$
By the associative law and 3) in {\bf Lemma \ref{lem5.5}} and 9)
in {\bf Lemma \ref{lem5.7}} :
$$b_3(\overline{d}_3x_6)=(b_3\overline{d}_3)x_6=x_6y_9,\ $$
$$b_3y_{15}+b_3y_3=x_6y_9.\eqno  (8) $$
By 16) in {\bf Lemma \ref{lem5.5}} and 5) in {\bf Lemma
\ref{lem5.7}} and the associative law :
\begin{equation*}
x_6(b_3d_3) = b_3(x_6d_3),\ \\
x_6z_9 = x_{10}b_3+c_8b_3 ,\
\end{equation*}
then by 1) in {\bf Lemma \ref{lem5.3}} and 4) in {\bf Lemma
\ref{lem5.7}} :
$$x_6z_9=x_{15}+x_6+x_9+x_9+x_{15}=2x_{15}+2x_9+x_6.\eqno (9) $$
By (1) and 4) in {\bf Lemma \ref{lem5.6}} and 5) in {\bf Lemma
\ref{lem5.11}} and 7) in {\bf Lemma \ref{lem5.7}}:
$$(x_6c_8,\  x_{15})=(x_6,\  c_8x_{15})=2,\  (x_6c_8,\  x_6)=(x_6\overline{x}_6,\
c_8)=1,\ $$
$$(x_6c_8,\  x_9)=(x_6,\  c_8x_9)=1,\   (x_6c_8,\
\overline{d}_3)=(c_8d_3,\  \overline{x}_6)=1,\ $$ then $$
x_6c_8=2x_{15}+x_6+x_9+\overline{d}_3.\eqno(10)$$ By the
associative law and 9) in {\bf Lemma \ref{lem5.7}} and 1) in {\bf
Lemma \ref{lem5.9}} :
\begin{equation*}
(d_3\overline{x}_6)y_9 = (d_3y_9)\overline{x}_6,\ \\
y_{15}y_9+y_9y_3 =
b_3\overline{x}_6+x_{15}\overline{x}_6+x_9\overline{x}_6,\
\end{equation*}
then by (6) and 6),\  9) in {\bf Lemma \ref{lem5.9}} and 3) in
{\bf Lemma \ref{lem5.2}} :
\begin{eqnarray*}
y_{15}y_9+z_9+b_8+x_{10}&=&
b_8+x_{10}+2b_8+3x_{10}+2z_9+2c_8\\
&&+b_5+c_5+2z_9+x_{10}+b_8+c_8+b_5+c_5,\ \\
y_{15}y_9&=&4x_{10}+3b_8+3z_9+3c_8+2c_5+2b_5,\
\end{eqnarray*}
\end{proof}

\begin{lemma}\label{lem5.13}Let $(A,\  B)$ satisfy {\bf Hypothesis \ref{hyp2}} and $( b_8b_3,\
b_8b_3)=3$,\  then {\bf Lemma \ref{lem5.12}} hold and we get the
following :

1) $ y_3y_{15}=x_{10}+b_8+b_5+c_5+c_8+z_9,\ $

2) $y_{15}^2=1+3b_5+6x_{10}+6z_9+5b_8+3c_5+5c_8,\ $

3) $x_{10}y_{15}=6y_{15}+3b_6+4y_9+c_3+y_3,\ $

4) $z_9x_{10}=3x_{10}+b_5+c_5+2b_8+2z_9+2c_8,\ $

5) $ y_9x_{10}=4y_{15}+2y_9+y_3+b_6+c_3,\ $

6) $y_9c_8=3y_{15}+2y_9+b_6+c_3,\ $

7) $ z_9y_{15}=6y_{15}+2b_6+3y_9+c_3+y_3,\ $

8) $ x_{10}x_{15}=3x_6+6x_{15}+b_3+\overline{d}_3+4x_9,\ $

9) $ x_{10}x_9=4x_{15}+x_6+2x_9+b_3+\overline{d}_3,\ $

10) $b_6x_9=2\overline{x}_{15}+2\overline{x}_9+\overline{x}_6,\ $

11) $b_8y_9=3y_{15}+y_3,\ $

12) $ \overline{x}_6d_3=y_{15}+y_3,\ $

 13) $ b_8y_{15}=5y_{15}+2b_6+3y_9+c_3+y_3,\ $

 14) $x_{15}z_9=2x_6+6x_{15}+b_3+\overline{d}_3+3x_9,\ $

 15) $x_9y_{15}=6\overline{x}_{15}+2\overline{x}_6+3\overline{x}_9+\overline{b}_3+d_3,\ $

 16) $ y_{15}c_8=2b_6+5y_{15}+c_3+y_3+3y_9.$
 \end{lemma}

\begin{proof}: By the associative law and 9) in {\bf Lemma \ref{lem5.7}} and 15) in {\bf Lemma \ref{lem5.5}} and
 7) in {\bf Lemma \ref{lem5.3}} :
\begin{eqnarray*}
d_3^2\overline{x}_6^2 &=& y_{15}^2+2y_{15}y_3+y_3^2,\ \\
(b_6+y_3)(2b_6+y_{15}+y_9)&=&y_{15}^2+2y_{15}y_3+y_3^2,\ \\
2b_6^2+y_{15}b_6+b_6y_9+2b_6y_3+y_3y_{15}+y_3y_9&=&
y_{15}^2+2y_{15}y_3+y_3^2,\
\end{eqnarray*}
by 3) in {\bf Lemma \ref{lem5.6}} and 7) ,\  9) ,\ 10) in {\bf
Lemma \ref{lem5.8}} and 9) in {\bf Lemma \ref{lem5.9}} and 17) in
{\bf Lemma \ref{lem5.5}}:
\begin{eqnarray*}
y_{15}^2+1+c_8+y_{15}y_3 &=&
2+2b_8+2z_9+2c_8+2b_5+2c_5+2b_8+3x_{10}\\
 & &+2z_9+2c_8+b_5+c_5
+b_8+2z_9+b_5+c_5+x_{10}+c_8   (1) .
\end{eqnarray*}
By the associative law and 11),\  15) in {\bf Lemma \ref{lem5.5}}:
\begin{equation*}
(d_3^2)y_{15} = (d_3y_{15})d_3,\ \\
b_6y_{15}+y_3y_{15} = 2x_{15}d_3+x_6d_3+x_9d_3 ,\
\end{equation*}
so by 8) ,\  9) in {\bf Lemma \ref{lem5.8}} and 5) in {\bf Lemma
\ref{lem5.7}} and 8) in {\bf Lemma \ref{lem5.9}} :
\begin{eqnarray*}
2b_8+3x_{10}+2z_9+2c_8+b_5+c_5+y_3y_{15}&=&
2b_8+2z_9+2b_5+2c_5+2x_{10}+2c_8\\
& & +x_{10}+c_8+x_{10}+b_8+z_9 ,\
\end{eqnarray*}
$$y_3y_{15} = x_{10}+b_8+b_5+c_5+c_8+z_9.\eqno  (2)$$

then by (1),\ $$y_{15}^2=1+3b_5+6x_{10}+6z_9+5b_8+3c_5+5c_8.\eqno
(3)$$ By 9) in {\bf Lemma \ref{lem5.8}} and 11) in {\bf Lemma
\ref{lem5.12}} and 8) in {\bf Lemma \ref{lem5.6}} and (2) :
$$(x_{10}y_{15},\  y_{15})=(x_{10},\  y_{15}^2) =6,\  (x_{10}y_{15},\
b_6)=(x_{10},\  y_{15}b_6)=3,\  $$
$$(x_{10}y_{15},\  y_9)=(x_{10},\  y_{15}y_9)=4,\  (x_{10}y_{15},\
c_3)=(x_{10},\  y_{15}c_3) =1,\  $$
$$(y_{15}x_{10},\  y_3)= (x_{10},\  y_{15}y_3)=1.$$
So $$x_{10}y_{15}=6y_{15}+3b_6+4y_9+c_3+y_3.\eqno (4)$$
 By 16) in {\bf Lemma \ref{lem5.5}} and 7) in {\bf Lemma \ref{lem5.9}} :
 $$b_3(d_3x_{10})=(b_3d_3)x_{10}=z_9x_{10},\ $$
 $$\overline{x}_{15}b_3+\overline{x}_6b_3+\overline{x}_9b_3=z_9x_{10},\ $$
 then by 2) in {\bf Lemma \ref{lem5.6}} and 3) in {\bf Lemma \ref{lem5.2}} and 2) in {\bf Lemma \ref{lem5.7}}:
 $$x_{10}+b_8+c_8+z_9+b_5+c_5+b_8+x_{10}+z_9+x_{10}+c_8=z_9x_{10},\ $$
 $$ z_9x_{10}=3x_{10}+b_5+c_5+2b_8+2z_9+2c_8.\eqno (5)$$
 \end{proof}

\begin{lemma}\label{lem5.14}Let $(A,\ B)$ satisfy {\bf Hypothesis \ref{hyp2}} and
$(b_8b_3,\ b_8b_3)=3$,\   then $b_5\neq \overline{c}_5$.
\end{lemma}

\begin{proof}Assume $b_5=\overline{c}_5$.  By 2) in {\bf Lemma \ref{lem5.6}},\ $(
\overline{b}_3b_5,\ \overline{x}_{15})=(b_5,\
b_3\overline{x}_{15})=1,\ (b_3b_5,\
x_{15})=(b_3\overline{x}_{15},\ \overline{b}_5)=1$,\  then
$$ \overline{b}_3b_5=\overline{x}_{15},\ \eqno(1)$$
$$b_3b_5=x_{15}.\eqno(2)$$
By 9),\ 10) in {\bf Lemma \ref{lem5.8}} and 3) in {\bf Lemma
\ref{lem5.6}}:
$$(b_6b_5,\ y_{15})=(b_5,\ b_6y_{15})=1,\ $$
$$(b_6b_5,\ y_9)=(b_5,\ b_6y_9)=1,\ $$
$$(b_6b_5,\ b_6)=(b_5,\ b_{6}^2)=1,\ $$
then
$$b_6b_5=y_{15}+y_9+b_6.\eqno(3)$$
By 4) in {\bf Lemma \ref{lem5.6}} and 6) in {\bf Lemma
\ref{lem5.9}} and 6) in {\bf Lemma \ref{lem5.12}}:
$$(x_6b_5,\ x_6)=(b_5,\ \overline{x}_6x_6)=1,\ $$
$$(x_6b_5,\ x_{15})=(b_5,\ \overline{x}_6x_{15})=1,\ $$
$$(x_6b_5,\ x_9)=(b_5,\ \overline{x}_6x_9)=1$$
then $$x_6b_5=x_6+x_{15}+x_9.\eqno(4)$$ By 6) in {\bf Lemma
\ref{lem5.9}} and 6) in {\bf Lemma \ref{lem5.12}} and 4) in {\bf
Lemma \ref{lem5.6}}:
$$(
\overline{x}_6b_5,\ \overline{x}_{15})=(b_5,\
x_6\overline{x}_{15})=1,\ $$
$$( \overline{x}_6b_5,\ \overline{x}_9)=(
\overline{x}_6x_9,\ \overline{b}_5)=1,\ $$
$$( \overline{x}_6b_5,\ \overline{x}_6)=(
\overline{x}_6x_6,\ \overline{b}_5)=1$$ then
$$
\overline{x}_6b_5=\overline{x}_{15}+\overline{x}_9+\overline{x}_6.\eqno(5)$$
By 6) in {\bf Lemma \ref{lem5.9}} and 8),\ 11) in {\bf Lemma
\ref{lem5.8}} and 2) in {\bf Lemma \ref{lem5.6}} and 3) in {\bf
Lemma \ref{lem5.12}}:
$$(x_{15}b_5,\ x_6)=(b_5,\ \overline{x}_{15}x_6)=1,\ $$
$$(x_{15}b_5,\ \overline{d}_3)=(x_{15}d_3,\ \overline{b}_5)=1,\ $$
$$(x_{15}b_5,\ x_{15})=(b_5,\ \overline{x}_{15}x_{15})=3,\ $$
$$(x_{15}b_5,\ x_9)=(b_5,\ \overline{x}_{15}x_9)=2,\ $$
$$(x_{15}b_5,\ b_3)=(x_{15}\overline{b}_3,\ \overline{b}_5)=1,\ $$
then $$x_{15}b_5=x_6+3x_{15}+2x_9+b_3+\overline{d}_3.\eqno(6)$$ By
6) in {\bf Lemma \ref{lem5.9}} and 8),\ 11) in {\bf Lemma
\ref{lem5.8}} and 3) in {\bf Lemma \ref{lem5.12}} and 2) in {\bf
Lemma \ref{lem5.6}}
$$(
\overline{x}_{15}b_5,\ \overline{x}_6)=(b_5,\
x_{15}\overline{x}_6)=1,\ $$
$$( \overline{x}_{15}b_5,\ d_3)=(b_5,\ x_{15}d_3)=1,\ $$
$$( \overline{x}_{15}b_5,\ \overline{x}_{15})=(
\overline{x}_{15}x_{15},\ \overline{b}_5)=3,\ $$
$$( \overline{x}_{15}b_5,\ \overline{x}_9)=(
\overline{x}_{15}x_9,\ \overline{b}_5)=1,\ $$
$$( \overline{x}_{15}b_5,\ \overline{b}_3)=(
\overline{x}_{15}b_3,\ \overline{b}_5)=1,\ $$ then
$$
\overline{x}_{15}b_5=\overline{x}_6+3\overline{x}_{15}+2\overline{x}_9+\overline{b}_3+d_3.\eqno(7)$$
By 3),\  6) in {\bf Lemma \ref{lem5.12}} and 7) in {\bf Lemma
\ref{lem5.11}}:
$$(x_9b_5,\ x_{15})=(x_9\overline{x}_{15},\ \overline{b}_5)=2,\ $$
$$(x_9b_5,\ x_6)=(b_5,\ \overline{x}_9x_6)=1,\ $$
$$(x_9b_5,\ x_9)=(x_9\overline{x}_9,\ \overline{b}_5)=1,\ $$
then $$x_9b_5=2x_{15}+x_6+x_9.\eqno(8)$$ By 3),\ 6) in {\bf Lemma
\ref{lem5.12}} and 7) in {\bf Lemma \ref{lem5.11}}:
$$( \overline{x}_9b_5,\ \overline{x}_{15})=(
\overline{x}_9x_{15},\ \overline{b}_5)=2,\ $$
$$( \overline{x}_9b_5,\ \overline{x}_6)=(
b_5,\ \overline{x}_9x_6)=1,\ $$
$$(x_9b_5,\ x_9)=(x_9\overline{x}_9x_9,\ \overline{b}_5)=1,\ $$
then $$
\overline{x}_9b_5=2\overline{x}_{15}+\overline{x}_6+\overline{x}_9.\eqno(9)$$
By 9) in {\bf Lemma \ref{lem5.8}} and 1),\ 2) in {\bf Lemma
\ref{lem5.13}} and 11) in {\bf Lemma \ref{lem5.12}} and 8) in {\bf
Lemma \ref{lem5.6}}
$$(y_{15}b_5,\ b_6)=(b_5,\ y_{15}b_6)=1,\ $$
$$(y_{15}b_5,\ y_{15})=(b_5,\ y_{15}^2)=3,\ $$
$$(y_{15}b_5,\ y_3)=(y_{15}y_3,\ \overline{b}_5)=1,\ $$
$$(y_{15}b_5,\ c_3)=(b_5,\ y_{15}c_3)=1,\ $$
$$(y_{15}b_5,\ y_9)=(b_5,\ y_{15}y_9)=2,\ $$
then $$y_{15}b_5=b_6+3y_{15}+2y_9+c_3+y_3.\eqno(10)$$ By 11) in
{\bf Lemma \ref{lem5.12}} and 10) in {\bf Lemma \ref{lem5.8}} and
4) in {\bf Lemma \ref{lem5.11}}:
$$(y_9b_5,\ y_{15})=(b_5,\ y_9y_{15})=2,\ $$
$$(y_9b_5,\ b_6)=(b_5,\ y_9b_6)=1,\ $$
$$(y_9b_5,\ y_9)=(b_5,\ y_{9}^2)=1,\ $$
then $$y_9b_5=2y_{15}+b_6+y_9.\eqno(11)$$ By 8) in {\bf Lemma
\ref{lem5.8}} ,\
$$(d_3b_5,\ \overline{x}_{15})=(d_3x_{15},\ \overline{b}_5)=1,\ $$
$$( \overline{d}_3b_5,\ x_{15})=(b_5,\ d_3x_{15})=1,\ $$
then $$d_3b_5=\overline{x}_{15},\ \eqno(12)$$
$$ \overline{d}_3b_5=x_{15}.\eqno(13)$$
Then by 16) in {\bf Lemma \ref{lem5.5}} and associative law:
$$b_3(b_5d_3)=(b_3d_3)b_5,\ \overline{x}_{15}b_3=z_9b_5,\ $$  then by
2) in {\bf Lemma \ref{lem5.6}},\
$$z_9b_5=x_{10}+b_8+c_8+z_9+b_5+c_5.\eqno(14)$$
By 1) in {\bf Lemma \ref{lem5.13}},\ $(y_3b_5,\ y_{15})=(b_5,\
y_3y_{15})=1$,\   then
$$y_3b_5=y_{15}.\eqno(15)$$
Then $1=(y_3b_5,\ y_3\overline{b}_5)=(b_{5}^2,\ y_{3}^2)$,\   then
by 17) in {\bf Lemma \ref{lem5.5}},\
$$(b_{5}^2,\ c_8)=1,\ \eqno(16)$$
By (1),\ (2) we get $1=(b_3\overline{b}_5,\
b_3b_5)=(b_3\overline{b}_3,\ b_{5}^2)$,\   then by {\bf Hypothesis
\ref{hyp2}},\
$$(b_{5}^2,\ b_8)=1.\eqno(17)$$
By (2),\ (12) $1=(\overline{d}_3\overline{b}_5,\
b_3b_5)=(\overline{d}_3\overline{b}_3,\ b_{5}^2)$,\ then by 16) in
{\bf Lemma \ref{lem5.5}} $(b_{5}^2,\ z_9)=1$,\  so by (16),\
(17),\ we get
$$b_{5}^2=c_8+b_8+z_9.\eqno(18)$$
By (1) $1=(\overline{b}_3b_5,\
\overline{b}_3b_5)=(\overline{b}_3b_3,\ \overline{b}_5b_5)$,\
then by {\bf Hypothesis \ref{hyp2}},\
$$(b_5\overline{b}_5,\ b_8)=0.\eqno(19)$$
By (2),\ (4) ,\ $1=(x_6b_5,\ b_3b_5)=(x_6\overline{b}_3,\
\overline{b}_5b_5)$,\   then by 3) in {\bf Lemma \ref{lem5.2}} and
(19):
$$( \overline{b}_5b_5,\ x_{10})=1.$$
By (13),\  $3=(b_6b_5,\ b_6b_5)=(b_{6}^2,\ \overline{b}_5b_5)$,\
so by 3) in {\bf Lemma \ref{lem5.6}},\ $m_5=b_5$,\   or
$m_5=\overline{b}_5$,\  in the two cases $( \overline{b}_5b_5,\
b_5) \neq 0$  or  $( \overline{b}_5b_5,\ \overline{b}_5)\neq 0$, a
contradiction to (18).
\end{proof}

\begin{lemma}\label{lem5.15} Let $(A,\ B)$ satisfy {\bf Hypothesis \ref{hyp2}} and
$(b_8b_3,\ b_8b_3)=3$,\ then {\bf Lemma \ref{lem5.13}} hold and we
get following:

1) $b_3b_5=x_{15}$;

2) $b_3b_5=x_{15}$;

3) $b_6b_5=b_6+y_{15}+y_9$;

4) $b_6c_5=b_6+y_{15}+y_9$;

5) $b_5x_6=x_{15}+x_9+x_6$;

6) $c_5x_6=x_{15}+x_9+x_6$;

7) $x_{15}b_5=x_6+\overline{d}_3+3x_{15}+b_3+2x_9$;

8) $x_{15}c_5=x_6+\overline{d}_3+3x_{15}+b_3+2x_9$;

9) $x_9b_5=2x_{15}+x_6+x_9$;

10) $x_9c_5=2x_{15}+x_6+x_9$;

11) $y_{15}b_5=b_6+y_3+3y_{15}+2y_9+c_3$;

12) $y_{15}c_5=b_6+y_3+3y_{15}+2y_9+c_3$;

13) $y_9b_5=2y_{15}+b_6+y_9$;

14) $y_9c_5=2y_{15}+b_6+y_9$;

15) $d_3b_5=\overline{x}_{15}$;

16) $d_3c_5=\overline{x}_{15}$;

17) $z_9b_5=x_{10}+b_8+c_8+z_9+b_5+c_5$;

18) $z_9c_5=x_{10}+b_8+c_8+z_9+b_5+c_5$;

19) $y_3b_5=y_{15}$;

20) $y_3c_5=y_{15}$;

21) $c_5b_5=c_8+b_8+z_9$;

22) $b_8b_5=x_{10}+b_8+c_8+c_5+z_9$;

23) $b_8c_5=x_{10}+b_8+c_8+c_5+z_9$;

24) $b_5^{2}=1+b_5+z_9+x_{10}$;

25) $c_5^{2}=1+c_5+z_9+x_{10}$;

26) $b_5x_{10}=2x_{10}+b_8+z_9+b_5+c_8$;

27) $c_5x_{10}=2x_{10}+b_8+z_9+c_5+c_8$;

28) $c_8b_5=b_8+z_9+c_5+x_{10}+c_8$;

29) $c_8c_5=b_8+z_9+x_{10}+c_8+b_5$.

\end{lemma}
 \begin{proof} By 2)\ in {\bf Lemma \ref{lem5.6}} ,\ $ (b_3b_5,\ x_{15})=(b_3\overline{x}_{15},\ b_5)=1,\
  (b_3c_5,\ x_{15})=(b_3\overline{x}_{15},\ c_5)=1 $,\
then
     $$ b_3b_5=x_{15},\ \eqno(1)$$ $$ b_3c_5=x_{15}. \eqno(2)$$

 By 9),\ \ 10) in {\bf Lemma \ref{lem5.8}} and 3)\ in {\bf Lemma \ref{lem5.6}}:

$(b_6b_5,\ b_6)=(b_{6}^{2},\ b_5)=1,\ (b_6c_5,\ b_6)=(b_{6}^{2},\
c_5)=1$.

$(b_6b_5,\ y_{15})=(b_{5},\ b_6y_{15})=1,\ (b_6c_5,\
y_{15})=(c_{5},\ b_6y_{15})=1$.

$(b_6b_5,\ y_{9})=(b_{5},\ b_6y_{9})=1,\ (b_6c_5,\ y_{9})=(c_{5},\
b_6y_{9})=1$.

then
$$b_5b_6=b_6+y_{15}+y_{9}\eqno(3),$$
 $$b_6c_5=y_{15}+y_9+b_{6}.\eqno(4)$$

By 4) in {\bf Lemma \ref{lem5.6}} and 6) in {\bf Lemma
\ref{lem5.9}} and 6) in {\bf Lemma \ref{lem5.12}}:

$(b_5x_6,\ x_6)=(b_5,\ \overline{x_{6}}x_6)=1,\ (c_5x_6,\
x_6)=(c_5,\ \overline{x_{6}}x_6)=1$.

$(b_5x_6,\ x_{15})=(x_6\overline{x}_{15},\ b_5)=1,\ (c_5x_6,\
x_{15})=(x_6\overline{x}_{15},\ c_5)=1$.

$(b_5x_6,\ x_{9})=(x_6\overline{x_{9}},\ b_{5})=1,\ (c_5x_6,\
x_{9})=(x_6\overline{x_{9}},\ c_{5})=1$.

then
$$b_5x_6=x_6+x_{15}+x_{9}\eqno(5),$$
 $$c_5x_6=x_{15}+x_{9}+x_{6}.\eqno(6)$$

By 6) in {\bf Lemma \ref{lem5.9}} and 8),\ 11) in {\bf Lemma
\ref{lem5.8}} and 2) in {\bf Lemma \ref{lem5.6}} and 3) in {\bf
Lemma \ref{lem5.12}}:

$(x_{15}b_5,\ x_6)=(b_5,\ \overline{x}_{15}x_6)=1,\ (x_{15}c_5,\
x_6)= (c_5,\ \overline{x}_{15}x_6)=1 $.

$(x_{15}b_5,\ \overline{d}_3)= (x_{15}d_3,\ b_5)=1,\ (x_{15}c_5,\
\overline{d}_3)= (x_{15}d_3,\ c_5)=1$.

$(x_{15}b_5,\ x_{15})=(x_{15}\overline{x}_{15},\ b_5)=3,\
(x_{15}c_5,\ x_{15})=(x_{15}\overline{x}_{15},\ c_5)=3$.

$(x_{15}b_5,\ b_3)=(b_5,\ \overline{x}_{15}b_3)=1,\ (x_{15}c_5,\
b_3)=(c_5,\ \overline{x}_{15}b_3)=1$.

$(x_{15}b_5,\ x_9)=(b_5,\ \overline{x}_{15}x_9)=2,\ (x_{15}c_5,\
x_9)=(x_{15}\overline{x}_9,\ c_5)=2$.

then
$$x_{15}b_5=x_6+\overline{d}_3+3x_{15}+b_3+2x_9\eqno(7),$$\
$$x_{15}c_5=x_6+\overline{d}_3+3x_{15}+b_3+2x_9.
\eqno(8) $$

By 3),\ 6) in {\bf Lemma \ref{lem5.12}} and 7) in {\bf Lemma
\ref{lem5.11}}:

$(x_9b_5,\ x_{15})=(b_5,\ \overline{x}_9x_{15})=2,\ (x_9c_5,\
x_{15})=(c_5,\ \overline{x}_9x_{15})=2$.

$(x_9b_5,\ x_{6})=(b_5,\ \overline{x}_9x_{6})=1,\ (x_9c_5,\
x_{6})=(c_5,\ \overline{x}_9x_{6})=1$.

$(x_9b_5,\ x_{9})=(b_5,\ \overline{x}_9x_{9})=1,\ (x_9c_5,\
x_{9})=(c_5,\ \overline{x}_9x_{9})=1$.

then$$x_{9}b_5=2x_{15}+x_6+x_9, \eqno(9)$$\
$$x_9c_5=2x_{15}+x_6+x_9. \eqno(10)$$

By 9) in {\bf Lemma \ref{lem5.8}} and 1),\ 2) in {\bf Lemma
\ref{lem5.13}} and 11) in {\bf Lemma \ref{lem5.12}} and 8) in {\bf
Lemma \ref{lem5.6}}:

$(y_{15}b_5,\ b_6)=(y_{15}b_6,\ b_5)=1,\ (y_{15}c_5,\
b_6)=(y_{15}b_6,\ c_5)=1$.

$(y_{15}b_5,\ y_3)=(y_{15}y_3,\ b_5)=1,\ (y_{15}c_5,\
y_3)=(y_{15}y_3,\ c_5)=1$.

$(y_{15}b_5,\ y_{15})=(b_5,\ y_{15}^{2})=3,\ (y_{15}c_5,\
y_{15})=(y_{15}^{2},\ c_5)=3$.

$(y_{15}b_6,\ c_3)=(y_{15}c_3,\ b_6)=1,\ (y_{15}c_5,\
c_3)=(y_{15}c_3,\ c_5)=1$.

$(y_{15}b_5,\ y_9)=(y_{15}y_9,\ b_5)=2,\ (y_{15}c_5,\ y_9)=(c_5,\
y_{15}y_9)=2$.

Then $$y_{15}b_5=b_6+y_3+3y_{15}+2y_9+c_3, \eqno(11)$$ \
$$y_{15}c_5=b_6+y_3+3y_{15}+2y_9+c_3. \eqno(12)$$

By 11) in {\bf Lemma \ref{lem5.12}} and 10) in {\bf Lemma
\ref{lem5.8}} and 4) in {\bf Lemma \ref{lem5.11}}:

$(y_9b_5,\ y_{15})=(b_5,\ y_9y_{15})=2,\ (y_9c_5,\ y_{15})=(c_5,\
y_9y_{15})=2$.

$(y_9b_5,\ b_{6})=(y_9b_{6},\ b_5)=1,\ (y_9c_5,\
b_{6})=(y_9b_{6},\ c_5)=1$.

$(y_9b_5,\ y_{9})=(b_5,\ y_9^{2})=1,\ (y_9c_5,\ y_{9})=(c_5,\
y_9^{2})=1$.

then $$y_{9}b_5=2y_{15}+b_6+y_9 ,\eqno(13)$$\
$$y_{9}c_5=2y_{15}+b_6+y_9. \eqno(14) $$

By 8) in {\bf Lemma \ref{lem5.8}},\ $(d_3b_5,\
\overline{x}_{15})=(d_3x_{15},\ b_5)=1,\ (d_3c_5,\
\overline{x}_{15})=(d_3x_{15},\ c_5)=1$.

then $$d_3b_5=\overline{x}_{15}, \eqno(15)$$\
$$d_3c_5=\overline{x}_{15}. \eqno(16) $$

So by the associative law and 16) in {\bf Lemma \ref{lem5.5}}, we
get $(b_3d_3)b_5=b_3(d_3b_5),\ (b_3d_3)c_5=b_3(d_3c_5)$,\
$z_9b_5=b_3\overline{x}_{15},\ z_9c_5=b_3\overline{x}_{15}$,\ then
by 2) in {\bf Lemma \ref{lem5.6}} ,\
$$z_9b_5=x_{10}+b_8+c_8+z_9+b_5+c_5, \eqno(17)$$ \
$$z_9c_5=x_{10}+b_8+c_8+z_9+b_5+c_5. \eqno(18)$$

By 1) in {\bf Lemma \ref{lem5.13}} ,\ $(y_3b_5,\ y_{15})=(b_5,\
y_3y_{15})=1,\ (y_3c_5,\ y_{15})=(c_5,\ y_3y_{15})=1$,\ then
$$y_3b_5=\overline{y}_{15}, \eqno(19)$$\
$$y_3c_5=\overline{y}_{15}. \eqno(20) $$

So $(y_3^{2},\ b_5^{2})=(y_3b_5,\ y_3b_5)=1$,\ $(y_3^{2},\
c_5b_5)=(y_3b_5,\ y_3c_5)=1$,\ then by 17) in {\bf Lemma
\ref{lem5.5}} ,\ $(c_5b_5,\ c_8)=1$. By (1),\ (2),\
$(\overline{b}_3b_3,\ b_5c_5)=(b_3b_5,\ b_3c_5)=1$. then by {\bf
Hypothesis \ref{hyp2}},\ $(c_5b_5,\ b_8)=1$. By (17),\
$$c_5b_5=c_8+b_8+z_9. \eqno(21)$$

By the associative and {\bf Hypothesis \ref{hyp2}} and (1),\ (2),
it follows that$(b_3\overline{b}_3)b_5=(b_3b_5)\overline{b}_3,\
(b_3\overline{b}_3)c_5=(b_3c_5)\overline{b}_3$,\ $b_5 + b_8b_5=
x_{15}\overline{b}_3,\ c_5 + b_8c_5= x_{15}\overline{b}_3$,\  then
by 2) in {\bf Lemma \ref{lem5.6}},\

$$b_8b_5=x_{10}+b_8+c_8+c_5+z_9, \eqno(22)$$
\ $$b_8c_5=x_{10}+b_8+c_8+b_5+z_9. \eqno(23)$$

By 5) in {\bf Lemma \ref{lem5.6}},\ $(x_{10}b_5,\ x_{10})=(b_5,\
x_{10}^{2})=2,\ (x_{10}c_5,\ x_{10})=(c_5,\ x_{10}^{2})=2$. By
(22) and (23),\ $(x_{10}b_5,\ b_{8})=(x_{10},\ b_5b_8)=1,\
(x_{10}c_5,\ b_{8})=(x_{10},\ c_5b_8)=1$. By (17),\ (18),\
$(x_{10}b_5,\ z_{9})=(x_{10},\ z_9b_5)=1,\ (x_{10}c_5,\
b_{8})=(x_{10},\ c_5b_8)=1$, $(b_5^{2},\ z_9)=(b_5,\ b_5z_9)=1$
and $ (c_5^{2},\ z_9)=(c_5,\ c_5z_9)=1$.

By (19),\ (20), we come to
 $(y_3^{2},\ b_5^{2})=(y_3b_5,\ y_3b_5)=1,\ (y_3^{2},\ c_5^{2})=(y_3c_5,\
 y_3c_5)=1$,\ then by 17) in {\bf Lemma \ref{lem5.5}} $(b_5^{2},\ c_8)=0,\ (c_5^{2},\ c_8)=0$.

 By
 (1) and (2), we see that $(\overline{b}_3b_3,\ b_5^{2})=(b_3b_5,\ b_3b_5)=1,\ (\overline{b}_3b_3,\ c_5^{2})=(b_3c_5,\ b_3c_5)=1$,\
 then by {\bf Hypothesis \ref{hyp2}},\  $(b_5^{2},\ b_8)=0,\ (c_5^{2},\ b_8)=0$ holds.

By (5), (6), (15) and (16),\ $(b_5^{2},\
d_3x_6)=(b_5\overline{d}_3,\ b_5x_6)=1,\ (c_5^{2},\
d_3x_6)=(c_5\bar{d}_3,\ c_5x_6)=1$,\ so
 by 5) in {\bf Lemma \ref{lem5.7}},\  $(b_5^{2},\ x_{10})=1,\ (c_5^{2},\ x_{10})=1$.

 By (5),\ (21) and (6) we have the following inner-products: $(b_5^{2},\ \overline{x}_6x_6)=(b_5x_6,\ b_5x_6)=3,\ (c_5^{2},\ \overline{x}_6x_6)=(c_5x_6,\ c_5x_6)=3$.
$(c_5^{2},\ b_5)=(c_5,\ c_5b_5)=0$ and $(b_5^{2},\ c_5)=(b_5,\
b_5c_5)=0$,\ then
 by 4) in {\bf Lemma \ref{lem5.6}}  $(b_5^{2},\ b_5)=1,\ (c_5^{2},\ c_5)=1$,\ hence
 $$b_5^{2}=1+b_5+z_9+x_{10},\ c_5^{2}=1+b_5+z_9+x_{10}.$$
So
 $(b_5x_{10},\ b_5)=(x_{10},\ b_5^{2})=1,\ (c_5x_{10},\ c_5)=(x_{10},\ c_5^{2})=1$.

By 14) in {\bf Lemma \ref{lem5.5}} and (15),\ (16), one has
$b_5(\overline{d}_3d_3)=\overline{d}_3(b_5d_3),\
c_5(\overline{d}_3d_3)=\overline{d}_3(c_5d_3)$. Hence
$b_5+c_8b_5=\overline{x}_{15}\overline{d}_3,\
c_5+c_8c_5=\overline{d}_3\overline{x}_{15}$.

Therefore by 8) in {\bf Lemma \ref{lem5.8}}, we have equations
$c_8b_5=b_8+z_9+c_5+x_{10}+c_8$ and
$c_8c_5=b_8+z_9+x_{10}+c_8+b_5$. Consequently $(b_5x_{10},\
c_8)=(x_{10},\ b_5c_8)=(x_{10},\ c_5c_8)=1$.

At last we get $b_5x_{10}=c_8+2x_{10}+b_8+z_9+b_8+z_9+b_5$ and $
c_5x_{10}=2x_{10}+b_8+z_9+c_5+c_8$.

\end{proof}

\begin{theorem}Let $(A,\ B)$  be a $NITA$ generated by an non-real
element $b_{3}\in B$  of degree 3 and without non-identity basis
element of degree 1 or 2,\ such that: $\overline{b}_{3}b_3=1+b_{8}$
and $b_{3}^{2}=c_3+b_{6},\ c_3=\overline{c}_3,\ (b_8b_3,\
b_8b_3)=3$. Then $(A,\ B)$ is a table algebra of dimension 22:
 $$B=\{b_1,\ \overline{b}_3,\ b_3,\ c_3,\ b_6,\ x_6,\ \overline{x}_6,\ b_8,\ x_{15},\ \overline{x}_{15},\ x_{10},\ x_9,\ \overline{x}_9,\ y_{15},\
 y_9,\ d_3,\ \overline{d}_3,\ z_9,\ c_8,\ b_5,\ c_5,\ y_3\}$$
 and $B$ has increasing series of table subsets $\{b_1\}\subseteq
 \{b_1,\ b_8,\ x_{10},\ z_9,\ c_8,\ b_5,\ c_5\}\subseteq B$ defined by:

 1)  $b_3\overline{b}_3=1+b_8$;

 2)  $b_{3}^2=c_3+b_6$;

 3)  $b_3c_3=\overline{b}_3+\overline{x}_6$;

 4)  $b_3b_6=\overline{b}_3+\overline{x}_{15}$;

 5)  $b_3x_6=c_3+y_{15}$;

 6)  $b_3\overline{x}_6=b_8+x_{10}$;

 7)  $b_3b_8=x_6+b_3+x_{15}$;

 8)  $b_3x_{15}=2y_{15}+b_6+y_9$;

 9)  $b_3x_{10}=x_{15}+x_6+x_9$;

 10)  $b_3\overline{x}_{15}=x_{10}+b_8+z_9+c_8+b_5+c_5$;

 11)  $b_3x_9=y_9+y_{15}+y_3$;

 12)  $b_3\overline{x}_9=x_{10}+z_9+c_8$;

 13)
 $b_3y_{15}=2\overline{x}_{15}+\overline{x}_6+\overline{x}_9$;

 14)  $b_3y_9=d_3+\overline{x}_{15}+\overline{x}_9$;

 15)  $b_3d_3=z_9$;

 16)  $b_3\overline{d}_3=y_9$;

 17)  $b_3z_9=x_9+\overline{d}_3+x_{15}$;

 18)  $b_3c_8=x_{15}+x_9$;

 19)  $b_3b_5=x_{15}$;

 20)  $b_3c_5=x_{15}$;

 21)  $b_3y_3=\overline{x}_9$;

 22)  $c_{3}^2=1+b_8$;

 23)  $c_3b_6=b_8+x_{10}$;

 24)  $c_3x_6=\overline{b}_3+\overline{x}_{15}$;

 25)  $c_3b_8=b_6+c_3+y_{15}$;

 26)
 $c_3x_{15}=2\overline{x}_{15}+\overline{x}_6+\overline{x}_9$;

 27)  $c_3x_{10}=b_6+y_{15}+y_9$;

 28)  $c_3x_9=d_3+\overline{x}_{15}+\overline{x}_9$;

 29)  $c_3y_{15}=x_{10}+b_8+z_9+c_8+b_5+c_5$;

 30)  $c_3y_9=z_9+x_{10}+c_8$;

 31)  $c_3d_3=x_9$;

 32)  $c_3z_9=y_3+y_9+y_{15}$;

 33)  $c_3c_8=y_9+y_{15}$;

 34)  $c_3b_5=y_{15}$;

 35)  $c_3c_5=y_{15}$;

 36)  $c_3y_3=z_9$;

 37)  $b_{6}^2=1+b_8+z_9+c_8+b_5+c_5$;

 38)  $b_6x_6=2\overline{x}_6+\overline{x}_{15}+\overline{x}_9$;

 39)  $b_6b_8=2y_{15}+b_6+y_9+c_3$;

 40)
 $b_6x_{15}=4\overline{x}_{15}+\overline{x}_6+2\overline{x}_9+\overline{b}_3+d_3$;

 41)  $b_6x_{10}=c_3+3y_{15}+y_9+y_3$;

 42)  $b_6x_9=2\overline{x}_{15}+2\overline{x}_9+\overline{x}_6$;

 43)  $b_6y_{15}=2b_8+3x_{10}+2z_9+2c_8+b_5+c_5$;

 44)  $b_6y_9=b_8+2z_9+c_8+x_{10}+b_5+c_5$;

 45)  $b_6d_3=\overline{d}_3+x_{15}$;

 46)  $b_6z_9=2y_{15}+2y_9+b_6$;

 47)  $b_6c_8=b_6+y_3+y_9+2y_{15}$;

 48)  $b_6b_5=y_{15}+y_9+b_6$;

 49)  $b_6c_5=y_{15}+b_6+y_9$;

 50)  $b_6y_3=x_{10}+c_8$;

 51)  $x_{6}^2=2b_6+y_{15}+y_9$;

 52)  $x_6\overline{x}_6=1+b_8+z_9+c_8+b_5+c_5$;

 53)  $x_6b_8=b_3+2x_{15}+x_6+x_9$;

 54)  $x_6x_{15}=4y_{15}+b_6+2y_9+c_3+y_3$;

 55)  $x_6x_{10}=b_3+3x_{15}+\overline{d}_3+x_9$;

 56)  $x_6\overline{x}_{15}=2b_8+3x_{10}+2z_9+2c_8+b_5+c_5$;

 57)  $x_6x_9=2y_{15}+2y_9+b_6$;

 58)  $x_6\overline{x}_9=x_{10}+b_8+2z_9+c_8+b_5+c_5$;

 59)
 $x_6y_{15}=4\overline{x}_{15}+2\overline{x}_9+\overline{x}_6+\overline{b}_3+d_3$;

 60)  $x_6y_9=2\overline{x}_{15}+2\overline{x}_9+\overline{x}_6$;

 61)  $x_6d_3=x_{10}+c_8$;

 62)  $x_6\overline{d}_3=y_{15}+y_3$;

 63)  $x_6z_9=2x_9+2x_{15}+x_6$;

 64)  $x_6c_8=2x_{15}+x_6+x_9+\overline{d}_3$;

 65)  $x_6b_5=x_6+x_{15}+x_9$;

 66)  $x_6c_5=x_{15}+x_9+x_6$;

 67)  $x_6y_3=d_3+\overline{x}_{15}$;

 68)  $b_{8}^2=1+2b_8+2x_{10}+z_9+c_8+b_5+c_5$;

 69)  $b_8x_{10}=2x_{10}+2z_9+2c_8+2b_8+b_5+c_5$;

 70)  $b_8x_9=3x_{15}+2x_9+x_6+\overline{d}_3$;

 71)  $b_8y_{15}=5y_{15}+2b_6+3y_9+c_3+y_3$;

 72)  $b_8y_9=3y_{15}+2y_9+b_6+y_3$;

 73)  $b_8d_3=\overline{x}_{15}+\overline{x}_9$;

 74)  $b_8z_9=2x_{10}+b_8+2c_8+2z_9+b_5+c_5$;

 75)  $b_8c_8=2x_{10}+2z_9+b_8+c_8+b_5+c_5$;

 76)  $b_8b_5=b_8+z_9+c_5+c_8+x_{10}$;

 77)  $b_8c_5=b_8+x_{10}+c_8+z_9+b_5$;

 78)  $b_8y_3=y_9+y_{15}$;

 79)  $b_8x_{15}=5x_{15}+2x_6+3x_9+b_3+\overline{d}_3$;

 80)  $x_{15}^2=9y_{15}+4b_6+6y_9+2c_3+2y_3$;

 81)  $x_{15}x_{10}=3x_6+6x_{15}+b_3+\overline{d}_3+4x_9$;

 82)
 $x_{15}\overline{x}_{15}=5b_8+6x_{10}+6z_9+5c_8+3b_5+3c_5+1$;

 83)  $x_{15}x_9=6y_{15}+3y_9+2b_6+c_3+y_3$;

 84)  $x_{15}\overline{x}_9=3b_8+4x_{10}+3z_9+3c_8+2b_5+2c_5$;

 85)  $x_{15}y_{15}=4\overline{x}_6+9\overline{x}_{15}+6\overline{x}_9+2\overline{b}_3+2d_3$;

 86)  $x_{15}y_9=6\overline{x}_{15}+2\overline{x}_6+3\overline{x}_9+\overline{b}_3+d_3$;

 87)  $x_{15}d_3=b_8+z_9+b_5+c_5+x_{10}+c_8$;

 88)  $x_{15}\overline{d}_3=2y_{15}+y_9+b_6$;

 89)  $x_{15}z_9=2x_6+6x_{15}+b_3+\overline{d}_3+3x_9$;

 90)  $x_{15}c_8=5x_{15}+2x_6+3x_9+b_3+\overline{d}_3$;

 91)  $x_{15}b_5=x_6+3x_{15}+2x_9+b_3+\overline{d}_3$;

 92)  $x_{15}c_5=x_6+3x_{15}+2x_9+b_3+\overline{d}_3$;

 93)
 $x_{15}y_3=\overline{x}_6+2\overline{x}_{15}+\overline{x}_9$;

 94)  $x_{10}^2=1+2b_8+2x_{10}+2c_8+2b_5+2c_5+3z_9$;

 95)  $x_{10}x_9=b_3+4x_{15}+\overline{d}_3+2x_9+x_6$;

 96)  $x_{10}y_{15}=6y_{15}+3b_6+4y_9+c_3+y_3$;

 97)  $x_{10}y_9=4y_{15}+2y_9+y_3+b_6+c_3$;

 98)  $x_{10}d_3=\overline{x}_{15}+\overline{x}_6+\overline{x}_9$;

 99)  $x_{10}z_9=3x_{10}+b_5+c_5+2b_8+2z_9+2c_8$;

 100)  $x_{10}c_8=2x_{10}+2c_8+2b_8+2z_9+b_5+c_5$;

 101)  $x_{10}b_5=b_8+2x_{10}+z_9+c_8+b_5$;

 102)  $x_{10}c_5=b_8+z_9+c_5+2x_{10}+c_8$;

 103)  $x_{10}y_3=y_{15}+y_9+b_6$;

 104)  $x_{9}^2=y_3+c_3+2b_6+3y_{15}+2y_9$;

 105)  $x_9\overline{x}_9=1+2b_8+2x_{10}+2c_8+2z_9+b_5+c_5$;

 106)
 $x_9y_{15}=6\overline{x}_{15}+2\overline{x}_6+3\overline{x}_9+\overline{b}_3+d_3$;

 107)  $x_9y_9=d_3+3\overline{x}_{15}+2\overline{x}_9+\overline{b}_3+2\overline{x}_6$;

 108)  $x_9d_3=x_{10}+b_8+z_9$;

 109)  $x_9\overline{d}_3=y_9+c_3+y_{15}$;

 110)  $x_9z_9=3x_{15}+2x_9+2x_6+b_3+\overline{d}_3$;

 111)  $x_9b_5=2x_{15}+x_6+x_9$;

 112)  $x_9c_5=2x_{15}+x_6+x_9$;

 113)  $x_9y_3=\overline{b}_3+\overline{x}_9+\overline{x}_{15}$;

 114)  $y_{15}^2=1+3b_5+6x_{10}+6z_{9}+5b_8+3c_5+5c_8$;

 115)  $y_{15}y_9=4x_{10}+3b_8+3z_9+3c_8+2c_5+2b_5$;

 116)  $y_{15}d_3=2x_{15}+x_6+x_9$;

 117)  $y_{15}z_9=6y_{15}+2b_6+3y_9+c_3+y_3$;

 118)  $y_{15}c_8=2b_6+5y_{15}+c_3+y_3+3y_9$;

 119)  $y_{15}b_5=b_6+3y_{15}+2y_9+c_3+y_3$;

 120)  $y_{15}c_5=b_6+c_3+3y_{15}+y_3+2y_9$;

 121)  $y_{15}y_3=b_5+b_8+z_9+c_5+x_{10}+c_8$;

 122)  $y_{9}^2=1+2z_9+2b_8+2c_8+2x_{10}+b_5+c_5$;

 123)  $y_9d_3=b_3+x_{15}+x_9$;

 124)  $y_9z_9=3y_{15}+2b_6+2y_9+c_3+y_3$;

 125)  $y_9c_8=3y_{15}+2y_9+b_6+c_3$;

 126)  $y_9b_5=2y_{15}+b_6+y_9$;

 127)  $y_9c_5=2y_{15}+y_9+b_6$;

 128)  $y_9y_3=z_9+b_8+x_{10}$;

 129)  $d_3\overline{d}_3=1+c_8$;

 130)  $d_3^2=b_6+y_3$;

 131)  $d_3z_9=\overline{b}_3+\overline{x}_{15}+\overline{x}_9$;

 132)  $d_3c_8=\overline{x}_{15}+d_3+\overline{x}_6$;

 133)  $d_3b_5=\overline{x}_{15}$;

 134)  $d_3c_5=\overline{x}_{15}$;

 135)  $d_3y_3=\overline{d}_3+x_6$;

 136)  $z_{9}^2=1+2x_{10}+2z_9+2b_8+2c_8+b_5+c_5$;

 137)  $z_9c_8=2x_{10}+2b_8+2z_9+b_5+c_5+c_8$;

 138)  $z_9b_5=x_{10}+b_8+z_9+c_8+b_5+c_5$;

 139)  $z_9c_5=x_{10}+b_8+z_9+c_8+b_5+c_5$;

 140)  $z_9y_3=y_9+c_3+y_{15}$;

 141)  $c_{8}^2=1+2c_8+2x_{10}+b_8+z_9+b_5+c_5$;

 142)  $c_8b_5=b_8+z_9+c_5+x_{10}+c_8$;

 143)  $c_8c_5=b_8+z_9+x_{10}+c_8+b_5$;

 144)  $c_8y_3=b_6+y_3+y_{15}$;

 145)  $b_{5}^2=1+b_5+x_{10}+z_9$;

 146)  $b_5c_5=b_8+z_9+c_8$;

 147)  $b_5y_3=y_{15}$;

 148)  $c_{5}^2=1+z_9+x_{10}+c_5$;

 149)  $c_5y_3=y_{15}$;

 150)  $y_{3}^2=1+c_8$.

\end{theorem}

\begin{proof} The equations hold by lemmas 8.1-8.13 and 8.15.\end{proof}

\end{document}